\documentclass[reqno,11pt]{amsart}
\usepackage[T1]{fontenc}
\usepackage{lmodern,amsmath,amssymb,amsthm,mathrsfs,geometry,microtype,tikz}
\usepackage[hidelinks]{hyperref}
\usepackage[backend=biber,url=false,eprint=false,maxbibnames=99]{biblatex}
\usepackage{seqsplit}
\newcommand*\hidekey[1]{\expandafter\def\csname quietkey@#1\endcsname{}}
\hidekey{lem:elliptic-high-frequency}
\hidekey{zero:constraints}
\hidekey{r2:compact-damping}
\numberwithin{equation}{section}
\newtheorem{theorem}{Theorem}[section]
\newtheorem{lemma}[theorem]{Lemma}
\newtheorem{proposition}[theorem]{Proposition}
\newtheorem{corollary}[theorem]{Corollary}
\newtheorem{definition}{Definition}[section]
\newtheorem{remark}[theorem]{Remark}
\newcommand{\T}{\mathbb T}
\newcommand{\R}{\mathbb R}
\newcommand{\Z}{\mathbb Z}
\newcommand{\eps}{\epsilon}

\newcommand{\bv}{\mathbf{V}}
\allowdisplaybreaks[1]
\title[Linearized Euler around dipoles]{Dynamics of the linearized two-dimensional Euler equation around the symmetric dipole}
\author{Joel Dahne}
\address{University of Minnesota}
\email{jdahne@umn.edu}
\author{Hao Jia}
\thanks{HJ was partially supported by NSF DMS-2453270, NSF DMS-2245021 and a grant
from the Simons Foundation.}
\address{University of Minnesota}
\email{jia@umn.edu}
\begin{document}
\begin{abstract}
We study the two-dimensional Euler equation linearized around the symmetric dipole with stream function $\cos x+\cos y$ on the square torus. For data with
zero streamline average and orthogonal to the (generalized) zero modes, we prove optimal inviscid damping estimates, and sharp vorticity depletion at the elliptic critical points. The proof is based on a limiting absorption principle and regularity estimates for the spectral density function. The main difficulty is to control the singularities near the elliptic critical points and the separatrix in a genuinely two-dimensional setting. 
\end{abstract}

\maketitle
\tableofcontents
\section{Introduction and main results}

Consider the two-dimensional incompressible Euler equation in the vorticity formulation
\begin{equation}\label{Int1}
    \begin{split}
        &\partial_t\omega+u\cdot\nabla \omega=0,\\
        &u=\nabla^\perp \psi=(-\partial_y\psi,\partial_x\psi),\quad \Delta\psi=\omega,\\
        &\omega|_{t=0}=\omega_0,
    \end{split}
\end{equation}
for $(x,y)\in\T^2,\quad t\ge0$.

\begin{figure}[t]
\centering
\begin{tikzpicture}[x=.8cm,y=.8cm]
\draw[gray] (-3.14159,-3.14159) rectangle (3.14159,3.14159);
\draw[gray] (2.0944,-3.1416) -- (2.0967,-3.0788) -- (2.1036,-3.0159) -- (2.1150,-2.9531) -- (2.1252,-2.9112) -- (2.1375,-2.8693) -- (2.1572,-2.8135) -- (2.1776,-2.7646) -- (2.1975,-2.7227) -- (2.2197,-2.6808) -- (2.2410,-2.6443) -- (2.2619,-2.6112) -- (2.2859,-2.5761) -- (2.3170,-2.5342) -- (2.3457,-2.4987) -- (2.3693,-2.4714) -- (2.4083,-2.4295) -- (2.4295,-2.4083) -- (2.4714,-2.3693) -- (2.4987,-2.3457) -- (2.5342,-2.3170) -- (2.5761,-2.2859) -- (2.6112,-2.2619) -- (2.6443,-2.2410) -- (2.6808,-2.2197) -- (2.7227,-2.1975) -- (2.7646,-2.1776) -- (2.8135,-2.1572) -- (2.8693,-2.1375) -- (2.9112,-2.1252) -- (2.9531,-2.1150) -- (3.0159,-2.1036) -- (3.0788,-2.0967) -- (3.1416,-2.0944) -- (3.1416,-2.0944);
\draw[gray] (-3.1416,-2.0944) -- (-3.0997,-2.0954) -- (-3.0369,-2.1008) -- (-2.9740,-2.1107) -- (-2.9322,-2.1198) -- (-2.8732,-2.1363) -- (-2.8274,-2.1519) -- (-2.7855,-2.1685) -- (-2.7437,-2.1873) -- (-2.7018,-2.2083) -- (-2.6599,-2.2317) -- (-2.6180,-2.2575) -- (-2.5804,-2.2829) -- (-2.5515,-2.3038) -- (-2.5133,-2.3337) -- (-2.4744,-2.3667) -- (-2.4504,-2.3884) -- (-2.4086,-2.4293) -- (-2.3876,-2.4513) -- (-2.3511,-2.4923) -- (-2.3248,-2.5244) -- (-2.3011,-2.5552) -- (-2.2714,-2.5970) -- (-2.2443,-2.6389) -- (-2.2201,-2.6801) -- (-2.1991,-2.7195) -- (-2.1782,-2.7633) -- (-2.1599,-2.8065) -- (-2.1445,-2.8484) -- (-2.1311,-2.8903) -- (-2.1153,-2.9515) -- (-2.1069,-2.9950) -- (-2.0985,-3.0578) -- (-2.0946,-3.1206) -- (-2.0944,-3.1416);
\draw[gray] (3.1416,2.0944) -- (3.0788,2.0967) -- (3.0159,2.1036) -- (2.9531,2.1150) -- (2.9112,2.1252) -- (2.8693,2.1375) -- (2.8135,2.1572) -- (2.7646,2.1776) -- (2.7227,2.1975) -- (2.6808,2.2197) -- (2.6443,2.2410) -- (2.6112,2.2619) -- (2.5761,2.2859) -- (2.5342,2.3170) -- (2.4987,2.3457) -- (2.4714,2.3693) -- (2.4295,2.4083) -- (2.4083,2.4295) -- (2.3693,2.4714) -- (2.3457,2.4987) -- (2.3170,2.5342) -- (2.2859,2.5761) -- (2.2619,2.6112) -- (2.2410,2.6443) -- (2.2197,2.6808) -- (2.1975,2.7227) -- (2.1776,2.7646) -- (2.1572,2.8135) -- (2.1375,2.8693) -- (2.1252,2.9112) -- (2.1150,2.9531) -- (2.1036,3.0159) -- (2.0967,3.0788) -- (2.0944,3.1416) -- (2.0944,3.1416);
\draw[gray] (-2.0944,3.1416) -- (-2.0954,3.0997) -- (-2.1008,3.0369) -- (-2.1107,2.9740) -- (-2.1198,2.9322) -- (-2.1363,2.8732) -- (-2.1519,2.8274) -- (-2.1685,2.7855) -- (-2.1873,2.7437) -- (-2.2083,2.7018) -- (-2.2317,2.6599) -- (-2.2575,2.6180) -- (-2.2829,2.5804) -- (-2.3038,2.5515) -- (-2.3337,2.5133) -- (-2.3667,2.4744) -- (-2.3884,2.4504) -- (-2.4293,2.4086) -- (-2.4513,2.3876) -- (-2.4923,2.3511) -- (-2.5244,2.3248) -- (-2.5552,2.3011) -- (-2.5970,2.2714) -- (-2.6389,2.2443) -- (-2.6801,2.2201) -- (-2.7195,2.1991) -- (-2.7633,2.1782) -- (-2.8065,2.1599) -- (-2.8484,2.1445) -- (-2.8903,2.1311) -- (-2.9515,2.1153) -- (-2.9950,2.1069) -- (-3.0578,2.0985) -- (-3.1206,2.0946) -- (-3.1416,2.0944);
\draw[gray] (1.5708,-3.1416) -- (1.5728,-3.0788) -- (1.5787,-3.0159) -- (1.5885,-2.9531) -- (1.5972,-2.9112) -- (1.6127,-2.8513) -- (1.6265,-2.8065) -- (1.6411,-2.7646) -- (1.6574,-2.7227) -- (1.6755,-2.6803) -- (1.6965,-2.6356) -- (1.7174,-2.5943) -- (1.7387,-2.5552) -- (1.7630,-2.5133) -- (1.7888,-2.4714) -- (1.8163,-2.4295) -- (1.8431,-2.3907) -- (1.8640,-2.3617) -- (1.8917,-2.3248) -- (1.9247,-2.2829) -- (1.9478,-2.2548) -- (1.9773,-2.2201) -- (2.0106,-2.1825) -- (2.0338,-2.1572) -- (2.0735,-2.1156) -- (2.0944,-2.0944) -- (2.1363,-2.0535) -- (2.1596,-2.0316) -- (2.1991,-1.9957) -- (2.2300,-1.9687) -- (2.2619,-1.9418) -- (2.3038,-1.9080) -- (2.3337,-1.8850) -- (2.3667,-1.8604) -- (2.4086,-1.8306) -- (2.4504,-1.8023) -- (2.4851,-1.7802) -- (2.5195,-1.7593) -- (2.5557,-1.7383) -- (2.5970,-1.7159) -- (2.6389,-1.6948) -- (2.6808,-1.6753) -- (2.7297,-1.6546) -- (2.7853,-1.6336) -- (2.8274,-1.6198) -- (2.8693,-1.6076) -- (2.9322,-1.5927) -- (2.9740,-1.5848) -- (3.0369,-1.5763) -- (3.0997,-1.5717) -- (3.1416,-1.5708) -- (3.1416,-1.5708);
\draw[gray] (-3.1416,-1.5708) -- (-3.0997,-1.5717) -- (-3.0369,-1.5763) -- (-2.9740,-1.5848) -- (-2.9322,-1.5927) -- (-2.8693,-1.6076) -- (-2.8274,-1.6198) -- (-2.7853,-1.6336) -- (-2.7297,-1.6546) -- (-2.6808,-1.6753) -- (-2.6389,-1.6948) -- (-2.5970,-1.7159) -- (-2.5557,-1.7383) -- (-2.5195,-1.7593) -- (-2.4851,-1.7802) -- (-2.4504,-1.8023) -- (-2.4086,-1.8306) -- (-2.3667,-1.8604) -- (-2.3337,-1.8850) -- (-2.3038,-1.9080) -- (-2.2619,-1.9418) -- (-2.2300,-1.9687) -- (-2.1991,-1.9957) -- (-2.1596,-2.0316) -- (-2.1363,-2.0535) -- (-2.0944,-2.0944) -- (-2.0735,-2.1156) -- (-2.0338,-2.1572) -- (-2.0106,-2.1825) -- (-1.9773,-2.2201) -- (-1.9478,-2.2548) -- (-1.9247,-2.2829) -- (-1.8917,-2.3248) -- (-1.8640,-2.3617) -- (-1.8431,-2.3907) -- (-1.8163,-2.4295) -- (-1.7888,-2.4714) -- (-1.7630,-2.5133) -- (-1.7387,-2.5552) -- (-1.7174,-2.5943) -- (-1.6965,-2.6356) -- (-1.6755,-2.6803) -- (-1.6574,-2.7227) -- (-1.6411,-2.7646) -- (-1.6265,-2.8065) -- (-1.6127,-2.8513) -- (-1.5972,-2.9112) -- (-1.5885,-2.9531) -- (-1.5787,-3.0159) -- (-1.5728,-3.0788) -- (-1.5708,-3.1416) -- (-1.5708,-3.1416);
\draw[gray] (3.1416,1.5708) -- (3.0788,1.5728) -- (3.0159,1.5787) -- (2.9531,1.5885) -- (2.9112,1.5972) -- (2.8513,1.6127) -- (2.8065,1.6265) -- (2.7646,1.6411) -- (2.7227,1.6574) -- (2.6803,1.6755) -- (2.6356,1.6965) -- (2.5943,1.7174) -- (2.5552,1.7387) -- (2.5133,1.7630) -- (2.4714,1.7888) -- (2.4295,1.8163) -- (2.3907,1.8431) -- (2.3617,1.8640) -- (2.3248,1.8917) -- (2.2829,1.9247) -- (2.2548,1.9478) -- (2.2201,1.9773) -- (2.1825,2.0106) -- (2.1572,2.0338) -- (2.1156,2.0735) -- (2.0944,2.0944) -- (2.0535,2.1363) -- (2.0316,2.1596) -- (1.9957,2.1991) -- (1.9687,2.2300) -- (1.9418,2.2619) -- (1.9080,2.3038) -- (1.8850,2.3337) -- (1.8604,2.3667) -- (1.8306,2.4086) -- (1.8023,2.4504) -- (1.7802,2.4851) -- (1.7593,2.5195) -- (1.7383,2.5557) -- (1.7159,2.5970) -- (1.6948,2.6389) -- (1.6753,2.6808) -- (1.6546,2.7297) -- (1.6336,2.7853) -- (1.6198,2.8274) -- (1.6076,2.8693) -- (1.5927,2.9322) -- (1.5848,2.9740) -- (1.5763,3.0369) -- (1.5717,3.0997) -- (1.5708,3.1416);
\draw[gray] (-1.5708,3.1416) -- (-1.5717,3.0997) -- (-1.5763,3.0369) -- (-1.5848,2.9740) -- (-1.5927,2.9322) -- (-1.6076,2.8693) -- (-1.6198,2.8274) -- (-1.6336,2.7853) -- (-1.6546,2.7297) -- (-1.6753,2.6808) -- (-1.6948,2.6389) -- (-1.7159,2.5970) -- (-1.7383,2.5557) -- (-1.7593,2.5195) -- (-1.7802,2.4851) -- (-1.8023,2.4504) -- (-1.8306,2.4086) -- (-1.8604,2.3667) -- (-1.8850,2.3337) -- (-1.9080,2.3038) -- (-1.9418,2.2619) -- (-1.9687,2.2300) -- (-1.9957,2.1991) -- (-2.0316,2.1596) -- (-2.0535,2.1363) -- (-2.0944,2.0944) -- (-2.1156,2.0735) -- (-2.1572,2.0338) -- (-2.1825,2.0106) -- (-2.2201,1.9773) -- (-2.2548,1.9478) -- (-2.2829,1.9247) -- (-2.3248,1.8917) -- (-2.3617,1.8640) -- (-2.3907,1.8431) -- (-2.4295,1.8163) -- (-2.4714,1.7888) -- (-2.5133,1.7630) -- (-2.5552,1.7387) -- (-2.5943,1.7174) -- (-2.6356,1.6965) -- (-2.6803,1.6755) -- (-2.7227,1.6574) -- (-2.7646,1.6411) -- (-2.8065,1.6265) -- (-2.8513,1.6127) -- (-2.9112,1.5972) -- (-2.9531,1.5885) -- (-3.0159,1.5787) -- (-3.0788,1.5728) -- (-3.1416,1.5708) -- (-3.1416,1.5708);
\draw[gray] (1.0472,-3.1416) -- (1.0495,-3.0788) -- (1.0562,-3.0159) -- (1.0675,-2.9531) -- (1.0774,-2.9112) -- (1.0892,-2.8693) -- (1.1100,-2.8073) -- (1.1265,-2.7646) -- (1.1444,-2.7227) -- (1.1639,-2.6808) -- (1.1850,-2.6389) -- (1.2075,-2.5970) -- (1.2314,-2.5552) -- (1.2566,-2.5133) -- (1.2776,-2.4801) -- (1.2985,-2.4480) -- (1.3252,-2.4086) -- (1.3546,-2.3667) -- (1.3823,-2.3287) -- (1.4032,-2.3007) -- (1.4329,-2.2619) -- (1.4660,-2.2201) -- (1.4870,-2.1941) -- (1.5174,-2.1572) -- (1.5499,-2.1188) -- (1.5708,-2.0944) -- (1.6075,-2.0525) -- (1.6336,-2.0233) -- (1.6642,-1.9897) -- (1.6965,-1.9549) -- (1.7228,-1.9268) -- (1.7593,-1.8888) -- (1.7834,-1.8640) -- (1.8221,-1.8248) -- (1.8459,-1.8012) -- (1.8850,-1.7630) -- (1.9105,-1.7383) -- (1.9478,-1.7031) -- (1.9774,-1.6755) -- (2.0106,-1.6451) -- (2.0467,-1.6127) -- (2.0735,-1.5890) -- (2.1153,-1.5528) -- (2.1435,-1.5289) -- (2.1782,-1.5000) -- (2.2200,-1.4661) -- (2.2464,-1.4451) -- (2.2829,-1.4168) -- (2.3248,-1.3852) -- (2.3574,-1.3614) -- (2.3876,-1.3398) -- (2.4295,-1.3109) -- (2.4714,-1.2831) -- (2.5133,-1.2566) -- (2.5479,-1.2357) -- (2.5842,-1.2147) -- (2.6223,-1.1938) -- (2.6628,-1.1729) -- (2.7064,-1.1519) -- (2.7540,-1.1310) -- (2.8065,-1.1103) -- (2.8484,-1.0958) -- (2.8903,-1.0831) -- (2.9504,-1.0681) -- (2.9950,-1.0595) -- (3.0578,-1.0512) -- (3.1206,-1.0474) -- (3.1416,-1.0472);
\draw[gray] (-3.1416,-1.0472) -- (-3.0997,-1.0482) -- (-3.0369,-1.0535) -- (-2.9740,-1.0633) -- (-2.9322,-1.0722) -- (-2.8698,-1.0891) -- (-2.8274,-1.1028) -- (-2.7855,-1.1182) -- (-2.7437,-1.1353) -- (-2.7018,-1.1540) -- (-2.6599,-1.1743) -- (-2.6180,-1.1961) -- (-2.5761,-1.2193) -- (-2.5342,-1.2438) -- (-2.4923,-1.2697) -- (-2.4504,-1.2969) -- (-2.4169,-1.3195) -- (-2.3867,-1.3404) -- (-2.3457,-1.3698) -- (-2.3038,-1.4008) -- (-2.2733,-1.4242) -- (-2.2410,-1.4493) -- (-2.1991,-1.4829) -- (-2.1686,-1.5080) -- (-2.1363,-1.5350) -- (-2.0944,-1.5708) -- (-2.0704,-1.5917) -- (-2.0316,-1.6262) -- (-2.0002,-1.6546) -- (-1.9687,-1.6835) -- (-1.9326,-1.7174) -- (-1.9059,-1.7428) -- (-1.8672,-1.7802) -- (-1.8431,-1.8040) -- (-1.8040,-1.8431) -- (-1.7802,-1.8672) -- (-1.7428,-1.9059) -- (-1.7174,-1.9326) -- (-1.6835,-1.9687) -- (-1.6546,-2.0002) -- (-1.6262,-2.0316) -- (-1.5917,-2.0704) -- (-1.5708,-2.0944) -- (-1.5350,-2.1363) -- (-1.5080,-2.1686) -- (-1.4829,-2.1991) -- (-1.4493,-2.2410) -- (-1.4242,-2.2733) -- (-1.4008,-2.3038) -- (-1.3698,-2.3457) -- (-1.3404,-2.3867) -- (-1.3195,-2.4169) -- (-1.2969,-2.4504) -- (-1.2697,-2.4923) -- (-1.2438,-2.5342) -- (-1.2193,-2.5761) -- (-1.1961,-2.6180) -- (-1.1743,-2.6599) -- (-1.1540,-2.7018) -- (-1.1353,-2.7437) -- (-1.1182,-2.7855) -- (-1.1028,-2.8274) -- (-1.0891,-2.8698) -- (-1.0722,-2.9322) -- (-1.0633,-2.9740) -- (-1.0535,-3.0369) -- (-1.0482,-3.0997) -- (-1.0472,-3.1416) -- (-1.0472,-3.1416);
\draw[gray] (3.1416,1.0472) -- (3.0788,1.0495) -- (3.0159,1.0562) -- (2.9531,1.0675) -- (2.9112,1.0774) -- (2.8693,1.0892) -- (2.8073,1.1100) -- (2.7646,1.1265) -- (2.7227,1.1444) -- (2.6808,1.1639) -- (2.6389,1.1850) -- (2.5970,1.2075) -- (2.5552,1.2314) -- (2.5133,1.2566) -- (2.4801,1.2776) -- (2.4480,1.2985) -- (2.4086,1.3252) -- (2.3667,1.3546) -- (2.3287,1.3823) -- (2.3007,1.4032) -- (2.2619,1.4329) -- (2.2201,1.4660) -- (2.1941,1.4870) -- (2.1572,1.5174) -- (2.1188,1.5499) -- (2.0944,1.5708) -- (2.0525,1.6075) -- (2.0233,1.6336) -- (1.9897,1.6642) -- (1.9549,1.6965) -- (1.9268,1.7228) -- (1.8888,1.7593) -- (1.8640,1.7834) -- (1.8248,1.8221) -- (1.8012,1.8459) -- (1.7630,1.8850) -- (1.7383,1.9105) -- (1.7031,1.9478) -- (1.6755,1.9774) -- (1.6451,2.0106) -- (1.6127,2.0467) -- (1.5890,2.0735) -- (1.5528,2.1153) -- (1.5289,2.1435) -- (1.5000,2.1782) -- (1.4661,2.2200) -- (1.4451,2.2464) -- (1.4168,2.2829) -- (1.3852,2.3248) -- (1.3614,2.3574) -- (1.3398,2.3876) -- (1.3109,2.4295) -- (1.2831,2.4714) -- (1.2566,2.5133) -- (1.2357,2.5479) -- (1.2147,2.5842) -- (1.1938,2.6223) -- (1.1729,2.6628) -- (1.1519,2.7064) -- (1.1310,2.7540) -- (1.1103,2.8065) -- (1.0958,2.8484) -- (1.0831,2.8903) -- (1.0681,2.9504) -- (1.0595,2.9950) -- (1.0512,3.0578) -- (1.0474,3.1206) -- (1.0472,3.1416);
\draw[gray] (-1.0472,3.1416) -- (-1.0482,3.0997) -- (-1.0535,3.0369) -- (-1.0633,2.9740) -- (-1.0722,2.9322) -- (-1.0891,2.8698) -- (-1.1028,2.8274) -- (-1.1182,2.7855) -- (-1.1353,2.7437) -- (-1.1540,2.7018) -- (-1.1743,2.6599) -- (-1.1961,2.6180) -- (-1.2193,2.5761) -- (-1.2438,2.5342) -- (-1.2697,2.4923) -- (-1.2969,2.4504) -- (-1.3195,2.4169) -- (-1.3404,2.3867) -- (-1.3698,2.3457) -- (-1.4008,2.3038) -- (-1.4242,2.2733) -- (-1.4493,2.2410) -- (-1.4829,2.1991) -- (-1.5080,2.1686) -- (-1.5350,2.1363) -- (-1.5708,2.0944) -- (-1.5917,2.0704) -- (-1.6262,2.0316) -- (-1.6546,2.0002) -- (-1.6835,1.9687) -- (-1.7174,1.9326) -- (-1.7428,1.9059) -- (-1.7802,1.8672) -- (-1.8040,1.8431) -- (-1.8431,1.8040) -- (-1.8672,1.7802) -- (-1.9059,1.7428) -- (-1.9326,1.7174) -- (-1.9687,1.6835) -- (-2.0002,1.6546) -- (-2.0316,1.6262) -- (-2.0704,1.5917) -- (-2.0944,1.5708) -- (-2.1363,1.5350) -- (-2.1686,1.5080) -- (-2.1991,1.4829) -- (-2.2410,1.4493) -- (-2.2733,1.4242) -- (-2.3038,1.4008) -- (-2.3457,1.3698) -- (-2.3867,1.3404) -- (-2.4169,1.3195) -- (-2.4504,1.2969) -- (-2.4923,1.2697) -- (-2.5342,1.2438) -- (-2.5761,1.2193) -- (-2.6180,1.1961) -- (-2.6599,1.1743) -- (-2.7018,1.1540) -- (-2.7437,1.1353) -- (-2.7855,1.1182) -- (-2.8274,1.1028) -- (-2.8698,1.0891) -- (-2.9322,1.0722) -- (-2.9740,1.0633) -- (-3.0369,1.0535) -- (-3.0997,1.0482) -- (-3.1416,1.0472);
\draw[black,semithick] (-0.0000,-3.1416) -- (0.0419,-3.0997) -- (0.0628,-3.0788) -- (0.1047,-3.0369) -- (0.1257,-3.0159) -- (0.1676,-2.9740) -- (0.1885,-2.9531) -- (0.2304,-2.9112) -- (0.2513,-2.8903) -- (0.2932,-2.8484) -- (0.3142,-2.8274) -- (0.3560,-2.7855) -- (0.3770,-2.7646) -- (0.4189,-2.7227) -- (0.4398,-2.7018) -- (0.4817,-2.6599) -- (0.5027,-2.6389) -- (0.5445,-2.5970) -- (0.5655,-2.5761) -- (0.6074,-2.5342) -- (0.6283,-2.5133) -- (0.6702,-2.4714) -- (0.6912,-2.4504) -- (0.7330,-2.4086) -- (0.7540,-2.3876) -- (0.7959,-2.3457) -- (0.8168,-2.3248) -- (0.8587,-2.2829) -- (0.8796,-2.2619) -- (0.9215,-2.2201) -- (0.9425,-2.1991) -- (0.9844,-2.1572) -- (1.0053,-2.1363) -- (1.0472,-2.0944) -- (1.0681,-2.0735) -- (1.1100,-2.0316) -- (1.1310,-2.0106) -- (1.1729,-1.9687) -- (1.1938,-1.9478) -- (1.2357,-1.9059) -- (1.2566,-1.8850) -- (1.2985,-1.8431) -- (1.3195,-1.8221) -- (1.3614,-1.7802) -- (1.3823,-1.7593) -- (1.4242,-1.7174) -- (1.4451,-1.6965) -- (1.4870,-1.6546) -- (1.5080,-1.6336) -- (1.5499,-1.5917) -- (1.5708,-1.5708) -- (1.6127,-1.5289) -- (1.6336,-1.5080) -- (1.6755,-1.4661) -- (1.6965,-1.4451) -- (1.7383,-1.4032) -- (1.7593,-1.3823) -- (1.8012,-1.3404) -- (1.8221,-1.3195) -- (1.8640,-1.2776) -- (1.8850,-1.2566) -- (1.9268,-1.2147) -- (1.9478,-1.1938) -- (1.9897,-1.1519) -- (2.0106,-1.1310) -- (2.0525,-1.0891) -- (2.0735,-1.0681) -- (2.1153,-1.0263) -- (2.1363,-1.0053) -- (2.1782,-0.9634) -- (2.1991,-0.9425) -- (2.2410,-0.9006) -- (2.2619,-0.8796) -- (2.3038,-0.8378) -- (2.3248,-0.8168) -- (2.3667,-0.7749) -- (2.3876,-0.7540) -- (2.4295,-0.7121) -- (2.4504,-0.6912) -- (2.4923,-0.6493) -- (2.5133,-0.6283) -- (2.5552,-0.5864) -- (2.5761,-0.5655) -- (2.6180,-0.5236) -- (2.6389,-0.5027) -- (2.6808,-0.4608) -- (2.7018,-0.4398) -- (2.7437,-0.3979) -- (2.7646,-0.3770) -- (2.8065,-0.3351) -- (2.8274,-0.3142) -- (2.8693,-0.2723) -- (2.8903,-0.2513) -- (2.9322,-0.2094) -- (2.9531,-0.1885) -- (2.9950,-0.1466) -- (3.0159,-0.1257) -- (3.0578,-0.0838) -- (3.0788,-0.0628) -- (3.1206,-0.0209) -- (3.1206,0.0209) -- (3.0997,0.0419) -- (3.0578,0.0838) -- (3.0369,0.1047) -- (2.9950,0.1466) -- (2.9740,0.1676) -- (2.9322,0.2094) -- (2.9112,0.2304) -- (2.8693,0.2723) -- (2.8484,0.2932) -- (2.8065,0.3351) -- (2.7855,0.3560) -- (2.7437,0.3979) -- (2.7227,0.4189) -- (2.6808,0.4608) -- (2.6599,0.4817) -- (2.6180,0.5236) -- (2.5970,0.5445) -- (2.5552,0.5864) -- (2.5342,0.6074) -- (2.4923,0.6493) -- (2.4714,0.6702) -- (2.4295,0.7121) -- (2.4086,0.7330) -- (2.3667,0.7749) -- (2.3457,0.7959) -- (2.3038,0.8378) -- (2.2829,0.8587) -- (2.2410,0.9006) -- (2.2201,0.9215) -- (2.1782,0.9634) -- (2.1572,0.9844) -- (2.1153,1.0263) -- (2.0944,1.0472) -- (2.0525,1.0891) -- (2.0316,1.1100) -- (1.9897,1.1519) -- (1.9687,1.1729) -- (1.9268,1.2147) -- (1.9059,1.2357) -- (1.8640,1.2776) -- (1.8431,1.2985) -- (1.8012,1.3404) -- (1.7802,1.3614) -- (1.7383,1.4032) -- (1.7174,1.4242) -- (1.6755,1.4661) -- (1.6546,1.4870) -- (1.6127,1.5289) -- (1.5917,1.5499) -- (1.5499,1.5917) -- (1.5289,1.6127) -- (1.4870,1.6546) -- (1.4661,1.6755) -- (1.4242,1.7174) -- (1.4032,1.7383) -- (1.3614,1.7802) -- (1.3404,1.8012) -- (1.2985,1.8431) -- (1.2776,1.8640) -- (1.2357,1.9059) -- (1.2147,1.9268) -- (1.1729,1.9687) -- (1.1519,1.9897) -- (1.1100,2.0316) -- (1.0891,2.0525) -- (1.0472,2.0944) -- (1.0263,2.1153) -- (0.9844,2.1572) -- (0.9634,2.1782) -- (0.9215,2.2201) -- (0.9006,2.2410) -- (0.8587,2.2829) -- (0.8378,2.3038) -- (0.7959,2.3457) -- (0.7749,2.3667) -- (0.7330,2.4086) -- (0.7121,2.4295) -- (0.6702,2.4714) -- (0.6493,2.4923) -- (0.6074,2.5342) -- (0.5864,2.5552) -- (0.5445,2.5970) -- (0.5236,2.6180) -- (0.4817,2.6599) -- (0.4608,2.6808) -- (0.4189,2.7227) -- (0.3979,2.7437) -- (0.3560,2.7855) -- (0.3351,2.8065) -- (0.2932,2.8484) -- (0.2723,2.8693) -- (0.2304,2.9112) -- (0.2094,2.9322) -- (0.1676,2.9740) -- (0.1466,2.9950) -- (0.1047,3.0369) -- (0.0838,3.0578) -- (0.0419,3.0997) -- (0.0209,3.1206) -- (-0.0209,3.1206) -- (-0.0628,3.0788) -- (-0.0838,3.0578) -- (-0.1257,3.0159) -- (-0.1466,2.9950) -- (-0.1885,2.9531) -- (-0.2094,2.9322) -- (-0.2513,2.8903) -- (-0.2723,2.8693) -- (-0.3142,2.8274) -- (-0.3351,2.8065) -- (-0.3770,2.7646) -- (-0.3979,2.7437) -- (-0.4398,2.7018) -- (-0.4608,2.6808) -- (-0.5027,2.6389) -- (-0.5236,2.6180) -- (-0.5655,2.5761) -- (-0.5864,2.5552) -- (-0.6283,2.5133) -- (-0.6493,2.4923) -- (-0.6912,2.4504) -- (-0.7121,2.4295) -- (-0.7540,2.3876) -- (-0.7749,2.3667) -- (-0.8168,2.3248) -- (-0.8378,2.3038) -- (-0.8796,2.2619) -- (-0.9006,2.2410) -- (-0.9425,2.1991) -- (-0.9634,2.1782) -- (-1.0053,2.1363) -- (-1.0263,2.1153) -- (-1.0681,2.0735) -- (-1.0891,2.0525) -- (-1.1310,2.0106) -- (-1.1519,1.9897) -- (-1.1938,1.9478) -- (-1.2147,1.9268) -- (-1.2566,1.8850) -- (-1.2776,1.8640) -- (-1.3195,1.8221) -- (-1.3404,1.8012) -- (-1.3823,1.7593) -- (-1.4032,1.7383) -- (-1.4451,1.6965) -- (-1.4661,1.6755) -- (-1.5080,1.6336) -- (-1.5289,1.6127) -- (-1.5708,1.5708) -- (-1.5917,1.5499) -- (-1.6336,1.5080) -- (-1.6546,1.4870) -- (-1.6965,1.4451) -- (-1.7174,1.4242) -- (-1.7593,1.3823) -- (-1.7802,1.3614) -- (-1.8221,1.3195) -- (-1.8431,1.2985) -- (-1.8850,1.2566) -- (-1.9059,1.2357) -- (-1.9478,1.1938) -- (-1.9687,1.1729) -- (-2.0106,1.1310) -- (-2.0316,1.1100) -- (-2.0735,1.0681) -- (-2.0944,1.0472) -- (-2.1363,1.0053) -- (-2.1572,0.9844) -- (-2.1991,0.9425) -- (-2.2201,0.9215) -- (-2.2619,0.8796) -- (-2.2829,0.8587) -- (-2.3248,0.8168) -- (-2.3457,0.7959) -- (-2.3876,0.7540) -- (-2.4086,0.7330) -- (-2.4504,0.6912) -- (-2.4714,0.6702) -- (-2.5133,0.6283) -- (-2.5342,0.6074) -- (-2.5761,0.5655) -- (-2.5970,0.5445) -- (-2.6389,0.5027) -- (-2.6599,0.4817) -- (-2.7018,0.4398) -- (-2.7227,0.4189) -- (-2.7646,0.3770) -- (-2.7855,0.3560) -- (-2.8274,0.3142) -- (-2.8484,0.2932) -- (-2.8903,0.2513) -- (-2.9112,0.2304) -- (-2.9531,0.1885) -- (-2.9740,0.1676) -- (-3.0159,0.1257) -- (-3.0369,0.1047) -- (-3.0788,0.0628) -- (-3.0997,0.0419) -- (-3.1416,-0.0000) -- (-3.0997,-0.0419) -- (-3.0788,-0.0628) -- (-3.0369,-0.1047) -- (-3.0159,-0.1257) -- (-2.9740,-0.1676) -- (-2.9531,-0.1885) -- (-2.9112,-0.2304) -- (-2.8903,-0.2513) -- (-2.8484,-0.2932) -- (-2.8274,-0.3142) -- (-2.7855,-0.3560) -- (-2.7646,-0.3770) -- (-2.7227,-0.4189) -- (-2.7018,-0.4398) -- (-2.6599,-0.4817) -- (-2.6389,-0.5027) -- (-2.5970,-0.5445) -- (-2.5761,-0.5655) -- (-2.5342,-0.6074) -- (-2.5133,-0.6283) -- (-2.4714,-0.6702) -- (-2.4504,-0.6912) -- (-2.4086,-0.7330) -- (-2.3876,-0.7540) -- (-2.3457,-0.7959) -- (-2.3248,-0.8168) -- (-2.2829,-0.8587) -- (-2.2619,-0.8796) -- (-2.2201,-0.9215) -- (-2.1991,-0.9425) -- (-2.1572,-0.9844) -- (-2.1363,-1.0053) -- (-2.0944,-1.0472) -- (-2.0735,-1.0681) -- (-2.0316,-1.1100) -- (-2.0106,-1.1310) -- (-1.9687,-1.1729) -- (-1.9478,-1.1938) -- (-1.9059,-1.2357) -- (-1.8850,-1.2566) -- (-1.8431,-1.2985) -- (-1.8221,-1.3195) -- (-1.7802,-1.3614) -- (-1.7593,-1.3823) -- (-1.7174,-1.4242) -- (-1.6965,-1.4451) -- (-1.6546,-1.4870) -- (-1.6336,-1.5080) -- (-1.5917,-1.5499) -- (-1.5708,-1.5708) -- (-1.5289,-1.6127) -- (-1.5080,-1.6336) -- (-1.4661,-1.6755) -- (-1.4451,-1.6965) -- (-1.4032,-1.7383) -- (-1.3823,-1.7593) -- (-1.3404,-1.8012) -- (-1.3195,-1.8221) -- (-1.2776,-1.8640) -- (-1.2566,-1.8850) -- (-1.2147,-1.9268) -- (-1.1938,-1.9478) -- (-1.1519,-1.9897) -- (-1.1310,-2.0106) -- (-1.0891,-2.0525) -- (-1.0681,-2.0735) -- (-1.0263,-2.1153) -- (-1.0053,-2.1363) -- (-0.9634,-2.1782) -- (-0.9425,-2.1991) -- (-0.9006,-2.2410) -- (-0.8796,-2.2619) -- (-0.8378,-2.3038) -- (-0.8168,-2.3248) -- (-0.7749,-2.3667) -- (-0.7540,-2.3876) -- (-0.7121,-2.4295) -- (-0.6912,-2.4504) -- (-0.6493,-2.4923) -- (-0.6283,-2.5133) -- (-0.5864,-2.5552) -- (-0.5655,-2.5761) -- (-0.5236,-2.6180) -- (-0.5027,-2.6389) -- (-0.4608,-2.6808) -- (-0.4398,-2.7018) -- (-0.3979,-2.7437) -- (-0.3770,-2.7646) -- (-0.3351,-2.8065) -- (-0.3142,-2.8274) -- (-0.2723,-2.8693) -- (-0.2513,-2.8903) -- (-0.2094,-2.9322) -- (-0.1885,-2.9531) -- (-0.1466,-2.9950) -- (-0.1257,-3.0159) -- (-0.0838,-3.0578) -- (-0.0628,-3.0788) -- (-0.0209,-3.1206) -- (-0.0000,-3.1416);
\draw[gray] (-0.1885,-2.0741) -- (-0.1257,-2.0853) -- (-0.0628,-2.0921) -- (-0.0000,-2.0944) -- (0.0628,-2.0921) -- (0.1257,-2.0853) -- (0.1885,-2.0741) -- (0.2304,-2.0642) -- (0.2723,-2.0524) -- (0.3343,-2.0316) -- (0.3770,-2.0151) -- (0.4189,-1.9972) -- (0.4608,-1.9777) -- (0.5027,-1.9566) -- (0.5445,-1.9341) -- (0.5864,-1.9102) -- (0.6283,-1.8850) -- (0.6615,-1.8640) -- (0.6936,-1.8431) -- (0.7330,-1.8164) -- (0.7749,-1.7870) -- (0.8129,-1.7593) -- (0.8409,-1.7383) -- (0.8796,-1.7087) -- (0.9215,-1.6756) -- (0.9475,-1.6546) -- (0.9844,-1.6242) -- (1.0228,-1.5917) -- (1.0472,-1.5708) -- (1.0891,-1.5341) -- (1.1183,-1.5080) -- (1.1519,-1.4774) -- (1.1867,-1.4451) -- (1.2147,-1.4188) -- (1.2528,-1.3823) -- (1.2776,-1.3582) -- (1.3168,-1.3195) -- (1.3404,-1.2957) -- (1.3786,-1.2566) -- (1.4032,-1.2311) -- (1.4385,-1.1938) -- (1.4661,-1.1642) -- (1.4965,-1.1310) -- (1.5289,-1.0949) -- (1.5525,-1.0681) -- (1.5888,-1.0263) -- (1.6127,-0.9981) -- (1.6416,-0.9634) -- (1.6755,-0.9216) -- (1.6965,-0.8952) -- (1.7248,-0.8587) -- (1.7564,-0.8168) -- (1.7802,-0.7842) -- (1.8018,-0.7540) -- (1.8307,-0.7121) -- (1.8585,-0.6702) -- (1.8850,-0.6283) -- (1.9059,-0.5936) -- (1.9268,-0.5574) -- (1.9478,-0.5193) -- (1.9687,-0.4788) -- (1.9897,-0.4352) -- (2.0106,-0.3876) -- (2.0313,-0.3351) -- (2.0458,-0.2932) -- (2.0585,-0.2513) -- (2.0735,-0.1912) -- (2.0821,-0.1466) -- (2.0904,-0.0838) -- (2.0941,-0.0209) -- (2.0944,0.0000) -- (2.0921,0.0628) -- (2.0853,0.1257) -- (2.0741,0.1885) -- (2.0642,0.2304) -- (2.0524,0.2723) -- (2.0316,0.3343) -- (2.0151,0.3770) -- (1.9972,0.4189) -- (1.9777,0.4608) -- (1.9566,0.5027) -- (1.9341,0.5445) -- (1.9102,0.5864) -- (1.8850,0.6283) -- (1.8640,0.6615) -- (1.8431,0.6936) -- (1.8164,0.7330) -- (1.7870,0.7749) -- (1.7593,0.8129) -- (1.7383,0.8409) -- (1.7087,0.8796) -- (1.6756,0.9215) -- (1.6546,0.9475) -- (1.6242,0.9844) -- (1.5917,1.0228) -- (1.5708,1.0472) -- (1.5341,1.0891) -- (1.5080,1.1183) -- (1.4774,1.1519) -- (1.4451,1.1867) -- (1.4188,1.2147) -- (1.3823,1.2528) -- (1.3582,1.2776) -- (1.3195,1.3168) -- (1.2957,1.3404) -- (1.2566,1.3786) -- (1.2311,1.4032) -- (1.1938,1.4385) -- (1.1642,1.4661) -- (1.1310,1.4965) -- (1.0949,1.5289) -- (1.0681,1.5525) -- (1.0263,1.5888) -- (0.9981,1.6127) -- (0.9634,1.6416) -- (0.9216,1.6755) -- (0.8952,1.6965) -- (0.8587,1.7248) -- (0.8168,1.7564) -- (0.7842,1.7802) -- (0.7540,1.8018) -- (0.7121,1.8307) -- (0.6702,1.8585) -- (0.6283,1.8850) -- (0.5936,1.9059) -- (0.5574,1.9268) -- (0.5193,1.9478) -- (0.4788,1.9687) -- (0.4352,1.9897) -- (0.3876,2.0106) -- (0.3351,2.0313) -- (0.2932,2.0458) -- (0.2513,2.0585) -- (0.1912,2.0735) -- (0.1466,2.0821) -- (0.0838,2.0904) -- (0.0209,2.0941) -- (-0.0000,2.0944) -- (-0.0628,2.0921) -- (-0.1257,2.0853) -- (-0.1885,2.0741) -- (-0.2304,2.0642) -- (-0.2723,2.0524) -- (-0.3343,2.0316) -- (-0.3770,2.0151) -- (-0.4189,1.9972) -- (-0.4608,1.9777) -- (-0.5027,1.9566) -- (-0.5445,1.9341) -- (-0.5864,1.9102) -- (-0.6283,1.8850) -- (-0.6615,1.8640) -- (-0.6936,1.8431) -- (-0.7330,1.8164) -- (-0.7749,1.7870) -- (-0.8129,1.7593) -- (-0.8409,1.7383) -- (-0.8796,1.7087) -- (-0.9215,1.6756) -- (-0.9475,1.6546) -- (-0.9844,1.6242) -- (-1.0228,1.5917) -- (-1.0472,1.5708) -- (-1.0891,1.5341) -- (-1.1183,1.5080) -- (-1.1519,1.4774) -- (-1.1867,1.4451) -- (-1.2147,1.4188) -- (-1.2528,1.3823) -- (-1.2776,1.3582) -- (-1.3168,1.3195) -- (-1.3404,1.2957) -- (-1.3786,1.2566) -- (-1.4032,1.2311) -- (-1.4385,1.1938) -- (-1.4661,1.1642) -- (-1.4965,1.1310) -- (-1.5289,1.0949) -- (-1.5525,1.0681) -- (-1.5888,1.0263) -- (-1.6127,0.9981) -- (-1.6416,0.9634) -- (-1.6755,0.9216) -- (-1.6965,0.8952) -- (-1.7248,0.8587) -- (-1.7564,0.8168) -- (-1.7802,0.7842) -- (-1.8018,0.7540) -- (-1.8307,0.7121) -- (-1.8585,0.6702) -- (-1.8850,0.6283) -- (-1.9059,0.5936) -- (-1.9268,0.5574) -- (-1.9478,0.5193) -- (-1.9687,0.4788) -- (-1.9897,0.4352) -- (-2.0106,0.3876) -- (-2.0313,0.3351) -- (-2.0458,0.2932) -- (-2.0585,0.2513) -- (-2.0735,0.1912) -- (-2.0821,0.1466) -- (-2.0904,0.0838) -- (-2.0941,0.0209) -- (-2.0934,-0.0419) -- (-2.0881,-0.1047) -- (-2.0783,-0.1676) -- (-2.0694,-0.2094) -- (-2.0525,-0.2718) -- (-2.0388,-0.3142) -- (-2.0234,-0.3560) -- (-2.0063,-0.3979) -- (-1.9876,-0.4398) -- (-1.9673,-0.4817) -- (-1.9455,-0.5236) -- (-1.9223,-0.5655) -- (-1.8978,-0.6074) -- (-1.8719,-0.6493) -- (-1.8447,-0.6912) -- (-1.8221,-0.7247) -- (-1.8012,-0.7549) -- (-1.7718,-0.7959) -- (-1.7407,-0.8378) -- (-1.7174,-0.8683) -- (-1.6922,-0.9006) -- (-1.6587,-0.9425) -- (-1.6336,-0.9730) -- (-1.6066,-1.0053) -- (-1.5708,-1.0472) -- (-1.5499,-1.0712) -- (-1.5154,-1.1100) -- (-1.4870,-1.1414) -- (-1.4581,-1.1729) -- (-1.4242,-1.2090) -- (-1.3988,-1.2357) -- (-1.3614,-1.2744) -- (-1.3376,-1.2985) -- (-1.2985,-1.3376) -- (-1.2744,-1.3614) -- (-1.2357,-1.3988) -- (-1.2090,-1.4242) -- (-1.1729,-1.4581) -- (-1.1414,-1.4870) -- (-1.1100,-1.5154) -- (-1.0712,-1.5499) -- (-1.0472,-1.5708) -- (-1.0053,-1.6066) -- (-0.9730,-1.6336) -- (-0.9425,-1.6587) -- (-0.9006,-1.6922) -- (-0.8683,-1.7174) -- (-0.8378,-1.7407) -- (-0.7959,-1.7718) -- (-0.7549,-1.8012) -- (-0.7247,-1.8221) -- (-0.6912,-1.8447) -- (-0.6493,-1.8719) -- (-0.6074,-1.8978) -- (-0.5655,-1.9223) -- (-0.5236,-1.9455) -- (-0.4817,-1.9673) -- (-0.4398,-1.9876) -- (-0.3979,-2.0063) -- (-0.3560,-2.0234) -- (-0.3142,-2.0388) -- (-0.2718,-2.0525) -- (-0.2094,-2.0694) -- (-0.1885,-2.0741);
\draw[gray] (-0.1885,-1.5531) -- (-0.1257,-1.5629) -- (-0.0628,-1.5688) -- (-0.0000,-1.5708) -- (0.0628,-1.5688) -- (0.1257,-1.5629) -- (0.1885,-1.5531) -- (0.2304,-1.5444) -- (0.2903,-1.5289) -- (0.3351,-1.5151) -- (0.3770,-1.5005) -- (0.4189,-1.4842) -- (0.4613,-1.4661) -- (0.5060,-1.4451) -- (0.5473,-1.4242) -- (0.5864,-1.4029) -- (0.6283,-1.3786) -- (0.6702,-1.3528) -- (0.7121,-1.3253) -- (0.7509,-1.2985) -- (0.7799,-1.2776) -- (0.8168,-1.2498) -- (0.8587,-1.2169) -- (0.8868,-1.1938) -- (0.9215,-1.1643) -- (0.9591,-1.1310) -- (0.9844,-1.1078) -- (1.0260,-1.0681) -- (1.0472,-1.0472) -- (1.0881,-1.0053) -- (1.1100,-0.9820) -- (1.1459,-0.9425) -- (1.1729,-0.9116) -- (1.1997,-0.8796) -- (1.2336,-0.8378) -- (1.2566,-0.8079) -- (1.2812,-0.7749) -- (1.3110,-0.7330) -- (1.3392,-0.6912) -- (1.3614,-0.6565) -- (1.3823,-0.6221) -- (1.4032,-0.5859) -- (1.4256,-0.5445) -- (1.4468,-0.5027) -- (1.4663,-0.4608) -- (1.4870,-0.4119) -- (1.5080,-0.3562) -- (1.5218,-0.3142) -- (1.5339,-0.2723) -- (1.5489,-0.2094) -- (1.5568,-0.1676) -- (1.5653,-0.1047) -- (1.5699,-0.0419) -- (1.5706,0.0209) -- (1.5673,0.0838) -- (1.5601,0.1466) -- (1.5499,0.2048) -- (1.5394,0.2513) -- (1.5281,0.2932) -- (1.5080,0.3560) -- (1.4926,0.3979) -- (1.4755,0.4398) -- (1.4567,0.4817) -- (1.4364,0.5236) -- (1.4145,0.5655) -- (1.3910,0.6074) -- (1.3659,0.6493) -- (1.3404,0.6893) -- (1.3195,0.7207) -- (1.2963,0.7540) -- (1.2657,0.7959) -- (1.2357,0.8350) -- (1.2147,0.8613) -- (1.1822,0.9006) -- (1.1519,0.9357) -- (1.1271,0.9634) -- (1.0891,1.0043) -- (1.0679,1.0263) -- (1.0263,1.0679) -- (1.0043,1.0891) -- (0.9634,1.1271) -- (0.9357,1.1519) -- (0.9006,1.1822) -- (0.8613,1.2147) -- (0.8350,1.2357) -- (0.7959,1.2657) -- (0.7540,1.2963) -- (0.7207,1.3195) -- (0.6893,1.3404) -- (0.6493,1.3659) -- (0.6074,1.3910) -- (0.5655,1.4145) -- (0.5236,1.4364) -- (0.4817,1.4567) -- (0.4398,1.4755) -- (0.3979,1.4926) -- (0.3560,1.5080) -- (0.2932,1.5281) -- (0.2513,1.5394) -- (0.2048,1.5499) -- (0.1466,1.5601) -- (0.0838,1.5673) -- (0.0209,1.5706) -- (-0.0419,1.5699) -- (-0.1047,1.5653) -- (-0.1676,1.5568) -- (-0.2094,1.5489) -- (-0.2723,1.5339) -- (-0.3142,1.5218) -- (-0.3562,1.5080) -- (-0.4119,1.4870) -- (-0.4608,1.4663) -- (-0.5027,1.4468) -- (-0.5445,1.4256) -- (-0.5859,1.4032) -- (-0.6221,1.3823) -- (-0.6565,1.3614) -- (-0.6912,1.3392) -- (-0.7330,1.3110) -- (-0.7749,1.2812) -- (-0.8079,1.2566) -- (-0.8378,1.2336) -- (-0.8796,1.1997) -- (-0.9116,1.1729) -- (-0.9425,1.1459) -- (-0.9820,1.1100) -- (-1.0053,1.0881) -- (-1.0472,1.0472) -- (-1.0681,1.0260) -- (-1.1078,0.9844) -- (-1.1310,0.9591) -- (-1.1643,0.9215) -- (-1.1938,0.8868) -- (-1.2169,0.8587) -- (-1.2498,0.8168) -- (-1.2776,0.7799) -- (-1.2985,0.7509) -- (-1.3253,0.7121) -- (-1.3528,0.6702) -- (-1.3786,0.6283) -- (-1.4029,0.5864) -- (-1.4242,0.5473) -- (-1.4451,0.5060) -- (-1.4661,0.4613) -- (-1.4842,0.4189) -- (-1.5005,0.3770) -- (-1.5151,0.3351) -- (-1.5289,0.2903) -- (-1.5444,0.2304) -- (-1.5531,0.1885) -- (-1.5629,0.1257) -- (-1.5688,0.0628) -- (-1.5708,-0.0000) -- (-1.5688,-0.0628) -- (-1.5629,-0.1257) -- (-1.5531,-0.1885) -- (-1.5444,-0.2304) -- (-1.5289,-0.2903) -- (-1.5151,-0.3351) -- (-1.5005,-0.3770) -- (-1.4842,-0.4189) -- (-1.4661,-0.4613) -- (-1.4451,-0.5060) -- (-1.4242,-0.5473) -- (-1.4029,-0.5864) -- (-1.3786,-0.6283) -- (-1.3528,-0.6702) -- (-1.3253,-0.7121) -- (-1.2985,-0.7509) -- (-1.2776,-0.7799) -- (-1.2498,-0.8168) -- (-1.2169,-0.8587) -- (-1.1938,-0.8868) -- (-1.1643,-0.9215) -- (-1.1310,-0.9591) -- (-1.1078,-0.9844) -- (-1.0681,-1.0260) -- (-1.0472,-1.0472) -- (-1.0053,-1.0881) -- (-0.9820,-1.1100) -- (-0.9425,-1.1459) -- (-0.9116,-1.1729) -- (-0.8796,-1.1997) -- (-0.8378,-1.2336) -- (-0.8079,-1.2566) -- (-0.7749,-1.2812) -- (-0.7330,-1.3110) -- (-0.6912,-1.3392) -- (-0.6565,-1.3614) -- (-0.6221,-1.3823) -- (-0.5859,-1.4032) -- (-0.5445,-1.4256) -- (-0.5027,-1.4468) -- (-0.4608,-1.4663) -- (-0.4119,-1.4870) -- (-0.3562,-1.5080) -- (-0.3142,-1.5218) -- (-0.2723,-1.5339) -- (-0.2094,-1.5489) -- (-0.1885,-1.5531);
\draw[gray] (-0.1885,-1.0266) -- (-0.1257,-1.0380) -- (-0.0628,-1.0449) -- (-0.0000,-1.0472) -- (0.0628,-1.0449) -- (0.1257,-1.0380) -- (0.1885,-1.0266) -- (0.2304,-1.0164) -- (0.2723,-1.0041) -- (0.3281,-0.9844) -- (0.3770,-0.9640) -- (0.4189,-0.9441) -- (0.4608,-0.9219) -- (0.4973,-0.9006) -- (0.5304,-0.8796) -- (0.5655,-0.8557) -- (0.6074,-0.8246) -- (0.6429,-0.7959) -- (0.6702,-0.7723) -- (0.7121,-0.7332) -- (0.7332,-0.7121) -- (0.7723,-0.6702) -- (0.7959,-0.6429) -- (0.8246,-0.6074) -- (0.8557,-0.5655) -- (0.8796,-0.5304) -- (0.9006,-0.4973) -- (0.9219,-0.4608) -- (0.9441,-0.4189) -- (0.9640,-0.3770) -- (0.9844,-0.3281) -- (1.0041,-0.2723) -- (1.0164,-0.2304) -- (1.0266,-0.1885) -- (1.0380,-0.1257) -- (1.0449,-0.0628) -- (1.0472,-0.0000) -- (1.0469,0.0209) -- (1.0431,0.0838) -- (1.0347,0.1466) -- (1.0263,0.1901) -- (1.0105,0.2513) -- (0.9971,0.2932) -- (0.9817,0.3351) -- (0.9634,0.3783) -- (0.9425,0.4221) -- (0.9215,0.4615) -- (0.8973,0.5027) -- (0.8702,0.5445) -- (0.8405,0.5864) -- (0.8168,0.6172) -- (0.7905,0.6493) -- (0.7540,0.6903) -- (0.7330,0.7123) -- (0.6912,0.7532) -- (0.6672,0.7749) -- (0.6283,0.8079) -- (0.5901,0.8378) -- (0.5612,0.8587) -- (0.5236,0.8841) -- (0.4817,0.9099) -- (0.4398,0.9333) -- (0.3979,0.9543) -- (0.3560,0.9731) -- (0.3142,0.9897) -- (0.2684,1.0053) -- (0.2094,1.0218) -- (0.1676,1.0309) -- (0.1047,1.0408) -- (0.0419,1.0462) -- (-0.0000,1.0472) -- (-0.0419,1.0462) -- (-0.1047,1.0408) -- (-0.1676,1.0309) -- (-0.2094,1.0218) -- (-0.2684,1.0053) -- (-0.3142,0.9897) -- (-0.3560,0.9731) -- (-0.3979,0.9543) -- (-0.4398,0.9333) -- (-0.4817,0.9099) -- (-0.5236,0.8841) -- (-0.5612,0.8587) -- (-0.5901,0.8378) -- (-0.6283,0.8079) -- (-0.6672,0.7749) -- (-0.6912,0.7532) -- (-0.7330,0.7123) -- (-0.7540,0.6903) -- (-0.7905,0.6493) -- (-0.8168,0.6172) -- (-0.8405,0.5864) -- (-0.8702,0.5445) -- (-0.8973,0.5027) -- (-0.9215,0.4615) -- (-0.9425,0.4221) -- (-0.9634,0.3783) -- (-0.9817,0.3351) -- (-0.9971,0.2932) -- (-1.0105,0.2513) -- (-1.0263,0.1901) -- (-1.0347,0.1466) -- (-1.0431,0.0838) -- (-1.0469,0.0209) -- (-1.0462,-0.0419) -- (-1.0408,-0.1047) -- (-1.0309,-0.1676) -- (-1.0218,-0.2094) -- (-1.0053,-0.2684) -- (-0.9897,-0.3142) -- (-0.9731,-0.3560) -- (-0.9543,-0.3979) -- (-0.9333,-0.4398) -- (-0.9099,-0.4817) -- (-0.8841,-0.5236) -- (-0.8587,-0.5612) -- (-0.8378,-0.5901) -- (-0.8079,-0.6283) -- (-0.7749,-0.6672) -- (-0.7532,-0.6912) -- (-0.7123,-0.7330) -- (-0.6903,-0.7540) -- (-0.6493,-0.7905) -- (-0.6172,-0.8168) -- (-0.5864,-0.8405) -- (-0.5445,-0.8702) -- (-0.5027,-0.8973) -- (-0.4615,-0.9215) -- (-0.4221,-0.9425) -- (-0.3783,-0.9634) -- (-0.3351,-0.9817) -- (-0.2932,-0.9971) -- (-0.2513,-1.0105) -- (-0.1901,-1.0263) -- (-0.1885,-1.0266);
\node[below] at (-3.14159,-3.14159) {$-\pi$};
\node[below] at (0,-3.14159) {$0$};
\node[below] at (3.14159,-3.14159) {$\pi$};
\node[left] at (-3.14159,3.14159) {$\pi$};
\node[left] at (-3.14159,0) {$0$};
\end{tikzpicture}
\caption{Level sets of $\Psi(x,y)=\cos x+\cos y$. The separatrix $\Psi=0$ is drawn in black.}
\label{fig:stream-function}
\end{figure}
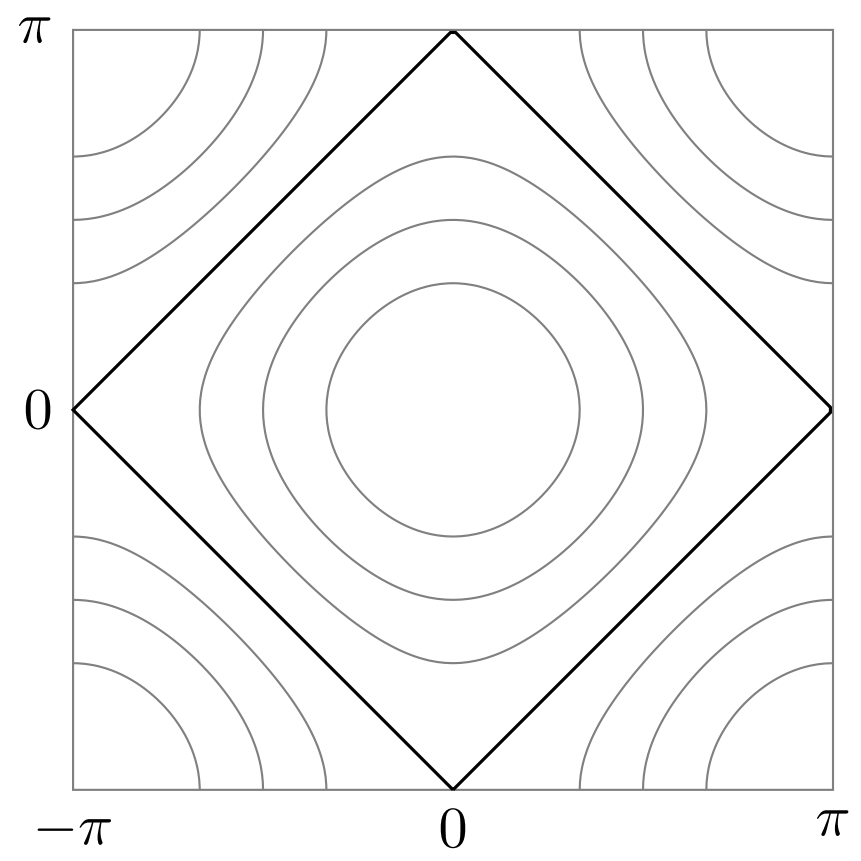

The vorticity field
\begin{equation}\label{Int1.1}
\Omega=-\cos x - \cos y
\end{equation}
is a steady solution of \eqref{Int1}, called the symmetric dipole.
Its velocity and stream function are $(\sin y,-\sin x)$ and
$\Psi(x,y)=\cos x+\cos y$, respectively; see
Figure~\ref{fig:stream-function}.

We study the Euler equation linearized around $\Omega$:
\begin{equation}\label{Int2}
    \begin{split}
        \partial_t\omega+\Big[\sin y \,\partial_x - \sin x \,\partial_y\Big](\omega+\psi)&=0,\\
        \Delta\psi &=\omega,
    \end{split}
\end{equation}
for $(x,y)\in\T^2,\,
\,t\ge0$ with the initial condition $\omega(0,x,y)=\omega_0(x,y)\in C^\infty(\T^2)$.

Our main goal in this paper is to understand the precise dynamics of \eqref{Int2}. Prior to this paper, precise decay estimates for the linearized Euler equation have only been obtained around shear flows and vortices; see Section~\ref{sec:background} for a more detailed discussion. A common feature is that in both cases the linearized flow decouples via a Fourier transform to one-dimensional problems. This reduction in dimension is essential in controlling the singularities of the resolvent when the spectrum is near the essential spectrum. As far as we know, the only results on the asymptotic stability of linearized Euler around dipole solutions are qualitative and obtained as a consequence of the RAGE theorem. In a pioneering application of the RAGE theorem, Constantin, Kiselev,
Ryzhik, and Zlato\v{s}~\cite{CKRZ} obtained time-averaged decay of
low-frequency components of solutions to an advection-diffusion equation. Lin and Xu~\cite{Zhiwu3}
adapted this spectral approach to establish qualitative inviscid
damping for linearized Euler equations. The RAGE theorem controls
suitable time averages on the continuous spectral subspace, but
does not by itself imply decay at each time or provide a quantitative
rate. Quantitative decay estimates generally require finer information
on the resolvent near the continuous spectrum.

In this paper we initiate the study of quantitative asymptotic stability of genuinely two-dimensional incompressible inviscid flows, by developing robust new methods and proving explicit decay rates of the stream function and the velocity field for the linearized flow around dipoles.

To state our main results precisely, we define the operators
\begin{equation}\label{Intm1}
L = \Big[\sin y \,\partial_x - \sin x \,\partial_y\Big](1+\Delta^{-1}), \qquad L_0=\sin y \,\partial_x - \sin x \,\partial_y.
\end{equation}
Let
\begin{equation}\label{Intm2}
    X = \Big\{g\in L^2(\T^2):\, \int_{\T^2}g\,dxdy=0\Big\},\quad Y=\overline{{\rm Ran}\,L_0}
\end{equation}
and
\begin{equation}\label{Intm3}
    \mathcal{D}=\Big\{g\in X:\, L_0g\in L^2(\T^2)\Big\}.
\end{equation}
We use $\mathcal D$ as the domain of both $L$ and $L_0$ on $X$.
The closure in the definition of $Y$ is taken in $L^2$, and
$\Delta^{-1}$ denotes the inverse on mean-zero functions. We have
\begin{equation}\label{Intm4}
    L = L_0 (1+\Delta^{-1}).
\end{equation}
Thus $Y$ is invariant under the linearized flow. For data
$\omega_{\rm ini}\in\ker L_0$, the $\ker L_0$ component remains
constant, while the component $\eta$ in $Y$ satisfies
\begin{equation}\label{sol8}
    \partial_t\eta+ L\eta +L_0\big(1+\Delta^{-1}\big)\omega_{\rm ini}=0, \quad \eta|_{t=0}\equiv0,
\end{equation}
and can be studied using the homogeneous equation and the Duhamel's principle.

We will show below that
\begin{equation}\label{Intm5}
    X=\ker L_0 \oplus Y \quad {\rm and} \quad \ker L \cap Y={\rm span}\big\{\sin x, \sin y, \cos x - \cos y\big\}.
\end{equation}
\begin{definition}\label{defPQ}
Let $P$ be the $L^2$ orthogonal projection onto
\[
 E:=\ker L\cap Y
   =\operatorname{span}\{\sin x,\sin y,\cos x-\cos y\},
 \qquad Q:=I-P.
\]
\end{definition}
Besides the smooth zero modes above, there are also three {\it adjoint} zero modes $g_1, g_2$ and $g_3$ for $L$, which are given explicitly as
\begin{equation}\label{IntzeroDual}
 g_1=-y_c,\qquad g_2=x_c,\qquad
 g_3=\log\left|\frac{\cos((x-y)/2)}{\cos((x+y)/2)}\right|.
\end{equation}
Here $(x_c,y_c)$ are Cartesian coordinates measured from the
elliptic point in each cell; see Figure~\ref{fig:stream-function}. More precisely, for $\Psi(x,y)\ne0$,
they are the unique pair satisfying $|x_c|+|y_c|<\pi$ and
\[
 (x,y)\equiv
 \begin{cases}
  (x_c,y_c),&\Psi(x,y)>0,\\
  (\pi+x_c,\pi+y_c),&\Psi(x,y)<0,
 \end{cases}
 \qquad \pmod{2\pi}.
\]
Their values on the separatrix $\{\Psi=0\}$ are not important.
Although $x_c$ and $y_c$ jump across the separatrix, $L_0$ is tangent to it, and as a result these jumps introduce no distributional boundary terms. The functions $g_j$ belong to $L^2(\T^2)$ and satisfy
\[
 L_0g_1=\sin x,\qquad L_0g_2=\sin y,\qquad
 L_0g_3=\cos x-\cos y
\]
in the sense of distributions. See also Proposition~\ref{zero:root} for the generalized zero modes of $L$.

Our main result is the following theorem.

\begin{theorem}\label{th1:Int}
Suppose that $\omega_0\in H^5(\T^2)\cap Y$ and
\begin{equation}\label{IntzeroMoments}
 \langle\omega_0,g_j\rangle_{L^2}=0,\qquad j=1,2,3.
\end{equation}
Let $\omega\in C([0,\infty),H^5(\T^2))$ be the solution to
\eqref{Int2} with initial data $\omega_0$. Then, for every $t\ge0$,
\begin{equation}\label{Intm7}
 \|\nabla\Delta^{-1}Q\omega(t)\|_{L^2(\T^2)}
 \lesssim (1+t)^{-1}\|\omega_0\|_{H^5(\T^2)},
\end{equation}
\begin{equation}\label{Intquadratic}
 \|\Delta^{-1}Q\omega(t)\|_{L^2(\T^2)}
 \lesssim (1+t)^{-2}\|\omega_0\|_{H^5(\T^2)},
\end{equation}
and
\begin{equation}\label{Intm9}
 |Q\omega(t,x,y)|\lesssim
 \left[d(x,y)^3+\frac{d(x,y)}{1+t}\right]
 \|\omega_0\|_{H^5(\T^2)},
\end{equation}
where $d$ is the torus distance to $\{(0,0),(\pi,\pi)\}$.
The implied constants are independent of $t$ and $\omega_0$.
\end{theorem}

The following result on the absence of nonzero embedded eigenvalues for $L$ provides a key spectral ingredient for the proof of the above theorem. The proof of Theorem \ref{Intm6} was found by GPT-6 Astra and then formalized in LEAN; see  \cite{NoEmbeddedEigenvalues}. 
\begin{theorem}\label{Intm6}
There is no nonzero embedded eigenvalue for $L$. That is, for every $\lambda\in\R\setminus\{0\}$,
\[
 \ker(L-i\lambda)=\{0\}.
\]
\end{theorem}
For the purpose of this paper, we work on the basis of Theorem \ref{Intm6}.

\medskip

\begin{remark} We make the following remarks on Theorem \ref{th1:Int}.

\begin{enumerate}
\item The operator $L$ has three generalized zero modes $b_j$ in addition to the zero modes in \eqref{Intm5}; see Proposition~\ref{zero:root}. $(1+\Delta^{-1})b_j=g_j$ and the orthogonality condition \eqref{IntzeroMoments} is natural projection off these generalized zero modes with respect to the conserved energy form
\begin{equation}\label{int:conserv}
    (f,g):=\langle (1+\Delta^{-1})f, g\rangle_{L^2}
\end{equation}
for \eqref{Int2}.

\item Simple calculation shows that $Q\Delta=\Delta Q$ and therefore the order in which the projection $Q$ is taken with respect to $\Delta^{-1}$ does not matter.

    \item More generally, one can study the family of ``dipoles'' given by $$A\cos(x-a)+B\cos(y-b),\quad A, a, B, b\in\R, \,\,AB\neq0,$$ which are steady states of the two-dimensional Euler equation. The four-dimensional manifold of steady states has tangent space at $\cos x+\cos y$ given by $${\rm span}\{\sin x, \,\sin y, \,\cos x-\cos y, \,\cos x + \cos y\}.$$ It is easy to check that $\cos x+\cos y\in {\rm Ker} \,L_0$ and hence is perpendicular to $Y$. Therefore $Q$ exactly projects away from the three tangent directions in $Y$.

 \item Write
\begin{equation}\label{eq:decomposition-p-eta}
 \omega(t)=p(t)+\eta(t),\qquad p=P\omega,\quad\eta=Q\omega.
\end{equation}
Since $Lp=0$, projecting the equation onto $E$ implies
\begin{equation}\label{eq:evolution-p}
 \partial_t p=-PL\eta.
\end{equation}
For $\phi\in E$, a direct calculation shows that
$\Delta L_0\phi=-2L_0\phi$. Integration by parts therefore yields
\begin{equation}\label{eq:PLeta-bound}
\begin{aligned}
 \langle L\eta,\phi\rangle_{L^2}
 &=-\langle\eta+\Delta^{-1}\eta,L_0\phi\rangle_{L^2}\\
 &=-\langle\Delta^{-1}\eta,(\Delta+1)L_0\phi\rangle_{L^2}
 =\langle\Delta^{-1}\eta,L_0\phi\rangle_{L^2}.
\end{aligned}
\end{equation}
Since $E$ is finite dimensional, it follows from
\eqref{eq:evolution-p}--\eqref{eq:PLeta-bound} and
\eqref{Intquadratic} that, for every $m\ge0$,
\[
 \|\partial_t p(t)\|_{H^m}
 \le C_m\|\Delta^{-1}\eta(t)\|_{L^2}
 \le C_m(1+t)^{-2}\|\omega_0\|_{H^5}.
\]
In particular, the integral defining
$p_\infty:=p(0)-\int_0^\infty PL\eta(s)\,ds\in E$
converges, and
\begin{equation}\label{eq:p-convergence}
\begin{aligned}
 \|p(t)-p_\infty\|_{H^m}\le C_m\|\omega_0\|_{H^5}\int_t^\infty(1+s)^{-2}\,ds
 =\frac{C_m}{1+t}\|\omega_0\|_{H^5}.
\end{aligned}
\end{equation}
The same rate holds in every fixed $C^m$ norm. Finally, since
$\Delta^{-1}(p-p_\infty)=-(p-p_\infty)$, it follows from
\eqref{Intm7} and \eqref{eq:p-convergence} that the full velocity
converges at the rate
\[
 \|\nabla\Delta^{-1}\omega(t)-\nabla\Delta^{-1}p_\infty\|_{L^2}
 \le C(1+t)^{-1}\|\omega_0\|_{H^5}.
\]
The rate in \eqref{eq:p-convergence} can be improved using the later
analysis. Indeed, it follows from \eqref{quad:stream-decay} in
Theorem~\ref{r2:global-damping} and the identity
$p-p_\infty=-P\Delta^{-1}(\omega-p_\infty)$ that
\[
 \|p(t)-p_\infty\|_{H^m}
 \le C_m\|\Delta^{-1}(\omega(t)-p_\infty)\|_{L^2}
 \le C_m(1+t)^{-2}\|\omega_0\|_{H^5}.
\]

    \item The estimates \eqref{Intquadratic} and \eqref{Intm7} imply
the same decay rates as those for monotonic shear flows: $t^{-2}$
for the projected stream function and $t^{-1}$ for the projected
velocity. The stream-function estimate follows by taking a second
spectral derivative and estimating its singularity after integration
in the spectral parameter. The vorticity depletion estimate \eqref{Intm9} is also expected to be sharp. This can be seen from the comparison with the case of monotonic vortices in~\cite{JiaVortex} as well as~\cite{Bed2}, which provides a good model for the analysis near the elliptic centers. Notice that since for the dipoles the angular modes are closely coupled it is generally not possible to obtain enhanced rate of depletion near the elliptic centers for higher angular modes. We use the same $H^5$ assumption for all
estimates and do not optimize the regularity of the initial data.

    \item In numerical simulations of slightly viscous fluids with viscosity $\nu>0$, metastable ``dipole states'' with vorticity $e^{-\nu t}(A\cos (x-a)+B\cos (y-b))$ ($AB\neq0$) appear to emerge from generic smooth initial data and remain stable for a long time. To prove this for the full nonlinear perturbed equation is a challenging problem. See~\cite{CJWZ} for a result in this direction for the Kolmogorov flow, and~\cite{Elgindi} for an orbital stability result for dipoles. We expect that the methods developed in this paper can also be applied to study similar problems for general dipoles ($A\neq B$).

    \item Lastly we note that the presence of elliptic centers and separatrices is a generic feature for smooth divergence-free vector fields in two dimensions. We therefore expect the methods introduced here to be useful for the analysis of the precise long-time behavior of the linearized flow of general steady states to the two-dimensional incompressible Euler equation.
\end{enumerate}
\end{remark}

\subsection{Background and some related results}\label{sec:background}
The equation \eqref{Int1} is globally well-posed for sufficiently smooth initial data, and smooth solutions exist for all times. By a result of Yudovich~\cite{Yudovich1}, there exists a unique solution if the initial vorticity $\omega_0\in L^\infty(\T^2)$. In comparison with regularity,
the long-time behavior of solutions to \eqref{Int1} in general is much harder to study, due to the lack of global relaxation mechanisms.

Over the last decade, there has been tremendous progress in understanding the precise dynamics of smooth solutions to \eqref{Int1} in perturbative regimes. The main stabilizing mechanisms are \emph{inviscid damping} and \emph{vorticity depletion}. Although the vorticity field is only transported and there is no decay in any $L^p$ spaces, the stream functions and velocity fields may still decay over time, leading to convergence of the vorticity field in the weak sense. This is the so-called ``inviscid damping''.  Moreover, near the elliptic centers of the vorticity field where the inviscid damping mechanism is weaker, remarkably the perturbation vorticity is dynamically depleted, which in turn recovers the full inviscid damping estimates. This is called vorticity depletion.

Inviscid damping was first studied by Orr~\cite{Orr} in 1907 near the Couette flow at the linearized level. However, the extension to the full nonlinear asymptotic stability has only been achieved relatively recently, by Bedrossian and Masmoudi in~\cite{BM}. Inspired by this remarkable work, intensive research was carried out to prove similar results for more general shear flows, see for example~\cite{Zillinger1,Zillinger2,dongyi,Dongyi2,Dongyi3,JiaG,JiaUM,JiaVortex,CWZ2} and references therein. In particular, Wei-Zhang-Zhao~\cite{dongyi} established sharp decay rates for the linearized Euler equation around general spectrally stable monotonic shear flows. The linearized equation around non-monotonic shear flows was rigorously studied by Wei-Zhang-Zhao~\cite{Dongyi2,Dongyi3}, following the numerical work by Bouchet and Morita~\cite{Bouchet}. Iyer, Ionescu and the second author then provided a more quantitative description of the vorticity depletion effects~\cite{Iyer}. In the important special case of Kolmogorov flows, the result~\cite{CJWZ} by Chen, Wei, Zhang and the second author contains sharp inviscid damping and vorticity depletion estimates. For more general non-monotonic shear, Bian and Grenier~\cite{BG26} obtained nearly sharp depletion estimates with a log loss; finally, Beekie, Chen and the second author established the sharp depletion estimates in~\cite{BCJ26}.  More recently, uniform-in-viscosity inviscid damping and vorticity depletion estimates for the linearized Navier--Stokes equation around shear flows with small viscosity have also been studied, see~\cite{CWZ2,JiaUM} for the case of monotonic shear flows, and~\cite{BeekieChenJia2026} for non-monotonic shear flows.

Besides shear flows, another important class of steady states for the two-dimensional Euler equation is the (radially monotonic) vortices. Quantitative inviscid damping and vorticity depletion estimates were first calculated in the physics literature, by Bassom and Gilbert~\cite{Bassom} using a matched asymptotics method. Bedrossian, Coti Zelati and Vicol~\cite{Bed2} obtained the first rigorous precise bounds on the linearized Euler equation around general radially decreasing vortices. Ionescu and the second author~\cite{JiaVortex} established optimal Gevrey estimates for the profiles of the vorticity for the linearized equation, which are expected to be important for nonlinear applications.

Extending the linearized stability estimates to full nonlinear asymptotic stability is a difficult problem. For spectrally stable monotonic shear flows in a periodic channel, nonlinear inviscid damping was proved independently, by Ionescu and the second author~\cite{IJacta}, and Masmoudi and Zhao~\cite{NaZh}. The important case of vortices, however, remains an important open problem.

\subsection{New difficulties and ingredients}
Our proof follows the spectral density function approach developed
in~\cite{JiaSIMA,JiaG,JiaVortex,Iyer}. We first explain the new difficulties that arise for the dipole and the key ingredients used to overcome them.

For shear flows and vortices, Fourier decomposition reduces the linearized equation to independent one-dimensional problems.
For the dipole, an important first step is to use streamline coordinates (see Subsection~\ref{sec:coordinates}) to diagonalize the transport operator $L_0$. However, the inverse Laplacian in $L=L_0(1+\Delta^{-1})$ still couples the angular modes.
Moreover, these coordinates degenerate both at the elliptic points and at the separatrix. The resolvent must therefore be estimated as a two-dimensional elliptic equation with a singular nonlocal term. This is the fundamental difference between the dipole problem and the cases of shear flows and vortices.

Near each elliptic point, we can show that the metric in the streamline coordinates is a quadratic perturbation of the Euclidean metric. This allows us to compare with an explicit radial operator while keeping the coupling in the elliptic energy form. By now it is well known that near the elliptic points the nonlocal term in the Rayleigh equation  \eqref{rayleigh1} cannot be treated as a compact perturbation; see~\cite{JiaVortex}. Instead, its limiting inverse-square potential is critically singular and favorable. As a result, the radial exponent increases from $|k|$ to $\sqrt{k^2+8}$ leading to enhanced vanishing of the vorticity amplitude that compensates for the weaker mixing near the critical streamline.

The presence of the separatrix associated with hyperbolic critical points is a new feature for the dipole, in comparison with shear flows and vortices. Near the separatrix, the difficulty is the accumulation of infinitely many angular resonances, in the sense that there is a spatial accumulation of critical layers for different angular modes at the same spectral parameter $\lambda\neq0$. Crucially, the physical area of streamlines with large period is exponentially small. Consequently, in the period variable, this overpowers the polynomial losses in the resolvent and its derivatives.


The analysis near the spectral parameter $\lambda=0$ requires
a separate treatment of the stationary and generalized zero
modes. The projection $Q$ removes the stationary modes,
while the moment conditions \eqref{IntzeroMoments} exclude
the generalized modes. Under these conditions, the projected
solution admits the spectral representation \eqref{sub1}.

\subsection{Main ideas of the proof}
Throughout this outline, $f=\omega_0\in H^5(\T^2)\cap Y$
satisfies the moment conditions \eqref{IntzeroMoments}.
We first specify the coordinate and parameter conventions.

\subsubsection{Coordinates and conventions}
We describe the coordinates on the positive cell
$\mathcal C=\{\Psi>0\}$; the other cell is treated by reflection. Set $r=2+\Psi$. The regular streamlines in $\mathcal C$ satisfy $2<r<4$, with $r=2$ at the separatrix and $r=4$ at the elliptic point $(0,0)$. Let $S(r)$ be the period of the flow generated by $L_0$ on a streamline. Define $\theta\in\R/(2\pi\Z)$ to be $2\pi$ times
the elapsed flow time divided by $S(r)$, measured from
$(\arccos(r-3),0)$. Then
\[
 L_0=b(r)\partial_\theta,\qquad
 b(r)=\frac{2\pi}{S(r)},\qquad
 2\pi<S(r)<\infty,\qquad 0<b(r)<1.
\]
The period is strictly decreasing, with $S'(r)<0$.
As $r\to4$, $S(r)\to2\pi$ and $b(r)\to1$; as $r\to2$,
$S(r)\to\infty$ and $b(r)\to0$. The construction and the
coordinate bounds are given in Subsection~\ref{sec:coordinates}.
For a function $h$, we write
\[
 h_k(r)=\frac1{2\pi}\int_0^{2\pi}h(r,\theta)e^{-ik\theta}\,d\theta,
 \qquad k\in\Z.
\]
The condition $h\in Y$ means that $h_0=0$ on each cell.
We note that the angular mean of the stream function need not vanish, since the Laplacian couples these modes.

Near an elliptic point we use the radius and logarithmic radius
\[
 q=\sqrt{8(1-b)},\qquad 0<q<q_0,\qquad
 v=\log q<\log q_0,
\]
where $q_0>0$ is fixed and sufficiently small. Lemma~\ref{geom:metric}
shows that $q$ is comparable to the distance from the center and that the coefficient matrix of the Dirichlet form in $(v,\theta)$ differs from the identity by $O(q^2)$. We write $b(q)$ for the same angular velocity in this coordinate; thus $b(q)=1-q^2/8$ exactly.

\subsubsection{The spectral representation and the resolvent equation}

Our starting point is the spectral representation of the projected
stream function:
\begin{equation}\label{sub1}
 \Delta^{-1}Q\omega(t)
 =\int_{\R}e^{-it\lambda}Q\Gamma_f(\lambda)\,d\lambda,
 \qquad t\ge0.
\end{equation}
Here $\Gamma_f$ is the spectral density associated with the initial
data $f$. To define it, set
\[
 \psi_f(z)=\Delta^{-1}(z-L)^{-1}f,
 \qquad z=\varepsilon+i\lambda,\quad
 \lambda\in\R,\quad \varepsilon\in\R\setminus\{0\},
\]
and write
\[
 \psi_{f,\pm}(\lambda)
 =\lim_{\varepsilon\to0+}\psi_f(i\lambda\pm\varepsilon),
 \qquad
 \Gamma_f(\lambda)
 =\frac{\psi_{f,+}(\lambda)-\psi_{f,-}(\lambda)}{2\pi},
 \qquad \lambda\ne0.
\]
The limiting absorption principle, Theorem~\ref{global:LAP},
justifies these boundary values.
As usual, once the singularities of the spectral density have been
resolved, integration by parts in $\lambda$ implies the decay
estimates.

The projection $Q$ and the moment conditions
\eqref{IntzeroMoments} remove different obstructions to decay.
The projection removes the stationary modes in $E$.
The generalized zero modes $b_j$, constructed in
\eqref{zero:C}, satisfy
$Lb_j=e_j$, $Le_j=0$, and $Qb_j=b_j$, where $e_j=L_0g_j$.
Thus $Qe^{-tL}b_j=b_j$, and these modes must also be excluded.
The identity
\[
 \langle(I+\Delta^{-1})Qf,b_j\rangle_{L^2}
 =\langle f,g_j\rangle_{L^2}=0
\]
shows that the moment conditions impose energy orthogonality
to the generalized zero modes. These conditions are preserved
by the evolution because $L^*g_j=0$; see
Proposition~\ref{zero:constraints}.

To estimate the spectral density, we study the equation for
the full resolvent $\psi_f$. Set
\[
 T(z)=(L_0-z)^{-1}L_0.
\]
The identity $(z-L)\Delta\psi_f=f$ implies
\begin{equation}\label{rayleigh1}
 (\Delta+T(z))\psi_f=-(L_0-z)^{-1}f.
\end{equation}
In the streamline coordinates introduced above, $T(z)$ is
diagonal in the angular variable:
\[
 T_k(i\lambda\pm0)
 =\frac{kb(r)}{kb(r)-\lambda\pm i0},
 \qquad k\ne0,
 \qquad T_0=0.
\]
The symbols $\pm i0$ denote distributional boundary values.
For the $k$th angular component, a critical streamline is a
level satisfying $kb(r)=\lambda$. When $k\ge1$, such a level
exists for $0<\lambda<k$ and contracts to the elliptic point
as $\lambda\to k-$. We call $\lambda=k$ an elliptic
threshold. Although $T(z)$ is diagonal in these coordinates,
the physical Laplacian couples the angular components.
Thus \eqref{rayleigh1} remains a two-dimensional equation.

\subsubsection{The elliptic thresholds and the comparison operator}
Consider $k\ge1$ and the resonant side of the threshold,
\[
 k(1-q_0^2/8)<\lambda<k,\qquad
 \delta=\sqrt{8(1-\lambda/k)}\in(0,q_0),\qquad
 s=\log(q/\delta).
\]
The critical streamline is $q=\delta$, or $s=0$. In the physical
center coordinate chart we have $s<\log(q_0/\delta)$; after localization we
compare with the operator on the whole line $s\in\R$,
\[
 -\partial_s^2+k^2+\frac{8e^{2s}}{e^{2s}-1\mp i0}.
\]
This rescaling removes the critical radius from the comparison
operator. At $\lambda=k$ we use $v=\log q$, since $s$ is no
longer defined; the threshold operator is then
$-\partial_v^2+k^2+8$.

The comparison operator can be inverted uniformly with respect to its parameters from the dual of the angular energy space with norm
\[
 \|u\|_{H^1_k}^2=
       \int_{\R}\bigl(|\partial_su|^2+k^2|u|^2\bigr)\,ds
\]
to the angular energy space.

By Lemma~\ref{geom:metric}, the metric errors are $O(q^2)$,
together with the corresponding derivative bounds. Their coupling of the
angular modes is estimated in divergence form and can be bounded perturbatively in a small neighborhood of elliptic centers. This rules out concentration at a collapsing
critical layer. Together with compactness away from the centers, the limiting absorption is reduced to a spectral condition that holds for the symmetric dipole.

To obtain the cubic depletion rate, we apply the comparison estimates
to $g=f+i\lambda\psi_f$. $g$ satisfies the equation:
\[
 (\Delta+T(i\lambda\pm0))g=(\Delta+1)f.
\]
The zero streamline average and $H^5$ regularity imply
$(\Delta+1)f=O(q)$ at the center. For $|k|\ge2$, the exponent
$\sqrt{k^2+8}$ is strictly larger than $3$, and the weighted estimate for the comparison operator near the elliptic centers implies a cubic bound at the critical circle.
For $k=\pm1$ the exponent
equals $3$, and this argument alone results in a logarithmic loss.
Here the algebraic structure of the dipole proves crucial: there is an exact local
homogeneous solution. Near $(0,0)$ it is given by  $(L_0-i\lambda)(x-iy)$; its first angular coefficient vanishes on
the critical circle, while its normal derivative does not.
Green's identity then implies a cubic bound for the critical trace;
see Subsection~\ref{subsec:cubic-first-trace}.

The two boundary values of the spectral density as $\epsilon\to0+$ have the same source term, which cancels in
the difference. Together with the cubic critical trace bound, this
implies the cubic estimate for the spectral jump without a logarithmic
loss.

\subsubsection{The separatrix}
In the region $S>S_*$, with $S_*$ fixed and sufficiently large,
we use the period as the radial coordinate. The area measure is
\[
 dxdy=\frac{S}{2\pi}\left|\frac{dr}{dS}\right|dS\,d\theta
       \sim S e^{-S/4}\,dS\,d\theta;
\]
see \eqref{geom:period} and \eqref{eq:weight}. For each fixed
$\lambda>0$, every integer $k>\lambda$ has a critical streamline $r_k$ satisfying $S(r_k)=2\pi k/\lambda$. As $k\to\infty$, $r_k\to2$. The infinitely many resonances are therefore critical streamlines accumulating spatially at the separatrix for a fixed spectral parameter.

We write each Cauchy multiplier as the derivative of a logarithm and integrate by parts in $S$. The area weight makes the resulting operator forms summable over the angular modes; Proposition~\ref{prop:tail} implies the basic estimate. As $\lambda\to0$, all nonzero-mode resonances satisfy
$S\ge2\pi/|\lambda|$. After separating the zero-mode singular
terms, this implies an exponentially small jump near zero;
see \eqref{tailD:zero-derivative}.

\subsubsection{Spectral derivatives and sharp damping}

The spectral representation \eqref{sub1} reduces inviscid damping
to regularity of the spectral density in $\lambda$. One integration
by parts implies $t^{-1}$ decay, while a second implies the optimal
$t^{-2}$ decay of the stream function. The difficulty is that
the critical streamlines move with $\lambda$: differentiating
the resolvent at a fixed spatial point produces increasingly
singular terms near these streamlines.

As in~\cite{BeekieChenJia2026}, we introduce an adapted spectral derivative \eqref{good-der} which does not worsen the singularity at the critical streamlines. As a result, we can apply the same limiting absorption principle to estimate the differentiated spectral density and the usual spectral derivative can then be recovered from the adapted spectral derivative of the spectral density. However, due to the degeneracy at the elliptic centers, this procedure does not remove all singularities. Near the elliptic centers, we separate
the explicitly identifiable logarithmic and principal-value
terms from a more regular remainder. The remainder can be
estimated by the limiting absorption principle, whereas the
singular terms are treated directly in the spectral integral.
It is essential to perform this integration before taking
norms, so that the cancellation in these terms is preserved.
The distinction becomes especially important for the second
spectral derivative: an absolute integrability bound for the
full derivative is not required. What we prove instead is
that its oscillatory integral is uniformly bounded in $L^2$,
which is sufficient for two integrations by parts.

The first spectral derivative estimates are established in
Proposition~\ref{prop:spectral-bounds}, and the second derivative
estimate in Proposition~\ref{quad:spectral}. Together
with the local estimates near the elliptic points, this closes
the second derivative argument for $f\in H^5$ and implies
\[
 \|\Delta^{-1}Q\omega(t)\|_{L^2}
 \lesssim (1+t)^{-2}\|f\|_{H^5}.
\]
The projected vorticity is uniformly bounded in $L^2$.
Interpolation therefore implies the velocity bound \eqref{Intm7}.

\subsection{Organization and notations}
In Section~\ref{sec:preliminaries} we introduce the coordinates and obtain
estimates on $\Delta^{-1}(\lambda-L_0)^{-1}$.
In Section~\ref{sec:lap} we prove the limiting absorption principle.
In Section~\ref{sec:spectral-preliminaries} we establish the spectral
representation and the preliminary bounds near the elliptic points,
including the quadratic critical-trace estimate in
Proposition~\ref{prop:weighted-center}. In Section~\ref{sec:higher} we treat the
first spectral derivative and the pointwise bounds, and in
Section~\ref{sec:second-spectral} we study the second spectral derivative.
In Section~\ref{sec:cubic-centers} we prove the cubic trace and spectral
jump bounds near elliptic centers, first for $|k|\ge2$ and then
for $k=\pm1$. Theorem~\ref{th1:Int} is
proved in Section~\ref{sec:main-proof}. Throughout the paper, $A\lesssim B$ means $A\le CB$, where $C$ depends only on fixed parameters; dependence on an additional
parameter is indicated by a subscript.

Unless another domain or measure is specified, we use the complex
$L^2$ inner product
\[
 \langle f,g\rangle=\langle f,g\rangle_{L^2(\T^2)}
 :=\int_{\T^2}f\overline g\,dx\,dy,
\]
which is linear in the first argument. The same notation denotes
its continuous extension to Sobolev duality, in particular the
$H^{-1}$--$H^1$ pairing. For one-dimensional model operators and
angular Fourier coefficients, the pairings use $dv$; for cylinder
operators, they use $dv\,d\theta$, on the domain under consideration.

\section{Bounds on coordinates and the compact operator}\label{sec:preliminaries}
In this section, we first establish coordinate estimates near the elliptic points
and the separatrix. We then use these estimates to control the
free resolvent composed with the inverse Laplacian.

\subsection{Bounds on coordinates}\label{sec:coordinates}
Recall that for $(x,y)\in\T^2$,
\[
 r=2+\cos x+\cos y,\qquad L_0=\sin y\,\partial_x-\sin x\,\partial_y.
\]
The flow generated by $L_0$ follows the level sets of $\Psi$:
$$\dot{x}(s)=\sin y(s),\quad \dot{y}(s)=-\sin x(s).$$
On $\{r>2\}$ we start the flow at $(\arccos(r-3),0)$ and set
$s=S(r)\theta/(2\pi)$, where $S(r)$ is the orbit period.
We set $S=\infty$ on the separatrix and treat the other cell by symmetry.
We normalize angular Fourier coefficients by
\[
 h_k(r)=\frac1{2\pi}\int_0^{2\pi}h(r,\theta)e^{-ik\theta}\,d\theta.
\]
In these coordinates,
\begin{equation}\label{geom:area}
 L_0=b(r)\partial_\theta,\qquad b=\frac{2\pi}{S},\qquad
 dxdy=\frac{S(r)}{2\pi}\,drd\theta.
\end{equation}
In the lower cell ($0<r<2$), $S(r)=S(4-r)$.
Unless otherwise indicated, Sobolev norms are on $\T^2$, and
$\Delta^{-1}$ acts on functions of integral zero.

\subsubsection{The period and the coordinates near an elliptic point}
We begin with the following properties of the period $S$.
\begin{lemma}\label{geom:period-lemma}
Let $S(r)$ be the minimal period of the $L_0$ flow on
$\{\Psi=r-2\}$, where $2<r<4$. Then
\begin{equation}\label{geom:period-formula}
 S(r)=4K\!\bigg(1-\frac{(r-2)^2}{4}\bigg),\qquad
 K(m)=\int_0^{\pi/2}\frac{du}{\sqrt{1-m\sin^2u}},
 \quad 0\le m<1.
\end{equation}
The function $S$ is real analytic on $(2,4)$, with $S'(r)<0$,
and maps this interval onto $(2\pi,\infty)$. It extends real
analytically across $r=4$, with
\begin{equation}\label{geom:period-center}
 S(4)=2\pi,\qquad S'(4)=-\frac\pi2.
\end{equation}
There is $d_0\in(0,1)$ such that, for $0<d\le d_0$,
\begin{equation}\label{geom:period}
 S(2+d)=-4\log d+12\log2+O(d^2|\log d|).
\end{equation}
More precisely,
$S(2+d)=S_1(d^2)\log d+S_2(d^2)$, where $S_1,S_2$ are
real analytic near zero, $S_1(0)=-4$ and $S_2(0)=12\log2$.
For every integer $j\ge0$ there is a constant $C_j$, independent
of $d$, such that
\begin{equation}\label{geom:period-remainder}
 \left|\partial_d^j\Bigl[S(2+d)+4\log d-12\log2\Bigr]\right|
 \le C_j d^{2-j}(1+|\log d|),\qquad 0<d\le d_0.
\end{equation}
In particular,
\begin{equation}\label{geom:period-derivatives}
 S'(2+d)=-\frac4d+O(d|\log d|),\qquad
 S''(2+d)=\frac4{d^2}+O(|\log d|).
\end{equation}
Consequently $b=2\pi/S$ is strictly increasing on $(2,4)$,
with $\lim_{r\to2^+}b(r)=0$, $b(4)=1$ and $b'(4)=1/4$.
\end{lemma}
\begin{proof}
Set $d=r-2$, $\xi=(x+y)/2$ and $\eta=(x-y)/2$.
On the positive cell the streamline equation and the flow imply
$\cos\xi\cos\eta=d/2$ and $\dot\xi=-\cos\xi\sin\eta$.
By symmetry, one quarter of the period is
\[
 \frac{S(r)}4
 =\int_0^{\arccos(d/2)}\frac{d\xi}{\sqrt{\cos^2\xi-d^2/4}}
 =K\!\left(1-\frac{d^2}{4}\right),
\]
where the last equality follows from
$\sin\xi=\sqrt{1-d^2/4}\sin u$.
Differentiation under the integral defining $K$, we obtain that
$K'(m)>0$, $K(0)=\pi/2$ and $K'(0)=\pi/8$.
Thus $S'(r)=-2dK'(1-d^2/4)<0$, and analyticity of $K$ at
$m=0$ implies \eqref{geom:period-center}.

For completeness, we note that the logarithmic expansion follows from the
regular singular equation for $F(t)=K(1-t)$:
\[
 t(1-t)F''+(1-2t)F'-\tfrac14F=0.
\]
This identity follows by differentiation under the integral and
one integration by parts. The indicial root at $t=0$ is a double
zero, so the Frobenius expansion has the form
$F(t)=F_1(t)\log t+F_2(t)$, with $F_1,F_2$ analytic near zero.
To identify the leading coefficients, write
\[
 F(t)=\int_0^1
 \frac{dx}{\sqrt{1-x^2}\sqrt{t+(1-t)x^2}}.
\]
Subtracting $\int_0^1(t+x^2)^{-1/2}\,dx
=\operatorname{arsinh}(t^{-1/2})$ and using dominated convergence, we obtain
\[
 \lim_{t\to0^+}\bigl(F(t)+\tfrac12\log t\bigr)
 =\log2+\int_0^1\frac{(1-x^2)^{-1/2}-1}{x}\,dx
 =\log4.
\]
Hence $F_1(0)=-1/2$ and $F_2(0)=\log4$.
Substituting $t=d^2/4$ proves \eqref{geom:period} and the stated
analytic form. The remainder in \eqref{geom:period-remainder}
can therefore be written as
$d^2[A(d^2)\log d+B(d^2)]$, with $A,B$ analytic near zero.
Choose $d_0>0$ so that these functions are analytic on a
neighborhood of $[0,d_0^2]$. For every integer $\ell\ge0$, we have
\[
 \left|\partial_d^\ell\big[d^2A(d^2)\big]\right|
 +\left|\partial_d^\ell\big[d^2B(d^2)\big]\right|
 \le C_\ell d^{2-\ell},\qquad 0<d\le d_0.
\]
For $\ell\le2$ this follows from the double zero at $d=0$;
for $\ell>2$ it follows from boundedness of the analytic
derivatives and $d_0<1$. \eqref{geom:period-remainder} then follows from Leibniz's rule and
$|\partial_d^m\log d|=(m-1)!d^{-m}$ for $m\ge1$. Its cases $j=1,2$ imply
\eqref{geom:period-derivatives}. The range of $S$ and the
conclusion about $b=2\pi/S$ follow easily.
\end{proof}

We next describe the coordinates near the elliptic points.
\begin{lemma}\label{geom:metric}
In a small neighborhood of $(0,0)$ set $\rho=\arccos(r-3)$ and $v=\log\rho$. Then
\[
 (x,y)=\rho(\cos\theta,-\sin\theta)+\rho^3F(\rho^2,\theta),
\]
for some smooth $F$. The exact Dirichlet form is
\[
 \int\nabla f\cdot\nabla\bar g\,dxdy
 =\int \mathsf G\nabla_{v,\theta}f\cdot
                \nabla_{v,\theta}\bar g\,dv\,d\theta,
 \qquad
 |\partial_v^a\partial_\theta^b(\mathsf G-I)|
 \lesssim_{a,b}e^{2v},
\]
where $\mathsf G$ is symmetric. We can replace the radial coordinate $\rho$ by $q$ and redefine
the logarithmic coordinate $v$ by
\begin{equation}\label{geom:q}
 q=\sqrt{8(1-b(r))}=\rho+O(\rho^3),\qquad v=\log q.
\end{equation}
The map $\rho\mapsto q$ is smooth and invertible near the center.
In the coordinates $(q,\theta)$, the preceding coordinate expansion
and Dirichlet form estimates remain valid with $\rho$ replaced by $q$
and with the remainder $F$ and coefficient matrix $\mathsf G$
expressed in the new coordinates.
\end{lemma}
\begin{proof}
We first notice that the velocity field is analytic and odd and the period is analytic in $\rho^2$. The expansion follows from analytic dependence on the initial point. The metric in
$(v,\theta)$ is $\rho^2(I+\rho^2h)$ with $h$ smooth, and therefore its
Dirichlet coefficient is
$\sqrt{\det(I+\rho^2h)}(I+\rho^2h)^{-1}=\mathsf{G}$.
The derivative bounds then follow from Taylor expansion and \eqref{geom:q} follows from the expansion of $S$
at $r=4$.
The change of radius from $\rho$ to $q$ preserves these bounds. By \eqref{geom:area},
the area measure has the form
\[
 dx\,dy=\xi_0(q)q\,dq\,d\theta,
\]
where $\xi_0$ is a smooth positive function of $q^2$, with
$\xi_0(q)=1+O(q^2)$ as $q\to0$.
\end{proof}

Let $X(r,\theta)=(x(r,\theta),y(r,\theta))$ denote the streamline
coordinate map on the positive cell, with $2<r<4$ and
$\theta\in\R/(2\pi\Z)$. Its derivatives
$\partial_\theta X=(x_\theta,y_\theta)$ and
$\partial_r X=(x_r,y_r)$ are functions of $(r,\theta)$.
The corresponding vector fields in physical coordinates are
$(\partial_\theta X)\circ X^{-1}$ and
$(\partial_r X)\circ X^{-1}$, respectively.
Since $|\det DX|=S(r)/(2\pi)$, we obtain from the coordinate formula for
divergence that
\[
 \begin{aligned}
 \left[\operatorname{div}_{x,y}
   \bigl((\partial_\theta X)\circ X^{-1}\bigr)\right]\circ X
 &=0,\\
 \left[\operatorname{div}_{x,y}
   \bigl((\partial_r X)\circ X^{-1}\bigr)\right]\circ X
 &=\frac{S'(r)}{S(r)}.
 \end{aligned}
\]

With a slight abuse of notation, we write
\begin{equation}\label{geom:identities}
 \operatorname{div}_{x,y}(x_\theta,y_\theta)=0,\qquad
 \operatorname{div}_{x,y}(x_r,y_r)=S'/S,\qquad
 [L_0,\Delta]=L_0+2(\cos x-\cos y)\partial_x\partial_y.
\end{equation}
Although $\omega\in Y$ implies $(\Delta\psi)_0=0$, in general
$\psi_0\ne0$. The metric couples the angular modes, so the elliptic
estimates below include the zero mode.

\subsubsection{Coordinates near the separatrix}
Near the separatrix it is convenient to use the period itself as
the radial coordinate. On the positive cell, set $\tau=S(r)$,
so that $2\pi<\tau<\infty$ and $\tau\to\infty$ as $r\to2^+$.
Write $X(\tau,\theta)=(x,y)$ for the coordinate map, with
$\theta\in\R/(2\pi\Z)$. Since $b(r)=2\pi/S(r)$, the resonance
condition for the $k$th angular mode becomes
\[
 kb(r)=\lambda\quad\Longleftrightarrow\quad 2\pi k=\lambda\tau.
\]
For each fixed $\lambda\ne0$, the resonances in $(2\pi,\infty)$
correspond to $k\lambda>0$ and $|k|>|\lambda|$, and are equally
spaced in $\tau$, with spacing $2\pi/|\lambda|$. Thus the infinitely
many critical streamlines accumulating at the separatrix become
separated critical values tending to infinity in this coordinate.
Moreover, the multiplier takes the form
\[
 \frac{kb(r)}{kb(r)-\lambda+i\eps}
 =\frac{2\pi k}{2\pi k-\lambda\tau+i\eps\tau}.
\]
Its denominator is linear in $\tau$, so its singularity can be
estimated by integration by parts using a logarithm;
see Proposition~\ref{prop:tail}.

We must also express physical area in the new coordinates.
It follows from \eqref{geom:area}, \eqref{geom:period} and
\eqref{geom:period-derivatives} that
\begin{equation}\label{eq:weight}
 w(\tau):=\tau\left|\frac{dr}{d\tau}\right|
 \sim\tau e^{-\tau/4},\qquad |w'|\lesssim w,\qquad
 dxdy=\frac{w(\tau)}{2\pi}\,d\tau d\theta.
\end{equation}
Thus $w$ is the area density, up to the angular normalization
$2\pi$, and is the natural weight for the physical norm estimates
below. Its exponential decay compensates for the polynomial
factors in $\tau$ that arise when summing the resonant contributions.

Let $X_\tau=(\partial_\tau X)\circ X^{-1}$ denote the corresponding
physical vector field. For any smooth test function $\phi$ compactly
supported in the image of $X$, \eqref{eq:weight} and integration
by parts imply
\[
\begin{aligned}
 \int\phi\,\operatorname{div}X_\tau\,dx\,dy
 &=-\frac1{2\pi}\int\partial_\tau(\phi\circ X)\,w(\tau)\,d\tau\,d\theta\\
 &=\frac1{2\pi}\int(\phi\circ X)\,w'(\tau)\,d\tau\,d\theta.
\end{aligned}
\]
Thus, using $w(\tau)=\tau/|S'(r(\tau))|$ and
$dr/d\tau=1/S'(r)$, we obtain
\begin{equation}\label{jxtb1}
 (\operatorname{div}X_\tau)\circ X
 =\frac{w'(\tau)}{w(\tau)}
 =\frac1\tau-\frac{S''(r(\tau))}{(S'(r(\tau)))^2}.
\end{equation}
It follows from \eqref{geom:period-derivatives}, and symmetry in the
other cell, that $|\operatorname{div}X_\tau|\le C$ for large $\tau$.

\begin{lemma}\label{lem:geometry}
For large $\tau$, $|\partial_\tau X|\lesssim1$. Consequently
\begin{align}
 \sum_k\int|\partial_\tau h_k|^2w\,d\tau
 &\lesssim\|\nabla h\|_2^2,\label{eq:derivative-L2}\\
 \|\partial_\tau h_k\|_{L^q(wd\tau)}
 &\lesssim_q\|\nabla h\|_q,\label{eq:derivative-Lq}\\
 \|h_k\|_{L^q(wd\tau)}&\lesssim_q\|h\|_q,
 \qquad 1\le q<\infty.\label{eq:coefficient-Lq}
\end{align}
Estimates \eqref{eq:derivative-Lq} and \eqref{eq:coefficient-Lq}
hold for every $k\in\Z$, with constants independent of $k$.
\end{lemma}
\begin{proof}
We first consider the positive cell. Fix $\tau_0$ sufficiently
large, and let $r=r(\tau)$ be determined by $\tau=S(r)$.
All integrals in $\tau$ below are over $[\tau_0,\infty)$.
It follows from \eqref{geom:period} and
\eqref{geom:period-derivatives} that
\[
 \tau=-4\log(r-2)+O(1),\qquad
 |S'(r)|\sim(r-2)^{-1}.
\]
Consequently,
\[
 \left|\frac{dr}{d\tau}\right|\lesssim r-2,
 \qquad
 e^{\tau/8}\lesssim(r-2)^{-1/2}.
\]

Let $\Phi_s$ denote the flow generated by
$U(x,y)=(\sin y,-\sin x)$. With the angular normalization
fixed above,
\[
 X(\tau,\theta)
 =\Phi_{\tau\theta/(2\pi)}
       \bigl(\arccos(r(\tau)-3),0\bigr).
\]
Since
\[
 DU(x,y)=
 \begin{pmatrix}
  0&\cos y\\
  -\cos x&0
 \end{pmatrix},
 \qquad \|DU(x,y)\|\le1,
\]
it follows from the variational equation and Gronwall's
inequality that
\[
 \|D\Phi_s(x,y)\|\le e^{|s|}.
\]
Moreover,
\[
 \left|\frac{d}{dr}\bigl(\arccos(r-3),0\bigr)\right|
 =\frac1{\sqrt{1-(r-3)^2}}
 \lesssim(r-2)^{-1/2}.
\]
Thus, with $s$ held fixed,
\[
 \left|\partial_r
 \bigl[\Phi_s(\arccos(r-3),0)\bigr]\right|
 \lesssim e^{|s|}(r-2)^{-1/2}.
\]

We next justify restricting attention to $|\theta|\le\pi/4$.
The clockwise quarter-turn $R(x,y)=(y,-x)$ preserves $\Psi$
and commutes with the flow. The four successive quadrant
arcs of each streamline therefore have equal travel times.
Their total travel time is $\tau$, so
\[
 \Phi_{\tau/4}\bigl(\arccos(r-3),0\bigr)
 =R\bigl(\arccos(r-3),0\bigr).
\]
It follows that
\[
 X(\tau,\theta+j\pi/2)=R^jX(\tau,\theta),
 \qquad j=0,1,2,3.
\]
The intervals $j\pi/2+[-\pi/4,\pi/4]$ cover the angular
circle. Since the rotations and angular shifts in this
identity are independent of $\tau$, it suffices to estimate
$\partial_\tau X$ for $|\theta|\le\pi/4$.

For such $\theta$, set $s=\tau\theta/(2\pi)$, so that
$|s|\le\tau/8$. Differentiating at fixed $\theta$, we obtain
\[
 \partial_\tau X
 =\frac{dr}{d\tau}\,
   \partial_r\bigl[\Phi_s(\arccos(r-3),0)\bigr]
   +\frac{\theta}{2\pi}U(X),
\]
where the derivative in the first term is taken at fixed $s$.
It follows from the preceding estimates and the boundedness
of $U$ that
\[
 |\partial_\tau X|
 \lesssim (r-2)e^{\tau/8}(r-2)^{-1/2}+1
 \lesssim1.
\]
By the rotation symmetry this bound holds for every $\theta$.
The other cell can be handled by translation and reversal of the flow.

We now prove the Fourier estimates. By approximation, it suffices to consider smooth $h$.
By definition,
\[
 h_k(\tau)=\frac1{2\pi}\int_0^{2\pi}
 h(X(\tau,\theta))e^{-ik\theta}\,d\theta,
 \qquad
 \partial_\tau(h\circ X)
 =\nabla h(X)\cdot\partial_\tau X.
\]
It follows from Parseval's identity and \eqref{eq:weight} that
\[
 \begin{aligned}
 \sum_{k\in\Z}\int_{\tau_0}^{\infty}
       |\partial_\tau h_k|^2w(\tau)\,d\tau
 &=\frac1{2\pi}\int_{\tau_0}^{\infty}\int_0^{2\pi}
       |\nabla h(X)\cdot\partial_\tau X|^2
       w(\tau)\,d\theta\,d\tau\lesssim\|\nabla h\|_{L^2(\T^2)}^2.
 \end{aligned}
\]
For $1\le q<\infty$, it follows from Hölder's inequality that
\[
 \begin{aligned}
 |h_k(\tau)|^q
 &\le\frac1{2\pi}\int_0^{2\pi}
          |h(X(\tau,\theta))|^q\,d\theta,\\
 |\partial_\tau h_k(\tau)|^q
 &\le\frac1{2\pi}\int_0^{2\pi}
          |\nabla h(X)\cdot\partial_\tau X|^q\,d\theta.
 \end{aligned}
\]
Integrating against $w(\tau)\,d\tau$ and using
\eqref{eq:weight}, we conclude that
\[
 \|h_k\|_{L^q(w\,d\tau)}
 \le\|h\|_{L^q(\T^2)},\qquad
 \|\partial_\tau h_k\|_{L^q(w\,d\tau)}
 \lesssim_q\|\nabla h\|_{L^q(\T^2)},
\]
with constants independent of $k$.

\end{proof}

\subsubsection{A vector field adapted to the spectral parameter}
For spectral derivatives we combine $\partial_\lambda$ with a radial field $\bv$, writing $\bv h=V\cdot\nabla h$. The field below satisfies
$\bv b(r)=b(r)/\lambda$, and hence
\[
 (\partial_\lambda+\bv)\frac{kb(r)}{kb(r)-\lambda\pm i0}=0
\]
away from the critical layer $kb(r)=\lambda$. We estimate $V$ in physical coordinates, including across the separatrix.
\begin{lemma}\label{geom:common-field}
For $\lambda$ in a compact subset $I\subset\R\setminus\{0\}$,
the field
\[
 \bv=-\frac{S}{\lambda S'}\partial_r
\]
has a continuous physical extension across the separatrix. After a smooth cutoff near the elliptic centers, $V\in W^{1,p}$ for every finite $p$, with $\|V\|_\infty+\|DV\|_p\le C_{I,p}$. There is no jump measure in its distributional derivatives.
\end{lemma}
\begin{proof}
Throughout the proof, $\partial_r$ in the definition of $\bv$ is taken
at fixed $\theta$. It follows from Lemma~\ref{geom:period-lemma}
and the symmetry $S(r)=S(4-r)$ that, for $0<|r-2|$ small,
\begin{equation}\label{geom:field-period-bounds}
 \left|\frac{S(r)}{S'(r)}\right|
 \lesssim |r-2|(1+|\log|r-2||),\qquad
 \left|\frac{d}{dr}\frac{S(r)}{S'(r)}\right|
 \lesssim 1+|\log|r-2||.
\end{equation}
Although $S/S'$ tends to zero at the separatrix, the coordinate field $\partial_r$ is singular there and as a result $V$ may not vanish on the separatrix. We first make its tangential component explicit away from the saddles.

We first consider a compact subarc of the separatrix away from
the saddle points. By symmetry, it suffices to work on the arc
joining $(\pi,0)$ to $(0,-\pi)$. Choose the transverse section
\[
 \gamma(r)=
 \left(\arccos\frac{r-2}{2},-\arccos\frac{r-2}{2}\right),
\]
and let $\sigma$ be the signed flow time from this section.
Since $\gamma$ is smooth and transverse to $L_0$ near $r=2$,
smooth dependence of the flow on its initial point yields
smooth coordinates $(r,\sigma)$ across the chosen subarc.
After shrinking the neighborhood, $\sigma$ ranges over a
bounded interval, and its physical derivatives through
order two are bounded.

The relation $L_0\sigma=1$ implies
\[
 \det\frac{\partial(r,\sigma)}{\partial(x,y)}
 =r_x\sigma_y-r_y\sigma_x=1,
\]
because $L_0=-r_y\partial_x+r_x\partial_y$.
Inverting the Jacobian, we obtain the coordinate vector fields
in physical variables:
\[
 \left.\partial_r\right|_\sigma
 =\sigma_y\partial_x-\sigma_x\partial_y,
 \qquad
 \partial_\sigma=-r_y\partial_x+r_x\partial_y.
\]
Their coefficients and first derivatives are bounded on this
neighborhood. In particular, physical derivatives and
derivatives in $(r,\sigma)$ are uniformly comparable:
\[
 |\nabla f|
 \sim
 \left|\left.\partial_r\right|_\sigma f\right|
 +|\partial_\sigma f|.
\]

By the symmetry of the streamlines, the section $\gamma$
has a constant angular phase $\theta_c$ in each cell.
Thus, with a local choice of $\theta$,
\[
 \sigma=\frac{S(r)}{2\pi}(\theta-\theta_c).
\]
The value of $\theta_c$ may differ between the cells, but
$\sigma$ is the same smooth physical flow time on both sides.
It follows that
\[
 \left.\partial_r\right|_\theta
 =
 \left.\partial_r\right|_\sigma
 +\frac{S'(r)}{S(r)}\sigma\partial_\sigma.
\]
Consequently, the physical vector field $V$ is given by
\[
 \lambda V
 =-\frac{S(r)}{S'(r)}(\sigma_y,-\sigma_x)
   -\sigma(-r_y,r_x).
\]
It follows from \eqref{geom:field-period-bounds} and the
bounded first and second derivatives of $r$ and $\sigma$ that
\[
 |V|\lesssim_I 1,\qquad
 |DV|\lesssim_I 1+|\log|r-2||,
\]
where the constants may depend on the chosen subarc.
Since $dx\,dy=dr\,d\sigma$, $DV$ belongs
locally to $L^p$ for every $1\le p<\infty$.

Finally, $S/S'$ tends to zero as $r$ tends to $2$ from either
side. The two traces of $V$ therefore agree and equal
\[
 V\big|_{r=2}
 =-\frac{\sigma}{\lambda}(-r_y,r_x)\big|_{r=2}.
\]
Hence $V$ extends continuously across this subarc, and
integration by parts on the two sides produces no interface
term. The preceding bounds therefore prove that
$V\in W^{1,p}$ locally across the regular part of the
separatrix for every finite $p$.

It remains to examine the saddles. By symmetry it suffices to
work near $(\pi,0)$. Introduce the local coordinates
\[
 a=\sqrt2\sin\frac{x-\pi+y}{2},\qquad
 b=\sqrt2\sin\frac{x-\pi-y}{2}.
\]
The Jacobian determinant of this change of variables equals $-1$
at $(\pi,0)$, so $(a,b)$ form a smooth coordinate chart on a full
neighborhood of the saddle. Moreover,
\[
 r-2=\Psi(x,y)=ab.
\]
Thus the two branches of the separatrix are the coordinate axes
$a=0$ and $b=0$, and nearby streamlines are the hyperbolas
$ab=\text{constant}$.
The identity \(r-2=ab\) holds on both sides of the separatrix.
Recall that \(S(r)=S(4-r)\) for \(r<2\), so
\(S^{\prime}(r)=-S^{\prime}(4-r)\) on that interval.
Writing
\[
 j(u)=(1-u^2/2)^{-1/2},\qquad J=j(a)j(b),
\]
we have $L_0=(L_0a)\partial_a+(L_0b)\partial_b=J^{-1}(a\partial_a-b\partial_b)$.
Fix a small $a_0>0$ and first consider
$0<|b|\le |a|\le a_0$. The remaining forward travel time
to the symmetry section $\gamma_2$ with $|a|=\sqrt2$ is
\begin{equation}\label{geom:field-travel-time}
 t(a,b)=\int_{|a|}^{\sqrt2}j(u)j(ab/u)\,\frac{du}{u}.
\end{equation}
Let $\theta_c$ denote the constant angular phase of $\gamma_2$.
Choose a local lift of $\theta$ so that $\theta_c-\theta$
is the forward angular displacement to this section.
Since $L_0\theta=2\pi/S(r)$, the remaining travel time is
\[
 t=\frac{S(r)}{2\pi}(\theta_c-\theta).
\]
For the sector of the positive cell adjacent to $(\pi,0)$,
the section $\gamma_2$ lies on $x+y=0$, and hence
\[
 \theta_c=\frac{\pi}{4},\qquad
 t=S(r)\left(\frac18-\frac{\theta}{2\pi}\right),
 \qquad 0\le\theta\le\frac{\pi}{4}.
\]
In particular, $t/S(r)$ is independent of $r$ at fixed $\theta$. Since $\lambda\bv S=-S$, it follows that
\[
 \lambda\bv t=-t,\qquad
 \lambda\bv(ab)=-\frac{S(r)}{S'(r)}.
\]
Also, $L_0t=-1$, so $a t_a-bt_b=-J$. Solving the preceding
two equations for the components of $\bv$ we obtain
\begin{equation}\label{geom:field-components}
 \lambda\bv a=\frac{a t-(S(r)/S'(r))t_b}{J},\qquad
 \lambda\bv b=\frac{(S(r)/S'(r))t_a-bt}{J}.
\end{equation}
These identities preserve the cancellation between the vanishing
coefficient $S/S'$ and the singular coordinate derivatives.

We establish the derivative bounds needed to use
\eqref{geom:field-components}. Since $|ab/u|\le |b|\le a_0$,
we may use $|j'(z)|\lesssim|z|$ and $|j''(z)|\lesssim1$
when differentiating \eqref{geom:field-travel-time} in $b$.
It follows from
\[
 \int_{|a|}^{\sqrt2}\frac{j(u)}{u}\,du
 \lesssim1+|\log|a||,\qquad
 \int_{|a|}^{\sqrt2}\frac{j(u)}{u^3}\,du
 \lesssim |a|^{-2}
\]
that
\[
 |t|\lesssim1+|\log|a||,\qquad
 |t_b|\lesssim|b|,\qquad |t_{bb}|\lesssim1.
\]
For example,
$|t_b|\lesssim |a|^2|b|\int_{|a|}^{\sqrt2}j(u)u^{-3}\,du$.
The identity $a t_a-bt_b=-J$ and its first derivatives imply
\[
 |a t_a|\lesssim1,\qquad
 |a t_{ab}|\lesssim|b|,\qquad
 |a^2t_{aa}|\lesssim1.
\]
Here we have used $J_a=O(a)$ and $J_b=O(b)$ in the local chart.

Substitute these bounds and \eqref{geom:field-period-bounds}
into \eqref{geom:field-components}, with $r=2+ab$.
When differentiating the period quotient, the chain rule is
\[
 \partial_a\!\left(\frac{S(r)}{S'(r)}\right)
 =b\left(\frac{S}{S'}\right)'(r),\qquad
 \partial_b\!\left(\frac{S(r)}{S'(r)}\right)
 =a\left(\frac{S}{S'}\right)'(r).
\]
In particular, the potentially singular terms satisfy
\[
 \left|b\left(\frac{S}{S'}\right)'t_a\right|
 +\left|\frac{S}{S'}t_{aa}\right|
 \lesssim\frac{|b|}{|a|}(1+|\log|ab||)
 \lesssim1+|\log|a||+|\log|b||.
\]
All the other terms are controlled directly by the above
bounds. Since the coordinate change is smooth and invertible,
we conclude in physical coordinates that
\begin{equation}\label{geom:field-saddle-bounds}
 \begin{split}
 |V|&\lesssim_I
 |a|(1+|\log|a||)+|b|(1+|\log|b||),\\
 |DV|&\lesssim_I1+|\log|a||+|\log|b||.
 \end{split}
\end{equation}
In the sector $0<|a|\le |b|\le a_0$, the same argument applies
after exchanging $a,b$ and using backward travel time to the
section $|b|=\sqrt2$.

For fixed $a\ne0$, the functions $t,t_a,t_b$ extend smoothly
through $b=0$. Since $S/S'$ tends to zero there, it follows
from \eqref{geom:field-components} that the traces from both
cells agree and satisfy
\[
 \bv a=\frac{a\,t(a,0)}{\lambda j(a)},\qquad \bv b=0.
\]
The traces across $a=0$ agree by the analogous backward-time
argument. The first bound in \eqref{geom:field-saddle-bounds}
also shows that $V$ tends to zero at the saddle. This proves
continuity across the entire separatrix.

Finally, the logarithmic bound for $DV$ is locally in $L^p$
for every $1\le p<\infty$.
\end{proof}

\subsection{Bounds on the compact operator}\label{sec:compact}
We estimate the free resolvent for the limiting absorption principle.
For $z=\eps+i\lambda$, $\eps\ne0$, define on the mean-zero
Sobolev spaces on $\T^2$
\begin{equation}\label{eq:T}
 \begin{gathered}
 K(z)=\Delta^{-1}T(z),\qquad T(z)=(L_0-z)^{-1}L_0,\\
 (T(z)h)_k=\frac{2\pi k}{2\pi k-\lambda S+i\eps S}h_k
 \quad(k\ne0),\qquad (T(z)h)_0=0.
 \end{gathered}
\end{equation}
We write $T(i\lambda\pm0)$ for the limits as $\eps\to0^\pm$.
The singularities $\lambda=kb(r)$ lie on regular streamlines when
$\lambda\notin\Z$, where $K$ has compact boundary values. As
$\lambda$ approaches an integer, a critical streamline can approach
a center. We first establish a uniform inverse estimate for the model operator obtained by rescaling the resolvent equation near an elliptic point.

\subsubsection{The comparison operator near an elliptic point}
Write $\|u\|_k^2=\int_{\R}(|u'|^2+k^2|u|^2)\,dv$ and
let $\|\cdot\|_{k,*}$ be the dual norm.
\begin{proposition}\label{compact:normal}
For $k\ge1$ and real $a,\eta$ with $|a|+|\eta|=1$,
\[
 \mathcal N_{k,a,\eta}=-\partial_v^2+k^2+
                   \frac{8e^{2v}}{e^{2v}-a-i\eta}
\]
is a bijection $H^1_k\to H^{-1}_k$, with
$\|u\|_k\le C\|\mathcal N_{k,a,\eta}u\|_{k,*}$ uniformly
in all parameters. At $(a,\eta)=(1,0)$ we use the boundary values:
\[
 \frac{8e^{2v}}{e^{2v}-1\mp i0}
 =8e^{2v}\operatorname{p.v.}\frac1{e^{2v}-1}
                                      \mathbin{\pm}4\pi i\delta_0.
\]
The same inverse bound holds for $-\partial_v^2+k^2+8$. Moreover,
for arbitrary real $a,\eta$, including the boundary values,
\begin{equation}\label{eq:cylinder-form}
 \left|\int c(v)\frac{8e^{2v}}{e^{2v}-a-i\eta}
                         u\bar w\,dv\right|
 \le C\|c\|_{W^{1,\infty}}k^{-1/2}\|u\|_k\|w\|_k.
\end{equation}
\end{proposition}
\begin{proof}
We first assume $|a|+|\eta|=1$. Set
\[
 \chi(v)=\frac{e^{2v}}{1+e^{2v}},
 \qquad
 \begin{aligned}
 F_{a,\eta}(v)
 ={}&\log(e^{2v}-a-i\eta)-\log(1+e^{2v})\\
    &-(1-\chi(v))\log(-a-i\eta).
 \end{aligned}
\]
At $(a,\eta)=(1,0)$, the logarithms are understood as limits as
$\eta\to0^\pm$. In particular,
\[
 \log(e^{2v}-1\mp i0)
 =\log|e^{2v}-1|\mp i\pi\mathbf1_{\{v<0\}}.
\]
By direct differentiation, we obtain, in the sense of distributions,
\begin{equation}\label{eq:model-log-identity}
 \frac{8e^{2v}}{e^{2v}-a-i\eta}
 =8\chi-4\chi'\log(-a-i\eta)+4F_{a,\eta}'.
\end{equation}
Moreover,
\[
 |F_{a,\eta}(v)|\le Ce^{-2|v|}\quad (|v|\ge1),
 \qquad
 \sup_{|a|+|\eta|=1}\|F_{a,\eta}\|_{L^4(-1,1)}\le C.
\]
For the local bound, the logarithm is uniformly bounded on
$[-1,1]$ when $a\le e^{-2}/2$. Otherwise, the mean value theorem
implies
\[
 \bigl|\log(e^{2v}-a-i\eta)\bigr|
 \le C+\left|\log\left|v-\frac12\log a\right|\right|.
\]
Since $\frac12\log a\in[-1-\frac12\log2,0]$, the right-hand side
has uniformly bounded $L^4(-1,1)$ norm.
These bounds and pointwise convergence show that $F_{a,\eta}$ is
uniformly bounded in $L^2(\mathbb R)$ and continuous in that space
on each closed parameter set
\[
 \{|a|+|\eta|=1,\ \eta\ge0\},
 \qquad
 \{|a|+|\eta|=1,\ \eta\le0\},
\]
with the respective boundary values.

For $u,w\in H^1(\mathbb R)$ and $c\in W^{1,\infty}(\mathbb R)$,
integration by parts in \eqref{eq:model-log-identity} yields
\[
\begin{aligned}
 \int_{\mathbb R}
 c\,\frac{8e^{2v}}{e^{2v}-a-i\eta}\,u\overline w\,dv
 ={}&\int_{\mathbb R}
 c\bigl(8\chi-4\chi'\log(-a-i\eta)\bigr)
 u\overline w\,dv\\
 &-4\int_{\mathbb R}F_{a,\eta}
 \bigl(c'u\overline w+cu'\overline w
                    +cu\overline{w'}\bigr)\,dv.
\end{aligned}
\]
Since
\[
 \|u\|_{L^2}\le k^{-1}\|u\|_k,
 \qquad
 \|u\|_{L^\infty}\le Ck^{-1/2}\|u\|_k,
\]
it follows that
\[
 \left|\int_{\mathbb R}
 c\,\frac{8e^{2v}}{e^{2v}-a-i\eta}\,u\overline w\,dv\right|
 \le C\|c\|_{W^{1,\infty}}k^{-1/2}\|u\|_k\|w\|_k.
\]
Translation by $\frac12\log(|a|+|\eta|)$ proves the same estimate
for every nonzero parameter pair. The case $a=\eta=0$ is immediate.

We next turn to the proof of the bounds on $\mathcal{N}_{k,a,\eta}$, which is based on injectivity and the Fredholm theory. Suppose $\mathcal N_{k,a,\eta}u=0$.
If $\eta\ne0$, then
\[
 0=\operatorname{Im}\langle\mathcal N_{k,a,\eta}u,u\rangle
 =\int_{\mathbb R}
 \frac{8\eta e^{2v}}{(e^{2v}-a)^2+\eta^2}|u|^2\,dv,
\]
hence $u=0$. At $(a,\eta)=(-1,0)$,
\[
 0=\langle\mathcal N_{k,-1,0}u,u\rangle
 =\|u\|_k^2+8\int_{\mathbb R}\chi|u|^2\,dv,
\]
with the same conclusion.

For either boundary value at $(1,0)$, taking imaginary parts, we obtain
\[
 0=\pm4\pi|u(0)|^2.
\]
Thus $u_+=\mathbf1_{(0,\infty)}u\in H^1(\mathbb R)$, and testing
against $u_+$ yields
\[
 0=\int_0^\infty
 \left(|u'|^2+k^2|u|^2+
              \frac{8e^{2v}}{e^{2v}-1}|u|^2\right)\,dv.
\]
Consequently $u=0$ on $[0,\infty)$. Since $u(0)=0$,
\[
 |u(v)|\le |v|^{1/2}\|u'\|_{L^2(-1,1)}
 \qquad (|v|<1).
\]
The distributional equation therefore implies
\[
 u''=\left(k^2+\frac{8e^{2v}}{e^{2v}-1}\right)u
 \in L^1_{\mathrm{loc}},
\]
so $u'$ is continuous and $u'(0)=0$. For $-h\le v\le0$,
\[
 u'(v)=-\int_v^0
 \left(k^2+\frac{8e^{2s}}{e^{2s}-1}\right)u(s)\,ds.
\]
Using $|u(s)|\le |s|\sup_{[-h,0]}|u'|$, we obtain
\[
 \sup_{[-h,0]}|u'|
 \le\left(\frac{k^2h^2}{2}+Ch\right)
                       \sup_{[-h,0]}|u'|.
\]
For sufficiently small $h$, this implies $u=0$ on $[-h,0]$.
Ordinary differential equation uniqueness then implies $u=0$
everywhere.

Finally, \eqref{eq:model-log-identity} expresses
$\mathcal N_{k,a,\eta}$ as a compact perturbation of the coercive
operator
\[
 -\partial_v^2+k^2+8\chi.
\]
Indeed, if $F_j\in C_c^\infty$ converges to $F_{a,\eta}$ in $L^2$,
then
\[
 \bigl|\langle(F_{a,\eta}'-F_j')u,w\rangle\bigr|
 \le C\|F_{a,\eta}-F_j\|_{L^2}
                   \|u\|_{H^1}\|w\|_{H^1}.
\]
Multiplication by $F_j'$ is compact from $H^1$ to $H^{-1}$, as is
multiplication by the exponentially decaying function $\chi'$.
By the Fredholm theory and injectivity we conclude that $\mathcal{N}_{k,a,\eta}$ is bijective. The continuity of the norm of the operator in the parameters follows from the identity \eqref{eq:model-log-identity}.
Hence compactness of the two closed parameter sets implies a uniform inverse bound for $1\le k\le K$, for every fixed $K$.

For large $k$, we have
\[
 \|u\|_k
 \le\|\mathcal N_{k,a,\eta}u\|_{k,*}
       +Ck^{-1/2}\|u\|_k,
\]
and the desired bounds follow directly. The required bound for
$-\partial_v^2+k^2+8$ follows directly from coercivity.
\end{proof}

To see the relation of this operator to \eqref{eq:T}, take the radius
\eqref{geom:q} and set
$\sigma_k(\lambda)=8(1-\lambda/k)$. The corresponding multiplier
in the Dirichlet form is
\[
 -\xi_0(q)b(q)\frac{8e^{2v}}{e^{2v}-\sigma_k-i\eta_k},\qquad
 \eta_k=8\eps/k.
\]
Rescaling $q=\delta e^v$, $\delta^2=|\sigma_k|+|\eta_k|$, we obtain
$\mathcal N_{k,\sigma_k/\delta^2,\eta_k/\delta^2}$ as the resonant
normal operator. The other nonzero modes have free normal operator
$-\partial_v^2+n^2$; the metric error is $O(q^2)$ in form norm.
For the zero mode we use the homogeneous derivative norm. At
$\lambda=k$, the characteristic exponents of the resonant equation are
$\pm\sqrt{k^2+8}$.

\subsubsection{Bounds near the separatrix}
We have the following bounds for $T(z),\, z=i\lambda+\epsilon$ localized near the separatrix.
\begin{proposition}\label{prop:tail}
Let $0\le\chi_R\le1$ satisfy
\[
 \chi_R(\tau)=0\quad(\tau\le R),\qquad
 \chi_R(\tau)=1\quad(\tau\ge R+1),\qquad
 \sup_R\|\chi_R'\|_{L^\infty}<\infty.
\]
There exist $R_0,c>0$, independent of $M$, such that, for every
$M>0$ and $R\ge R_0$,
\begin{equation}\label{eq:tail}
 \sup_{\substack{|\lambda|\le M\\0<|\varepsilon|\le1}}
 \|\chi_R(S)T(\varepsilon+i\lambda)\|
 _{H^1(\mathbb T^2)\to H^{-1}(\mathbb T^2)}
 \le C_Me^{-cR}.
\end{equation}
For each integer $m\ge0$ and every $h\in H^1(\mathbb T^2)$,
\begin{equation}\label{eq:tail-weighted}
\begin{aligned}
 &\sup_{\substack{|\lambda|\le M\\0<|\varepsilon|\le1}}
 \left\|\chi_R(S)\sum_{k\ne0}
 \frac{S^m}{k^m}
 \frac{2\pi k}{2\pi k-\lambda S+i\varepsilon S}
 h_ke^{ik\theta}\right\|_{H^{-1}(\mathbb T^2)}\le C_{M,m}e^{-c_mR}\|h\|_{H^1(\mathbb T^2)}.
\end{aligned}
\end{equation}
The Fourier series is interpreted separately on each cell, with
$h_k$ being the angular Fourier coefficients of $h$ on that cell.
The constants $C_{M,m},c_m>0$ are independent of
$R,\lambda,\varepsilon$ and $h$.
\end{proposition}

\begin{proof}
All physical-space norms below use the measure $dx\,dy$;
the notation $L^p(I,w\,d\tau)$ explicitly indicates the weighted
norm on an interval $I$.
Choose $R_0\ge\tau_0+2$, where the coordinate estimates of
Lemma~\ref{lem:geometry} hold for $\tau\ge\tau_0$.
By duality, it suffices to estimate the pairing \eqref{eq:tail-form} below on the positive cell
$\mathcal C=\{\Psi>0\}$; the other cell is treated identically.
Write $h_k,g_k$ for the angular Fourier coefficients on this cell.

Assume first that $0<\lambda\le M$. For $k>0$, set
\[
 \tau_k=\frac{2\pi k}{\lambda},\qquad
 d=\min\left(1,\frac{\pi}{4M}\right),\qquad
 I_k=(\tau_k-d,\tau_k+d).
\]
Choose $0\le\chi_k\le1$, compactly supported in $I_k$, equal to
one on $(\tau_k-d/2,\tau_k+d/2)$, with $|\chi_k'|\le C_M$.
Set $\chi_k=0$ for $k<0$. Since
$dx\,dy=w(\tau)\,d\tau\,d\theta/(2\pi)$, by Parseval's identity
and $1=(1-\chi_k)+\chi_k$ we have the decomposition
\begin{equation}\label{eq:tail-form}
\begin{aligned}
 \int_{\mathcal C}\chi_R(S)T(\varepsilon+i\lambda)h\,
                  \overline g\,dx\,dy&=\sum_{k\ne0}\int_{2\pi}^{\infty}
 \chi_R(1-\chi_k)
 \frac{2\pi k}{2\pi k-\lambda\tau+i\varepsilon\tau}
 h_k\overline{g_k}\,w\,d\tau\\
 &\qquad+\sum_{\substack{k\ge1\\\tau_k>R-d}}\int_{I_k}
 \chi_R\chi_k
 \frac{2\pi k}{2\pi k-\lambda\tau+i\varepsilon\tau}
 h_k\overline{g_k}\,w\,d\tau.
\end{aligned}
\end{equation}
The first sum is the regular part, and the second is the singular
part. Only intervals meeting $\{\tau>R\}$ contribute to the
second sum, which explains the restriction $\tau_k>R-d$.

For the regular part, we have
\begin{equation}\label{eq:nonres}
 \left|(1-\chi_k(\tau))
 \frac{2\pi k}{2\pi k-\lambda\tau+i\varepsilon\tau}\right|
 \le C_M(1+\tau).
\end{equation}
By Parseval's identity and
Hölder's inequality, the absolute value of the first sum in
\eqref{eq:tail-form} is bounded by
\begin{equation}\label{eq:tail-nonres-form}
\begin{aligned}
 &C_M\left(\int_{2\pi}^{\infty}
  (1+\tau)^2\chi_R(\tau)^2w(\tau)\,d\tau\right)^{1/2}
 \|h\|_{L^4(\mathbb T^2)}\|g\|_{L^4(\mathbb T^2)}\le C_Me^{-cR}\|h\|_{H^1(\mathbb T^2)}
                              \|g\|_{H^1(\mathbb T^2)}.
\end{aligned}
\end{equation}
Here we used $w(\tau)\sim\tau e^{-\tau/4}$.

We now estimate the second sum in \eqref{eq:tail-form}.
For every interval appearing there,
\[
 \tau_k>R-d,\qquad I_k\subset(\tau_0,\infty).
\]
Set
$\zeta_k=\frac{2\pi k}{\lambda-i\varepsilon}.$
Then $|\zeta_k|\le\tau_k$ and
\[
 \frac{2\pi k}{2\pi k-\lambda\tau+i\varepsilon\tau}
 =-\zeta_k\,\partial_\tau\log(\tau-\zeta_k).
\]
On $I_k$, we have $w(\tau)\sim\tau_ke^{-\tau_k/4}$.
Consequently, for $1\le p<\infty$,
\begin{equation}\label{eq:log}
 \|\log(\tau-\zeta_k)\|_{L^p(I_k,w\,d\tau)}
 \le C_{M,p}(\tau_ke^{-\tau_k/4})^{1/p}
       \bigl(1+\log(2+\tau_k)\bigr).
\end{equation}
Indeed, the negative part of $\log|\tau-\zeta_k|$ has uniformly
bounded local $L^p$ norm, while $|\tau-\zeta_k|\le2\tau_k+d$
and the argument of the logarithm is bounded by $\pi$.

Using integration by parts we get the identity
\begin{equation}\label{eq:tail-ibp}
\begin{aligned}
 \int_{I_k}\chi_R\chi_k
  \frac{2\pi k}{2\pi k-\lambda\tau+i\varepsilon\tau}
  h_k\overline{g_k}\,w\,d\tau &=\zeta_k\int_{I_k}\log(\tau-\zeta_k)
  \Bigg[\chi_R\chi_k
  \bigl(h_k'\overline{g_k}+h_k\overline{g_k'}\bigr)\\
 &\hspace{10mm}
  +\left(\chi_R'\chi_k+\chi_R\chi_k'
         +\chi_R\chi_k\frac{w'}{w}\right)
   h_k\overline{g_k}\Bigg]w\,d\tau.
\end{aligned}
\end{equation}

On every interval
with $\tau_k>R-d$, the cutoff and weight bounds imply
\[
 |\chi_R\chi_k|\le1,\qquad
 \left|\chi_R'\chi_k+\chi_R\chi_k'
       +\chi_R\chi_k\frac{w'}w\right|\le C_M.
\]
It follows from Hölder's inequality with exponents $4,2,4$ and
$|\zeta_k|\le\tau_k$ that
\[
\begin{aligned}
 &\left|\int_{I_k}\chi_R\chi_k
 \frac{2\pi k}{2\pi k-\lambda\tau+i\varepsilon\tau}
 h_k\overline{g_k}\,w\,d\tau\right|\\
 &\quad\le C_M\tau_k
 \|\log(\tau-\zeta_k)\|_{L^4(I_k,w\,d\tau)}\\
 &\qquad\times\Big[
 \big(\|h_k'\|_{L^2(I_k,w\,d\tau)}
      +\|h_k\|_{L^2(I_k,w\,d\tau)}\big)
       \|g_k\|_{L^4(I_k,w\,d\tau)}
 +\|h_k\|_{L^4(I_k,w\,d\tau)}
       \|g_k'\|_{L^2(I_k,w\,d\tau)}\Big].
\end{aligned}
\]
The norms on the right use the measure $w(\tau)\,d\tau$.
By Lemma~\ref{lem:geometry} and Parseval's identity,
\begin{equation}\label{eq:tail-coefficient-norms}
 \sum_{k\in\mathbb Z}\int_{\tau_0}^{\infty}
 \bigl(|h_k'|^2+|h_k|^2\bigr)w\,d\tau
 \le C\|h\|_{H^1(\mathbb T^2)}^2.
\end{equation}
Moreover, Hölder's inequality in $\theta$ and
$dx\,dy=w(\tau)\,d\tau\,d\theta/(2\pi)$ imply
\[
\begin{aligned}
 \|h_k\|_{L^4(I_k,w\,d\tau)}^4
 &\le\frac1{2\pi}\int_{I_k}\int_0^{2\pi}
 |h(X(\tau,\theta))|^4\,d\theta\,w(\tau)\,d\tau\\
 &\le\|h\|_{L^4(\mathbb T^2)}^4.
\end{aligned}
\]
The same estimates hold for $g$.

Using \eqref{eq:log} with $p=4$, then Cauchy--Schwarz in $k$,
\eqref{eq:tail-coefficient-norms}, and
$H^1(\mathbb T^2)\hookrightarrow L^4(\mathbb T^2)$, we obtain
\begin{equation}\label{eq:tail-resonant-sum}
\begin{aligned}
 &\sum_{\substack{k\ge1\\\tau_k>R-d}}
 \left|\int_{I_k}\chi_R\chi_k
 \frac{2\pi k}{2\pi k-\lambda\tau+i\varepsilon\tau}
 h_k\overline{g_k}\,w\,d\tau\right|\\
 &\quad\le C_M
 \bigg(\sum_{\substack{k\ge1\\\tau_k>R-d}}
 \tau_k^{5/2}e^{-\tau_k/8}
 (1+\log(2+\tau_k))^2\bigg)^{1/2}\\
 &\qquad\times\Bigg[
 \|g\|_{L^4(\mathbb T^2)}
 \bigg(\sum_{k\in\mathbb Z}\int_{\tau_0}^{\infty}
 (|h_k'|^2+|h_k|^2)w\,d\tau\bigg)^{1/2}\\
 &\hspace{39mm}
 +\|h\|_{L^4(\mathbb T^2)}
 \left(\sum_{k\in\mathbb Z}\int_{\tau_0}^{\infty}
 |g_k'|^2w\,d\tau\right)^{1/2}\Bigg]\\
 &\quad\le C_M
 \bigg(\sum_{\substack{k\ge1\\\tau_k>R-d}}
 \tau_k^{5/2}e^{-\tau_k/8}
 (1+\log(2+\tau_k))^2\bigg)^{1/2}
 \|h\|_{H^1(\mathbb T^2)}\|g\|_{H^1(\mathbb T^2)}.
\end{aligned}
\end{equation}
Finally, since $\tau_{k+1}-\tau_k=2\pi/\lambda\ge2\pi/M$ and
$d\le1$,
\[
\begin{aligned}
 \sum_{\substack{k\ge1\\\tau_k>R-d}}
 \tau_k^{5/2}e^{-\tau_k/8}
 (1+\log(2+\tau_k))^2
 &\le C\sum_{\substack{k\ge1\\\tau_k>R-d}}e^{-\tau_k/16}\\
 &\le\frac{Ce^{-(R-d)/16}}{1-e^{-\pi/(8M)}}
 \le C_Me^{-R/16}.
\end{aligned}
\]
Consequently,
\[
 \sum_{\substack{k\ge1\\\tau_k>R-d}}
 \left|\int_{I_k}\chi_R\chi_k
 \frac{2\pi k}{2\pi k-\lambda\tau+i\varepsilon\tau}
 h_k\overline{g_k}\,w\,d\tau\right|
 \le C_Me^{-R/32}
 \|h\|_{H^1(\mathbb T^2)}\|g\|_{H^1(\mathbb T^2)}.
\]

Finally, inserting $\tau^m/k^m$ into \eqref{eq:tail-form} and
\eqref{eq:tail-ibp} introduces only polynomial factors, since
\[
 \left|\frac{\tau^m}{k^m}\right|
 +\left|\partial_\tau\frac{\tau^m}{k^m}\right|
 \le C_m(1+\tau)^m,
 \qquad k\ne0,\quad \tau\ge\tau_0.
\]
These factors are absorbed by the exponential decay in
\eqref{eq:tail-nonres-form} and \eqref{eq:tail-resonant-sum}.
Summing the two cell contributions and using duality proves
\eqref{eq:tail-weighted} and completes the proof.
\end{proof}

\subsubsection{Large spectral parameters and compactness}
For large $|\lambda|$, resonant modes have large angular frequency.
Angular energy then implies decay of $K(z)$ in each coordinate region as the following proposition shows.
\begin{lemma}\label{prop:hyperbolic-high-frequency}
For a fixed streamline cutoff $\chi$ supported in the hyperbolic
region, away from the centers, and every $0<\alpha<1$,
\[
 \|\Delta^{-1}T(\eps+i\lambda)\chi h\|_{H^1}
 \le C_{\alpha,\chi}(1+|\lambda|)^{-\alpha}\|h\|_{H^1},
 \qquad 0<|\eps|\le1.
\]
\end{lemma}
\begin{proof}
All unweighted Sobolev and Lebesgue norms below are on $\mathbb T^2$.
Write $z=\varepsilon+i\lambda$. Since $L_0\chi=0$, we have
$T(z)\chi=\chi T(z)$. We estimate the
pairing with $g\in H^1$ on the positive cell $\mathcal C$; the other
cell is treated identically. Choose $\tau_*>2\pi$ so that
$\chi(\tau)=0$ for $\tau\le\tau_*$, and set
\[
 d=\min\left(1,\frac{\tau_*-2\pi}{4}\right).
\]
By Lemma~\ref{lem:geometry} and smoothness of the coordinates on compact
regular annuli, for $\tau\ge\tau_*-2d$ we have
\[
 w(\tau)\sim_\chi\tau e^{-\tau/4},\qquad
 |w'|\le C_\chi w,\qquad
 |\partial_\tau X|+|\chi|+|\chi'|\le C_\chi.
\]
It follows from Parseval's identity and
$\tau^{-1}\partial_\theta=(2\pi)^{-1}L_0$ that
\[
 \sum_{k\in\mathbb Z}\int_{\tau_*-2d}^{\infty}
 \bigg[|h_k'|^2+\bigg(1+\frac{k^2}{\tau^2}\bigg)|h_k|^2\bigg]
 w\,d\tau\le C_\chi\|h\|_{H^1}^2.
\]
The same bound holds for $g$. As in Proposition~\ref{prop:tail},
Hölder's inequality in $\theta$ also implies
\[
 \|h_k\|_{L^p(I,w\,d\tau)}\le\|h\|_{L^p},\qquad
 2\le p<\infty,
 \quad I\subset(2\pi,\infty).
\]

Suppose first that $\lambda>0$. Use the decomposition
\eqref{eq:tail-form}, with $\chi$ in place of $\chi_R$, and set
\[
 \tau_k=\frac{2\pi k}{\lambda},\qquad
 I_k=(\tau_k-d,\tau_k+d),\qquad
 \zeta_k=\frac{2\pi k}{\lambda-i\varepsilon},\qquad k\ge1.
\]
Choose $\chi_k$ supported in $I_k$, equal to one on its middle half,
with $0\le\chi_k\le1$ and $|\chi_k'|\le C_\chi$; set
$\chi_k=0$ for $k<0$. Only intervals with $\tau_k>\tau_*-d$
contribute to the singular part, and each such interval lies in
$(\tau_*-2d,\infty)$. Their overlap causes no loss, since the
integral on $I_k$ involves the $k$th angular coefficient.

For $0<\lambda\le1$, the regular and singular estimates
\eqref{eq:tail-nonres-form}--\eqref{eq:tail-resonant-sum}, with the
coordinate bounds above, imply
\[
\begin{aligned}
 \left|\int_{\mathcal C}\chi T(z)h\,\overline g\,dx\,dy\right|
 &\le C_\chi\Bigg[
 \left(\int_{\tau_*}^{\infty}(1+\tau)^2|\chi|^2w\,d\tau\right)^{1/2}\\
 &\qquad+
 \bigg(\sum_{\substack{k\ge1\\\tau_k>\tau_*-d}}
 \tau_k^{5/2}e^{-\tau_k/8}
 (1+\log(2+\tau_k))^2\bigg)^{1/2}\Bigg]
 \|h\|_{H^1}\|g\|_{H^1}\\
 &\le C_\chi\|h\|_{H^1}\|g\|_{H^1}.
\end{aligned}
\]
Indeed, the integral is finite, and the sum satisfies
\[
 \sum_{\substack{k\ge1\\\tau_k>\tau_*-d}}
 \tau_k^{5/2}e^{-\tau_k/8}(1+\log(2+\tau_k))^2
 \le C\sum_{\substack{k\ge1\\\tau_k>\tau_*-d}}e^{-\tau_k/16}
 \le\frac{Ce^{-(\tau_*-d)/16}}{1-e^{-\pi/8}},
\]
since $\tau_{k+1}-\tau_k\ge2\pi$.

Now let $\lambda\ge1$. On the support of $\chi$, the regular
multiplier satisfies
\[
 \left|(1-\chi_k)
 \frac{2\pi k}{2\pi k-\lambda\tau+i\varepsilon\tau}\right|
 \le C_\chi\frac{|k|}{\lambda},\qquad k\ne0.
\]
Thus Cauchy--Schwarz in $k$ and Hölder's inequality in $w\,d\tau$
with exponents $4,2,4$ imply
\[
\begin{aligned}
 &\left|\sum_{k\ne0}\int_{2\pi}^{\infty}
 \chi(1-\chi_k)
 \frac{2\pi k}{2\pi k-\lambda\tau+i\varepsilon\tau}
 h_k\overline{g_k}\,w\,d\tau\right|\\
 &\quad\le\frac{C_\chi}{\lambda}
 \int_{\tau_*}^{\infty}\tau
 \left(\sum_k\frac{k^2}{\tau^2}|h_k|^2\right)^{1/2}
 \left(\sum_k|g_k|^2\right)^{1/2}w\,d\tau\\
 &\quad\le\frac{C_\chi}{\lambda}
 \left(\int_{\tau_*}^{\infty}\tau^4w\,d\tau\right)^{1/4}
 \|L_0h\|_2\|g\|_4
 \le\frac{C_\chi}{\lambda}\|h\|_{H^1}\|g\|_{H^1}.
\end{aligned}
\]

To estimate the singular part, choose
\[
 s=\frac{2\alpha+1}{3}\in\left(\frac13,1\right),\qquad
 \frac1q=\frac{s}{2}+\frac{1-s}{8},\qquad
 \frac1r=\frac12-\frac1q=\frac{3(1-s)}8.
\]
For $\tau\in I_k$ with $\tau_k>\tau_*-d$, we have
$k/\tau\ge c_\chi\lambda$. Interpolation and
$H^1\hookrightarrow L^8$ therefore imply
\[
\begin{aligned}
 \|g_k\|_{L^q(I_k,w\,d\tau)}
 &\le \|g_k\|_{L^2(I_k,w\,d\tau)}^s
        \|g_k\|_{L^8(I_k,w\,d\tau)}^{1-s}\\
 &\le C_\chi\lambda^{-s}
 \left(\int_{I_k}\frac{k^2}{\tau^2}|g_k|^2w\,d\tau\right)^{s/2}
 \|g\|_{H^1}^{1-s}.
\end{aligned}
\]
The same estimate holds for $h$. Moreover, the proof of
\eqref{eq:log}, using $|\zeta_k|\le\tau_k$ and the fixed radius $d$,
implies, uniformly in $\lambda$ and $\varepsilon$,
\[
 \|\log(\tau-\zeta_k)\|_{L^r(I_k,w\,d\tau)}
 \le C_{\chi,r}(\tau_ke^{-\tau_k/4})^{1/r}
 (1+\log(2+\tau_k)).
\]
Applying Hölder's inequality with exponents $r,2,q$ to
\eqref{eq:tail-ibp}, with $\chi_R$ replaced by $\chi$, yields
\[
\begin{aligned}
 &\left|\int_{I_k}\chi\chi_k
 \frac{2\pi k}{2\pi k-\lambda\tau+i\varepsilon\tau}
 h_k\overline{g_k}\,w\,d\tau\right|\\
 &\quad\le C_{\chi,r}\tau_k(\tau_ke^{-\tau_k/4})^{1/r}
 (1+\log(2+\tau_k))\\
 &\qquad\times\Big[
 (\|h_k'\|_{L^2(I_k,w\,d\tau)}+\|h_k\|_{L^2(I_k,w\,d\tau)})
 \|g_k\|_{L^q(I_k,w\,d\tau)}
 +\|h_k\|_{L^q(I_k,w\,d\tau)}\|g_k'\|_{L^2(I_k,w\,d\tau)}\Big].
\end{aligned}
\]
In view of the interpolation bounds and applying Hölder in $k$ with
exponents $2,2/s,2/(1-s)$, we obtain from the energy bounds at the beginning
of the proof that
\[
\begin{aligned}
 &\sum_{\substack{k\ge1\\\tau_k>\tau_*-d}}
 \left|\int_{I_k}\chi\chi_k
 \frac{2\pi k}{2\pi k-\lambda\tau+i\varepsilon\tau}
 h_k\overline{g_k}\,w\,d\tau\right|\\
 &\quad\le C_{\chi,s}\lambda^{-s}
 \bigg(\sum_{\substack{k\ge1\\\tau_k>\tau_*-d}}
 \left[\tau_k(\tau_ke^{-\tau_k/4})^{1/r}
 (1+\log(2+\tau_k))\right]^{2/(1-s)}\bigg)^{(1-s)/2}
 \|h\|_{H^1}\|g\|_{H^1}.
\end{aligned}
\]
\begin{samepage}
Since $2/(r(1-s))=3/4$ and $\tau_k=2\pi k/\lambda$,
\[
\begin{aligned}
 &\sum_{\substack{k\ge1\\\tau_k>\tau_*-d}}
 \left[\tau_k(\tau_ke^{-\tau_k/4})^{1/r}
 (1+\log(2+\tau_k))\right]^{2/(1-s)}\\
 &\quad=\sum_{\substack{k\ge1\\\tau_k>\tau_*-d}}
 \tau_k^{2/(1-s)+3/4}e^{-3\tau_k/16}
 (1+\log(2+\tau_k))^{2/(1-s)}\\
 &\quad\le C_s\sum_{k=1}^{\infty}e^{-\tau_k/8}
 =\frac{C_s}{e^{\pi/(4\lambda)}-1}
 \le C_s\lambda.
\end{aligned}
\]
Thus the singular part is bounded by
\[
 C_{\chi,s}\lambda^{-s+(1-s)/2}\|h\|_{H^1}\|g\|_{H^1}
 =C_{\alpha,\chi}\lambda^{-\alpha}\|h\|_{H^1}\|g\|_{H^1}.
\]
\end{samepage}

For $\lambda=0$, the multiplier has modulus at most one, and
Parseval implies the bound directly. The case of negative $\lambda$ follows by
complex conjugation. Summing over the two cells and using duality,
we obtain
\[
 \|T(\varepsilon+i\lambda)\chi h\|_{H^{-1}}
 \le C_{\alpha,\chi}(1+|\lambda|)^{-\alpha}\|h\|_{H^1}.
\]
The output has zero angular average on each streamline, hence zero
spatial mean. Applying $\Delta^{-1}:H^{-1}_0\to H^1_0$ proves the
proposition.
\end{proof}

\begin{lemma}\label{lem:regular-high-frequency}\label{lem:elliptic-high-frequency}
If $\chi$ is supported either on a compact regular annulus or near
an elliptic center, then
\[
 \|\Delta^{-1}T(\eps+i\lambda)\chi h\|_{H^1}
       \lesssim_\chi(1+|\lambda|)^{-1/2}\|h\|_{H^1}.
\]
These operators have strong boundary values for fixed $h\in H^1$.
\end{lemma}
\begin{proof}
Write $z=\varepsilon+i\lambda$ and recall $K(z)=\Delta^{-1}T(z)$.
All physical-space norms are on
$\mathbb T^2$. Choose a smooth streamline cutoff $\widetilde\chi$
equal to one on $\operatorname{supp}\chi$. Since
\[
 T(z)(\chi h)=T(z)\widetilde\chi(\chi h),\qquad
 \|\chi h\|_{H^1}\le C_\chi\|h\|_{H^1},
\]
it suffices to consider streamline cutoffs. We henceforth assume
$L_0\chi=0$. In each local argument we may also multiply $h,g$
by a fixed larger streamline cutoff equal to one on
$\operatorname{supp}\chi$, and extend their angular coefficients
by zero at the outer endpoints of the coordinate interval.

We first estimate the modes $0<|k|<|\lambda|/2$. Since $0<b\le1$,
\[
 \left|\frac{kb}{kb-\lambda+i\varepsilon}\right|
 \le\frac{2|k|b}{|\lambda|}.
\]
Thus by Parseval and Cauchy--Schwarz we get
\[
\begin{aligned}
 \left|\sum_{0<|k|<|\lambda|/2}
 \int\chi\frac{kb}{kb-\lambda+i\varepsilon}
 h_k\overline{g_k}\,dx\,dy\right|
 &\le\frac{C_\chi}{|\lambda|}
       \|L_0h\|_2\|g\|_2\\
 &\le\frac{C_\chi}{|\lambda|}\|h\|_{H^1}\|g\|_{H^1}.
\end{aligned}
\]
This sum is empty when $|\lambda|\le2$.

On a regular annulus, use $b\in I\Subset(0,1)$ as radial coordinate.
It follows from smoothness of the coordinates and the one-dimensional Sobolev
inequality that
\[
\begin{aligned}
 &\sum_{k\ne0}\int_I
  (|\partial_bh_k|^2+k^2|h_k|^2)\,db
       \le C_\chi\|h\|_{H^1}^2,\\
 &\|h_k\|_{L^\infty(I)}
       +|k|^{1/2}\|h_k\|_{L^2(I)}
 \le C|k|^{-1/2}
 \left(\int_I(|\partial_bh_k|^2+k^2|h_k|^2)\,db\right)^{1/2}.
\end{aligned}
\]
The same estimates hold for $g$. Since
$dx\,dy=b^{-1}|dr/db|\,db\,d\theta$, by integration by parts,
\[
\begin{aligned}
 &2\pi\int_I\frac{\chi(b)|dr/db|}{b-\lambda/k+i\varepsilon/k}
                     h_k\overline{g_k}\,db\\
 &\quad=-2\pi\int_I\log(b-\lambda/k+i\varepsilon/k)
 \partial_b\left(\chi(b)\left|\frac{dr}{db}\right|
                         h_k\overline{g_k}\right)db.
\end{aligned}
\]
For $|k|\ge|\lambda|/2$, $k\ne0$, the logarithm has uniformly
bounded $L^2(I)$ norm, because $|\lambda/k|\le2$ and
$|\varepsilon/k|\le1$. Consequently,
\[
\begin{aligned}
 &\left|2\pi\int_I
 \frac{\chi(b)|dr/db|}{b-\lambda/k+i\varepsilon/k}
                      h_k\overline{g_k}\,db\right|\\
 &\quad\le C_\chi\Big[
 (\|\partial_bh_k\|_{L^2(I)}+\|h_k\|_{L^2(I)})
       \|g_k\|_{L^\infty(I)}
 +\|h_k\|_{L^\infty(I)}\|\partial_bg_k\|_{L^2(I)}\Big]\\
 &\quad\le C_\chi|k|^{-1/2}
 \left(\int_I(|\partial_bh_k|^2+k^2|h_k|^2)\,db\right)^{1/2}
 \left(\int_I(|\partial_bg_k|^2+k^2|g_k|^2)\,db\right)^{1/2}.
\end{aligned}
\]
Cauchy--Schwarz in $k$ bounds the sum of these modes by
$C_\chi(1+|\lambda|)^{-1/2}\|h\|_{H^1}\|g\|_{H^1}$.

Near an elliptic center, use $q$ and $v=\log q$ from
\eqref{geom:q}, so that
\[
 b(q)=1-\frac{q^2}{8},\qquad
 dx\,dy=\xi_0(q)q\,dq\,d\theta.
\]
The Dirichlet form in Lemma~\ref{geom:metric} implies
\[
 \sum_{k\ne0}\|h_k\|_{|k|}^2\le C_\chi\|h\|_{H^1}^2,
 \qquad
 \|h_k\|_{|k|}^2
 :=\int_{\mathbb R}(|\partial_vh_k|^2+k^2|h_k|^2)\,dv,
\]
and likewise for $g$. Recalling $v=\log q$, the exact form identity is
\[
\begin{aligned}
 &\int\chi\frac{kb}{kb-\lambda+i\varepsilon}
                      h_k\overline{g_k}\,dx\,dy\\
 &\quad=-2\pi\int_{\mathbb R}
 \xi_0(q)b(q)\chi(q)
 \frac{8e^{2v}}{e^{2v}-8(1-\lambda/k)-8i\varepsilon/k}
 h_k\overline{g_k}\,dv.
\end{aligned}
\]
Since $\xi_0(q)b(q)\chi(q)$ has bounded $W^{1,\infty}$ norm as a
function of $v$, the form estimate \eqref{eq:cylinder-form} implies
\[
 \left|\int\chi\frac{kb}{kb-\lambda+i\varepsilon}
                      h_k\overline{g_k}\,dx\,dy\right|
 \le C_\chi|k|^{-1/2}\|h_k\|_{|k|}\|g_k\|_{|k|}.
\]
In particular,
\[
\begin{aligned}
 &\sum_{\substack{k\ne0\\|k|\ge|\lambda|/2}}
 \left|\int\chi\frac{kb}{kb-\lambda+i\varepsilon}
                      h_k\overline{g_k}\,dx\,dy\right|\\
 &\quad\le C_\chi(1+|\lambda|)^{-1/2}
 \left(\sum_{k\ne0}\|h_k\|_{|k|}^2\right)^{1/2}
 \left(\sum_{k\ne0}\|g_k\|_{|k|}^2\right)^{1/2}\\
 &\quad\le C_\chi(1+|\lambda|)^{-1/2}
                       \|h\|_{H^1}\|g\|_{H^1}.
\end{aligned}
\]
Combining the two ranges of $k$ proves, in both regions,
\[
 \|T(z)\chi h\|_{H^{-1}}
 \le C_\chi(1+|\lambda|)^{-1/2}\|h\|_{H^1}.
\]
The zeroth angular mode is annihilated, so the output has zero
spatial mean. Applying $\Delta^{-1}:H^{-1}_0\to H^1_0$ proves the
claimed estimate, including bounded $\lambda$.

We now prove the existence of the one-sided limits
\[
 K(i\lambda\pm0)\chi h
 :=\lim_{\varepsilon\to0^\pm}K(\varepsilon+i\lambda)\chi h
 \quad\text{in }H^1,
 \qquad K(z)=\Delta^{-1}T(z),
\]
for every fixed $\lambda\in\mathbb R$ and $h\in H^1$.
We first prove operator-norm convergence when $\chi$ is supported
on a regular annulus. In the other case when $\chi$ is supported near a center, we remove a small disk and show
that its contribution tends to zero as the radius of the small disk goes to zero, uniformly in $\varepsilon$ for
each fixed $h$.

\smallskip\noindent
\emph{Convergence on a regular annulus.}
Use $b\in I\Subset(0,1)$ as radial coordinate. For each $k\ne0$,
the principal logarithm satisfies
\[
 \log(b-\lambda/k+i\varepsilon/k)
 \longrightarrow
 \log|b-\lambda/k|
 \mathbin{\pm}i\pi\operatorname{sgn}(k)
                 \mathbf1_{\{b<\lambda/k\}}
 \quad\text{in }L^2(I)
\]
as $\varepsilon\to0^\pm$. The logarithmic singularity is square
integrable.
Convergence in operator norm for any fixed finite set of angular modes follows from integration by parts arguments as above.
The remaining modes are uniformly small: for
$N>\max(1,|\lambda|)$, their bilinear form is bounded by
\[
 \left|2\pi\sum_{|k|>N}\int_I
 \frac{\chi(b)|dr/db|}{b-\lambda/k+i\varepsilon/k}
                  h_k(b)\overline{g_k(b)}\,db\right|
 \le C_\chi N^{-1/2}\|h\|_{H^1}\|g\|_{H^1}.
\]
Each nonzero angular mode has zero spatial mean. By duality and
boundedness of $\Delta^{-1}:H^{-1}_0\to H^1_0$,
these two estimates yield, for $\varepsilon,\varepsilon'$ of the same sign,
\[
\begin{aligned}
 &\|[K(\varepsilon+i\lambda)-K(\varepsilon'+i\lambda)]\chi
                                  \|_{H^1\to H^1}\\
 &\quad\le C_{\chi,N}\max_{0<|k|\le N}
 \|\log(b-\lambda/k+i\varepsilon/k)
       -\log(b-\lambda/k+i\varepsilon'/k)\|_{L^2(I)}
       +C_\chi N^{-1/2}.
\end{aligned}
\]
First let $\varepsilon,\varepsilon'\to0^\pm$ with $N$ fixed,
and then let $N\to\infty$. This proves operator-norm convergence
on the annulus.

\smallskip\noindent
\emph{Removing a small disk about the center.}
Fix $h\in H^1$. Choose $\eta\in C_c^\infty([0,\infty))$ with
$\eta=1$ on $[0,1]$ and $\eta=0$ on $[2,\infty)$, and write
\[
 \chi h=\chi\eta(q/\rho)h+\chi[1-\eta(q/\rho)]h.
\]
For fixed $\rho>0$, the second term is supported on a regular
annulus. It remains to control the first term uniformly in the
absorption parameter $\epsilon$. Only its nonzero angular modes matter,
because $T$ annihilates the zeroth mode and the cutoffs are radial.
This distinction is essential: the full norm
$\|\eta(q/\rho)h\|_{H^1}$ need not tend to zero, even for a
smooth function equal to one near the center.

Recall that $v=\log q$ and
\[
 \|u\|_{|k|}^2
 =\int_{\mathbb R}(|\partial_vu|^2+k^2|u|^2)\,dv,
 \qquad
 \partial_v\eta(q/\rho)=\frac q\rho\eta'(q/\rho).
\]
The last expression is bounded independently of $\rho$ and
supported where $\rho<q<2\rho$. Since $|k|\ge1$, it follows from Parseval's identity and the center Dirichlet-form estimate that
\[
\begin{aligned}
 \sum_{k\ne0}\|\eta(q/\rho)h_k\|_{|k|}^2
 &\le C\sum_{k\ne0}\int_{v<\log(2\rho)}
       (|\partial_vh_k|^2+k^2|h_k|^2)\,dv\\
 &\le C\int_{q<2\rho}|\nabla h|^2\,dx\,dy
 \to0\qquad(\rho\to0).
\end{aligned}
\]
In particular, the derivative of the cutoff is controlled by the
angular energy; no factor $\rho^{-1}$ remains in this estimate.

Choose a fixed real smooth streamline cutoff $\widetilde\chi$, equal
to one on $\operatorname{supp}\chi$ and supported in $\{q<q_0\}$.
For $g\in H^1(\mathbb T^2)$, extend $\widetilde\chi g_k$ by zero
for $v\ge\log q_0$.
By Parseval's identity and the comparison of the Dirichlet forms,
\[
\begin{aligned}
 2\pi\sum_{k\ne0}\|\widetilde\chi g_k\|_{|k|}^2
 &\le \int_{\mathbb R\times\mathbb T}
 \bigl(|\partial_v(\widetilde\chi g)|^2
       +|\partial_\theta(\widetilde\chi g)|^2\bigr)\,dv\,d\theta\\
 &\le C\int_{\mathbb T^2}|\nabla(\widetilde\chi g)|^2\,dx\,dy\\
 &\le C\int_{\mathbb T^2}
 \bigl(\widetilde\chi^2|\nabla g|^2
       +|\nabla\widetilde\chi|^2|g|^2\bigr)\,dx\,dy
 \le C_\chi\|g\|_{H^1}^2.
\end{aligned}
\]
Since $\chi\widetilde\chi=\chi$ and $T$ annihilates the zeroth
angular mode, for $0<2\rho<q_0$ the exact form identity and the
cylinder form bound imply
\[
\begin{aligned}
 &\bigl|\langle T(\varepsilon+i\lambda)
                    \chi\eta(q/\rho)h,g\rangle\bigr|\\
 &\quad=2\pi\left|\sum_{k\ne0}\int_{\mathbb R}
 \xi_0(q)b(q)\chi(q)
 \frac{8e^{2v}}{e^{2v}-8(1-\lambda/k)-8i\varepsilon/k}
 \eta(q/\rho)h_k\overline{\widetilde\chi g_k}\,dv\right|\\
 &\quad\le C_\chi\sum_{k\ne0}|k|^{-1/2}
             \|\eta(q/\rho)h_k\|_{|k|}
             \|\widetilde\chi g_k\|_{|k|}\\
 &\quad\le C_\chi
 \left(\sum_{k\ne0}\|\eta(q/\rho)h_k\|_{|k|}^2\right)^{1/2}
 \left(\sum_{k\ne0}\|\widetilde\chi g_k\|_{|k|}^2\right)^{1/2}\\
 &\quad\le C_\chi
 \left(\int_{q<2\rho}|\nabla h|^2\,dx\,dy\right)^{1/2}
 \|g\|_{H^1}.
\end{aligned}
\]
Here $q=e^v$ in the integral. The coefficient $\xi_0b\chi$
and the test-function cutoff are fixed, so the constant is independent
of $\rho,\lambda,\varepsilon$. The output has zero spatial mean;
hence $K(z)=\Delta^{-1}T(z)$ satisfies
\[
\begin{aligned}
 \|K(\varepsilon+i\lambda)\chi\eta(q/\rho)h\|_{H^1}&\quad\le C\|T(\varepsilon+i\lambda)\chi\eta(q/\rho)h\|_{H^{-1}}\\
 &\quad=C\sup_{\|g\|_{H^1}\le1}
 \bigl|\langle T(\varepsilon+i\lambda)\chi\eta(q/\rho)h,g\rangle\bigr|\\
 &\quad\le C_\chi
 \left(\int_{q<2\rho}|\nabla h|^2\,dx\,dy\right)^{1/2},
 \qquad \lambda\in\mathbb R,\quad 0<|\varepsilon|\le1.
\end{aligned}
\]
For fixed $h\in H^1$, the right-hand side tends to zero as
$\rho\to0$, uniformly in $\lambda$ and $\varepsilon$.

For fixed $\rho$, the annular part converges by the first step.
The preceding uniform estimate therefore implies
\[
\begin{aligned}
 &\limsup_{\varepsilon,\varepsilon'\to0^\pm}
 \|[K(\varepsilon+i\lambda)-K(\varepsilon'+i\lambda)]\chi h\|_{H^1}\le C_\chi
 \left(\int_{q<2\rho}|\nabla h|^2\,dx\,dy\right)^{1/2}.
\end{aligned}
\]
Letting $\rho\to0$ proves the Cauchy property in $H^1$ and hence
the asserted boundary values.
\end{proof}

We will use the following notation for the full resolvent:
\[
 \psi_f(z)=\Delta^{-1}(z-L)^{-1}f,\qquad
 \psi_{f,\eps}^{\iota}(\lambda)=\psi_f(\iota\eps+i\lambda),
 \quad \eps>0,\quad\iota\in\{+,-\}.
\]
Here $\iota\eps=\pm\eps$, and we denote the boundary values by
$\psi_{f,+}$ and $\psi_{f,-}$. A subscript $k$ denotes an angular
coefficient of the full solution, whose modes are coupled by $\Delta$.

\begin{corollary}\label{cor:global-high-frequency}
For sufficiently large $|\lambda|$ and $0<|\varepsilon|\le1$,
$\|K(\eps+i\lambda)\|_{H^1\to H^1}\lesssim|\lambda|^{-1/2}$.
Consequently $I+K$ is uniformly invertible there, with strong boundary
values, and for mean-zero $f\in H^3$,
\begin{align}
 \|\Delta^{-1}(\eps+i\lambda-L)^{-1}f\|_{H^1}
 &\lesssim |\lambda|^{-1}\|f\|_{H^1},\label{eq:large-resolvent}\\
 \|\Delta^{-1}[(i\lambda+0-L)^{-1}-(i\lambda-0-L)^{-1}]f\|_{H^1}
 &\lesssim |\lambda|^{-3}\|f\|_{H^3}.\label{eq:large-jump}
\end{align}
\end{corollary}
\begin{proof}
Write $z=\varepsilon+i\lambda$, with $0<|\varepsilon|\le1$.
Recall that $L=L_0(I+\Delta^{-1})$. Thus the resolvent equation
$(z-L)\omega=f$, with $\psi=\Delta^{-1}\omega$, takes the form
\[
 (z-L_0)\Delta\psi-L_0\psi=f.
\]
Applying $(L_0-z)^{-1}$ we obtain
\[
 \Delta\psi+T(z)\psi=-(L_0-z)^{-1}f.
\]
Since $T(z)-I=z(L_0-z)^{-1}$, taking $\Delta^{-1}$ we get
\[
 (I+K(z))\psi
 =-\Delta^{-1}(L_0-z)^{-1}f
 =z^{-1}\bigl(\Delta^{-1}-K(z)\bigr)f.
\]

By Lemma~\ref{prop:hyperbolic-high-frequency} and Lemma~\ref{lem:regular-high-frequency}, we obtain that
\[
 \|K(z)\|_{H^1\to H^1}\le C|\lambda|^{-1/2}.
\]
For sufficiently large $|\lambda|$, the right-hand side of the above inequality is at
most $1/2$, so the Neumann series implies
\[
 \|(I+K(z))^{-1}\|_{H^1\to H^1}\le2.
\]
The strong boundary convergence of $K(z)$ and the identity
\[
 (I+K(z))^{-1}-(I+K(w))^{-1}
 =(I+K(z))^{-1}\bigl(K(w)-K(z)\bigr)(I+K(w))^{-1}
\]
imply the existence of strong boundary values for the inverse. Moreover,
\[
 \|\psi\|_{H^1}
 \le \frac{2}{|z|}
 \left(\|\Delta^{-1}f\|_{H^1}+\|K(z)f\|_{H^1}\right)
 \le C|\lambda|^{-1}\|f\|_{H^1},
\]
which proves \eqref{eq:large-resolvent}.

We finally prove the jump estimate. Let $f\in H^3$ have zero spatial
mean, and let $z=\varepsilon+i\lambda$,
$0<|\varepsilon|\le1$, and $|\lambda|$ is sufficiently large.
Applying the resolvent identity first to $f$ and then to $Lf$, we obtain
\[
\begin{aligned}
 (z-L)^{-1}f
 &=z^{-1}f+z^{-1}(z-L)^{-1}Lf\\
 &=z^{-1}f+z^{-2}Lf+z^{-2}(z-L)^{-1}L^2f.
\end{aligned}
\]
Thus $\psi_f(z)=\Delta^{-1}(z-L)^{-1}f$ satisfies
\[
 \psi_f(z)=z^{-1}\Delta^{-1}f
          +z^{-2}\Delta^{-1}Lf+z^{-2}\psi_{L^2f}(z).
\]
The first two terms have the same limits as
$\varepsilon\to0^+$ and $\varepsilon\to0^-$, since $\lambda\ne0$.
Taking the strong $H^1$ boundary values and subtracting, we obtain
\[
 \psi_{f,+}(\lambda)-\psi_{f,-}(\lambda)
 =(i\lambda)^{-2}
 \bigl(\psi_{L^2f,+}(\lambda)-\psi_{L^2f,-}(\lambda)\bigr).
\]
Only boundary convergence of the stream functions in $H^1$ is used here.
Since $L=L_0(I+\Delta^{-1})$ and $L_0$ is a divergence-free transport
operator with smooth coefficients, $L^2f$ has zero spatial mean and
\[
 \|L^2f\|_{H^1}\le C\|Lf\|_{H^2}\le C\|f\|_{H^3}.
\]
The resolvent estimate \eqref{eq:large-resolvent} holds for every mean-zero $H^1$
source. Applying it to $L^2f$,
we conclude that
\[
\begin{aligned}
 &\|\psi_{f,+}(\lambda)-\psi_{f,-}(\lambda)\|_{H^1}\\
 &\quad\le |\lambda|^{-2}
 \bigl(\|\psi_{L^2f,+}(\lambda)\|_{H^1}
      +\|\psi_{L^2f,-}(\lambda)\|_{H^1}\bigr)\\
 &\quad\le C|\lambda|^{-3}\|L^2f\|_{H^1}
 \le C|\lambda|^{-3}\|f\|_{H^3}.
\end{aligned}
\]

\end{proof}

For bounded spectral parameters away from the integers we have
compactness in operator norm.
\begin{proposition}\label{prop:compact}
For compact $I\subset\R\setminus\Z$, $K(\eps+i\lambda)$ extends
to two norm-continuous compact families on mean-zero $H^1$ for
$\lambda\in I$, $\eps\to0^\pm$.
\end{proposition}
\begin{proof}
Fix a compact set $I\subset\mathbb R\setminus\mathbb Z$ and write
$z=\varepsilon+i\lambda$, with $\lambda\in I$ and
$0<|\varepsilon|\le1$. Recall that
\[
 K(z)=\Delta^{-1}T(z),\qquad
 T(z)h=\sum_{k\ne0}\frac{kb(r)}{kb(r)-\lambda+i\varepsilon}
                         h_k(r)e^{ik\theta},
 \qquad b(r)=\frac{2\pi}{S(r)}.
\]

\smallskip\noindent
\emph{Step 1: Removing the separatrix and the large angular modes.}
Choose smooth functions $0\le\chi_R\le1$ such that
\[
 \chi_R(s)=0\quad(s\le R),\qquad
 \chi_R(s)=1\quad(s\ge R+1),\qquad
 \sup_R\|\chi_R'\|_{L^\infty}<\infty.
\]
For sufficiently large $R$, Proposition~\ref{prop:tail} and the
boundedness of $\Delta^{-1}:H^{-1}_0\to H^1_0$ imply
\[
 \sup_{\substack{\lambda\in I\\0<|\varepsilon|\le1}}
 \|\Delta^{-1}\chi_R(S)T(\varepsilon+i\lambda)\|_{H^1\to H^1}
 \le C_Ie^{-cR}.
\]
On the support of $1-\chi_R(S)$, we have $b\ge2\pi/(R+1)$.
Choose an integer
\[
 N\ge1+\frac{R+1}{\pi}\max_{\lambda\in I}|\lambda|.
\]
Then, for $|k|>N$,
$|kb-\lambda+i\varepsilon|\ge |kb|/2$ and the multiplier has
absolute value at most $2$. Since
$\partial_\theta=b^{-1}L_0$, Parseval's identity yields
\[
\begin{aligned}
 &\left\|(1-\chi_R(S))\sum_{|k|>N}
 \frac{kb}{kb-\lambda+i\varepsilon}h_ke^{ik\theta}\right\|_{L^2}\\
 &\quad\le\frac2N\|(1-\chi_R(S))\partial_\theta h\|_{L^2}
 \le\frac{R+1}{\pi N}\|L_0h\|_{L^2}
 \le C\frac{R+1}{N}\|h\|_{H^1}.
\end{aligned}
\]
Define the truncated operator by
\[
 K_{R,N}(z)h
 =\Delta^{-1}\left[(1-\chi_R(S))\sum_{0<|k|\le N}
 \frac{kb}{kb-\lambda+i\varepsilon}h_ke^{ik\theta}\right].
\]
The preceding estimates imply
\begin{equation}\label{compact:explicit-approximation}
 \sup_{\substack{\lambda\in I\\0<|\varepsilon|\le1}}
 \|K(\varepsilon+i\lambda)-K_{R,N}(\varepsilon+i\lambda)\|_{H^1\to H^1}
 \le C_Ie^{-cR}+C\frac{R+1}{N}.
\end{equation}

\pagebreak[3]
\smallskip\noindent
\emph{Step 2: The finitely many remaining modes.}
Fix $R,N$. Near either center we use the radius $q$ defined by
$b=1-q^2/8$. Choose $q_0>0$ sufficiently small such that
\[
 \frac{Nq_0^2}{8}\le\frac12\operatorname{dist}(I,\mathbb Z),
\]
and a smooth streamline cutoff $0\le\beta(q)\le1$, equal to one
on $q\le q_0/2$ around both centers and zero outside $q<q_0$.
For $q\le q_0$ and $0<|k|\le N$,
\[
 |kb-\lambda|
 \ge |k-\lambda|-\frac{|k|q^2}{8}
 \ge\frac12\operatorname{dist}(I,\mathbb Z).
\]
Consequently, for $\lambda,\lambda'\in I$ and
$|\varepsilon|,|\varepsilon'|\le1$,
\[
\begin{aligned}
 \left|\frac{kb}{kb-\lambda+i\varepsilon}\right|
 &\le\frac{2N}{\operatorname{dist}(I,\mathbb Z)},\\
 \left|\frac{kb}{kb-\lambda+i\varepsilon}
       -\frac{kb}{kb-\lambda'+i\varepsilon'}\right|
 &\le\frac{4N(|\lambda-\lambda'|+|\varepsilon-\varepsilon'|)}
              {\operatorname{dist}(I,\mathbb Z)^2}.
\end{aligned}
\]
By Parseval and elliptic regularity,
\begin{equation}\label{jc0.1}
 \left\|\Delta^{-1}\left[\beta(1-\chi_R(S))
 \sum_{0<|k|\le N}\frac{kb}{kb-\lambda+i\varepsilon}
 h_ke^{ik\theta}\right]\right\|_{H^2}
 \le\frac{CN}{\operatorname{dist}(I,\mathbb Z)}\|h\|_{L^2}.
\end{equation}
The multiplier difference implies the corresponding Lipschitz bound
in $L^2\to H^2$. Thus the center contribution is uniformly
bounded and norm-continuous up to $\varepsilon=0$.

The cutoff $(1-\chi_R(S))(1-\beta)$ is supported on two regular
annuli. We work on one annulus, denote the restricted cutoff by
$\chi$, and sum the estimates over the two cells. Here $b$ is a smooth
radial coordinate in a compact subinterval of $(0,1)$; the angular
coefficients below are regarded as functions of $b$.
We use the following explicit
kernel estimates. With Fourier transform
$\mathcal{F}(u)(\xi):=\widehat u(\xi)=\int_{\mathbb R}e^{-ib\xi}u(b)\,db$, for
$a\in\mathbb R$ and $\eta\ne0$,
\[
 \mathcal{F}\Big(\frac1{b-a+i\eta}\Big)(\xi)
 =-2\pi i\operatorname{sgn}(\eta)\mathbf1_{\{\eta\xi>0\}}
                         e^{-ia\xi-|\eta||\xi|}.
\]
In particular,
\[
 \left\|\frac1{b-a+i\eta}\right\|_{H^{-3/4}(\mathbb R)}^2
 \le C\int_{\mathbb R}(1+\xi^2)^{-3/4}\,d\xi\le C.
\]
For $\eta,\eta'$ of the same sign, the same formula implies
\[
\begin{aligned}\label{j0.1}
 &\left\|\frac1{b-a+i\eta}-\frac1{b-a'+i\eta'}
                                      \right\|_{H^{-3/4}(\mathbb R)}^2\\
 &\quad\le C\int_{\mathbb R}(1+\xi^2)^{-3/4}
 \min\!\left\{1,(|a-a'|+|\eta-\eta'|)^2\xi^2\right\}\,d\xi\\
 &\quad\le C(|a-a'|+|\eta-\eta'|)^{1/2}.
\end{aligned}
\]

Choose a fixed smooth cutoff $\zeta(b)$ equal to one on
$\operatorname{supp}\chi$ and supported in a larger regular annulus.
Extend $\zeta h_k$ and $\zeta g_k$ by zero to the real $b$-axis.
It follows that
\[
 \|\zeta h_k\|_{H^1(\mathbb R)}\le C_\chi\|h\|_{H^1},\qquad
 \|\zeta g_k\|_{H^{3/4}(\mathbb R)}\le C_\chi\|g\|_{H^{3/4}}.
\]
Since $H^{3/4}(\mathbb R)$ is an algebra,
\[
 \left\|\chi\left|\frac{dr}{db}\right|h_k\overline{g_k}
                                               \right\|_{H^{3/4}(\mathbb R)}
 \le C_\chi\|\zeta h_k\|_{H^1(\mathbb R)}
                 \|\zeta g_k\|_{H^{3/4}(\mathbb R)}
 \le C_\chi\|h\|_{H^1}\|g\|_{H^{3/4}}.
\]
Using $dx\,dy=b^{-1}|dr/db|\,db\,d\theta$, duality with the
Cauchy kernel we obtain
\[
\begin{aligned}
 &\bigg|\bigg\langle\chi\sum_{0<|k|\le N}
 \frac{kb}{kb-\lambda+i\varepsilon}h_ke^{ik\theta},g
 \bigg\rangle\bigg|=2\pi\bigg|\sum_{0<|k|\le N}\int_{\mathbb R}
 \frac{\chi(b)|dr/db|}{b-\lambda/k+i\varepsilon/k}
 h_k\overline{g_k}\,db\bigg|\le C_{\chi,N}\|h\|_{H^1}\|g\|_{H^{3/4}}.
\end{aligned}
\]
Applying \eqref{j0.1} we get,
for $\varepsilon,\varepsilon'$ on the same side of zero,
\[
\begin{aligned}
 &\left\|\chi\sum_{0<|k|\le N}
 \left(\frac{kb}{kb-\lambda+i\varepsilon}
       -\frac{kb}{kb-\lambda'+i\varepsilon'}\right)
 h_ke^{ik\theta}\right\|_{H^{-3/4}}\le C_{\chi,N}
 (|\lambda-\lambda'|+|\varepsilon-\varepsilon'|)^{1/4}
 \|h\|_{H^1}.
\end{aligned}
\]
The case of $k<0$ can be treated identically.
Thus the annular contribution maps $H^1$ into $H^{-3/4}$,
continuously in operator norm up to either boundary.

\smallskip\noindent
\emph{Step 3: Compactness and passage to the full operator.}
Applying $\Delta^{-1}:H^{-3/4}\to H^{5/4}$ over zero-mean functions and adding the
center contribution \eqref{jc0.1}, we obtain on each closed parameter set
$\lambda\in I$, $0\le\pm\varepsilon\le1$,
\[
\begin{aligned}
 \|K_{R,N}(\varepsilon+i\lambda)\|_{H^1\to H^{5/4}}
 &\le C_{I,R,N},\\
 \|K_{R,N}(\varepsilon+i\lambda)
      -K_{R,N}(\varepsilon'+i\lambda')\|_{H^1\to H^{5/4}}
 &\le C_{I,R,N}
 (|\lambda-\lambda'|+|\varepsilon-\varepsilon'|)^{1/4}.
\end{aligned}
\]
For fixed $R,N$, the image of the $H^1$ unit ball under $K(z)$ in the prescribed parameter range is bounded in $H^{5/4}$, hence relatively compact in $H^1$.
By \eqref{compact:explicit-approximation}, first choosing $R$ large
and then $N$ large makes the approximation error arbitrarily small.
The continuous extensions of the approximants are uniformly Cauchy
on either closed parameter range. Their uniform limit
therefore defines a norm-continuous extension of $K$.
The compactness of $K(z)$ then follows.

\end{proof}

\subsubsection{Bounds in $W^{1,4}$ away from the elliptic points}
The spectral derivative estimates below involve products of the form $DV\,\nabla h$, where $h$ is at the level of regularity of the spectral density of the stream function;
since for fixed $\lambda\ne0$ Lemma~\ref{geom:common-field} shows
$DV\in L^4(\Omega)$ for every region $\Omega$ whose closure avoids
the elliptic centers, we establish exterior $W^{1,4}$ bounds for $h$
in order to apply
\[
 \|DV\,\nabla h\|_{L^2(\Omega)}
 \le \|DV\|_{L^4(\Omega)}\|\nabla h\|_{L^4(\Omega)}.
\]

\label{sec:exterior-W14}
For an open set $U\subset\T^2$, define
\[
 \mathcal B(U)=
 \{\operatorname{div}F+f:
 F\in L^4(U;\mathbb C^2),\ f\in L^{4/3}(U)\},
\]
and set
\[
 \|u\|_{\mathcal B(U)}
 =\inf_{u=\operatorname{div}F+f}
       \bigl(\|F\|_{L^4(U)}+\|f\|_{L^{4/3}(U)}\bigr).
\]
All these norms and derivatives are in physical coordinates.
We write $\mathcal B=\mathcal B(\T^2)$. By Sobolev embedding,
\begin{equation}\label{hg:B-elliptic}
 \|u\|_{W^{1,4}(\T^2)}
 \lesssim\|\Delta u\|_{\mathcal B}+\|u\|_2.
\end{equation}

Fix $c>0$ small enough so that $q<4c$ lies in the two center charts.
For $\lambda\ge1$, set $q_0=c\lambda^{-1/2}$ and define the exterior
region by
\[
 \Omega_\lambda
 =\T^2\setminus
   \{q\le2q_0\text{ in either center chart}\}.
\]
Choose a smooth streamline cutoff $\chi_\lambda$ which vanishes
for $q\le q_0$, equals one for $q\ge2q_0$ and outside the center
charts, and satisfies
$|\nabla\chi_\lambda|\lesssim q_0^{-1}$ and
$|\Delta\chi_\lambda|\lesssim q_0^{-2}$.
The derivatives are supported on annular region with area $O(q_0^2)$.

We allow a fixed streamline weight $w_0=w_0(r)$, smooth up to the
centers on each side of $r=2$, such that, in the coordinate
$\tau=S(r)$ near the separatrix,
\begin{equation}\label{hg:weight-class}
 |w_0|+|\partial_\tau w_0|\le C(1+\tau)^m
\end{equation}
for some fixed $m\ge0$. Constants below may depend on $c$ and
this weight, but not on $\lambda$. Besides $w_0=1$, this includes
the weight
\[
 w_0(r)=\frac{S(r)S''(r)}{(S'(r))^2}-2
       =\lambda\operatorname{div}V
\]
from Lemma~\ref{geom:common-field}. The period estimates imply
that this weight is smooth at the centers and satisfies
\eqref{hg:weight-class}. Polynomial factors in $S$ are also allowed.

\begin{lemma}\label{lem:exterior-W14}
Let $\lambda\ge1$ and let $T=T(i\lambda+0)$ or $T(i\lambda-0)$.
For every $h\in H^1(\T^2)$,
\begin{equation}\label{hg:exterior-B}
 \|\chi_\lambda w_0Th\|_{\mathcal B}
 \le C\lambda\|h\|_{H^1}.
\end{equation}
In particular, the same bound holds for
$\|w_0Th\|_{\mathcal B(\Omega_\lambda)}$.
If $(\Delta+T)h=H\in L^2$, then
\begin{equation}\label{hg:exterior-W14}
 \|h\|_{W^{1,4}(\Omega_\lambda)}
 \le C\bigl(\lambda\|h\|_{H^1}+\|H\|_2\bigr).
\end{equation}
If $L_0^jg\in H^1$ for $0\le j\le2$, then
\begin{equation}\label{hg:T-improved}
 \|\chi_\lambda w_0Tg\|_{\mathcal B}
 \le C\lambda^{-1}\sum_{j=0}^2\|L_0^jg\|_{H^1}.
\end{equation}
Consequently, if also $(\Delta+T)g=H\in L^2$, then on any
fixed exterior region $\Omega_*$ whose closure does not intersect with the centers,
\begin{equation}\label{hg:fixed-exterior-W14}
 \|g\|_{W^{1,4}(\Omega_*)}
 \le C_{\Omega_*}\left(\|H\|_2+
                       \sum_{j=0}^2\|L_0^jg\|_{H^1}\right).
\end{equation}
The analogous estimates hold for $\lambda\le-1$ by conjugation.
\end{lemma}

\begin{proof}
The main point is to write each singular multiplier as a derivative of a
logarithm, as in~\cite{JiaSIMA}. Integrating by parts once, we obtain a decomposition
$\operatorname{div}F+f$, with $F\in L^4$ and $f\in L^{4/3}$.
We estimate these terms on regular annuli, near the centers, and near the
separatrix; the remaining conclusions follow from the elliptic equation.

It suffices first to prove, for smooth $h$,
\begin{equation}\label{w14:basic}
 \|\chi_\lambda w_0Th\|_{\mathcal B}\le C\lambda\|h\|_{H^1}.
\end{equation}
We use a fixed smooth streamline partition, by smooth cutoff functions
$\zeta$, into regular annuli, the regions $q<3c$ around the two centers (with $0<c\ll1$),
and $\tau=S(r)>\tau_0+1$. Here $\tau_0$ is fixed and sufficiently large;
outside the center regions the cutoffs vanish on $q\le2c$, so
$\chi_\lambda=1$. Fix $\eta\in C_c^\infty((-1,1))$, equal to one
on $[-1/2,1/2]$, with $0\le\eta\le1$.

We record the integration by parts used below. For a radial coordinate $\varrho$,
let $X_{\varrho}=(\partial_{\varrho}X)\circ X^{-1}$ be the physical vector field induced
by the coordinate map $X(\varrho,\theta)$. For a smooth radial coefficient $A$
supported inside this coordinate neighborhood,
\begin{equation}\label{w14:ibp}
\begin{aligned}
 \frac{Ah_ke^{ik\theta}}{\varrho-a\pm i0}
 &=\operatorname{div}\!\left(A\log(\varrho-a\pm i0)h_ke^{ik\theta}X_{\varrho}\right)\\
 &\quad-\log(\varrho-a\pm i0)e^{ik\theta}
             \left[\operatorname{div}(AX_{\varrho})h_k+A\partial_{\varrho}h_k\right].
\end{aligned}
\end{equation}
Here $\log(\varrho-a\pm i0)=\log|\varrho-a|\pm i\pi\mathbf1_{\{\varrho<a\}}$,
and the terms are extended by zero outside the support of $A$. Since the area density is radial,
$\|h_k\|_{L^p(U)}\le\|h\|_{L^p(U)}$ on every region $U$ consisting of
complete streamlines. It follows from Hölder's inequality and the Sobolev embedding $H^1\hookrightarrow L^8\cap L^4$ that
\begin{equation}\label{w14:holder}
\begin{aligned}
 \left\|\frac{Ah_ke^{ik\theta}}{\varrho-a\pm i0}\right\|_{\mathcal B}
 &\le C\left(\|AX_{\varrho}\log(\varrho-a\pm i0)\|_8
       +\|\operatorname{div}(AX_{\varrho})\log(\varrho-a\pm i0)\|_2\right)\|h\|_{H^1}\\
 &\quad+\|A\log(\varrho-a\pm i0)\|_4
                    \|\partial_{\varrho}h_k\|_{L^2(\operatorname{supp}A)}.
\end{aligned}
\end{equation}

\smallskip\noindent
\emph{Step 1: regular annuli.}
Choose $b_*>0$ so that $\operatorname{supp}\zeta\subset\{b_*\le b\le1-b_*\}$.
For $k>0$ with $\lambda/k\in[b_*/2,1-b_*/2]$, write
\[
 \zeta w_0(Th)_ke^{ik\theta}
 =\frac{\zeta w_0b\,h_ke^{ik\theta}}{b-\lambda/k+i0}.
\]
In \eqref{w14:holder}, take $\varrho=b$, $a=\lambda/k$, and $A=\zeta w_0b$.
The coordinate coefficients and their first derivatives are bounded, and
\[
 \|\log(b-\lambda/k+i0)\|_{L^p(\operatorname{supp}\zeta)}\le C_p,
 \qquad \|\partial_bh_k\|_{L^2(\operatorname{supp}\zeta)}\le C\|\nabla h\|_2.
\]
Each term therefore has $\mathcal B$ norm at most $C\|h\|_{H^1}$.
There are at most $C_{b_*}\lambda$ such integers. For all remaining modes,
$|kb/(kb-\lambda)|\le C_{b_*}$ on $\operatorname{supp}\zeta$,
so Parseval bounds their sum in $L^2$ by $C\|h\|_2$. This proves
\eqref{w14:basic} on the regular annuli.

\smallskip\noindent
\emph{Step 2: elliptic centers.}
Near either center, $b(q)=1-q^2/8$. For $k>\lambda$, set
\[
 \varrho_{k}^2=8(1-\lambda/k),\qquad v=\log(q/\varrho_{k}),\qquad
 A_k(v)=-\frac{8}{\varrho_{k}^2}\zeta(q)\chi_\lambda(q)w_0(q)\eta(4v)b(q)
                         \frac{v}{e^{2v}-1}.
\]
The exact decomposition is
\begin{equation}\label{w14:center}
 \zeta\chi_\lambda w_0Th
 =\sum_{k\in\mathcal J_\lambda}
       \frac{A_k(v)h_k(q)e^{ik\theta}}{v-i0}
   +\zeta\chi_\lambda w_0\sum_k r_k(q)h_k(q)e^{ik\theta},
\end{equation}
where
\[
\begin{gathered}
 \mathcal J_\lambda
 =\{k\in\mathbb Z:k>\lambda,\ e^{-1/4}q_0<\varrho_{k}<3e^{1/4}c\},\\
 r_k(q)=
 \begin{cases}
  \displaystyle\frac{kb(q)}{kb(q)-\lambda},&k\le\lambda,\\[5pt]
  \displaystyle[1-\eta(4\log(q/\varrho_{k}))]
                         \frac{kb(q)}{kb(q)-\lambda},&k>\lambda .
 \end{cases}
\end{gathered}
\]
Indeed, $kb/(kb-\lambda+i0)=-8b/[\varrho_{k}^2(e^{2v}-1-i0)]$.
The regular multipliers satisfy $|r_k(q)|\le Cq^{-2}$:
for $0<k\le\lambda$, $|kb-\lambda|=(\lambda-k)+kq^2/8$;
for $k>\lambda$, the cutoff removes $|q^2-\varrho_{k}^2|<C^{-1}q^2$.
This includes $k=\lambda$ when $\lambda$ is an integer. Hence the
second term of \eqref{w14:center} has $L^{4/3}$ norm at most
$Cq_0^{-2}\|h\|_2\le C\lambda\|h\|_2$.

On $U_k=\{q<3c:|\log(q/\varrho_{k})|<1/4\}$, Lemma~\ref{geom:metric} implies
\[
\begin{gathered}
 dx\,dy=\xi_0(q)q\,dq\,d\theta,\qquad
 |X_v|\le C\varrho_{k},\qquad
 |\operatorname{div}X_v|=\left|2+\frac{q\xi_0'}{\xi_0}\right|\le C,\\
 |A_k|+|\partial_vA_k|\le C\varrho_{k}^{-2},\qquad
 \|\log(v-i0)\|_{L^p(U_k)}\le C_p \varrho_{k}^{2/p}.
\end{gathered}
\]
The last estimate follows by substituting $q=\varrho_{k}e^v$.
Since $\partial_v=q\partial_q$, \eqref{w14:holder} implies
\[
 \left\|\frac{A_kh_ke^{ik\theta}}{v-i0}\right\|_{\mathcal B}
 \le C\varrho_{k}^{-1}\|h\|_{H^1}
       +C\varrho_{k}^{-1/2}\|\partial_qh_k\|_{L^2(U_k)}.
\]
To sum, note that $k\le C\lambda$ on $\mathcal J_\lambda$ and
\[
 \varrho_{k}^2\ge\frac{e^{-1/2}c^2}{\lambda},\qquad
 \varrho_{k+1}^2-\varrho_{k}^2=\frac{8\lambda}{k(k+1)}\ge\frac1{C\lambda}.
\]
Thus, enumerating these roots increasingly,
$\varrho_{k_n}^2\ge n/(C_c\lambda)$, and
\[
\begin{aligned}
 \sum_{k\in\mathcal J_\lambda}\varrho_{k}^{-1}
 &\le C_c\sqrt\lambda\sum_{n\le C\lambda}n^{-1/2}\le C_c\lambda,\\
 \sum_{k\in\mathcal J_\lambda}\varrho_{k}^{-1/2}
                         \|\partial_qh_k\|_{L^2(U_k)}
 &\le\left(\sum_{k\in\mathcal J_\lambda}\varrho_{k}^{-1}\right)^{1/2}
       \left(\sum_k\|\partial_qh_k\|_{L^2(q<3c)}^2\right)^{1/2}
 \le C_c\sqrt\lambda\,\|\nabla h\|_2.
\end{aligned}
\]
The last inequality uses Parseval and the boundedness of $X_q$.
Together with \eqref{w14:center}, this proves the center estimate.

\smallskip\noindent
\emph{Step 3: the separatrix.}
In $\tau>\tau_0+1$, it follows from Lemma~\ref{lem:geometry} and \eqref{jxtb1} that
\[
 dx\,dy=\frac{w(\tau)}{2\pi}\,d\tau\,d\theta,\qquad
 w(\tau)\sim\tau e^{-\tau/4},\qquad
 |X_\tau|+|\operatorname{div}X_\tau|\le C.
\]
For $k>0$, put $\tau_k=2\pi k/\lambda$ and
$A_k(\tau)=-\tau_k\zeta(\tau)w_0(\tau)\eta(\tau-\tau_k)$. Then
\begin{equation}\label{w14:separatrix}
 \zeta w_0Th
 =\sum_{k>0}\frac{A_k(\tau)h_ke^{ik\theta}}{\tau-\tau_k-i0}
      +\zeta w_0\sum_k r_k(\tau)h_ke^{ik\theta},
\end{equation}
where
\[
 r_k(\tau)=
 \begin{cases}
  \displaystyle-\frac{[1-\eta(\tau-\tau_k)]\tau_k}{\tau-\tau_k},
                   &k>0,\\[5pt]
  \displaystyle\frac{2\pi k}{2\pi k-\lambda\tau},&k\le0.
 \end{cases}
\]
Only $\tau_k>\tau_0$ contribute to the first sum. On
$U_k=\{\tau>\tau_0:|\tau-\tau_k|<1\}$,
\[
 |A_k|+|A_k'|\le C(1+\tau_k)^{m+1},\qquad
 \|\log(\tau-\tau_k-i0)\|_{L^p(U_k)}
 \le C_p(\tau_k e^{-\tau_k/4})^{1/p}.
\]
It follows from \eqref{w14:holder} that
\[
 \left\|\frac{A_kh_ke^{ik\theta}}{\tau-\tau_k-i0}\right\|_{\mathcal B}
 \le C(1+\tau_k)^{m+2}e^{-\tau_k/32}\|h\|_{H^1}.
\]
Since $\#\{k:n\le\tau_k<n+1\}\le1+\lambda/(2\pi)\le C\lambda$,
the sum is bounded by $C\lambda\|h\|_{H^1}$.
For the remainder, $|r_k|\le1+2\tau$, and Parseval and Hölder imply
\[
 \left\|\zeta w_0\sum_k r_kh_ke^{ik\theta}\right\|_{4/3}
 \le C\|(1+\tau)^{m+1}\|_{L^4(\tau>\tau_0)}\|h\|_2
 \le C\|h\|_2.
\]
Here every polynomial in $\tau$ is integrable against $w(\tau)\,d\tau$.

These estimates also justify extension across the separatrix.
For $R\ge\tau_0+2$, choose $\nu_R=0$ on $\tau\le R$ and $\nu_R=1$
on $\tau\ge R+1$, with $0\le\nu_R\le1$ and $|\nu_R'|\le C$.
Applying the preceding estimates with $\zeta$ replaced by $\nu_R\zeta$ we obtain
\[
\begin{aligned}
 \|\nu_R\zeta w_0Th\|_{\mathcal B}
 &\le C\|h\|_{H^1}
       \sum_{\tau_k>R-1}(1+\tau_k)^{m+2}e^{-\tau_k/32}
       +C\|(1+\tau)^{m+1}\|_{L^4(\tau>R)}\|h\|_2\\
 &\le C\lambda e^{-c_1R}\|h\|_{H^1}
\end{aligned}
\]
for a fixed $c_1>0$. Each individual divergence in \eqref{w14:ibp}
is supported away from the separatrix; the estimate proves convergence
of their zero extensions and of the truncations in $\mathcal B$.
Adding the three regions proves \eqref{w14:basic}.

\smallskip\noindent
\emph{Step 4: consequences for the elliptic equation.}
For any smooth cutoff $\xi$ and $(\Delta+T)h=H$,
\[
 \Delta(\xi h)=\xi H-\xi Th
        +2\operatorname{div}(h\nabla\xi)-h\Delta\xi.
\]
The elliptic estimate therefore yields
\begin{equation}\label{w14:localized}
 \|\xi h\|_{W^{1,4}}
 \le C\left(\|\xi H\|_2+\|\xi Th\|_{\mathcal B}
       +\|h\nabla\xi\|_4+\|h\Delta\xi\|_{4/3}+\|\xi h\|_2\right).
\end{equation}
For $\xi=\chi_\lambda$, its derivative supports have area $O(q_0^2)$, and
\[
\begin{aligned}
 \|h\nabla\chi_\lambda\|_4
 &\le\|h\|_8\|\nabla\chi_\lambda\|_8
       \le Cq_0^{-3/4}\|h\|_{H^1},\\
 \|h\Delta\chi_\lambda\|_{4/3}
 &\le\|h\|_4\|\Delta\chi_\lambda\|_2
       \le Cq_0^{-1}\|h\|_{H^1}.
\end{aligned}
\]
Using \eqref{w14:basic} with $w_0=1$ we get
$\|h\|_{W^{1,4}(\Omega_\lambda)}
 \le C(\|H\|_2+\lambda\|h\|_{H^1})$.

For $L_0^jg\in H^1$, $0\le j\le2$, the exact multiplier identity
\[
 Tg=\frac{i}{\lambda}L_0g+\frac1{\lambda^2}L_0^2g
                           -\frac1{\lambda^2}TL_0^2g
\]
holds in distributions. Since $w_0\in L^4$, \eqref{w14:basic} implies
\begin{equation}\label{w14:improved}
\begin{aligned}
 \|\chi_\lambda w_0Tg\|_{\mathcal B}
 &\le\frac{\|w_0\|_4}{\lambda}\|L_0g\|_2
       +\frac{\|w_0\|_4}{\lambda^2}\|L_0^2g\|_2
       +\frac C\lambda\|L_0^2g\|_{H^1}
 \le\frac C\lambda\sum_{j=0}^2\|L_0^jg\|_{H^1}.
\end{aligned}
\end{equation}
Finally, choose a fixed smooth $\xi=1$ near $\overline{\Omega}_\ast$,
supported away from the centers. For this estimate choose $c$ smaller
if necessary, depending on $\Omega_*$, so that $\chi_\lambda=1$
on $\operatorname{supp}\xi$ for all $\lambda\ge1$.
Multiplication by $\xi$ is bounded on $\mathcal B$, as
\[
 \xi(\operatorname{div}F+f)
 =\operatorname{div}(\xi F)+\xi f-\nabla\xi\cdot F,\qquad
 \|\nabla\xi\cdot F\|_{4/3}\le\|\nabla\xi\|_2\|F\|_4.
\]
Thus \eqref{w14:improved}, with $w_0=1$, bounds $\|\xi Tg\|_{\mathcal B}$ by
$C_\xi\lambda^{-1}\sum_{j=0}^2\|L_0^jg\|_{H^1}$.
If $(\Delta+T)g=H$, it follows from \eqref{w14:localized} that
\[
 \|g\|_{W^{1,4}(\Omega_*)}
 \le C_{\Omega_*}\left(\|H\|_2+\sum_{j=0}^2\|L_0^jg\|_{H^1}\right).
\]
The cases for $T(i\lambda-0)$ and $\lambda\le-1$ follow similarly.

\end{proof}

\section{Limiting absorption principle}\label{sec:lap}
We identify the generalized kernel and prove the limiting absorption
principle, using the absence of nonzero eigenvalues in
Theorem~\ref{Intm6}. The final subsection describes the resolvent
near zero and the effect of the moment conditions
\eqref{IntzeroMoments}.

\subsection{The generalized kernel at zero}
Write
\[
 E=\operatorname{span}\{e_1,e_2,e_3\},\qquad
 (e_1,e_2,e_3)=(\sin x,\sin y,\cos x-\cos y).
\]
Recall the
functions in Theorem~\ref{th1:Int}:
\begin{equation}\label{zero:primitives}
 g_1=-y_c,\qquad g_2=x_c,\qquad
 g_3=\log\left|\frac{\cos((x-y)/2)}{\cos((x+y)/2)}\right|.
\end{equation}
Here $(x_c,y_c)$ are the centered representatives in the two cells; see the discussion below \eqref{IntzeroDual} for the details.

\begin{lemma}\label{zero:primitives-lemma}
We have $g_j\in\mathcal D\cap Y\cap E^\perp$, $L_0g_j=e_j$, and
$g_j\in L^p$ for every finite $p$.
\end{lemma}
\begin{proof}
In each cell, $L_0(-y_c)=\sin x$ and $L_0x_c=\sin y$.
With $a=(x+y)/2$, $b=(x-y)/2$,
\[
 L_0=-\cos a\sin b\,\partial_a+\sin a\cos b\,\partial_b,
 \qquad L_0(\log|\cos b|-\log|\cos a|)=-2\sin a\sin b=e_3.
\]
The rest of the Lemma is straightforward.
\end{proof}

For $u\in Y$, the Fourier expansion implies
\begin{equation}\label{zero:energy-coercivity}
 \frac12\|Qu\|_2^2
 \le \langle(I+\Delta^{-1})u,u\rangle
 \le \|Qu\|_2^2.
\end{equation}
Indeed, $u\perp\Psi=\cos x+\cos y$, and removing $E$ eliminates
the remaining first harmonics. For $j=1,2,3$, define $b_j$ as the
unique element of $Y\cap E^\perp$ satisfying
\begin{equation}\label{zero:C}
 \langle(I+\Delta^{-1})b_j,h\rangle=\langle g_j,h\rangle
 \qquad\text{for every }h\in Y\cap E^\perp.
\end{equation}
Existence and uniqueness follow from \eqref{zero:energy-coercivity}
and the Riesz representation theorem.

\begin{proposition}\label{zero:root}\label{zero:constraints}
The generalized kernel of $L$ in $Y$ consists of three chains:
\begin{equation}\label{zero:full-root}
 Lb_j=e_j,\qquad Le_j=0,\qquad
 \ker L^n\cap Y
 =E\oplus\operatorname{span}\{b_1,b_2,b_3\},\quad n\ge2.
\end{equation}
On $Y\cap E^\perp$, the operator $QL$, with domain
$\mathcal D\cap Y\cap E^\perp$, is skew-adjoint for the inner product
$\langle(I+\Delta^{-1})\cdot,\cdot\rangle$.
Its kernel is $\operatorname{span}\{b_1,b_2,b_3\}$, whose
orthogonal complement in this inner product is
\begin{equation}\label{zero:continuous-space}
\mathcal{Z}:= \{z\in Y\cap E^\perp:\ \langle z,g_j\rangle=0,
       \ j=1,2,3\}.
\end{equation}
This subspace is invariant under the projected evolution
$z\mapsto Qe^{-tL}z$, and $QL$ has trivial kernel on $\mathcal{Z}$.
Moreover,
\begin{equation}\label{zero:moments}
 \langle\omega(t),g_j\rangle=\langle\omega_0,g_j\rangle,
 \qquad j=1,2,3.
\end{equation}
Thus vanishing of the initial moments implies that $Q\omega(t)$
belongs to $\mathcal{Z}$ for all $t$.
\end{proposition}

\begin{proof}
If $Lu=0$ with $u\in\mathcal D\cap Y$, then
$(I+\Delta^{-1})u\in\ker L_0=Y^\perp$. Hence
$\langle(I+\Delta^{-1})u,u\rangle=0$, and
\eqref{zero:energy-coercivity} implies $u\in E$.
Conversely, $(I+\Delta^{-1})e_j=0$, so $Le_j=0$.
Therefore $\ker L\cap Y=E$.

By \eqref{zero:C}, the difference
$(I+\Delta^{-1})b_j-g_j$ is perpendicular to $Y\cap E^\perp$.
It is also perpendicular to $E$, since
$(I+\Delta^{-1})E=0$ and $g_j\perp E$. Consequently,
\[
 (I+\Delta^{-1})b_j-g_j\in Y^\perp=\ker L_0.
\]
Every function in $\ker L_0$ is even in $x$ and $y$, since each
regular level set of $\Psi$ is a single closed streamline.
Reflections preserve $Y\cap E^\perp$ and commute with
$I+\Delta^{-1}$, so uniqueness in \eqref{zero:C} implies that $b_j$
has the same reflection symmetries as $g_j$. Since $g_1$ is odd
in $y$ and $g_2,g_3$ are odd in $x$, the preceding difference
vanishes:
\begin{equation}\label{zero:full-duality}
 (I+\Delta^{-1})b_j=g_j.
\end{equation}
Lemma~\ref{zero:primitives-lemma} shows $g_j\in\mathcal D$ and
$L_0g_j=e_j$. Since $\Delta^{-1}b_j\in H^2\cap X\subset\mathcal D$,
we obtain $b_j\in\mathcal D$ and
\[
 Lb_j=L_0(I+\Delta^{-1})b_j=L_0g_j=e_j.
\]
If $L^2u=0$ with $u\in Y$, write $Lu=\sum_j a_je_j$.
Then $L(u-\sum_j a_jb_j)=0$, which proves
\eqref{zero:full-root} for $n=2$. The sum is direct because
$b_j\perp E$, and the $b_j$ are independent because the $e_j$ are.

To exclude longer chains, observe that for every
$u\in\mathcal D\cap Y$,
\begin{equation}\label{zero:adjoint-pairing}
 \langle Lu,g_j\rangle_{L^2}
 =-\langle(I+\Delta^{-1})u,L_0g_j\rangle_{L^2}
 =-\langle u,(I+\Delta^{-1})e_j\rangle_{L^2}=0.
\end{equation}
Suppose $u\in Y$ and $L^3u=0$. The description of $\ker L^2\cap Y$ implies
$Lu=e+v$, where $e\in E$ and $v=\sum_j a_jb_j$.
Since $(I+\Delta^{-1})v=\sum_j a_jg_j$ and
$v\in Y\cap E^\perp$, it follows that
\[
 0=\Big\langle Lu,\sum_j a_jg_j\Big\rangle_{L^2}
  =\left\langle v,(I+\Delta^{-1})v\right\rangle_{L^2}
  \ge\frac12\|v\|_2^2.
\]
Thus $Lu\in E$ and $L^2u=0$. Induction now proves
\eqref{zero:full-root} for every $n\ge2$.

The rest of the conclusions follow from properties of $L$ and \eqref{zero:full-root}-\eqref{zero:continuous-space}.
\end{proof}

In the following, the spectral calculus of $-iQL$ is taken on
$Y\cap E^\perp$ with the inner product
$\langle(I+\Delta^{-1})\cdot,\cdot\rangle_{L^2}$.
For $f$ in the range of $QL$, we denote by $(QL)^{-1}f$ the unique
solution $v\in\mathcal D\cap Y\cap E^\perp$ of
$QLv=f$ satisfying $\langle v,g_j\rangle=0$ for $j=1,2,3$.
This inverse is generally unbounded; its composition with
$\chi(-iQL)$ is bounded for a bounded spectral cutoff $\chi$
supported away from zero.

The matrix of pairings of the generalized modes and the adjoint
zero modes is positive definite. Indeed, by \eqref{zero:full-duality},
for $a\in\mathbb C^3$,
\begin{equation}\label{zero:projection}
\begin{aligned}
 \sum_{i,j=1}^3a_j\overline{a_i}\langle b_j,g_i\rangle_{L^2}
 &=\Big\langle(I+\Delta^{-1})\sum_{j=1}^3a_jb_j,
                         \sum_{j=1}^3a_jb_j\Big\rangle_{L^2}\\
 &\ge\tfrac12\left\|\sum_{j=1}^3a_jb_j\right\|_2^2.
\end{aligned}
\end{equation}
We impose the moments directly, since the corresponding energy
projection subtracts generalized modes and need not preserve $H^3$.
In particular, $b_3=g_3-\Delta^{-1}b_3$ is logarithmically unbounded,
whereas $\Delta^{-1}b_3\in H^2\subset L^\infty$.

For nonzero imaginary eigenvalues, Theorem~\ref{Intm6} implies
\begin{equation}\label{zero:nonzero-spectrum}
 \ker(L-i\lambda)=\{0\}\quad(\lambda\ne0),\qquad
 \mathcal D(L)=\{h\in X:L_0h\in L^2\}.
\end{equation}
Nonzero eigenvectors of $L$ belong to $Y$ and project to
eigenvectors of $QL$. Conversely, if $QLz=i\lambda z$ with
$z\in\mathcal D\cap Y\cap E^\perp$ and $\lambda\ne0$, then
\[
 L\bigl(z+(i\lambda)^{-1}PLz\bigr)
 =i\lambda\bigl(z+(i\lambda)^{-1}PLz\bigr).
\]
Thus $QL$ has no nonzero eigenvalues. There is also no nontrivial kernel as excluded by
the moment conditions in \eqref{zero:continuous-space}.

\subsection{Local elliptic estimates and homogeneous solutions}
Recall that, for $f\in H^1\cap Y$ and
$z=\varepsilon+i\lambda$, $\varepsilon\ne0$, the stream
resolvent satisfies
\[
 \psi_f(z)=\Delta^{-1}(z-L)^{-1}f,\qquad
 (\Delta+T(z))\psi_f(z)=-(L_0-z)^{-1}f,
 \qquad T(z)=(L_0-z)^{-1}L_0.
\]
We compare this equation near an elliptic point with the
operator in Proposition~\ref{compact:normal}.
In the coordinates
\[
 q=\sqrt{8(1-2\pi/S)},\qquad v=\log q,\qquad
 0<q<2q_0,\qquad \theta\in\mathbb T,
\]
it follows from Lemma~\ref{geom:metric} that
\[
 L_0=(1-q^2/8)\partial_\theta,\qquad
 dx\,dy=\xi_0(q)q^2\,dv\,d\theta,\qquad
 \mathsf G=I+O(q^2),\quad \xi_0=1+O(q^2).
\]
For $n\ne0$, set
\[
 \sigma_n=8(1-\lambda/n),\qquad
 \eta_n=8\varepsilon/n,\qquad
 V_n(v)=\frac{8e^{2v}}{e^{2v}-\sigma_n-i\eta_n}.
\]
The exact multiplier identity
\[
 (T(z)h)_n
 =-\frac{8(1-q^2/8)}{q^2-\sigma_n-i\eta_n}h_n
\]
implies
\[
 -\xi_0q^2(T(z)h)_n=cV_nh_n,\qquad
 c(v)=\xi_0(q)(1-q^2/8)=1+O(q^2).
\]
Consequently, for functions supported in the coordinate disk,
\[
 \langle(-\Delta-T(z))h,\varphi\rangle
 =\int\mathsf G\nabla_{v,\theta}h\cdot
                  \nabla_{v,\theta}\bar\varphi\,dv\,d\theta
   +2\pi\sum_{n\ne0}\langle cV_nh_n,\varphi_n\rangle.
\]

Let $\mathscr E$ be the cylinder energy space modulo constants,
with
\[
 \|h\|_{\mathscr E}^2
 =\int_{\mathbb R\times\mathbb T}
       \bigl(|\partial_vh|^2+|\partial_\theta h|^2\bigr)
       \,dv\,d\theta,
\]
and define the comparison operator by
\[
 (\mathscr Nh)_0=-\partial_v^2h_0,\qquad
 (\mathscr Nh)_n=(-\partial_v^2+n^2+V_n)h_n
 \quad(n\ne0).
\]
When $|\sigma_n|+|\eta_n|>0$, the translation
$v\mapsto v+\frac12\log(|\sigma_n|+|\eta_n|)$ transforms the
$n$th operator into
\[
 \mathcal N_{|n|,\,
       \sigma_n/(|\sigma_n|+|\eta_n|),\,
       \eta_n/(|\sigma_n|+|\eta_n|)}.
\]
When $\sigma_n=\eta_n=0$, it is
$-\partial_v^2+n^2+8$. By Proposition~\ref{compact:normal}
and the isomorphism
$-\partial_v^2:\dot H^1(\mathbb R)\to\dot H^{-1}(\mathbb R)$
for the zero mode, we obtain
\[
 \mathscr N:\mathscr E\longrightarrow\mathscr E^*
 \quad\text{bijective},\qquad
 \|\mathscr N^{-1}\|_{\mathscr E^*\to\mathscr E}\le C.
\]

Multiply $\mathsf G-I$ and $c-1$ by a smooth radial cutoff
equal to one on $q<q_0$ and zero on $q>2q_0$, and extend
the resulting coefficients to the whole cylinder.
Keeping the same notation, we have
\[
 \|\mathsf G-I\|_{L^\infty}
       +\|c-1\|_{W^{1,\infty}(\mathbb R_v)}
 \le Cq_0^2.
\]
The gradient error satisfies
\[
 \left|\int(\mathsf G-I)\nabla_{v,\theta}h\cdot
                 \nabla_{v,\theta}\bar\varphi\,dv\,d\theta\right|
 \le Cq_0^2\|h\|_{\mathscr E}\|\varphi\|_{\mathscr E}.
\]
For the potential error, it follows from \eqref{eq:cylinder-form} and Cauchy--Schwarz that
\[
\begin{aligned}
 \Big|2\pi\sum_{n\ne0}
       \langle(c-1)V_nh_n,\varphi_n\rangle\Big|
 &\le C\|c-1\|_{W^{1,\infty}}
       \sum_{n\ne0}|n|^{-1/2}
          \|h_n\|_{|n|}\|\varphi_n\|_{|n|}\\
 &\le Cq_0^2\|h\|_{\mathscr E}\|\varphi\|_{\mathscr E}.
\end{aligned}
\]
Thus the extended operator differs from $\mathscr N$ by
$O(q_0^2)$ in $\mathcal L(\mathscr E,\mathscr E^*)$.
For sufficiently small $q_0$, we obtain a uniform inverse
by a Neumann series. All these estimates hold for both boundary
values as $\varepsilon\to0^\pm$. By localizing this inverse,
we obtain the following estimate.

\begin{lemma}\label{global:residual}
Fix a sufficiently small $q_0>0$ and let
$z=\varepsilon+i\lambda$, where $\lambda\in\mathbb R$ and
$\varepsilon\ne0$.
If $h\in H^1$ satisfies
\[
 (\Delta+T(z))h=r\in H^{-1},
\]
and $r$ has zero angular mean on $q<q_0$, then
\begin{equation}\label{global:residual-estimate}
 \|h\|_{H^1(q<q_0/2)}
 \le C\bigl(\|r\|_{H^{-1}}
            +\|h\|_{H^1(q_0/2<q<q_0)}\bigr).
\end{equation}
If $f\in H^1\cap Y$ and $\psi\in H^1$ satisfies the resolvent
equation
\[
 (\Delta+T(z))\psi=-(L_0-z)^{-1}f,
\]
then
\begin{equation}\label{global:resolvent-source-estimate}
 \|\psi\|_{H^1(q<q_0/2)}
 \le C\bigl(\|f\|_{H^1(q<2q_0)}
            +\|\psi\|_{H^1(q_0/2<q<q_0)}\bigr).
\end{equation}
The estimates hold near either elliptic point, with $C$
depending on $q_0$ but independent of $\lambda$ and $\varepsilon$.
The same estimates hold for both boundary equations, with the
resolvent multipliers interpreted as the limits $\varepsilon\to0^\pm$.
All Sobolev norms are taken in physical coordinates with
measure $dx\,dy$.
\end{lemma}

\begin{remark}
Limiting absorption bounds for arbitrary $H^{-1}$ sources of
zero spatial mean, without the zero angular mean condition,
are established in Lemma~\ref{r2:full-source-LAP}, uniformly
on compact subsets of $\mathbb R\setminus\{0\}$.
\end{remark}

\begin{proof}
Choose a radial cutoff $\chi$ equal to one on $q<q_0/2$
and supported in $q<q_0$.
The zero angular mean assumption means that $r$ vanishes
against radial test functions supported in this disk. Throughout this proof, $\nabla$ and $\operatorname{div}$
are taken with respect to $(v,\theta)$, where $v=\log q$.

Since $c$ and $\xi_0$ depend only on $q$, averaging the equation
\[
 -\operatorname{div}(\mathsf G\nabla h)
       +\sum_{n\ne0}cV_nh_ne^{in\theta}
 =-\xi_0q^2r
\]
implies
\[
 \partial_v\int_{\mathbb T}(\mathsf G\nabla h)_v\,d\theta=0,
 \qquad
 \int_{-\infty}^{\log q_0}
 \left|\int_{\mathbb T}(\mathsf G\nabla h)_v\,d\theta\right|^2dv
 \le C\int_{q<q_0}|\nabla h|^2<\infty.
\]
Thus this constant flux is zero. In particular,
\[
 \int_{\mathbb T}\chi'(v)(\mathsf G\nabla h)_v\,d\theta=0.
\]

For a smooth test function $\varphi$, write
$\varphi_\perp=\varphi-(2\pi)^{-1}\int_{\mathbb T}\varphi\,d\theta$.
Since $T(z)$ commutes with $\chi$, integration by parts
and the two zero-mean identities imply
\[
\begin{aligned}
 &\int\mathsf G\nabla(\chi h)\cdot\nabla\bar\varphi
       +2\pi\sum_{n\ne0}\langle cV_n\chi h_n,\varphi_n\rangle\\
 &\qquad=-\langle r,\chi\varphi_\perp\rangle
       +\int h\mathsf G\nabla\chi\cdot\nabla\bar\varphi
       -\int\chi'(v)(\mathsf G\nabla h)_v\bar\varphi_\perp.
\end{aligned}
\]
Here the pairing with $r$ is in physical coordinates.
By angular Poincar\'e and metric equivalence,
\[
 \|\varphi_\perp\|_{L^2(dv\,d\theta)}
       \le\|\partial_\theta\varphi\|_{L^2(dv\,d\theta)},
 \qquad
 \|\chi\varphi_\perp\|_{H^1(\mathbb T^2)}
       \le C\|\varphi\|_{\mathscr E}.
\]
Moreover, the derivatives of $\chi$ are supported in the fixed
annulus $q_0/2<q<q_0$, where
\[
 \|h\mathsf G\nabla\chi\|_{L^2(dv\,d\theta)}
 +\|\chi'(v)(\mathsf G\nabla h)_v\|_{L^2(dv\,d\theta)}
 \le C\|h\|_{H^1(q_0/2<q<q_0)}.
\]
It follows from these bounds and the uniform inverse established
above that
\[
 \|\chi h\|_{\mathscr E}
 \le C\bigl(\|r\|_{H^{-1}}+
                    \|h\|_{H^1(q_0/2<q<q_0)}\bigr).
\]
Finally, \eqref{global:residual-estimate} follows from
$\chi h\in H^1_0(\{q<q_0\})$, Poincar\'e, and metric equivalence:
\[
 \|h\|_{H^1(q<q_0/2)}
 \le\|\chi h\|_{H^1(q<q_0)}
 \le C\|\chi h\|_{\mathscr E}.
\]

For \eqref{global:resolvent-source-estimate}, choose a radial cutoff
$\widetilde\chi$ equal to one on $\operatorname{supp}\chi$
and supported in $q<2q_0$.
Since $f\in Y$, its zero angular coefficient vanishes, and
\[
 \bigl(\xi_0q^2(L_0-z)^{-1}f\bigr)_n
       =\frac{i\xi_0}{n}V_nf_n\quad(n\ne0),
 \qquad
 \bigl(\xi_0q^2(L_0-z)^{-1}f\bigr)_0=0.
\]
It follows from \eqref{eq:cylinder-form}, with coefficient
$\chi\xi_0$, and Cauchy--Schwarz that
\[
\begin{aligned}
 \bigl|\langle(L_0-z)^{-1}f,\chi\varphi\rangle\bigr|
 &\le C\sum_{n\ne0}|n|^{-3/2}
       \|\widetilde\chi f_n\|_{|n|}\|\varphi_n\|_{|n|}\\
 &\le C\|\widetilde\chi f\|_{\mathscr E}\|\varphi\|_{\mathscr E}.
\end{aligned}
\]
The pairing on the left is in physical coordinates.
Using $f_0=0$, angular Poincar\'e, and metric equivalence,
we further obtain
\[
\begin{aligned}
 \|\widetilde\chi f\|_{\mathscr E}^2
 &\le C\int_{q<2q_0}
       \bigl(|\nabla f|^2+|\partial_v\widetilde\chi|^2|f|^2\bigr)\le C\int_{q<2q_0}|\nabla f|^2
 \le C\|f\|_{H^1(q<2q_0)}^2.
\end{aligned}
\]
Replacing the estimate for the residual pairing by these bounds
proves \eqref{global:resolvent-source-estimate}. The same identities and estimates
hold for both boundary values of the multipliers.
\end{proof}

Vanishing of the numerator at a critical streamline removes the singularity at the critical layers. The next lemma summarizes this point at the separatrix.
\begin{lemma}\label{lem:bootstrap}
Fix $\lambda\ne0$. If $h\in H^1$ and its coefficients vanish at
the separatrix resonances $\tau_k=2\pi k/\lambda$, then near the
separatrix,
\begin{align}
 \|T_\pm(i\lambda)h\|_p&\le C_{\lambda,p}\|h\|_{H^1},
                 &&1<p<2,\label{eq:Lpbootstrap}\\
 \|T_\pm(i\lambda)h\|_2&\le C_{\lambda,q}
                   (\|h\|_{H^1}+\|\nabla h\|_q),
                 &&q>2.\label{eq:L2bootstrap}
\end{align}
The second estimate assumes $h\in W^{1,q}$. Both boundary values
are strong limits in the above spaces and are the same.
\end{lemma}
\begin{proof}
It suffices to consider $\lambda>0$ and the positive cell; complex
conjugation and symmetry imply the other cases.
We use the decomposition of Proposition~\ref{prop:tail}.

Choose $R$ sufficiently large, halfway between two consecutive
resonances, and set
\[
 d=\min\left(1,\frac{\pi}{4\lambda}\right),\qquad
 I_k=(\tau_k-d,\tau_k+d),\qquad \tau_k>R.
\]
These intervals are disjoint and contained in $(R,\infty)$.
Choose $0\le\chi_k\le1$ in $C_c^\infty(I_k)$, with $\chi_k=1$
when $|\tau-\tau_k|\le d/2$, and set $\chi_k=0$ for the remaining
indices. Recall that
\[
 dx\,dy=\frac{w(\tau)}{2\pi}\,d\tau\,d\theta,\qquad
 w(\tau)\sim\tau e^{-\tau/4},\qquad
 w(\tau)\sim_\lambda w(\tau_k)\quad(\tau\in I_k).
\]
On $S>R$, the exact decomposition is
\[
\begin{aligned}
 T(\varepsilon+i\lambda)h
 &=\sum_{\tau_k>R}\chi_k
       \frac{2\pi k}{2\pi k-\lambda\tau+i\varepsilon\tau}
       h_ke^{ik\theta}\\
 &\quad+\sum_{k\ne0}(1-\chi_k)
       \frac{2\pi k}{2\pi k-\lambda\tau+i\varepsilon\tau}
       h_ke^{ik\theta}.
\end{aligned}
\]

\smallskip\noindent
\emph{Step 1: The regular part.}
For $k<0$, the multiplier has absolute value at most one.
For $k>0$, use $\tau_k=\tau+(\tau_k-\tau)$ to obtain
\[
 \left|(1-\chi_k)
       \frac{2\pi k}{2\pi k-\lambda\tau+i\varepsilon\tau}\right|
 \le (1-\chi_k)\left(1+\frac{\tau}{|\tau-\tau_k|}\right)
 \le C_\lambda(1+\tau).
\]
Consequently, Parseval's identity and Hölder's inequality imply
\[
\begin{aligned}
 &\left\|\sum_{k\ne0}(1-\chi_k)
       \frac{2\pi k}{2\pi k-\lambda\tau+i\varepsilon\tau}
       h_ke^{ik\theta}\right\|_{L^2(S>R)}^2\\
 &\quad\le C_\lambda\int_R^\infty
       (1+\tau)^2\sum_k|h_k(\tau)|^2w(\tau)\,d\tau\\
 &\quad\le C_\lambda\|(1+S)h\|_2^2
 \le C_\lambda\|1+S\|_4^2\|h\|_4^2
 \le C_\lambda\|h\|_{H^1}^2.
\end{aligned}
\]
Here $\|1+S\|_4<\infty$ follows from the exponential decay of $w$.
Since the physical area is finite, this also controls the regular
part in $L^p$ for $1<p<2$.

\smallskip\noindent
\emph{Step 2: The resonant part.}
For $\tau\in I_k$ and either sign of $\varepsilon$,
\[
 \left|\frac{2\pi k}{2\pi k-\lambda\tau+i\varepsilon\tau}\right|
 \le\frac{\tau_k}{|\tau-\tau_k|},\qquad
 \left\|\frac{h_k}{\tau-\tau_k}\right\|_{L^2(I_k,d\tau)}
 \le2\|\partial_\tau h_k\|_{L^2(I_k,d\tau)}.
\]
The second inequality follows from $h_k(\tau_k)=0$.
Weight comparison and Hölder's inequality on $I_k$ therefore imply
\[
 \left\|\chi_k
       \frac{2\pi k h_k}{2\pi k-\lambda\tau+i\varepsilon\tau}
       \right\|_{L^p(I_k,w\,d\tau)}
 \le C_{\lambda,p}\tau_k w(\tau_k)^{1/p-1/2}
       \|\partial_\tau h_k\|_{L^2(I_k,w\,d\tau)}.
\]
At each $\tau$, at most one summand in the resonant part is nonzero.
Taking the $p$th power and applying Hölder's inequality in $k$,
we obtain
\[
\begin{aligned}
 &\left\|\sum_{\tau_k>R}\chi_k
       \frac{2\pi k h_k e^{ik\theta}}
            {2\pi k-\lambda\tau+i\varepsilon\tau}
       \right\|_{L^p(S>R)}^p\\
 &\quad\le C_{\lambda,p}
 \left(\sum_{\tau_k>R}\tau_k^{2p/(2-p)}w(\tau_k)\right)^{1-p/2}
 \left(\sum_{\tau_k>R}\int_{I_k}
             |\partial_\tau h_k|^2w\,d\tau\right)^{p/2}
 \le C_{\lambda,p}\|\nabla h\|_2^p.
\end{aligned}
\]
The first sum is finite because $\tau_k=2\pi k/\lambda$ and
$w(\tau_k)\sim\tau_k e^{-\tau_k/4}$; the derivative sum is
controlled by Lemma~\ref{lem:geometry}.

If $h\in W^{1,q}$, $q>2$, the same Hardy inequality followed by
Hölder's inequality implies
\[
\begin{aligned}
 \left\|\chi_k
       \frac{2\pi k h_k}{2\pi k-\lambda\tau+i\varepsilon\tau}
       \right\|_{L^2(I_k,w\,d\tau)}
 &\le C_\lambda\tau_k w(\tau_k)^{1/2}
                 \|\partial_\tau h_k\|_{L^2(I_k,d\tau)}\\
 &\le C_{\lambda,q}\tau_k w(\tau_k)^{1/2-1/q}
                 \|\partial_\tau h_k\|_{L^q(I_k,w\,d\tau)}.
\end{aligned}
\]
Using Lemma~\ref{lem:geometry} for each coefficient and summing
the squares, we conclude that
\[
\begin{aligned}
 \left\|\sum_{\tau_k>R}\chi_k
       \frac{2\pi k h_k e^{ik\theta}}
            {2\pi k-\lambda\tau+i\varepsilon\tau}
       \right\|_{L^2(S>R)}^2
 &\le C_{\lambda,q}\sum_{\tau_k>R}
       \tau_k^2w(\tau_k)^{1-2/q}
       \|\partial_\tau h_k\|_{L^q(I_k,w\,d\tau)}^2\\
 &\le C_{\lambda,q}\|\nabla h\|_q^2
       \sum_{\tau_k>R}\tau_k^2w(\tau_k)^{1-2/q}
 \le C_{\lambda,q}\|\nabla h\|_q^2.
\end{aligned}
\]
The last series converges because $1-2/q>0$.
Together with Step~1, these estimates prove both bounds uniformly
for $\varepsilon\ne0$.

\smallskip\noindent
\emph{Step 3: Strong boundary limits.}
Away from the resonances, every multiplier converges as
$\varepsilon\to0^\pm$ to the same quotient
$2\pi k/(2\pi k-\lambda\tau)$.
The resonant part is dominated pointwise by
\[
 \sum_{\tau_k>R}\mathbf1_{I_k}(\tau)
       \frac{\tau_k|h_k(\tau)|}{|\tau-\tau_k|}.
\]
The estimates in Step~2 also apply with $\chi_k$ replaced by
$\mathbf1_{I_k}$, so this majorant belongs to $L^p$, and to $L^2$
under the additional $W^{1,q}$ assumption. Dominated convergence
therefore proves strong convergence of the resonant part.
For the regular part, dominated convergence on the product of
counting measure in $k$ and $w(\tau)\,d\tau$ implies
\[
 \lim_{\varepsilon\to0^\pm}\sum_{k\ne0}\int_R^\infty
 \left|(1-\chi_k)
 \left(\frac{2\pi k}{2\pi k-\lambda\tau+i\varepsilon\tau}
       -\frac{2\pi k}{2\pi k-\lambda\tau}\right)h_k\right|^2
 w\,d\tau=0,
\]
since the squared coefficient is bounded by
$C_\lambda(1+\tau)^2|h_k|^2$, whose weighted sum and integral are
finite by Step~1. Parseval's identity implies strong $L^2$ convergence,
and hence strong $L^p$ convergence for $p<2$.
The same argument applies in the other cell. Combining the two
contributions proves strong convergence in the claimed spaces on
$\{S>R\}$. Testing against $\varphi\in C_c^\infty(\{S>R\})$
identifies the limit with the already defined $H^{-1}$ boundary
value $T_\pm(i\lambda)h$. Both signs have the same limit.
\end{proof}

We now identify eigenfunctions
in $\mathcal D(L)$ for $\Delta+T(i\lambda\pm0)$.
\begin{proposition}\label{prop:nullvector}
For every $\lambda\ne0$, a mean-zero $H^1$ solution of
\begin{equation}\label{eq:null}
 (\Delta+T(i\lambda\pm0))h=0
\end{equation}
belongs to $H^2$, and $\Delta h\in\mathcal D(L)$ satisfies
$L\Delta h=i\lambda\Delta h$.
\end{proposition}
\begin{proof}
By complex conjugation, it suffices to consider $\lambda>0$.
We first remove the singularities on regular critical streamlines,
then treat the elliptic points, and finally use Lemma~\ref{lem:bootstrap}
at the separatrix. Sobolev norms without a specified domain are
taken on $\mathbb T^2$ with measure $dx\,dy$.

\smallskip\noindent
\emph{Step 1: Regular critical streamlines.}
The boundary convergence established in Section~2 and the equation imply
\[
 \lim_{\varepsilon\to0^\pm}
 \operatorname{Im}\langle T(\varepsilon+i\lambda)h,h\rangle
 =\operatorname{Im}\langle-\Delta h,h\rangle=0.
\]
Recall that $dx\,dy=S(r)\,dr\,d\theta/(2\pi)$ and
\[
 (T(\varepsilon+i\lambda)h)_n
 =\frac{nb(r)}{nb(r)-\lambda+i\varepsilon}h_n.
\]
Parseval's identity therefore implies
\[
\begin{aligned}
 &\sum_{\mathrm{cells}}\sum_{n>0}
 \int\frac{|\varepsilon|S(r)nb(r)}
              {(nb(r)-\lambda)^2+\varepsilon^2}|h_n(r)|^2\,dr\\
 &\quad=-\operatorname{sgn}(\varepsilon)
           \operatorname{Im}\langle T(\varepsilon+i\lambda)h,h\rangle
 +\sum_{\mathrm{cells}}\sum_{n<0}
 \int\frac{|\varepsilon|S(r)|nb(r)|}
              {(nb(r)-\lambda)^2+\varepsilon^2}|h_n(r)|^2\,dr\\
 &\quad\le
 \bigl|\operatorname{Im}\langle T(\varepsilon+i\lambda)h,h\rangle\bigr|
 +\frac{|\varepsilon|}{4\lambda}\|h\|_2^2
 \longrightarrow0.
\end{aligned}
\]
Here $a/(a+\lambda)^2\le(4\lambda)^{-1}$ for $a\ge0$.
At every regular resonance $nb(r_n)=\lambda$,
$h_n$ is locally $H^1$ in $r$. Restricting the nonnegative sum
to a small interval about $r_n$ and taking the limit, we obtain
\[
 \frac{\pi\lambda S(r_n)}{|nb'(r_n)|}|h_n(r_n)|^2=0,
 \qquad h_n(r_n)=0.
\]
On such an interval, Hardy's inequality implies
\[
 \left\|\frac{h_n}{r-r_n}\right\|_{L^2(dr)}
 \le2\|\partial_rh_n\|_{L^2(dr)},
 \qquad
 \left|\frac{nb(r)}{nb(r)-\lambda+i\varepsilon}\right|
 \le\frac{C_{\lambda,n}}{|r-r_n|}.
\]
On any fixed regular annulus, $b\ge b_*>0$, and the multiplier
is bounded by $2$ when $|n|>2\lambda/b_*$.
There are only finitely many remaining modes. Thus, for regular
annuli $A\Subset A'$, Parseval's identity and the preceding bounds imply
\[
 \|T(i\lambda\pm0)h\|_{L^2(A)}
 \le C_{\lambda,A,A'}\|h\|_{H^1(A')}.
\]
In particular, the two boundary multipliers agree with the ordinary
quotient on these annuli.

\smallskip\noindent
\emph{Step 2: The elliptic points.}
Use $q=\sqrt{8(1-b)}$, $v=\log q$, and
$c=\xi_0(q)(1-q^2/8)$, so that
\[
 dx\,dy=\xi_0(q)q^2\,dv\,d\theta,
 \qquad \mathsf G-I=O(q^2),\qquad c-1=O(q^2).
\]
If $\lambda\notin\mathbb Z$, all angular multipliers are uniformly
bounded on a sufficiently small disk, and hence $Th\in L^2$ there.
Suppose now that $\lambda=k\in\mathbb N$, and choose $q_0>0$ so
small that all multipliers except the $k$th are bounded on $q<q_0$.
Since $c$ is radial, the $k$th coefficient of the equation satisfies
\[
 (-\partial_v^2+k^2+8)h_k
 =\bigl[\operatorname{div}((\mathsf G-I)\nabla h)\bigr]_k
       +8(1-c)h_k
 \quad (v<\log q_0).
\]
In this step, $\nabla$ and $\operatorname{div}$ are taken in $(v,\theta)$.
In the remainder of this step, one-dimensional norms are on
$(-\infty,\log q_0)$ with measure $dv$, and $H^{-1}=(H_0^1)^*$.
By Parseval's identity and the Dirichlet-form comparison,
\[
 \int_{-\infty}^{\log q_0}(|h_k'|^2+k^2|h_k|^2)\,dv
 \le\frac1{2\pi}\int_{v<\log q_0}|\nabla h|^2\,dv\,d\theta
 \le C\|h\|_{H^1(q<q_0)}^2.
\]
For $\varphi\in C_c^\infty(-\infty,\log q_0)$, integration by parts implies
\[
\begin{aligned}
 &\left\langle e^{-v}
       \bigl[\operatorname{div}((\mathsf G-I)\nabla h)\bigr]_k,
       \varphi\right\rangle\\
 &\qquad=-\frac1{2\pi}\int_{v<\log q_0}\int_{\mathbb T}
       (\mathsf G-I)\nabla h\cdot
       \nabla\bigl(e^{-v}\overline\varphi e^{-ik\theta}\bigr)
       \,d\theta\,dv.
\end{aligned}
\]
It follows from $\mathsf G-I,1-c=O(e^{2v})$ and the Dirichlet-form
comparison that
\[
 \left\|e^{-v}
       \bigl[\operatorname{div}((\mathsf G-I)\nabla h)\bigr]_k
       \right\|_{H^{-1}}
 +\|8e^{-v}(1-c)h_k\|_{L^2}
 \le C_{k,q_0}\|h\|_{H^1(q<q_0)}.
\]
Moreover, for $\varphi\in H_0^1(-\infty,\log q_0)$,
\[
 \operatorname{Re}\int
 \bigl(|\varphi'|^2-2\varphi'\overline\varphi
                  +(k^2+7)|\varphi|^2\bigr)\,dv
 =\|\varphi'\|_2^2+(k^2+7)\|\varphi\|_2^2.
\]
After subtracting a fixed lift of the boundary value, the
Lax--Milgram theorem therefore implies the existence of $u\in H^1$ solving
\[
\begin{cases}
 (-\partial_v^2-2\partial_v+k^2+7)u
 =e^{-v}\bigl[\operatorname{div}((\mathsf G-I)\nabla h)\bigr]_k
       +8e^{-v}(1-c)h_k,\\[2pt]
 u(\log q_0)=q_0^{-1}h_k(\log q_0).
\end{cases}
\]
The boundary trace is controlled by the $H^1$ norm on the fixed
annulus $q_0/2<q<q_0$. Hence
$\|u\|_{H^1}\le C_{k,q_0}\|h\|_{H^1(q<q_0)}$.
Both $h_k$ and $e^vu$ belong to $H^1(-\infty,\log q_0)$,
solve the same equation, and have the same boundary value.
Testing their difference against itself implies $h_k=e^vu$.
Consequently,
\[
\begin{aligned}
 \|q^{-2}h_ke^{ik\theta}\|_{L^2(q<q_0;\,dx\,dy)}^2
 &=2\pi\int_{-\infty}^{\log q_0}
       \xi_0(e^v)|e^{-v}h_k(v)|^2\,dv\\
 &\le C_{k,q_0}\|h\|_{H^1(q<q_0)}^2.
\end{aligned}
\]
Finally,
\[
 \left|\frac{k(1-q^2/8)}{-kq^2/8+i\varepsilon}h_k\right|
 \le8q^{-2}|h_k|.
\]
Dominated convergence therefore implies that the $k$th mode converges
in physical $L^2$ to $(1-8q^{-2})h_k$. The other multipliers are
uniformly bounded on this disk, so $Th\in L^2$ near either elliptic
point, also at integer thresholds.

\smallskip\noindent
\emph{Step 3: The separatrix and the operator domain.}
Every separatrix resonance has vanishing trace by Step~1.
Combining Lemma~\ref{lem:bootstrap}, with $p=3/2$, and the bounds
away from the separatrix, we obtain
\[
 \|T(i\lambda\pm0)h\|_{L^{3/2}}
 \le C_\lambda\|h\|_{H^1}.
\]
Since $\Delta h=-T(i\lambda\pm0)h$ and $h$ has zero mean,
elliptic regularity and Sobolev embedding imply
\[
 \|\nabla h\|_{L^6}
 \le C\|h\|_{W^{2,3/2}}
 \le C\|T(i\lambda\pm0)h\|_{L^{3/2}}
 \le C_\lambda\|h\|_{H^1}.
\]
The second estimate of Lemma~\ref{lem:bootstrap}, with exponent $6$,
then implies
\[
 \|T(i\lambda\pm0)h\|_2
 \le C_\lambda\bigl(\|h\|_{H^1}+\|\nabla h\|_6\bigr),
 \qquad
 \|h\|_{H^2}\le C\|T(i\lambda\pm0)h\|_2.
\]
Finally, pass to the limit in the resolvent identity
\[
 (L_0-i\lambda)T(\varepsilon+i\lambda)h
 =L_0h+\varepsilon T(\varepsilon+i\lambda)h.
\]
Strong $H^{-1}$ convergence implies convergence in distributions,
and therefore
\[
 L_0\Delta h=i\lambda\Delta h-L_0h\in L^2.
\]
Thus $\Delta h\in\mathcal D(L)$. Since $\Delta^{-1}\Delta h=h$,
\[
 L\Delta h=L_0\Delta h+L_0h=i\lambda\Delta h.
\]

\end{proof}

\subsection{Limiting absorption away from zero}
We collect the free-operator estimates needed in the compactness
argument.
\begin{lemma}\label{global:free}
For every $M>0$, the operators $T(\varepsilon+i\lambda)$ and
their $L^2$ adjoints are uniformly bounded from $H^1$ to $H^{-1}$
for $|\lambda|\le M$ and $0<|\varepsilon|\le1$.
Both families extend strongly continuously to each closed side
$0\le\varepsilon\le1$ and $-1\le\varepsilon\le0$, including
integer thresholds. For every smooth cutoff $\chi$ supported away
from the elliptic points, the family
$\chi T(\varepsilon+i\lambda):H^1\to H^{-1}$ is collectively compact
on these closed parameter sets.
\end{lemma}
\begin{proof}
All Sobolev norms are in physical coordinates unless another
domain is indicated. We estimate the regular annuli first, then
remove small neighborhoods of the elliptic points and the separatrix.

\smallskip\noindent
\emph{Step 1: Regular annuli.}
Let $\alpha$ be a smooth streamline cutoff supported in regular
annuli, and let $\zeta$ be a larger such cutoff, equal to one on
$\operatorname{supp}\alpha$. In each cell use $b$ as radial
coordinate and extend $\zeta h_k,\zeta g_k$ by zero to
$\mathbb R_b$. The estimates for the Cauchy-type multiplier in the proof of
Proposition~\ref{prop:compact} are
\[
\begin{aligned}
 &\left\|\frac1{b-a+i\eta}\right\|_{H^{-3/4}(\mathbb R)}\le C,\\
 &\left\|\frac1{b-a+i\eta}-\frac1{b-a'+i\eta'}\right\|_{H^{-3/4}(\mathbb R)}
 \le C(|a-a'|+|\eta-\eta'|)^{1/4},
\end{aligned}
\]
where $\eta,\eta'$ have the same sign; the inequalities include
the corresponding one-sided boundary values.
Using $dx\,dy=b^{-1}|dr/db|\,db\,d\theta$, we obtain
\[
 \langle\alpha T(\varepsilon+i\lambda)h,g\rangle
 =2\pi\sum_{k\ne0}\int_{\mathbb R}
 \frac{\alpha(b)|dr/db|}{b-\lambda/k+i\varepsilon/k}
                    h_k\overline{g_k}\,db
\]
on each cell. The algebra property of $H^{3/4}(\mathbb R)$,
followed by Cauchy--Schwarz and Parseval, implies
\[
\begin{aligned}
 |\langle\alpha T(z)h,g\rangle|
 &\le C_\alpha\sum_{k\ne0}
       \|\zeta h_k\|_{H^1(\mathbb R)}
       \|\zeta g_k\|_{H^{3/4}(\mathbb R)}\\
 &\le C_\alpha
       \left(\sum_k\|\zeta h_k\|_{H^1(\mathbb R)}^2\right)^{1/2}
       \left(\sum_k\|\zeta g_k\|_{H^{3/4}(\mathbb R)}^2\right)^{1/2}
 \le C_\alpha\|h\|_{H^1}\|g\|_{H^{3/4}}.
\end{aligned}
\]
Applying the second kernel estimate in the same calculation, with
$z=\varepsilon+i\lambda$ and $z'=\varepsilon'+i\lambda'$, we conclude that
\[
\begin{aligned}
 \|\alpha T(z)\|_{H^1\to H^{-3/4}}&\le C_\alpha,\\
 \|\alpha[T(z)-T(z')]\|_{H^1\to H^{-3/4}}
 &\le C_\alpha(|\lambda-\lambda'|+|\varepsilon-\varepsilon'|)^{1/4}.
\end{aligned}
\]
The factor $|k|^{-1/4}$ in the difference estimate is at most one.
Thus all angular modes are included in these bounds.

\smallskip\noindent
\emph{Step 2: The elliptic points and the separatrix.}
Choose $\eta\in C_c^\infty([0,2))$, equal to one on $[0,1]$.
The small-disk estimate proved in Lemma~\ref{lem:regular-high-frequency}
states, at either elliptic point and for sufficiently small $\rho>0$, that
\[
 \sup_{\lambda\in\mathbb R,\,0<|\varepsilon|\le1}
 \|T(\varepsilon+i\lambda)\eta(q/\rho)h\|_{H^{-1}}
 \le C\left(\int_{q<2\rho}|\nabla h|^2\,dx\,dy\right)^{1/2}.
\]
Only nonzero angular modes enter this estimate, since $T$ annihilates
the zeroth mode. With $\chi_R(S)$ as in Proposition~\ref{prop:tail},
we also have
\[
 \sup_{|\lambda|\le M,\,0<|\varepsilon|\le1}
 \|\chi_R(S)T(\varepsilon+i\lambda)\|_{H^1\to H^{-1}}
 \le C_Me^{-cR}.
\]
These cutoffs commute with $T$. Together with Step~1 and a fixed
partition into these regions, the two bounds imply
$\|T(z)\|_{H^1\to H^{-1}}\le C_M$.
For fixed $h\in H^1$, if $z,z'\to i\lambda_0$ through the same side,
the remaining regular-annulus contribution tends to zero, so
\[
\begin{aligned}
 \limsup_{z,z'\to i\lambda_0}\|[T(z)-T(z')]h\|_{H^{-1}}
 &\le C\sum_{\mathrm{centers}}
       \left(\int_{q<2\rho}|\nabla h|^2\,dx\,dy\right)^{1/2}
       +C_Me^{-cR}\|h\|_{H^1}.
\end{aligned}
\]
Letting $\rho\to0$ and then $R\to\infty$ proves the Cauchy
property and joint strong continuity up to either boundary side.
This argument includes $\lambda_0=0$ and every integer threshold.

\smallskip\noindent
\emph{Step 3: Compactness and adjoints.}
For a smooth cutoff $\chi$ supported away from the centers,
Step~1, with a larger streamline cutoff inserted, implies
\[
 \|\chi[1-\chi_R(S)]T(z)\|_{H^1\to H^{-3/4}}
 \le C_{\chi,R},\qquad
 \|\chi\chi_R(S)T(z)\|_{H^1\to H^{-1}}
 \le C_{\chi,M}e^{-cR}.
\]
The first family has relatively compact image of the $H^1$ unit ball
in $H^{-1}$, since $H^{-3/4}\Subset H^{-1}$ on the torus.
Uniform approximation as $R\to\infty$ proves the asserted collective
compactness. Finally, for $\varepsilon\ne0$,
\[
 T(z)^*=T(-\overline z)=T(-\varepsilon+i\lambda).
\]
The uniform bounds and strong continuity for the adjoints follow
from the same results for the opposite boundary side.
\end{proof}

\begin{proposition}[Limiting absorption principle]\label{global:LAP}
For every compact $I\subset\mathbb R\setminus\{0\}$,
\begin{equation}\label{global:bound}
 \sup_{\lambda\in I,\,0<|\varepsilon|\le1}
 \|\Delta^{-1}(\varepsilon+i\lambda-L)^{-1}f\|_{H^1}
 \le C_I\|f\|_{H^1},\qquad f\in H^1\cap Y.
\end{equation}
Both boundary values exist strongly in $H^1$ and are strongly
continuous in $\lambda\in I$ with the boundary side fixed.
\end{proposition}
\begin{proof}
Write
$z=\varepsilon+i\lambda$, with $\lambda\in I$ and
$0\le|\varepsilon|\le1$.

\smallskip\noindent
\emph{Step 1: Injectivity.}
On mean-zero $L^2$, for $\operatorname{Re}z\ne0$,
\[
 (z-L)(z-L_0)^{-1}
 =I-L_0\Delta^{-1}(z-L_0)^{-1}
\]
is a compact perturbation of the identity, since
$L_0\Delta^{-1}:L^2\to H^1$.
If $L\omega=z\omega$, skew-adjointness of $L_0$ implies
\[
 (\operatorname{Re}z)\langle(I+\Delta^{-1})\omega,\omega\rangle
 =\operatorname{Re}\langle L_0(I+\Delta^{-1})\omega,
                              (I+\Delta^{-1})\omega\rangle=0.
\]
As $I+\Delta^{-1}\ge0$, we obtain $(I+\Delta^{-1})\omega=0$,
and then $z\omega=L\omega=0$. The Fredholm alternative therefore
implies invertibility of $z-L$.

It follows that $\Delta+T(z)$ has no nonzero mean-zero $H^1$
nullvector when $\operatorname{Re}z\ne0$. Indeed, $T(z)$ is bounded
on $L^2$, so such a nullvector belongs to $H^2$ and satisfies
\[
 L_0\Delta h=z\Delta h-L_0h\in L^2,
 \qquad (z-L)\Delta h=0.
\]
For $z=i\lambda\pm0$, $\lambda\ne0$, the same injectivity follows
from Proposition~\ref{prop:nullvector} and Theorem~\ref{Intm6}.

\smallskip\noindent
\emph{Step 2: A uniform residual estimate.}
Fix the two center disks from Lemma~\ref{global:residual}.
We claim that, on either closed boundary side (for $\epsilon$),
\begin{equation}\label{global:coercive-residual}
 \|h\|_{H^1}\le C_I\|(\Delta+T(z))h\|_{H^{-1}},
 \qquad \lambda\in I,\quad 0\le|\varepsilon|\le1,
\end{equation}
provided $h$ has zero spatial mean and $(\Delta+T(z))h$ has zero
angular mean in both center disks. At $\varepsilon=0$, $T(z)$
denotes the boundary value on the chosen side.

Otherwise there are $z_j\to z_*$ on that side and mean-zero $h_j$ with
\[
 \|h_j\|_{H^1}=1,\qquad
 r_j=(\Delta+T(z_j))h_j\to0\quad\hbox{in }H^{-1},
 \qquad h_j\rightharpoonup h_*\text{ in }H^1,
 \quad h_j\to h_*\text{ in }L^2.
\]
Here each $r_j$ has zero angular mean in the two disks.
For every $\varphi\in H^1$, Lemma~\ref{global:free} implies
\[
 \langle T(z_j)h_j,\varphi\rangle
 =\langle h_j,T(z_j)^*\varphi\rangle
 \longrightarrow\langle h_*,T(z_*)^*\varphi\rangle.
\]
Thus $(\Delta+T(z_*))h_*=0$, and Step~1 implies $h_*=0$.
Choose a fixed real smooth cutoff $\chi$, zero near the two centers
and equal to one outside their disks $q<q_0/2$.
Collective compactness and the preceding weak limit imply
\[
 \|\chi T(z_j)h_j\|_{H^{-1}}\longrightarrow0.
\]
Using the exact localization identity
\[
 \Delta(\chi h_j)=\chi r_j-\chi T(z_j)h_j
                  +2\operatorname{div}(h_j\nabla\chi)-h_j\Delta\chi,
\]
and the elliptic estimate on the torus, we obtain
\[
\begin{aligned}
 \|\chi h_j\|_{H^1}
 &\le C\bigl(\|\Delta(\chi h_j)\|_{H^{-1}}+\|\chi h_j\|_2\bigr)\\
 &\le C_\chi\bigl(\|r_j\|_{H^{-1}}
           +\|\chi T(z_j)h_j\|_{H^{-1}}+\|h_j\|_2\bigr)
 \longrightarrow0.
\end{aligned}
\]
In particular, $h_j\to0$ in $H^1$ on both annuli
$q_0/2<q<q_0$. Lemma~\ref{global:residual} now implies, at each center,
\[
 \|h_j\|_{H^1(q<q_0/2)}
 \le C\bigl(\|r_j\|_{H^{-1}}
            +\|h_j\|_{H^1(q_0/2<q<q_0)}\bigr)\longrightarrow0.
\]
This contradicts $\|h_j\|_{H^1}=1$ and proves
\eqref{global:coercive-residual}.

\smallskip\noindent
\emph{Step 3: The resolvent and its boundary values.}
For $\operatorname{Re}z\ne0$, let $\psi=\Delta^{-1}(z-L)^{-1}f$.
The resolvent equation and $T(z)-I=z(L_0-z)^{-1}$ imply
\[
 (\Delta+T(z))\psi=-\frac{T(z)-I}{z}f,
 \qquad
 \left\|\frac{T(z)-I}{z}f\right\|_{H^{-1}}
 \le C_I\|f\|_{H^1}.
\]
The source has zero angular mean because $f\in Y$ and every
$T(z)$ has zero angular mean in its range.
Applying \eqref{global:coercive-residual} proves \eqref{global:bound}.

Fix $f\in H^1\cap Y$ and let $z_j\to z_*=i\lambda_*$ through one
side, allowing $\lambda_j\to\lambda_*\in I$.
By the uniform bound, every subsequence of $\psi_j=\Delta^{-1}(z_j-L)^{-1}f$
has a weak $H^1$ limit along a further subsequence.
The adjoint continuity argument from Step~2 shows that each such limit
satisfies
\[
 (\Delta+T(z_*))\psi_*=-\frac{T(z_*)-I}{z_*}f.
\]
Step~1 implies uniqueness, so the limit is independent of the
subsequence. Subtracting the equations, we obtain
\[
\begin{aligned}
 (\Delta+T(z_j))(\psi_j-\psi_*)
 ={}&-\left(\frac{T(z_j)-I}{z_j}
            -\frac{T(z_*)-I}{z_*}\right)f\\
    &+[T(z_*)-T(z_j)]\psi_*\longrightarrow0
       \quad\hbox{in }H^{-1}.
\end{aligned}
\]
Every term on the right has zero angular mean. Hence
\eqref{global:coercive-residual} implies
$\|\psi_j-\psi_*\|_{H^1}\to0$.
After the boundary values have been constructed, the same argument
applies to $z_j=i\lambda_j\pm0$ on a fixed side. This proves their
strong continuity in $\lambda$.
\end{proof}

Away from integers, norm compactness permits a Fredholm formulation.
\begin{proposition}\label{thm:LAP}
Let $I\subset\mathbb R\setminus\mathbb Z$ be compact and
$K(z)=\Delta^{-1}T(z)$. There exist $\varepsilon_I,C_I>0$ such that
$I+K(\varepsilon+i\lambda)$ is invertible on mean-zero $H^1$ and
\[
 \sup_{\lambda\in I,\,0\le|\varepsilon|\le\varepsilon_I}
 \bigl\|[I+K(\varepsilon+i\lambda)]^{-1}\bigr\|_{H^1\to H^1}
 \le C_I.
\]
At $\varepsilon=0$, both boundary values are included.
\end{proposition}
\begin{proof}
By Proposition~\ref{prop:compact}, $K(i\lambda\pm0)$ is compact
and depends continuously on $\lambda$ in operator norm.
For mean-zero $h\in H^1$, Proposition~\ref{prop:nullvector}
and Theorem~\ref{Intm6} imply
\[
 [I+K(i\lambda\pm0)]h=0
 \ \Longrightarrow\ (\Delta+T(i\lambda\pm0))h=0
 \ \Longrightarrow\ (L-i\lambda)\Delta h=0
 \ \Longrightarrow\ h=0.
\]
The Fredholm alternative therefore implies invertibility. Continuity
of inversion and compactness of $I$ imply
\[
 \sup_{\lambda\in I,\,\pm}
 \bigl\|[I+K(i\lambda\pm0)]^{-1}\bigr\|_{H^1\to H^1}
 \le C_I.
\]
By the norm continuity in Proposition~\ref{prop:compact}, we may
choose $\varepsilon_I>0$ so that, on either fixed side,
\[
 \sup_{\lambda\in I,\,0<\pm\varepsilon\le\varepsilon_I}
 \bigl\|[I+K(i\lambda\pm0)]^{-1}
       [K(\varepsilon+i\lambda)-K(i\lambda\pm0)]\bigr\|_{H^1\to H^1}
 \le\tfrac12.
\]
The factorization
\[
\begin{aligned}
 I+K(\varepsilon+i\lambda)
 ={}&[I+K(i\lambda\pm0)]\\
 &\quad\cdot\left\{I+[I+K(i\lambda\pm0)]^{-1}
       [K(\varepsilon+i\lambda)-K(i\lambda\pm0)]\right\}
\end{aligned}
\]
and a Neumann series then imply the inverse bound $2C_I$, proving
the assertion after enlarging $C_I$.
\end{proof}

Differentiating the resolvent produces arbitrary $H^{-1}$ sources.
The residual estimate in Proposition~\ref{global:LAP} extends to
these sources by subtracting a function in $\ker L_0$.

\begin{lemma}\label{r2:full-source-LAP}
For every compact $I\subset\mathbb R\setminus\{0\}$, the operators
\[
 \Delta+T(\varepsilon+i\lambda):H^1(\mathbb T^2)
                  \longrightarrow H^{-1}(\mathbb T^2)
\]
are isomorphisms between the spaces of zero spatial mean, and
\begin{equation}\label{r2:full-source-estimate}
 \|h\|_{H^1}
 \le C_I\|(\Delta+T(\varepsilon+i\lambda))h\|_{H^{-1}},
 \qquad \lambda\in I,\quad 0\le|\varepsilon|\le1.
\end{equation}
At $\varepsilon=0$, the assertion holds for both boundary values,
including integer thresholds.
\end{lemma}
\begin{proof}
All global spaces in this proof have zero spatial mean. Fix
$z=\varepsilon+i\lambda$ and set $r=(\Delta+T(z))h$.
The subspace $H^1\cap\ker L_0$ is closed, since
$L_0:H^1\to L^2$ is bounded. By Poincar\'e and Lax--Milgram,
there is a unique $w\in H^1\cap\ker L_0$ satisfying
\[
 \int_{\mathbb T^2}\nabla w\cdot\nabla\overline\varphi\,dx\,dy
 =-\langle r,\varphi\rangle,
 \qquad \varphi\in H^1\cap\ker L_0.
\]
In particular,
\[
\begin{aligned}
 \|\nabla w\|_{L^2}^2
 &=-\operatorname{Re}\langle r,w\rangle
 \le\|r\|_{H^{-1}}\|w\|_{H^1}
 \le C\|r\|_{H^{-1}}\|\nabla w\|_{L^2},\\
 \|w\|_{H^1}&\le C\|r\|_{H^{-1}},\qquad T(z)w=0.
\end{aligned}
\]
The last identity also holds at the boundary by
Lemma~\ref{global:free}.

The variational equation implies
\[
 \langle r-\Delta w,\varphi\rangle=0
 \qquad(\varphi\in H^1\cap\ker L_0).
\]
For an $H^1$ test function $\varphi$ supported in either center disk
and depending only on $q$, subtracting its spatial mean makes it
admissible in this identity. Since $r-\Delta w$ annihilates constants,
\[
 \langle r-\Delta w,\varphi\rangle
 =\left\langle r-\Delta w,
 \varphi-\frac1{|\mathbb T^2|}\int_{\mathbb T^2}\varphi\,dx\,dy
 \right\rangle=0.
\]
Thus $r-\Delta w$ has zero angular mean in both disks. Consequently,
\[
 (\Delta+T(z))(h-w)=r-\Delta w,
\]
and the residual estimate \eqref{global:coercive-residual} applies
to $h-w$. We obtain
\[
\begin{aligned}
 \|h\|_{H^1}
 &\le\|h-w\|_{H^1}+\|w\|_{H^1}\\
 &\le C_I\|r-\Delta w\|_{H^{-1}}+\|w\|_{H^1}\\
 &\le C_I\bigl(\|r\|_{H^{-1}}+\|w\|_{H^1}\bigr)
 \le C_I\|r\|_{H^{-1}}.
\end{aligned}
\]
This proves \eqref{r2:full-source-estimate}.

The operator is bounded by Lemma~\ref{global:free}, and
\eqref{r2:full-source-estimate} implies that its range is closed.
If a mean-zero $u\in H^1$ annihilates this range, then
\[
 0=\langle(\Delta+T(z))h,u\rangle
   =\langle h,(\Delta+T(-\overline z))u\rangle
 \qquad\text{for every mean-zero }h\in H^1.
\]
At the boundary, $T(i\lambda\pm0)^*=T(i\lambda\mp0)$.
Hence $(\Delta+T(-\overline z))u=0$, and the same estimate on the
opposite side implies $u=0$. The range is thus dense and closed
in mean-zero $H^{-1}$, proving surjectivity and the uniform inverse bound.
\end{proof}

\subsection{The resolvent near zero}
Near $z=0$ we expand $T(z)$. We denote by $L_0^{-1}$ the inverse with zero
streamline average, extended by zero on $\ker L_0$. Thus, for
$m\ge1$,
\[
 (L_0^{-m}h)_k=\left(\frac{S}{2\pi ik}\right)^mh_k
 \quad(k\ne0),\qquad (L_0^{-m}h)_0=0.
\]
The projection $\Pi_{\ker L_0}$ is the orthogonal projection onto
$\ker L_0$, equivalently the streamline average. All Sobolev spaces
in this subsection are taken with zero spatial mean.

\begin{proposition}\label{prop:zero-expansion}
For every fixed $m\ge1$,
\begin{equation}\label{zeroT:D-bound}
 \|L_0^{-m}h\|_2\le C_m\|h\|_{H^1},\qquad
 \|L_0^{-m}T(z)\|_{H^1\to H^{-1}}\le C_m
\end{equation}
for sufficiently small $|z|$ with $\operatorname{Re}z\ne0$;
the second bound also holds for $m=0$. Consequently, for every $N$,
\begin{equation}\label{zeroT:expansion}
 T(z)=I-\Pi_{\ker L_0}+\sum_{j=1}^Nz^jL_0^{-j}
              +z^{N+1}L_0^{-(N+1)}T(z).
\end{equation}
Both boundary families have the same operator-norm asymptotic
coefficients at zero, and, for every $m\ge0$,
\begin{equation}\label{eq:flat-jump}
 \|L_0^{-m}(T(i\lambda+0)-T(i\lambda-0))\|_{H^1\to H^{-1}}
          \le C_m e^{-c_m/|\lambda|}.
\end{equation}
\end{proposition}
\begin{proof}
We first note that
\[
 \|L_0^{-m}h\|_2\lesssim_m\|S^mh\|_2
                         \lesssim_m\|h\|_{H^1}.
\]
On a fixed region of bounded $S$, the combined multiplier
$L_0^{-m}T(z)$ is nonresonant for small $|z|$; near the separatrix we use
the polynomially weighted version of Proposition~\ref{prop:tail}.
The finite geometric identity on each nonzero angular mode implies
\eqref{zeroT:expansion}. For the jump, choose the period cutoff
$R=\pi/(2|\lambda|)$. There are no poles on $S\le R+1$, while
the weighted tail is $O(e^{-cR})$. This proves
\eqref{eq:flat-jump}.
\end{proof}

\begin{corollary}\label{global:projected}
Let $f\in H^1\cap Y$ satisfy $\langle f,g_j\rangle=0$ for
$j=1,2,3$, and let
$\psi_f(z)=\Delta^{-1}(z-L)^{-1}f$.
For every $M>0$,
\[
 \sup_{\substack{|\lambda|\le M\\0<|\varepsilon|\le1}}
 \|Q\psi_f(\varepsilon+i\lambda)\|_{H^1}
 \le C_M\|f\|_{H^1}.
\]
The projected resolvent has strong $H^1$ boundary values including
zero, and the two boundary families have the same value at zero.
Moreover, there exist $\lambda_0,c,C>0$, independent of $f$, such that
\begin{equation}\label{zeroT:exponential-jump}
 \|\psi_{f,+}(\lambda)-\psi_{f,-}(\lambda)\|_{H^1}
 \le Ce^{-c/|\lambda|}\|f\|_{H^1},
 \qquad 0<|\lambda|<\lambda_0.
\end{equation}
\end{corollary}

\begin{proof}
All Sobolev spaces below have zero spatial mean. We first consider
$0<|z|\ll1$; the estimates are uniform for $\operatorname{Re}z\ne0$
and for both boundary families. The identities $Le_i=0$ and
$Lb_i=e_i$ imply
\[
 (z-L)(e_i+zb_i)=z^2b_i.
\]
This cancellation motivates us to write for $(\Delta+T(z))h=F$,
\begin{equation}\label{zeroT:decomposition}
 h=u+\sum_{i=1}^3a_i
       \left(\frac{e_i}{z}-\Delta^{-1}b_i\right),
 \qquad u\perp E.
\end{equation}
Since $\Delta^{-1}b_i\perp E$, this decomposition is unique.
The coefficients $a_i$ are $z$ times the coefficients of $Ph$.

\smallskip\noindent
\emph{Step 1: The coercive estimate and the three scalar equations.}
The only first harmonic in $E^\perp$ is a multiple of $\Psi$,
and $\Pi_{\ker L_0}\Psi=\Psi$. Hence
\begin{equation}\label{zeroT:coercivity}
\begin{aligned}
 -\langle(\Delta+I-\Pi_{\ker L_0})u,u\rangle
 &=\|\nabla u\|_2^2-\|(I-\Pi_{\ker L_0})u\|_2^2\\
 &\ge\tfrac12\|\nabla u\|_2^2,
 \qquad u\perp E.
\end{aligned}
\end{equation}
By \eqref{zero:full-duality},
\[
 (\Delta+I-\Pi_{\ker L_0})\Delta^{-1}b_i=g_i,
 \qquad g_i-\Delta^{-1}b_i=b_i.
\]
The expansion \eqref{zeroT:expansion} therefore implies
\begin{equation}\label{zeroT:corrected-modes}
 (\Delta+T(z))
 \left(\frac{e_i}{z}-\Delta^{-1}b_i\right)
 =zL_0^{-1}b_i+z^2L_0^{-2}b_i+O_{H^{-1}}(|z|^3).
\end{equation}
All terms are well defined: by \eqref{zeroT:D-bound},
\[
 L_0^{-m}b_i=L_0^{-(m+1)}e_i-L_0^{-m}\Delta^{-1}b_i
 \in L^2,\qquad m\ge1.
\]
Testing the equation for $h$ against $u$, and using
\eqref{zeroT:coercivity}--\eqref{zeroT:corrected-modes}, we obtain
\[
 \|u\|_{H^1}
 \le C\|F\|_{H^{-1}}+C|z|\sum_i|a_i|.
\]
Since $L_0^{-1}$ is skew-adjoint,
\[
 \langle(\Delta+T(z))u,e_j\rangle
 =-z\langle u,g_j\rangle+O(|z|^2\|u\|_{H^1}),
 \qquad
 \langle L_0^{-1}b_i,e_j\rangle=-\langle b_i,g_j\rangle.
\]
Testing the equation for $h$ against $e_j$ and dividing by $z$ thus implies
\begin{equation}\label{zeroT:coefficients}
\begin{aligned}
 \sum_i a_i\langle b_i,g_j\rangle
 ={}&-\frac{\langle F,e_j\rangle}{z}-\langle u,g_j\rangle\\
 &+O\!\left(|z|\left[\|u\|_{H^1}+\sum_i|a_i|\right]\right),
 \qquad j=1,2,3.
\end{aligned}
\end{equation}
The matrix $(\langle b_i,g_j\rangle)$ is positive definite by
\eqref{zero:projection}. Consequently,
\[
 \sum_i|a_i|
 \le C\left(|z|^{-1}\|PF\|_{H^{-1}}+\|u\|_{H^1}\right).
\]
Combining the last two estimates and absorbing the small terms,
we conclude that
\begin{equation}\label{zeroT:coefficient-bounds}
\begin{aligned}
 \|u\|_{H^1}&\le C\|F\|_{H^{-1}},\\
 \sum_i|a_i|&\le C\left(
       \|F\|_{H^{-1}}+|z|^{-1}\|PF\|_{H^{-1}}\right).
\end{aligned}
\end{equation}
The limiting system consists of the coercive equation for $u$
and the invertible scalar system for $a_i$. The exact equations,
with their $P$-components divided by $z$, differ from this system
by $O(|z|)$ as maps from $(H^1\cap E^\perp)\times\mathbb C^3$
to $QH^{-1}\times\mathbb C^3$.
They are therefore solvable for every $F\in H^{-1}$ when $|z|$
is sufficiently small. In particular,
\begin{equation}\label{zeroT:inverse}
 \|(\Delta+T(z))^{-1}\|_{H^{-1}\to H^1}\le C|z|^{-2}.
\end{equation}

\smallskip\noindent
\emph{Step 2: The moment conditions and the limit at zero.}
For $h=\psi_f(z)$, the source is
\begin{equation}\label{zeroT:stream-equation}
 F(z)=-L_0^{-1}T(z)f
 =-L_0^{-1}f-zL_0^{-2}f
       +O_{H^{-1}}(|z|^2\|f\|_{H^1}).
\end{equation}
The moment conditions imply
\begin{equation}\label{zeroT:moment-cancellation}
 \langle F(0),e_j\rangle
 =-\langle L_0^{-1}f,e_j\rangle
 =\langle f,g_j\rangle=0.
\end{equation}
Thus \eqref{zeroT:coefficient-bounds} implies
\[
 \|u\|_{H^1}+\sum_i|a_i|\le C\|f\|_{H^1}.
\]
Since $Q\psi_f=u-\sum_i a_i\Delta^{-1}b_i$, we obtain
\begin{equation}\label{zeroT:projected-LAP}
 \|Q\psi_f(z)\|_{H^1}+|z|\|P\psi_f(z)\|_{H^1}
 \le C\|f\|_{H^1}.
\end{equation}
Let $u_0\perp E$ and $a_i^0\in\mathbb C$ be determined successively by
\begin{equation}\label{zeroT:limiting-equations}
\begin{aligned}
 (\Delta+I-\Pi_{\ker L_0})u_0&=-L_0^{-1}f,\\
 \sum_i a_i^0\langle b_i,g_j\rangle
 &=\langle L_0^{-2}f,e_j\rangle-\langle u_0,g_j\rangle.
\end{aligned}
\end{equation}
The first equation is solvable by \eqref{zeroT:moment-cancellation}
and coercivity. Subtracting it from the $Q$-equation implies
\[
 \|u-u_0\|_{H^1}\le C|z|\|f\|_{H^1}.
\]
Then \eqref{zeroT:coefficients} implies
\[
 \sum_i|a_i-a_i^0|\le C|z|\|f\|_{H^1},\qquad
 Q\psi_f(z)\longrightarrow u_0-\sum_i a_i^0\Delta^{-1}b_i
 \quad\text{in }H^1.
\]
Both boundary families have this same limit.

We also keep the next term, which is used below. Let $u_1\perp E$
and $a_i^1\in\mathbb C$ solve
\[
\begin{aligned}
 (\Delta+I-\Pi_{\ker L_0})u_1
 &=-Q\left(L_0^{-2}f+L_0^{-1}u_0
                       +\sum_i a_i^0L_0^{-1}b_i\right),\\
 \sum_i a_i^1\langle b_i,g_j\rangle
 &=\langle L_0^{-3}f,e_j\rangle-\langle u_1,g_j\rangle
   +\langle L_0^{-2}u_0,e_j\rangle
   +\sum_i a_i^0\langle L_0^{-2}b_i,e_j\rangle.
\end{aligned}
\]
Their norms are bounded by $C\|f\|_{H^1}$. By
\eqref{zeroT:expansion}, the residual
\[
 R(z)=F(z)-(\Delta+T(z))
 \left[u_0+zu_1+\sum_i(a_i^0+za_i^1)
             \left(\frac{e_i}{z}-\Delta^{-1}b_i\right)\right]
\]
satisfies
\[
 \|R(z)\|_{H^{-1}}+|z|^{-1}\|PR(z)\|_{H^{-1}}
 \le C|z|^2\|f\|_{H^1}.
\]
Applying \eqref{zeroT:coefficient-bounds} to the difference, we obtain
\[
 \|u-u_0-zu_1\|_{H^1}
 +\sum_i|a_i-a_i^0-za_i^1|
 \le C|z|^2\|f\|_{H^1}.
\]
Hence, with
$a_{-1}=\sum_i a_i^0e_i$ and
$a_0=u_0+\sum_i a_i^1e_i-\sum_i a_i^0\Delta^{-1}b_i$,
\begin{equation}\label{zeroT:full-expansion}
\begin{aligned}
 \psi_f(z)&=\frac{a_{-1}}{z}+a_0
                  +O_{H^1}(|z|\|f\|_{H^1}),\qquad a_{-1}\in E,\\
 \|a_{-1}\|_{H^1}+\|a_0\|_{H^1}&\le C\|f\|_{H^1}.
\end{aligned}
\end{equation}
The coefficients are the same for both boundary families.

\smallskip\noindent
\emph{Step 3: The jump estimate.}
Subtracting the two boundary equations, we obtain
\[
\begin{aligned}
 &(\Delta+T(i\lambda+0))
                   (\psi_{f,+}-\psi_{f,-})\\
 &\qquad=-L_0^{-1}[T(i\lambda+0)-T(i\lambda-0)]f
          -[T(i\lambda+0)-T(i\lambda-0)]\psi_{f,-}.
\end{aligned}
\]
By \eqref{eq:flat-jump}, \eqref{zeroT:inverse}, and
\eqref{zeroT:projected-LAP},
\[
\begin{aligned}
 \|\psi_{f,+}-\psi_{f,-}\|_{H^1}
 &\le C|\lambda|^{-2}e^{-c_0/|\lambda|}
       \left(\|f\|_{H^1}+\|\psi_{f,-}\|_{H^1}\right)\\
 &\le C|\lambda|^{-3}e^{-c_0/|\lambda|}\|f\|_{H^1}
 \le Ce^{-c/|\lambda|}\|f\|_{H^1},\qquad 0<c<c_0.
\end{aligned}
\]
This proves \eqref{zeroT:exponential-jump}. Proposition~\ref{global:LAP}
near nonzero boundary parameters and resolvent continuity on the
remaining compact sets complete the proof.
\end{proof}

\section{Spectral representation formula and preliminary bounds}
\label{sec:spectral-preliminaries}
In this section we first express the evolution through the spectral density function and
identify its stationary part. We then prove the quadratic
critical-trace bound in Proposition~\ref{prop:weighted-center}.
The final subsection records a weighted energy estimate used later
in the spectral derivative bounds.

\subsection{The spectral density and representation formula}
\begin{proposition}\label{qual:theorem}
Let $f\in H^1\cap Y$ satisfy $\langle f,g_j\rangle=0$ for
$j=1,2,3$. Define
\begin{equation}\label{lap:density}
 \Gamma_f(\lambda)
 =\frac{\psi_{f,+}(\lambda)-\psi_{f,-}(\lambda)}{2\pi},
 \qquad \lambda\ne0.
\end{equation}
Then
\[
 \int_{\mathbb R}\|\Gamma_f(\lambda)\|_{H^1}\,d\lambda
 \le C\|f\|_{H^1}.
\]
There is a unique $p_\infty\in E$ such that, for every $t\ge0$,
\begin{equation}\label{lap:spectral-representation}
 \Delta^{-1}e^{-tL}f
 =-p_\infty+\int_{\mathbb R}e^{-it\lambda}
                       \Gamma_f(\lambda)\,d\lambda,
\end{equation}
where the integral converges absolutely in $H^1$. Moreover,
\begin{equation}\label{lap:residue}
\begin{aligned}
 &p_\infty=Pf+\int_{\mathbb R}P\Gamma_f(\lambda)\,d\lambda,\\
 &\|p_\infty\|_{H^s}\le C_s\|f\|_{H^1},\qquad s\ge0.
\end{aligned}
\end{equation}
For $f_*=f-p_\infty$, there is $\lambda_0>0$ such that
\begin{equation}\label{zp:bounded-full}
 \sup_{0<|\lambda|<\lambda_0}
 \|\psi_{f_*,\pm}(\lambda)\|_{H^1}\le C\|f\|_{H^1}.
\end{equation}
Finally,
\[
 \Delta^{-1}e^{-tL}f\longrightarrow-p_\infty
 \quad\text{in }H^1,\qquad
 e^{-tL}f\rightharpoonup p_\infty
 \quad\text{in }L^2,
 \qquad t\to\infty.
\]
\end{proposition}

\begin{proof}
\emph{Step 1: Integrability of the spectral density.}
By Corollary~\ref{cor:global-high-frequency}, for sufficiently large
$|\lambda|$ and $0<|\varepsilon|\le1$,
\[
 \|K(\varepsilon+i\lambda)\|_{H^1\to H^1}
 \le C|\lambda|^{-1/2},\qquad
 \|(I+K(\varepsilon+i\lambda))^{-1}\|_{H^1\to H^1}\le2.
\]
The resolvent equation implies
\[
 \psi_f(z)=\frac1z
 \left[(I+K(z))^{-1}(I+\Delta^{-1})f-f\right].
\]
Taking the jump and using the inverse identity, we obtain
\[
\begin{aligned}
 2\pi\Gamma_f(\lambda)
 ={}&\frac1{i\lambda}(I+K(i\lambda+0))^{-1}
       [K(i\lambda-0)-K(i\lambda+0)]\\
 &\hspace{35mm}\cdot(I+K(i\lambda-0))^{-1}(I+\Delta^{-1})f.
\end{aligned}
\]
Consequently,
\[
 \|\Gamma_f(\lambda)\|_{H^1}
 \le C|\lambda|^{-3/2}\|f\|_{H^1}
 \qquad (|\lambda|\gg1).
\]
Together with Proposition~\ref{global:LAP} on compact intervals
away from zero and the exponential jump estimate in
Corollary~\ref{global:projected}, this proves the desired
$L^1(\mathbb R;H^1)$ bound.

\smallskip\noindent
\emph{Step 2: The spectral representation.}
The projected resolvent equation is
\[
 Q\psi_f(z)=\Delta^{-1}(z-QL)^{-1}Qf.
\]
Since $-iQL$ is self-adjoint on $Y\cap E^\perp$ in the inner
product $\langle(I+\Delta^{-1})\cdot,\cdot\rangle$, Stone's
formula implies, for $\chi\in C_c^\infty(\mathbb R\setminus\{0\})$,
\begin{equation}\label{lap:Stone-formula}
\begin{aligned}
 \Delta^{-1}\chi(-iQL)Qf
 &=\lim_{\varepsilon\to0+}\frac1{2\pi}
   \int_{\mathbb R}\chi(\lambda)
    [Q\psi_f(\varepsilon+i\lambda)-Q\psi_f(-\varepsilon+i\lambda)]
    \,d\lambda\\
 &=\int_{\mathbb R}\chi(\lambda)Q\Gamma_f(\lambda)\,d\lambda.
\end{aligned}
\end{equation}
The last limit holds in $H^1$ by Proposition~\ref{global:LAP} on
the support of $\chi$. Choose
$\chi_n\in C_c^\infty(\mathbb R\setminus\{0\})$, $0\le\chi_n\le1$,
with $\chi_n(\lambda)\to1$ for every $\lambda\ne0$.
By the moment conditions and Proposition~\ref{zero:root},
$Qf$ has no spectral component at zero, and hence
$\chi_n(-iQL)Qf\to Qf$ in $L^2$. Applying
\eqref{lap:Stone-formula} to $e^{-it\lambda}\chi_n$ and
$\chi_n/(z-i\lambda)$, respectively, the spectral theorem and
the $L^1(H^1)$ bound imply
\begin{equation}\label{lap:projected-representations}
\begin{aligned}
 Q\Delta^{-1}e^{-tL}f
 &=\int_{\mathbb R}e^{-it\lambda}Q\Gamma_f(\lambda)\,d\lambda,\\
 Q\psi_f(z)
 &=\int_{\mathbb R}\frac{Q\Gamma_f(\lambda)}{z-i\lambda}\,d\lambda,
 \qquad \operatorname{Re}z\ne0.
\end{aligned}
\end{equation}

Since $LP=0$ and $P\Delta=-P$, projecting
$(z-L)\Delta\psi_f=f$ and then taking its jump implies
\begin{equation}\label{lap:P-resolvent}
 zP\psi_f(z)=-Pf-PL\Delta Q\psi_f(z),\qquad
 i\lambda P\Gamma_f(\lambda)=-PL\Delta Q\Gamma_f(\lambda).
\end{equation}
The operator $PL\Delta:H^1\to E$ is bounded. Indeed,
\[
 \langle PL\Delta h,e_j\rangle
 =-\langle(\Delta+I)h,L_0e_j\rangle,
 \qquad j=1,2,3,
\]
and $L_0e_j$ is smooth. In particular,
$\lambda P\Gamma_f\in L^1(\mathbb R;E)$.
Define
\[
 p_\infty=Pf+\int_{\mathbb R}P\Gamma_f(\lambda)\,d\lambda.
\]
The bound for $p_\infty$ follows from the $L^1(H^1)$ estimate
and equivalence of Sobolev norms on $E$.

The right-hand side of \eqref{lap:spectral-representation}
has the correct $Q$-component by
\eqref{lap:projected-representations}. Its $P$-component at
$t=0$ equals $-Pf=P\Delta^{-1}f$. Moreover,
\[
\begin{aligned}
 \frac{d}{dt}\int_{\mathbb R}e^{-it\lambda}
                         P\Gamma_f(\lambda)\,d\lambda
 &=\int_{\mathbb R}e^{-it\lambda}
                         PL\Delta Q\Gamma_f(\lambda)\,d\lambda\\
 &=PL\Delta Q\Delta^{-1}e^{-tL}f
 =\frac{d}{dt}P\Delta^{-1}e^{-tL}f.
\end{aligned}
\]
The differentiation is justified by
$\lambda P\Gamma_f\in L^1(\mathbb R;E)$.
This proves \eqref{lap:spectral-representation} and uniqueness
of $p_\infty$.

\smallskip\noindent
\emph{Step 3: Behavior at zero and as $t\to\infty$.}
By \eqref{zeroT:full-expansion},
\[
 \psi_f(z)=\frac{a_{-1}}{z}+O_{H^1}(1)\|f\|_{H^1},
 \qquad a_{-1}\in E,
\]
with the same coefficient for both boundary families. To identify
it, combine \eqref{lap:projected-representations} and
\eqref{lap:P-resolvent}:
\[
 zP\psi_f(z)
 =-Pf+\int_{\mathbb R}
       \frac{i\lambda}{z-i\lambda}P\Gamma_f(\lambda)\,d\lambda.
\]
For real $z>0$, $|i\lambda/(z-i\lambda)|\le1$. Dominated
convergence implies $a_{-1}=-Pf-\int P\Gamma_f\,d\lambda=-p_\infty$.
Since $\psi_{p_\infty}(z)=-p_\infty/z$, this proves
\eqref{zp:bounded-full}.
The $H^1$ convergence as $t\to\infty$ follows from
\eqref{lap:spectral-representation} and the vector-valued
Riemann--Lebesgue lemma. Also,
\[
 \sup_{t\ge0}\|Qe^{-tL}f\|_2\le C\|Qf\|_2,
 \qquad
 Pe^{-tL}f=-P\Delta^{-1}e^{-tL}f.
\]
Thus $e^{-tL}f$ is uniformly bounded in $L^2$. Its convergence
to $p_\infty$ in $H^{-1}$ implies weak convergence in $L^2$.
\end{proof}

\begin{lemma}\label{lap:vorticity-reconstruction}
Under the hypotheses of Proposition~\ref{qual:theorem}, there is
a unique $v\in\mathcal D\cap Y\cap E^\perp$ satisfying
\begin{equation}\label{lap:stationary}
 QLv=Qf,\qquad \langle v,g_j\rangle=0\quad(j=1,2,3),
 \qquad \|v\|_2\le C\|f\|_{H^1}.
\end{equation}
Moreover, $p_\infty=Pf-PLv$, and, for $t\ge0$,
\begin{equation}\label{lap:lift}
\begin{aligned}
 e^{-tL}f-p_\infty&=Qe^{-tL}Qf+PLQe^{-tL}v,\\
 \sup_{t\ge0}\|e^{-tL}f-p_\infty\|_2&\le C\|f\|_{H^1}.
\end{aligned}
\end{equation}
For every $\chi\in C_c^\infty(\mathbb R\setminus\{0\})$, we also have
\begin{equation}\label{lap:cutoff-representation}
\begin{aligned}
 \int_{\mathbb R}e^{-it\lambda}\chi(\lambda)
                  \Gamma_f(\lambda)\,d\lambda
 ={}&\Delta^{-1}\bigl[I+PL(QL)^{-1}\bigr]\\
 &\qquad\cdot Qe^{-tL}\chi(-iQL)Qf,\qquad t\ge0.
\end{aligned}
\end{equation}
Here $(QL)^{-1}\chi(-iQL)$ is the bounded spectral multiplier
with symbol $\chi(\lambda)/(i\lambda)$.
\end{lemma}
\begin{proof}
Let $\mu_f$ be the spectral measure of $Qf$ for $-iQL$ in the
energy inner product. Applying $\Delta+I$ to
\eqref{lap:Stone-formula} and pairing with $Qf\in H^1$ identifies
its density away from zero. The moment conditions exclude an
atom at zero, so $\mu_f$ is absolutely continuous and
\begin{equation}\label{lap:measure-density}
 \frac{d\mu_f}{d\lambda}
 =\langle(\Delta+I)Q\Gamma_f(\lambda),Qf\rangle,
 \qquad \lambda\ne0.
\end{equation}
The exponentially small jump near zero and
$\mu_f(\mathbb R)\le C\|f\|_2^2$ imply
\[
 \int_{\mathbb R}\lambda^{-2}\,d\mu_f(\lambda)
 \le C\|f\|_{H^1}^2.
\]
Spectral calculus with multiplier $(i\lambda)^{-1}$ therefore
defines $v$ and proves \eqref{lap:stationary}; uniqueness follows
from the description of $\ker(QL)$ in Proposition~\ref{zero:root}.

The exponential estimate and the $L^1(H^1)$ bound in
Proposition~\ref{qual:theorem} imply
$\Gamma_f/\lambda\in L^1(\mathbb R;H^1)$. Apply
\eqref{lap:Stone-formula} with exhausting cutoffs divided by
$i\lambda$, and use \eqref{lap:P-resolvent}, to obtain
\[
\begin{aligned}
 \int_{\mathbb R}\frac{Q\Gamma_f(\lambda)}{i\lambda}\,d\lambda
 &=\Delta^{-1}v,\\
 \int_{\mathbb R}P\Gamma_f(\lambda)\,d\lambda
 &=-PL\Delta\int_{\mathbb R}
           \frac{Q\Gamma_f(\lambda)}{i\lambda}\,d\lambda=-PLv.
\end{aligned}
\]
Thus \eqref{lap:residue} implies $p_\infty=Pf-PLv$, and
$f-p_\infty=Lv$. Commutation with the evolution and $LP=0$ imply
\[
 e^{-tL}f-p_\infty=Le^{-tL}v
 =Qe^{-tL}Qf+PLQe^{-tL}v.
\]
The bound in \eqref{lap:lift} follows from the energy estimate
for $Qe^{-tL}$ and boundedness of $PL:L^2\to E$.

Applying \eqref{lap:Stone-formula} to
$e^{-it\lambda}\chi$ and $e^{-it\lambda}\chi/(i\lambda)$,
and using \eqref{lap:P-resolvent}, we obtain
\[
\begin{aligned}
 \int e^{-it\lambda}\chi Q\Gamma_f\,d\lambda
 &=\Delta^{-1}Qe^{-tL}\chi(-iQL)Qf,\\
 \int e^{-it\lambda}\chi P\Gamma_f\,d\lambda
 &=-PL\Delta\int e^{-it\lambda}
                  \frac{\chi}{i\lambda}Q\Gamma_f\,d\lambda\\
 &=-PL(QL)^{-1}Qe^{-tL}\chi(-iQL)Qf.
\end{aligned}
\]
Since $\Delta^{-1}=-I$ on $E$, their sum proves
\eqref{lap:cutoff-representation}.
\end{proof}

\subsection{The quadratic bound at a critical circle}\label{sec:weighted-center}
The limiting absorption principle controls the resolvent in $H^1$.
Our purpose here is to prove quadratic vanishing of the modified
resolvent at a critical circle. For $f\in H^5\cap Y$, put
$M=\|f\|_{H^5}$ and
\[
 g_f(z)=f+z\psi_f(z),\qquad
 g_{f,\pm}(\lambda)=f+i\lambda\psi_{f,\pm}(\lambda).
\]
The trace $g_{f,-,k}(q,kb(q))$ is the amplitude in the vorticity
representation derived in Subsection~\ref{sec:vorticity-representation};
see \eqref{sd:amplitude}. The following proposition is the main
estimate of this subsection. Its cubic refinement is proved in
Lemmas~\ref{cubic:profiles} and~\ref{cubic:first-mode} of
Section~\ref{sec:cubic-centers}.

\begin{proposition}\label{prop:weighted-center}
Let $f\in H^5\cap Y$ satisfy $\langle f,g_j\rangle=0$ for
$1\le j\le3$. In either fixed small center disk, for $k\ne0$ and
$0<\delta<q_0/4$,
\begin{equation}\label{spectral:critical-amplitude}
 |g_{f,\pm,k}(\delta,kb(\delta))|
 \le C\delta^2|k|^{-5/2}M.
\end{equation}
\end{proposition}

We first estimate the source in cylinder measure. The proof of the
proposition then proceeds from angular estimates to a weighted
comparison estimate and the one-dimensional trace inequality.
Recall that $T-I=z(L_0-z)^{-1}$. The resolvent equation implies
\[
 (\Delta+T)(f+z\psi_f)
 =\Delta f+Tf-z(L_0-z)^{-1}f=(\Delta+I)f.
\]
Thus, writing $g=g_f$ when the datum is fixed, we have
\begin{equation}\label{r2:fixed-source}
 (\Delta+T(z))g=H,\qquad H=(\Delta+1)f.
\end{equation}
In particular, $H$ is independent of $z$.

\begin{lemma}\label{r2:cylinder-source}
For $f\in H^3\cap Y$, on a fixed small center disk,
\begin{equation}\label{r2:H-cylinder}
 \|H\|_{H^1(dv\,d\theta;q<q_0)}\le C\|f\|_{H^3},
 \qquad v=\log q.
\end{equation}
\end{lemma}
\begin{proof}
By Lemma~\ref{geom:metric}, $X_v,X_\theta=O(q)$ and
$dx\,dy=\xi_0(q)q^2dv\,d\theta$. Hence
\[
 \|H_v\|_{L^2(dv\,d\theta)}
 +\|H_\theta\|_{L^2(dv\,d\theta)}
 \le C\|\nabla H\|_{L^2(dx\,dy)}\le C\|f\|_{H^3}.
\]
Angular Poincar\'e controls $H-H_0$. For the zero mode, write
$\mathsf G=I+G_1$. Since $f_0=0$,
\[
 H_0=(\xi_0q^2)^{-1}
       \partial_v\bigl((G_1)_{vv}f_v+(G_1)_{v\theta}f_\theta\bigr)_0.
\]
The coefficients after division by $q^2$ are bounded, by
Lemma~\ref{geom:metric}. Since $f_v,f_\theta,f_{vv},f_{v\theta}$
are combinations of $O(q)\nabla f$ and $O(q^2)\nabla^2f$,
\[
 \|H_0\|_{L^2(dv)}
 \le C\bigl\|(f_v,f_\theta,f_{vv},f_{v\theta})
                         \bigr\|_{L^2(dv\,d\theta)}
 \le C\|f\|_{H^2}.
\]
These estimates prove \eqref{r2:H-cylinder}.
\end{proof}

\begin{proof}[Proof of Proposition~\ref{prop:weighted-center}]
\emph{Step 1: Angular estimates.}
Fix $\lambda_0>0$. Lemmas~\ref{global:free} and
\ref{r2:full-source-LAP}, together with
Corollary~\ref{cor:global-high-frequency}, imply, for
$|\lambda|\ge\lambda_0$ and $0<|\varepsilon|\le1$,
\[
 \|h\|_{H^1}\le C\|(\Delta+T(z))h\|_{H^{-1}},\qquad
 \|(L_0-z)^{-1}F\|_{H^{-1}}
       \le C|\lambda|^{-1}\|F\|_{H^1}.
\]
Here all functions have zero spatial mean. The second estimate
uses $(L_0-z)^{-1}=(T-I)/z$; at large $|\lambda|$, the first
uses $(\Delta+T)^{-1}=(I+K)^{-1}\Delta^{-1}$.
The first two commutators have divergence form with smooth
coefficients, so
\[
 \|[\Delta,L_0]h\|_{H^{-1}}
 +\|[[\Delta,L_0],L_0]h\|_{H^{-1}}
 \le C\|h\|_{H^1}.
\]
Since $L_0$ commutes with $T$ and $(L_0-z)^{-1}$,
\[
\begin{aligned}
 (\Delta+T)\psi_f
   &=-(L_0-z)^{-1}f,\\
 (\Delta+T)L_0\psi_f
   &=-(L_0-z)^{-1}L_0f+[\Delta,L_0]\psi_f,\\
 (\Delta+T)L_0^2\psi_f
   &=-(L_0-z)^{-1}L_0^2f+2[\Delta,L_0]L_0\psi_f
                          -[[\Delta,L_0],L_0]\psi_f.
\end{aligned}
\]
Applying the preceding estimates successively proves
\begin{equation}\label{r2:tangential-high}
 \|L_0^j\psi_f\|_{H^1}\le C|\lambda|^{-1}\|f\|_{H^{j+1}},
 \qquad j=0,1,2,\quad |\lambda|\ge\lambda_0>0.
\end{equation}
Each estimate is valid for mean-zero data in $H^{j+1}$, without
the moment conditions. Difference quotients along the smooth
area-preserving flow of $L_0$ justify the commutations in
$\mathcal L(H^1,H^{-1})$, and the uniform bounds pass to both
boundary values. Since $g=f+z\psi_f$ and
$L_0=b(q)\partial_\theta$, with $b$ bounded away from zero in
the center disk, it follows that
\begin{equation}\label{r2:angular-g}
 \sum_{j=0}^2\|L_0^jg\|_{H^1}
 +\sum_{j=0}^2\|\partial_\theta^jg\|_{\mathscr E(q<q_0)}
 \le CM,\qquad |\lambda|\ge\lambda_0.
\end{equation}
For the critical parameters below, $|\lambda|=|k|b(\delta)
\ge |k|/2$, so we may fix $\lambda_0=1/2$.

\emph{Step 2: The weighted inverse.}
It suffices to consider $k>0$. Fix $\lambda=kb(\delta)$, write
$g=g_{f,\pm}(\lambda)$, and use
the translated coordinate
\[
 v=\log(q/\delta),\qquad q=\delta e^v,\qquad
 W(v)=e^{v/2}(1+e^{2v})^{3/4}.
\]
Thus $\partial_v=q\partial_q$ and $dv=dq/q$, as for
$v=\log q$ in Lemma~\ref{r2:cylinder-source}. We retain the
norms from Proposition~\ref{compact:normal}:
\[
 \|u\|_k^2=\int_{\mathbb R}(|u'|^2+k^2|u|^2)\,dv,
 \qquad
 \|F\|_{k,*}=\sup_{\|\varphi\|_k=1}|\langle F,\varphi\rangle|.
\]
We claim that the conjugated operator is bijective from $H^1_k$
to $H^{-1}_k$, with
\begin{equation}\label{r2:weighted-normal}
 \|u\|_k\le C\|W^{-1}\mathcal N_{k,a,\eta}Wu\|_{k,*},
 \qquad k\ge1,\quad |a|+|\eta|=1,
\end{equation}
including both boundary values at $(a,\eta)=(1,0)$.
Here $\mathcal N_{k,a,\eta}$ is the operator in
Proposition~\ref{compact:normal}; explicitly,
\[
 W^{-1}\mathcal N_{k,a,\eta}Wu
 =-u''-2\frac{W'}W u'
  +\left(k^2+\frac{8e^{2v}}{e^{2v}-a-i\eta}
                     -\frac{W''}W\right)u.
\]
This expression defines a bounded map $H^1_k\to H^{-1}_k$
by \eqref{eq:cylinder-form}.

Since
\[
 \frac{W'}W=\frac12+\frac{3e^{2v}}{2(1+e^{2v})},\qquad
 \left(\frac{W'}W\right)^2
 \le\frac14+\frac{15e^{2v}}{4(1+e^{2v})},
\]
integration by parts implies
\[
\begin{aligned}
 &\operatorname{Re}\left\langle
 W^{-1}\left(-\partial_v^2+k^2+
                  \frac{8e^{2v}}{1+e^{2v}}\right)Wu,u\right\rangle\\
 &\quad=\int_{\mathbb R}\left[|u'|^2+
 \left(k^2+\frac{8e^{2v}}{1+e^{2v}}
                   -\left(\frac{W'}W\right)^2\right)|u|^2\right]dv
 \ge\frac34\|u\|_k^2.
\end{aligned}
\]
By \eqref{eq:model-log-identity}, the difference from
$W^{-1}\mathcal N_{k,a,\eta}W$ is the derivative of an
$L^2(dv)$ function plus a smooth exponentially decaying multiplier.
It is compact from $H^1_k$ to $H^{-1}_k$, as proved in
Proposition~\ref{compact:normal}: for $F_m\in C_c^\infty$
converging to $F$ in $L^2$,
\[
 \|(F'-F_m')u\|_{k,*}
 \le C_k\|F-F_m\|_{L^2(dv)}\|u\|_k\longrightarrow0.
\]
The conjugated operator therefore has index zero.

Suppose $W^{-1}\mathcal N_{k,a,\eta}Wu=0$.
The boundedness of $W,W'$ on $(-\infty,1]$ implies
$Wu\in H^1_k$ there. At $+\infty$, the potential is
$8+O(e^{-2v})$, so the two homogeneous solutions have asymptotics
\[
 e^{\pm\sqrt{k^2+8}\,v}(1+o(1)),\qquad v\to+\infty,
\]
also after one derivative. Since $W=e^{2v}(1+O(e^{-2v}))$
and $\sqrt{k^2+8}>2$, the condition $u\in L^2(dv)$ excludes
the growing solution. Thus $Wu\in H^1_k(\mathbb R)$, and
Proposition~\ref{compact:normal} implies $Wu=0$.
The Fredholm alternative proves bijectivity. Moreover,
\eqref{eq:cylinder-form} implies
\[
 \operatorname{Re}\langle W^{-1}\mathcal N_{k,a,\eta}Wu,u\rangle
 \ge(1-4k^{-2}-Ck^{-1/2})\|u\|_k^2.
\]
This proves \eqref{r2:weighted-normal} for large $k$.
For the finitely many remaining $k$, operator-norm continuity
on each of $\{|a|+|\eta|=1,\ \pm\eta\ge0\}$, established
by \eqref{eq:model-log-identity}, implies the same uniform bound.

\emph{Step 3: The source estimates and the critical trace.}
Choose radial cutoffs $\chi,\widetilde\chi$ with
\[
 \chi=1\ \text{on }q\le q_0/2,\qquad
 \operatorname{supp}\chi\subset\{q<3q_0/4\},\qquad
 \widetilde\chi=1\ \text{on }q\le3q_0/4,\qquad
 \operatorname{supp}\widetilde\chi\subset\{q<q_0\}.
\]
All derivatives in this step are in $(v,\theta)$; in particular,
$\chi'=q\partial_q\chi$. The cylinder energy in
\eqref{r2:angular-g} is
\[
 \|h\|_{\mathscr E(q<q_0)}^2
 =\int_{-\infty}^{\log(q_0/\delta)}\int_{\mathbb T}
       (|h_v|^2+|h_\theta|^2)\,d\theta\,dv.
\]
In particular, Parseval's identity and \eqref{r2:angular-g} imply
\begin{equation}\label{r2:localized-angular}
 k^2\|\widetilde\chi g_k\|_k
 \le C\|\partial_\theta^2g\|_{\mathscr E(q<q_0)}\le CM.
\end{equation}
All functions multiplied by these cutoffs are extended by zero
outside the coordinate disk.

The multiplier formula \eqref{eq:T},
$b(q)=1-q^2/8$, and the Dirichlet form in
Lemma~\ref{geom:metric} imply the exact identity
\begin{equation}\label{r2:localized-equation}
\begin{aligned}
 \mathcal N_{k,1,\pm0}(\chi g_k)
 ={}&-\chi\xi_0q^2H_k
      +\chi\bigl[\operatorname{div}((\mathsf G-I)\nabla g)\bigr]_k\\
 &+(1-\xi_0b)\frac{8e^{2v}}{e^{2v}-1\mp i0}\chi g_k
       -2\chi'g_k'-\chi''g_k.
\end{aligned}
\end{equation}
We estimate its right-hand side after multiplication by $W^{-1}$.
First, $q^2/W\le\delta^2$ and \eqref{r2:H-cylinder} imply
\[
\begin{aligned}
 \|W^{-1}\chi\xi_0q^2H_k\|_{k,*}
 &\le C\delta^2 k^{-1}\|H_k\|_{L^2(dv;q<q_0)}\\
 &\le C\delta^2 k^{-2}
       \|\partial_\theta H\|_{L^2(dv\,d\theta;q<q_0)}
 \le C\delta^2 k^{-2}M.
\end{aligned}
\]

For the metric term, the product rule implies
\[
\begin{aligned}
 W^{-1}\chi[\operatorname{div}((\mathsf G-I)\nabla g)]_k
 ={}&[\operatorname{div}(W^{-1}\chi(\mathsf G-I)\nabla g)]_k\\
 &-(W^{-1}\chi)'[((\mathsf G-I)\nabla g)_v]_k.
\end{aligned}
\]
Here the subscript $v$ on the vector field denotes its
$v$-component. Lemma~\ref{geom:metric} and the formula for $W$
imply, for $0\le j\le2$,
\[
 \|\partial_\theta^j[W^{-1}\chi(\mathsf G-I)]\|_{L^\infty}
 +\|\partial_\theta^j[(W^{-1}\chi)'(\mathsf G-I)]\|_{L^\infty}
 \le C\delta^2.
\]
Testing the divergence against $\varphi\in H^1_k$ and using
$\|\varphi\|_{L^2(dv)}\le k^{-1}\|\varphi\|_k$, we obtain
\[
\begin{aligned}
 &\|W^{-1}\chi[\operatorname{div}((\mathsf G-I)\nabla g)]_k\|_{k,*}\\
 &\quad\le\|[W^{-1}\chi(\mathsf G-I)\nabla g]_k\|_{L^2(dv)}
 +k^{-1}\|[(W^{-1}\chi)'((\mathsf G-I)\nabla g)_v]_k\|_{L^2(dv)}\\
 &\quad\le C\delta^2k^{-2}\sum_{j=0}^2
       \|\nabla\partial_\theta^jg\|_{L^2(dv\,d\theta;q<q_0)}
 \le C\delta^2k^{-2}M.
\end{aligned}
\]
The second inequality follows by taking two angular derivatives
before projecting onto the $k$th Fourier mode; the last is
\eqref{r2:angular-g}.

For the potential term, Lemma~\ref{geom:metric} implies
\[
 \|W^{-1}\chi(1-\xi_0b)\|_{W^{1,\infty}(\mathbb R_v)}
 \le C\delta^2.
\]
Apply \eqref{eq:cylinder-form} with
$c=W^{-1}\chi(1-\xi_0b)$, $(a,\eta)=(1,\pm0)$, and
$u=\widetilde\chi g_k$. Taking the supremum over
$\|\varphi\|_k=1$ and using \eqref{r2:localized-angular}, we obtain
\[
 \left\|W^{-1}\chi(1-\xi_0b)
       \frac{8e^{2v}}{e^{2v}-1\mp i0}g_k\right\|_{k,*}
 \le C\delta^2k^{-1/2}\|\widetilde\chi g_k\|_k
 \le C\delta^2k^{-5/2}M.
\]
Finally, on $\operatorname{supp}\chi'\cup\operatorname{supp}\chi''$,
$q_0/2\le q\le3q_0/4$ and $W^{-1}\le C\delta^2$.
Consequently,
\[
 \|W^{-1}(2\chi'g_k'+\chi''g_k)\|_{k,*}
 \le C\delta^2k^{-1}\|\widetilde\chi g_k\|_k
 \le C\delta^2k^{-3}M.
\]

By \eqref{r2:weighted-normal} and the estimates for
\eqref{r2:localized-equation}, the equation
\[
 W^{-1}\mathcal N_{k,1,\pm0}Wu
   =W^{-1}\mathcal N_{k,1,\pm0}(\chi g_k)
\]
has a unique solution $u\in H^1_k$ satisfying
$\|u\|_k\le C\delta^2k^{-2}M$.
As in Step~2, $Wu\in H^1_k$ at the left end.
The right-hand side vanishes for $q\ge3q_0/4$, so the
right-end asymptotics in Step~2 imply $Wu\in H^1_k$ there as
well. Therefore
\[
 Wu-\chi g_k\in H^1_k(\mathbb R),\qquad
 \mathcal N_{k,1,\pm0}(Wu-\chi g_k)=0.
\]
Proposition~\ref{compact:normal} implies $Wu=\chi g_k$, and hence
\begin{equation}\label{r2:weighted-input}
 \|W^{-1}\chi g_k\|_k\le C\delta^2k^{-2}M.
\end{equation}
The fixed cutoff is understood in subsequent whole-cylinder
uses of this estimate. Finally,
\[
 |u(0)|^2\le2\|u\|_{L^2(dv)}\|u'\|_{L^2(dv)}
       \le k^{-1}\|u\|_k^2.
\]
Since $\chi(\delta)=1$ and $W(0)=2^{3/4}$, we conclude that
\begin{equation}\label{r2:strong-onshell}
 |g_{f,\pm,k}(\delta,kb(\delta))|
 =|f_k(\delta)+ikb(\delta)\psi_{f,\pm,k}(\delta,kb(\delta))|
 \le C\delta^2k^{-5/2}M.
\end{equation}
Complex conjugation treats negative $k$, with the datum also
conjugated and $k$ replaced by $|k|$. The cell symmetries imply
the same estimate at the other center. This proves
\eqref{spectral:critical-amplitude}.
\end{proof}

\subsection{A weighted energy estimate}\label{subsec:weighted-energy}
We next record a weighted estimate for the full solution, including
its zeroth angular coefficient. It controls the commutator terms
on shrinking center annuli in
Lemma~\ref{r2:angular-regular-remainder}.

\begin{lemma}\label{lem:weighted-center-energy}
For $0<\alpha<1$, the solution of \eqref{r2:fixed-source} satisfies
\begin{equation}\label{r2:mild-center-weight}
 \|q^{-\alpha}(g_v,g_\theta)\|_{L^2(dv\,d\theta;q<q_0/2)}
 \le C_{\alpha,q_0}\bigl(\|H\|_{L^2(dv\,d\theta;q<q_0)}
                   +\|g\|_{H^1(q_0/2<q<q_0)}\bigr),
\end{equation}
uniformly in the spectral parameter.
\end{lemma}
\begin{proof}
Throughout the proof, $v=\log q$, $\nabla=(\partial_v,\partial_\theta)$,
and $L^2$ norms use $dv\,d\theta$ on the indicated set
(or $dv$ for Fourier coefficients). The $H^1$ norm on the
outer annulus is in physical coordinates. We also allow $\alpha=0$.
For $n\ne0$, set
\[
 \sigma_n=8(1-\lambda/n),\qquad \eta_n=8\varepsilon/n,
 \qquad z=\varepsilon+i\lambda.
\]
We use the formula for $\mathcal N_{|n|,\sigma_n,\eta_n}$ in
Proposition~\ref{compact:normal}, without normalizing its parameters.
The weighted inverse satisfies
\[
 \|u\|_{|n|}\le C_\alpha
 \|q^{-\alpha}\mathcal N_{|n|,\sigma_n,\eta_n}q^\alpha u\|_{|n|,*}.
\]
Indeed,
\[
 q^{-\alpha}\mathcal N_{|n|,\sigma_n,\eta_n}q^\alpha
 =-(\partial_v+\alpha)^2+n^2+
                \frac{8e^{2v}}{e^{2v}-\sigma_n-i\eta_n},
\]
and integration by parts implies
\[
 \operatorname{Re}\left\langle
 \left[-(\partial_v+\alpha)^2+n^2+
             \frac{8e^{2v}}{1+e^{2v}}\right]u,u\right\rangle
 \ge(1-\alpha^2)\|u\|_{|n|}^2.
\]
Thus the compactness and injectivity argument proving
\eqref{r2:weighted-normal} applies with $q^\alpha$ in place of $W$:
the endpoint exponents satisfy $\alpha<|n|$ and
$\alpha<\sqrt{n^2+8}$. Translation by
$\frac12\log(|\sigma_n|+|\eta_n|)$ reduces nonzero parameter pairs
to Proposition~\ref{compact:normal}; the zero pair has constant
potential $8$. The bound is therefore uniform in $n,z$, including
both boundary values.

Choose smooth radial cutoffs $\chi,\widetilde\chi$ such that
$\chi=1$ for $q\le q_0/2$, $\operatorname{supp}\chi\subset
\{q<3q_0/4\}$, $\widetilde\chi=1$ for $q\le3q_0/4$, and
$\operatorname{supp}\widetilde\chi\subset\{q<q_0\}$.
As in \eqref{r2:localized-equation},
\[
\begin{aligned}
 \mathcal N_{|n|,\sigma_n,\eta_n}(\chi g_n)
 ={}&-\chi\xi_0q^2H_n
   +\chi[\operatorname{div}((\mathsf G-I)\nabla g)]_n\\
 &+(1-\xi_0b)\frac{8e^{2v}}{e^{2v}-\sigma_n-i\eta_n}\chi g_n
   -2\chi'g_n'-\chi''g_n.
\end{aligned}
\]
Lemma~\ref{geom:metric} implies
\[
\begin{aligned}
 &\|q^{-\alpha}\chi(\mathsf G-I)\|_{L^\infty}
 +\|(q^{-\alpha}\chi)'(\mathsf G-I)\|_{L^\infty}
 +\|q^{-\alpha}\chi(1-\xi_0b)\|_{W^{1,\infty}(\mathbb R_v)}
 \le C_\alpha q_0^{2-\alpha}.
\end{aligned}
\]
For the divergence term we use the exact product rule
\[
\begin{aligned}
 q^{-\alpha}\chi\operatorname{div}((\mathsf G-I)\nabla g)
 ={}&\operatorname{div}(q^{-\alpha}\chi(\mathsf G-I)\nabla g)\\
 &-(q^{-\alpha}\chi)'((\mathsf G-I)\nabla g)_v.
\end{aligned}
\]
For the potential term, \eqref{eq:cylinder-form}, applied to
$\widetilde\chi g_n$, implies
\[
 \left\|q^{-\alpha}\chi(1-\xi_0b)
       \frac{8e^{2v}}{e^{2v}-\sigma_n-i\eta_n}g_n\right\|_{|n|,*}
 \le C_\alpha q_0^{2-\alpha}|n|^{-1/2}
                                  \|\widetilde\chi g_n\|_{|n|}.
\]
Using $\|F\|_{|n|,*}\le |n|^{-1}\|F\|_{L^2}$ for the
remaining terms and summing squares in $n\ne0$, we obtain
\[
\begin{aligned}
 &\left(\sum_{n\ne0}
 \|q^{-\alpha}\mathcal N_{|n|,\sigma_n,\eta_n}(\chi g_n)\|_{|n|,*}^2
                                                   \right)^{1/2}\\
 &\quad\le C_\alpha q_0^{2-\alpha}
 \bigl(\|H\|_{L^2(q<q_0)}+\|\nabla g\|_{L^2(q<q_0)}\bigr)
       +C_{\alpha,q_0}\|g\|_{H^1(q_0/2<q<q_0)}.
\end{aligned}
\]
Only the unweighted energy occurs on the right.
To justify the weighted estimate, solve
\[
 q^{-\alpha}\mathcal N_{|n|,\sigma_n,\eta_n}q^\alpha u_n
   =q^{-\alpha}\mathcal N_{|n|,\sigma_n,\eta_n}(\chi g_n),
 \qquad u_n\in H^1_{|n|}.
\]
Since $q^\alpha$ is bounded at $-\infty$ and the source vanishes
for $q\ge3q_0/4$, the endpoint argument in the proof of
\eqref{r2:weighted-normal} implies $q^\alpha u_n\in H^1_{|n|}$.
Proposition~\ref{compact:normal} then implies
$q^\alpha u_n=\chi g_n$. Since
$q^{-\alpha}(\chi g_n)'=u_n'+\alpha u_n$, the preceding bound
also controls $\|q^{-\alpha}\nabla(g-g_0)\|_{L^2(q<q_0/2)}$.

For the zero mode, recall that $v=\log q<\log q_0$ and
$\nabla=(\partial_v,\partial_\theta)$. For $h=g,H$, write
$h_0(v)=(2\pi)^{-1}\int_{\mathbb T}h(v,\theta)\,d\theta$.

\smallskip\noindent
\emph{(a) The averaged equation.}
Define
\[
 \mathcal F(v):=\frac1{2\pi}\int_{\mathbb T}
 \left(\mathsf G_{vv}(v,\theta)\partial_vg(v,\theta)
       +\mathsf G_{v\theta}(v,\theta)\partial_\theta g(v,\theta)
 \right)\,d\theta.
\]
Average \eqref{r2:fixed-source} in the cylinder coordinates of
Lemma~\ref{geom:metric}. The angular divergence integrates to
zero by periodicity, $\xi_0(q)q^2$ is radial, and
$[T(z)g]_0=0$. Thus, in distributions,
\begin{equation}\label{r2:zero-averaged-equation}
 \mathcal F'(v)=\xi_0(e^v)e^{2v}\bigl(H_0(v)-[T(z)g]_0(v)\bigr)
              =\xi_0(e^v)e^{2v}H_0(v).
\end{equation}

\smallskip\noindent
\emph{(b) The integration constant at the center.}
The finite physical energy of $g$, the source bound
\eqref{r2:H-cylinder}, metric equivalence, and
\eqref{r2:zero-averaged-equation} imply
\[
 \begin{aligned}
 \|\mathcal F\|_{L^2(dv;q<q_0)}
 &\le C\|\nabla g\|_{L^2(dv\,d\theta;q<q_0)}<\infty,\\
 \|\mathcal F'\|_{L^1(dv;q<q_0)}
 &\le C\left(\int_{-\infty}^{\log q_0}e^{4v}\,dv\right)^{1/2}
           \|H_0\|_{L^2(dv;q<q_0)}<\infty.
 \end{aligned}
\]
Thus $\mathcal F$ has a locally absolutely continuous
representative and a finite limit as $v\to-\infty$.
If this limit were nonzero, its $L^2(dv)$ norm on
$(-\infty,\log q_0)$ would be infinite. Hence
$\lim_{v\to-\infty}\mathcal F(v)=0$, and
\begin{equation}\label{r2:zero-integral-formula}
 \mathcal F(v)
 =\int_{-\infty}^{v}e^{2s}\xi_0(e^s)H_0(s)\,ds,
 \qquad v<\log q_0.
\end{equation}
This argument imposes no condition on the value of $g_0$ at
the center.

\smallskip\noindent
\emph{(c) The weighted estimate.}
Multiplying \eqref{r2:zero-integral-formula} by $e^{-2v}$,
we obtain
\[
 e^{-2v}\mathcal F(v)
 =\int_0^\infty e^{-2s}\xi_0(e^{v-s})H_0(v-s)\,ds.
\]
Extend $\xi_0(e^v)H_0(v)$ by zero for $v\ge\log q_0$.
Young's convolution inequality, $\int_0^\infty e^{-2s}\,ds=1/2$,
and $q^{2-\alpha}\le q_0^{2-\alpha}$ imply, for $0\le\alpha<1$,
\[
 \begin{aligned}
 \|q^{-2}\mathcal F\|_{L^2(dv;q<q_0)}
 &\le\tfrac12\|\xi_0H_0\|_{L^2(dv;q<q_0)},\\
 \|q^{-\alpha}\mathcal F\|_{L^2(dv;q<q_0)}
 &\le q_0^{2-\alpha}\|q^{-2}\mathcal F\|_{L^2(dv;q<q_0)}
 \le Cq_0^{2-\alpha}\|H_0\|_{L^2(dv;q<q_0)}.
 \end{aligned}
\]
Finally, the definition of $\mathcal F$ implies
\[
 \partial_vg_0
 =\mathcal F-\frac1{2\pi}\int_{\mathbb T}
 \left((\mathsf G_{vv}-1)\partial_vg
                  +\mathsf G_{v\theta}\partial_\theta g\right)
 \,d\theta.
\]
Since $|\mathsf G_{vv}-1|+|\mathsf G_{v\theta}|\le Cq^2$
by Lemma~\ref{geom:metric}, Cauchy--Schwarz in $\theta$ and
the preceding estimate for $\mathcal F$ imply
\begin{equation}\label{r2:zero-weighted-estimate}
 \begin{aligned}
 &\|q^{-\alpha}\partial_vg_0\|_{L^2(dv;q<q_0)}\\
 &\qquad\le Cq_0^{2-\alpha}
 \left(\|H\|_{L^2(dv\,d\theta;q<q_0)}
       +\|\nabla g\|_{L^2(dv\,d\theta;q<q_0)}\right).
 \end{aligned}
\end{equation}
The right-hand side involves only the unweighted energy of $g$.
The constant in this zero-mode estimate is independent of the
spectral parameter and of $q_0$ within a fixed sufficiently
small center chart.

\end{proof}

\section{The first spectral derivative and pointwise bounds}\label{sec:higher}
In this section we establish bounds on the first spectral derivative of the spectral density function. As is now well known, taking $\partial_\lambda$ to $\psi_f(i\lambda+\epsilon)$ creates logarithmic singularities. To reduce the strength of the resulting singularities, we use an adapted ``good derivative'' \eqref{good-der} instead, similarly to~\cite{BeekieChenJia2026}. Outside small neighborhoods around the elliptic points,
the adapted derivative does not create any singularities. Due to the degeneracy near the elliptic centers, a cutoff in the adapted derivative is necessary and as a result, there are still log singularities close to these centers. It is important to explicitly track the log singularity for our applications. After subtracting these log singularities, the resulting
source term belongs to physical $H^{-1}$, and therefore the limiting absorption
principle applies.
On the other hand, the explicit log singularities can be bounded in the
spectral representation using integration.

Fix $q_*>0$ sufficiently small. The set
\[
 \left\{q<q_*\right\},\qquad q=\sqrt{8(1-b(r))},
\]
has two connected components, containing $(0,0)$ and $(\pi,\pi)$,
respectively. Fix a sufficiently small $\kappa\in(0,1/256)$ and set
\[
 I_k=(k-\kappa,k+\kappa),\qquad
 J_k=[k+\kappa,k+1-\kappa],\qquad k\in\mathbb Z.
\]
These intervals cover $\mathbb R$, and the intervals $I_k$ are
pairwise disjoint. Choose $\eta\in C^\infty([0,\infty))$ with
$\eta=0$ on $[0,1]$ and $\eta=1$ on $[2,\infty)$.
Define the cutoff length scale for $I_k, J_k$ by
\[
 (q_{0,k}^{I})^2=
 \begin{cases}
  32\kappa/|k|,& k\ne0,\\[2pt]
  \kappa/2,& k=0,
 \end{cases}
 \qquad
 (q_{0,k}^{J})^2=\frac{\kappa}{2(1+|k|)},\qquad k\in\mathbb Z.
\]
For $\lambda\ne0$, set $V=-S(\lambda S')^{-1}\partial_r$ and
\begin{equation}\label{good-der}
 D_\lambda=\partial_\lambda+
 \begin{cases}
  \eta(q/q_{0,k}^{I})V,&\lambda\in I_k\setminus\{0\},\\[2pt]
  \eta(q/q_{0,k}^{J})V,&\lambda\in J_k.
 \end{cases}
\end{equation}
We choose
$\kappa$ sufficiently small such that $2q_{0,k}^{I}<q_*$ and $2q_{0,k}^{J}<q_*$
for every $k$.

Let $f\in H^5\cap Y$ satisfy
$\langle f,g_j\rangle_{L^2}=0$, $1\le j\le3$, and set
$M=\|f\|_{H^5}$. Let $\Gamma_f$ be defined by
\eqref{lap:density}. For $k\ne0$ and $\lambda\in I_k$ with $\lambda/k<1$, set
\[
 \sigma_k(\lambda)=8(1-\lambda/k)>0.
\]
Fix $\chi\in C_c^\infty((-1/2,1/2))$ with $\chi\equiv1$ on $[-1/8,1/8]$.
For either boundary value (as $\epsilon\to+$ or $\epsilon\to0-$), set
\[
 g_{f,\pm}=f+i\lambda\psi_{f,\pm},\qquad
 g_{f,\pm,k}(q,\lambda)
 =\frac1{2\pi}\int_0^{2\pi}
       g_{f,\pm}(q,\theta,\lambda)e^{-ik\theta}\,d\theta.
\]
Define
\begin{equation*}
\begin{aligned}
 \Lambda_f(\lambda,q,\theta)
 ={}&\sum_{\substack{k\in\mathbb Z\setminus\{0\}\\
                    \lambda\in I_k,\ \lambda/k<1}}
 \chi\!\left(\log\frac{q}{\sqrt{\sigma_k(\lambda)}}\right)\sum_{\pm}
 \log\!\left(
       \log\frac{q}{\sqrt{\sigma_k(\lambda)}}
           \mp i0\,\operatorname{sgn}k
      \right)b_{f,k}^{\pm}(\lambda,\theta),
\end{aligned}
\end{equation*}
with
\[
 b_{f,k}^{\pm}(\lambda,\theta)
 =
 \pm\frac{8\xi_0(\sqrt{\sigma_k(\lambda)})}
 {\pi i k^2\sigma_k(\lambda)\,
       \mathsf G_{vv}(\sqrt{\sigma_k(\lambda)},\theta)}
 g_{f,\pm,k}(\sqrt{\sigma_k(\lambda)},\lambda)e^{ik\theta}.
\]
Here $\sum_{\pm}$ denotes the two resolvent boundary values.
The local formula is applied in each center chart and both spatial
contributions are included.
The sum over $k$ has at most one term for each $\lambda$, since
the intervals $I_k$ are disjoint. In
particular, $\Lambda_f=0$ on $I_0$ and on every $J_k$.
The logarithmic contribution occurs only when
\[
 k-\kappa<\lambda<k\quad(k>0),\qquad
 k<\lambda<k+\kappa\quad(k<0),
\]
and its $k$th summand is supported where $\eta(q/q_{0,k}^{I})=0$.

All Sobolev norms below are on $\mathbb T^2$, unless an angular
or spectral variable is explicitly indicated.

\begin{proposition}\label{prop:spectral-bounds}
Let $f\in H^5\cap Y$ satisfy
$\langle f,g_j\rangle_{L^2}=0$, $1\le j\le3$, and set
$M=\|f\|_{H^5}$. Let $\Gamma_f$ be defined by
\eqref{lap:density}. Then we have the following conclusions.

\smallskip
\noindent
(i) The spectral density function, extended by zero at $\lambda=0$, satisfies
\begin{equation}\label{spectral:integrated}
 \|\Gamma_f\|_{L^1(\mathbb R;H^1)}
 +\|\Gamma_f\|_{H^1(\mathbb R;L^2)}
 \le CM.
\end{equation}

\smallskip
\noindent
(ii) For $|\lambda|\ge\kappa$, we have the decomposition
\begin{equation}\label{spectral:decomposition}
 D_\lambda\Gamma_f=\Lambda_f+R_f.
\end{equation}
There is a fixed $K\ge2$ such that for any $0<\alpha<1$,
\[
 \begin{aligned}
 \|R_f(\lambda)\|_{H^1}+\|L_0R_f(\lambda)\|_{H^1}
 &\le C_\alpha M,
 &&\kappa\le|\lambda|\le K,\\
 \|R_f(\lambda)\|_{H^1}
 +\|L_0R_f(\lambda)\|_{H^1}
 &\le C_\alpha|\lambda|^{-1-\alpha/2}M,
 &&|\lambda|\ge K.
 \end{aligned}
\]
In addition,
\begin{equation}\label{spectral:log-coefficient}
 \|b_{f,k}^{\pm}\|_{
 L^2(I_k\cap\{\lambda:\lambda/k<1\};H^1(\mathbb T_\theta))}
 \le C|k|^{-7/2}M,\qquad k\ne0.
\end{equation}

\smallskip
\noindent
(iii) For $0<|\lambda|<\kappa$,
\[
 \|D_\lambda\Gamma_f(\lambda)\|_{H^1}
 +\|\partial_\lambda\Gamma_f(\lambda)\|_{L^2}
 \le Ce^{-c/|\lambda|}M.
\]

\end{proposition}

\subsection{Proof of Proposition~\ref{prop:spectral-bounds}}
\label{subsec:proof-first-spectral}

We first treat $J_k$, then subtract the center logarithm on $I_k$,
$k\ne0$, and finally prove the small- and large-parameter bounds.
Throughout, $g=g_f$, $H=(\Delta+1)f$, and
$(\Delta+T)g=H$ by \eqref{r2:fixed-source}.

\subsubsection{The differentiated equation}
\label{subsec:critical-radius}

Fix $c\in\mathbb C$ with $\operatorname{Re}c=1$ and
$0<|c-1|\ll1$, and set $z=ic\lambda$.
On each spectral interval choose $|c-1|$ sufficiently small that
$|\operatorname{Re}(ic\lambda)|\le1$ and the limiting absorption
estimate \eqref{r2:full-source-estimate} applies.
All derivatives of $T(ic\lambda)$, $g_f(ic\lambda)$ and
$\psi_f(ic\lambda)$ below are derivatives with $c$ fixed. Since $V(kb)=kb/\lambda$,
\[
\begin{gathered}
 (\partial_\lambda+V)(\lambda S)=0,\qquad
 (\partial_\lambda+V)(kb-c\lambda)=\frac{kb-c\lambda}{\lambda},\\
 (\partial_\lambda+V)\frac{kb}{kb-c\lambda}=0.
\end{gathered}
\]
Equivalently,
\[
 \frac{d}{d\lambda}T(ic\lambda)+[V,T(ic\lambda)]=0,
 \qquad
 \frac{d}{d\lambda}\frac{kb}{kb-c\lambda}
       =\frac{c\,kb}{(kb-c\lambda)^2}.
\]
We obtain the boundary value with sign $\pm$ by sending
$c\to1$ with $\pm\operatorname{Re}(ic\lambda)>0$.

All divergences, gradients and Jacobians below are in the physical
coordinates $(x,y)$. For $\lambda\in I_k\setminus\{0\}$, we define the scalar function $\zeta(\lambda,x,y)$, matrix $K_\eta(\lambda,x,y)$ and operator $E_\eta(c,\lambda)$ as follows:
\[
\begin{aligned}
 \zeta&=\operatorname{div}\!\bigl(\eta(q/q_{0,k}^{I})V\bigr),\\
 K_\eta
 &=D_{x,y}\!\bigl(\eta(q/q_{0,k}^{I})V\bigr)
   +D_{x,y}\!\bigl(\eta(q/q_{0,k}^{I})V\bigr)^t
   -\zeta\,\mathrm{Id},\\
 E_\eta(c,\lambda)
 &=\frac{d}{d\lambda}T(ic\lambda)
   +[\eta(q/q_{0,k}^{I})V,T(ic\lambda)]=(1-\eta(q/q_{0,k}^{I}))
       \frac{d}{d\lambda}T(ic\lambda).
\end{aligned}
\]
On $J_k$, we use the same definitions and the identities below with
$q_{0,k}^{I}$ replaced throughout by $q_{0,k}^{J}$.
Using the product rule, we have the physical commutator identity:
\begin{equation}\label{spectral:physical-commutator}
 [\Delta,\eta(q/q_{0,k}^{I})V]u
 =\operatorname{div}(K_\eta\nabla u)+\zeta\Delta u.
\end{equation}
Applying this identity and \eqref{r2:fixed-source}, we obtain
\begin{equation}\label{eq:fixed-source-derivative}
\begin{aligned}
 (\Delta+T(ic\lambda))D_\lambda g_f(ic\lambda)
 ={}&\operatorname{div}\!\bigl(\eta(q/q_{0,k}^{I})VH\bigr)\\
 &+\operatorname{div}(K_\eta\nabla g_f(ic\lambda))
   -\zeta T(ic\lambda)g_f(ic\lambda)\\
 &-E_\eta(c,\lambda)g_f(ic\lambda).
\end{aligned}
\end{equation}
On $J_k$ we apply Lemma~\ref{r2:full-source-LAP} directly;
on the resonant side of $I_k$, $k\ne0$, we first subtract the
center logarithm.

We first consider $\lambda\in J_k$. For $q\le2q_{0,k}^{J}$,
\[
 0\le1-b(q)=\frac{q^2}{8}
 \le\frac{\kappa}{4(1+|k|)},\qquad
 |\lambda|\le1+|k|,\qquad
 \operatorname{dist}(\lambda,\mathbb Z)\ge\kappa.
\]
Splitting into $|n|\le2(1+|k|)$ and $|n|>2(1+|k|)$, we obtain
\[
 |nb(q)-c\lambda|
 \ge |nb(q)-\lambda|
 \ge\frac{\kappa}{2},
 \qquad n\ne0.
\]
Thus the cutoff region is nonresonant. By
Lemma~\ref{geom:common-field}, for every $1\le p<\infty$,
\[
 \|\eta(q/q_{0,k}^{J})V\|_{L^\infty}
 +\|\eta(q/q_{0,k}^{J})V\|_{W^{1,p}}
 \le C_{\kappa,k,p},
\]
with equal traces across the separatrix. In the following estimates, we allow the constants to depend on $\kappa$, $k$ and the indicated
exponents.

\begin{lemma}\label{aux:logregularity}
Let $1<p<2<q<\infty$. We assume $\operatorname{Re}c=1$ and $0<|c-1|$ sufficiently
small, depending on the compact spectral interval $I$.

\smallskip\noindent
\textup{(i)}
For every compact $I\subset\mathbb R\setminus\mathbb Z$,
$\lambda\in I$, and $h\in H^1(\mathbb T^2)$, there exist
$F\in L^q(\mathbb T^2;\mathbb C^2)$ and
$g\in L^p(\mathbb T^2)$ such that
\begin{equation}\label{spectral:log-regularity}
 T(ic\lambda)h=\operatorname{div}F+g,\qquad
 \|F\|_{L^q}+\|g\|_{L^p}
 \le C_{I,p,q}\|h\|_{H^1}.
\end{equation}
The same conclusion holds with $T(ic\lambda)$ replaced by
either boundary value $T_\pm(\lambda)$.

\smallskip\noindent
\textup{(ii)}
Let $\rho$ be a radial function, smooth in each cell and
through its elliptic center. Near the separatrix, write
$\rho$ as a function of $\tau=S(r)$, separately in each cell,
and assume that
\[
 |\rho(\tau)|+|\partial_\tau\rho(\tau)|
 \le A(1+\tau)^N,\qquad \tau\ge\tau_0,
\]
for some $A>0$, $N\ge0$, and $\tau_0>0$.
For every compact $I\subset\mathbb R\setminus\mathbb Z$,
$\lambda\in I$, and $h\in H^1(\mathbb T^2)$, there exist
$F_\rho\in L^q(\mathbb T^2;\mathbb C^2)$ and
$g_\rho\in L^p(\mathbb T^2)$ such that
\[
 \rho\,[T(ic\lambda)h]
 =\operatorname{div}F_\rho+g_\rho,\qquad
 \|F_\rho\|_{L^q}+\|g_\rho\|_{L^p}
 \le C_{I,p,q,\rho}\|h\|_{H^1}.
\]
 If $\rho$ vanishes in
neighborhoods of both elliptic centers, the same conclusion
holds for every compact $I\subset\mathbb R\setminus\{0\}$. The weighted boundary values
\[
 \rho T_\pm(\lambda)h
 :=
 \lim_{\substack{c\to1\\
          \pm\operatorname{Re}(ic\lambda)>0}}
 \rho\,[T(ic\lambda)h]
 \quad\text{in }\mathcal D'(\mathbb T^2)
\]
exist and satisfy the same decomposition and estimate.

\smallskip\noindent
\textup{(iii)}
For every compact $I\subset\mathbb R\setminus\mathbb Z$
and $f\in H^1(\mathbb T^2)\cap Y$,
\[
 \|\psi_f(ic\lambda)\|_{W^{1,q}}
 +\|\psi_{f,\pm}(\lambda)\|_{W^{1,q}}
 \le C_{I,q}\|f\|_{H^1},\qquad \lambda\in I.
\]
If $\chi\in C^\infty(\mathbb T^2)$ is a radial cutoff
vanishing in neighborhoods of both elliptic centers, then
for every compact $I\subset\mathbb R\setminus\{0\}$,
\[
 \|\chi\psi_f(ic\lambda)\|_{W^{1,q}}
 +\|\chi\psi_{f,\pm}(\lambda)\|_{W^{1,q}}
 \le C_{I,q,\chi}\|f\|_{H^1},\qquad \lambda\in I.
\]
All constants above are independent of $\lambda\in I$
and of $c$ as $c\to1$ from either side.
\end{lemma}
\begin{proof}
Near the centers, inside a fixed disk we have
$|nb-c\lambda|\ge c_I(1+|n|)$. Hence $T$ is bounded on $L^2$ in this case.
On a compact regular annulus, only finitely many indices resonate.
Using
\[
 \frac{nb}{nb-c\lambda}=1+\frac{c\lambda/n}{b-c\lambda/n},
 \qquad \partial_b\log(b-c\lambda/n)=(b-c\lambda/n)^{-1},
\]
apply \eqref{w14:ibp} with $X_b=\partial_b|_\theta$.
The flux contains $h_n\log(b-c\lambda/n)$ and the scalar terms
contain $h_n'\log(b-c\lambda/n)$ or bounded multiples of the flux.
It follows from
$\|h_n\|_\infty+\|h_n'\|_2\le C_I\|h\|_{H^1}$ and the uniform
finite-$L^a$ bounds for the logarithms that the desired norms are bounded.

Near the separatrix set $\tau_n=2\pi n/\lambda$, and choose
cutoffs $\eta_n$ supported on disjoint fixed-width intervals
$U_n$ about $\tau_n$, equal to one near their centers. Recall from \eqref{eq:weight} that
\[
 w(\tau)=\frac{\tau}{|S'(r(\tau))|},
 \qquad
 dx\,dy=\frac{w(\tau)}{2\pi}\,d\tau\,d\theta.
\]
It follows from \eqref{jxtb1} that
\[
 (\operatorname{div}X_\tau)\circ X
 =\frac{w'(\tau)}{w(\tau)},
\]
where the prime denotes differentiation with respect to $\tau$.
Since $\operatorname{div}X_\tau=w'/w$, the singular term is
\[
 -\frac{\tau_n\eta_nh_ne^{in\theta}}{c\tau-\tau_n}
       =\operatorname{div}F_n+f_n,
\]
where
\[
\begin{aligned}
 F_n&=-\frac{\tau_n}{c}\eta_n\log(\tau-\tau_n/c)h_ne^{in\theta}X_\tau,\\
 f_n&=\frac{\tau_n}{c}\log(\tau-\tau_n/c)e^{in\theta}
       (\eta_nh_n'+\eta_n'h_n+(w'/w)\eta_nh_n).
\end{aligned}
\]
Lemma~\ref{geom:common-field} bounds all geometric coefficients
by polynomials in $\tau$. With $w_n=\tau_ne^{-\tau_n/4}$,
\[
 |U_n|\le C_Iw_n,\qquad
 \|\log(\tau-\tau_n/c)\|_{L^a(U_n)}
       \le C_{I,a}w_n^{1/a}(1+\log(2+\tau_n)).
\]
Choose $P>q$, $1/a=1/q-1/P$, and $1/b=1/p-1/2$.
It follows from Hölder's inequality and $H^1\hookrightarrow L^P$ that
\[
 \|F_n\|_q+\|f_n\|_p
 \le C(1+\tau_n)^N(w_n^{1/a}+w_n^{1/b})\|h\|_{H^1}.
\]
The sum converges uniformly for $\lambda\in I$. Away from these
intervals the multiplier is bounded by a polynomial in $S$;
since $(1+S)^N\in L^a$ for every finite $a$, Hölder's inequality
controls the remaining terms. The weighted estimate  follows from
\[
 \rho\operatorname{div}F_n
 =\operatorname{div}(\rho F_n)
       -(\partial_\tau\rho)(F_n\cdot\nabla\tau),
\]
and only polynomial factors change. If we localize away from the centers the above analysis show that the same bounds hold for $I\subset\R\backslash\{0\}$ as claimed.

Finally,
\[
 \Delta\psi_f(ic\lambda)
 =-T(ic\lambda)\left(\psi_f(ic\lambda)+\frac f{ic\lambda}\right)
     +\frac f{ic\lambda}.
\]
Choose $2q/(q+2)<p<2$, so $L^p\hookrightarrow W^{-1,q}$.
Applying elliptic regularity, \eqref{spectral:log-regularity}, and
\eqref{r2:full-source-estimate} to
$(\Delta+T)\psi_f=i(T-I)f/(c\lambda)$, we obtain
\[
 \|\psi_f(ic\lambda)\|_{W^{1,q}}
 \le C_{I,q}(\|\psi_f(ic\lambda)\|_{H^1}+\|f\|_{H^1})
 \le C_{I,q}\|f\|_{H^1}.
\]
Weak compactness and Proposition~\ref{global:LAP} then imply all bounds up
to the boundary.
\end{proof}

\smallskip\noindent
For $f\in H^2\cap Y$, we first prove
\begin{equation}\label{cg:good-bound}
 \operatorname*{ess\,sup}_{\lambda\in J_k}
 \left(\|D_\lambda\psi_{f,\pm}(\lambda)\|_{H^1}
       +\|\partial_\lambda\psi_{f,\pm}(\lambda)\|_2\right)
 \le C_{\kappa,k}\|f\|_{H^2}.
\end{equation}
Write $\psi=\psi_f(ic\lambda)$, $T=T(ic\lambda)$ and
$\eta=\eta(q/q_{0,k}^J)$, with $c$ fixed under differentiation.
It follows from \eqref{spectral:log-regularity},
Lemma~\ref{geom:common-field}, and the separation of the multipliers
on $q\le2q_{0,k}^J$ that
\[
\begin{aligned}
 &\|Th\|_{H^{-1}}\le C_{\kappa,k}\|h\|_{H^1},\qquad
   \|\psi\|_{W^{1,4}}\le C_{\kappa,k}\|f\|_{H^1},\\
 &\|\eta V\|_{L^\infty\cap W^{1,4}}+\|K_\eta\|_4+\|\zeta\|_4
   \le C_{\kappa,k},\qquad
   \|E_\eta h\|_2\le C_{\kappa,k}\|h\|_2.
\end{aligned}
\]
We also use
\begin{equation}\label{cg:dM}
 \|\zeta Th\|_{H^{-1}}\le C_{\kappa,k}\|h\|_{H^1}.
\end{equation}
Indeed, $\zeta$ vanishes near the centers and, near the separatrix,
\[
 \zeta=-\lambda^{-1}(1+\tau w'/w),\qquad
 |\zeta|+|\partial_\tau\zeta|\le C_{\kappa,k}(1+\tau)^N.
\]
Thus \eqref{cg:dM} is the weighted bound of
Lemma~\ref{aux:logregularity}, with any $q>2$ and $p=4/3$.
For the cutoff on $I_k$, $k\ne0$, $\zeta$ still vanishes near
both centers and has this same formula near the separatrix;
applying
Lemma~\ref{aux:logregularity}, we obtain \eqref{cg:dM} in this case.

It follows from $(\Delta+T)\psi=i(T-I)f/(c\lambda)$ and
\eqref{spectral:physical-commutator} that
\begin{equation}\label{cg:equation}
\begin{aligned}
 (\Delta+T)D_\lambda\psi
 ={}&\frac{i}{c\lambda}(T-I)(\eta Vf)
       -\frac{i}{c\lambda^2}(T-I)f\\
 &+\frac{i}{c\lambda}\{E_\eta f+\zeta(T-I)f\}
       +\operatorname{div}(K_\eta\nabla\psi)
       -\zeta T\psi-E_\eta\psi.
\end{aligned}
\end{equation}
It follows from the product estimate
$\|\eta Vf\|_{H^1}\le C_{\kappa,k}\|f\|_{H^2}$ and
\[
 \|K_\eta\nabla\psi\|_2\le\|K_\eta\|_4\|\nabla\psi\|_4,
 \qquad \|\zeta f\|_2\le\|\zeta\|_4\|f\|_4
\]
that
\[
 \|(\Delta+T)D_\lambda\psi\|_{H^{-1}}
 \le C_{\kappa,k}\|f\|_{H^2},\qquad
 |\operatorname{mean}(D_\lambda\psi)|
 =|\operatorname{mean}(\zeta\psi)|
 \le C_{\kappa,k}\|f\|_{H^1}.
\]
At nonzero absorption, it follows from the elliptic equation
that $\psi,\partial_\lambda\psi\in H^2$, hence $D_\lambda\psi\in H^1$.
Applying Lemma~\ref{r2:full-source-LAP} to its mean-zero part and
using the bound on the mean value, we obtain
\[
 \|D_\lambda\psi\|_{H^1}+\|\partial_\lambda\psi\|_2
 \le C_{\kappa,k}\|f\|_{H^2},
 \qquad \partial_\lambda\psi=D_\lambda\psi-\eta V\psi.
\]
For $\varphi\in C_c^\infty(J_k^\circ)$, it follows from
Proposition~\ref{global:LAP} that
\[
 \int_{J_k}\partial_\lambda\psi_f(ic\lambda)\varphi\,d\lambda
 =-\int_{J_k}\psi_f(ic\lambda)\varphi'\,d\lambda
 \longrightarrow-\int_{J_k}\psi_{f,\pm}(\lambda)\varphi'\,d\lambda,
\]
where $c\to1$ with $\pm\operatorname{Re}(ic\lambda)>0$.
The spatial term converges in $L^1(J_k;L^2)$, so weak compactness
proves \eqref{cg:good-bound}. Taking the difference proves
\eqref{spectral:decomposition} on $J_k$ with $\Lambda_f=0$;
the uniform decay in $k$ is proved in
Lemma~\ref{r2:angular-regular-remainder}.

\subsubsection{The first derivative at an elliptic center}

We estimate the singular source in the equation for $D_\lambda g_f$
near an elliptic center. In either cell, use the coordinates of
Lemma~\ref{geom:metric}, with
\[
 0<q<q_*,\qquad \theta\in\mathbb T,\qquad v=\log q,
 \qquad dx\,dy=\xi_0(q)q^2\,dv\,d\theta.
\]
Here $\mathsf G_{vv}(q,\theta)$ denotes the radial entry of the
Dirichlet-form matrix in these coordinates, evaluated at $v=\log q$.
Fix
\[
 k\ge1,\qquad |\lambda-k|<5\kappa/4,\qquad
 \operatorname{Re}c=1,\qquad
 0<|\epsilon|\le\kappa/8,\qquad
 \epsilon=-\lambda\operatorname{Im}c,
\]
and set
\[
 g=g_f(ic\lambda),\qquad
 \sigma=8(1-\lambda/k),\qquad
 (q_{0,k}^{I})^2=32\kappa/k.
\]
Choose $q_*$ inside the disk where the cutoff in
\eqref{r2:weighted-input} equals one, and decrease $\kappa$ so that
$e^2\sqrt{10\kappa}<q_*/2$. In particular $2q_{0,k}^{I}<q_*$.
We first define the correction $\Lambda_{g,k}$.
If $\sigma>0$ and $0<8|\epsilon|/k\le\sigma/2$, set
$q_c=\sqrt{\sigma}$ and set
\[
\begin{aligned}
 \Lambda_{g,k}(c,\lambda;q,\theta)
 ={}&\chi\!\left(\log\frac q{q_c}\right)
 \frac{2c\,\xi_0(q_c)b(q_c)g_k(q_c)e^{ik\theta}}
      {k(1-c\lambda/k)\mathsf G_{vv}(q_c,\theta)}\\
 &\hspace{12mm}\times
 \log\!\left(\frac{q^2-8(1-c\lambda/k)}{\sigma}\right),
\end{aligned}
\]
where the logarithm is taken on the principal branch.
If $\sigma\le0$ or $8|\epsilon|/k>\sigma/2$, set
$\Lambda_{g,k}=0$.

Since $\chi\in C_c^\infty((-1/2,1/2))$, the nonzero correction
is supported in
\[
 e^{-1/2}q_c<q<e^{1/2}q_c.
\]
Moreover,
\[
 0<\sigma<\frac{10\kappa}{k},
 \qquad
 e\sigma<\frac{32\kappa}{k}=(q_{0,k}^{I})^2.
\]
Thus its support lies in $q<q_{0,k}^{I}$, strictly inside
the coordinate neighborhood. We regard $\Lambda_{g,k}$
as a function on $\mathbb T^2$ by setting it equal to zero
elsewhere.

With this correction fixed, define
\begin{equation}\label{adapted:first-center-source-definition}
 F_k=-\bigl(1-\eta(q/q_{0,k}^{I})\bigr)
       \sum_{n\ne0}\frac{c\,n b(q)}{(n b(q)-c\lambda)^2}
                         g_n(q)e^{in\theta}
       -(\Delta+T(ic\lambda))\Lambda_{g,k}.
\end{equation}
The first term is supported in $q<2q_{0,k}^{I}$.
Since $T(ic\lambda)$ acts separately on each streamline,
$F_k$ is supported in $q<2q_{0,k}^{I}$.

It follows from \eqref{r2:fixed-source} and
\eqref{eq:fixed-source-derivative}, with $H=(\Delta+1)f$,
that in this coordinate neighborhood,
\[
\begin{aligned}
 &(\Delta+T(ic\lambda))(D_\lambda g-\Lambda_{g,k})\\
 &\qquad=F_k+\eta(q/q_{0,k}^{I})VH
             +[\Delta,\eta(q/q_{0,k}^{I})V]g.
\end{aligned}
\]
On $\operatorname{supp}\Lambda_{g,k}$, we have $\eta=0$,
and therefore
\[
 D_\lambda g=\frac{d}{d\lambda}g_f(ic\lambda).
\]
We perform this construction separately at the two centers.
The global correction and source are the sums of the two
local corrections and sources.

In $q<2q_{0,k}^{I}$, only the angular index $n=k$ can
produce a singular denominator as $\epsilon\to0$.
Indeed, for $j=0,1$,
\[
 \sup_{\substack{0<q<2q_{0,k}^{I}\\n\ne0,k}}
 \frac{|c|\,|n b(q)|^{1+j}}
      {|n b(q)-c\lambda|^2|n|^3}
 \le C_\kappa k^{j-2}.
\]
For $f\in H^5\cap Y$ and $M=\|f\|_{H^5}$, it follows
from Parseval's identity and \eqref{r2:angular-g} that
\[
\begin{aligned}
 &\left\|\bigl(1-\eta(q/q_{0,k}^{I})\bigr)L_0^j
   \sum_{n\ne0,k}\frac{c\,n b(q)}{(n b(q)-c\lambda)^2}
                 g_n(q)e^{in\theta}\right\|_{L^2(dx\,dy)}\\
 &\qquad\le C_\kappa q_{0,k}^{I}k^{j-2}
       \|\partial_\theta^3g\|_{L^2(\{q<q_*\},\,dv\,d\theta)}
 \le C_\kappa k^{j-5/2}M,\qquad j=0,1.
\end{aligned}
\]
The following lemma estimates the full source $F_k$,
including the contribution from $n=k$.

\begin{lemma}\label{adapted:first-center-source}
Let $\Lambda_{g,k}$ be the correction defined above, and let
$F_k$ be defined by \eqref{adapted:first-center-source-definition}.
For the parameters fixed above and every $f\in H^5\cap Y$,
\[
 \langle F_k,1\rangle=0,
 \qquad
 \|F_k\|_{H^{-1}(\mathbb T^2)}
 +k^{-1}\|L_0F_k\|_{H^{-1}(\mathbb T^2)}
 \le C_\kappa k^{-3/2}\|f\|_{H^5}.
\]
The constant is independent of $k$, $\lambda$, $\epsilon$,
and the choice of center.

For $\lambda\ne k$, the same conclusions hold for either
boundary value with the following choice of correction.
For $\lambda>k$, set $\Lambda_{g,k}=0$.
For $\lambda<k$, put $q_c=\sqrt{\sigma}$ and set
\[
 \Lambda_{g,k}
 =\chi\!\left(\log\frac q{q_c}\right)
 \frac{2\,\xi_0(q_c)b(q_c)g_{f,\pm,k}(q_c,\lambda)e^{ik\theta}}
      {k(1-\lambda/k)\mathsf G_{vv}(q_c,\theta)}
 \log\!\left(\log(q/q_c)\mp i0\right).
\]
Define the corresponding boundary source $F_k$ by
\eqref{adapted:first-center-source-definition}, with
$c=1$, $g=g_{f,\pm}(\lambda)$, and $T(ic\lambda)$
replaced by $T_\pm(\lambda)$.
The singular expressions are interpreted in the corresponding
one-sided distributional sense. The double-pole terms in
\eqref{adapted:first-center-source-definition} cancel,
and the resulting source $F_k$ belongs to
$H^{-1}(\mathbb T^2)$ and satisfies the estimates above.
\end{lemma}
\begin{proof}
First assume $\sigma>0$ and $0<8|\epsilon|/k\le\sigma/2$.
Set $v=\log(q/q_c)$ and
\[
 \ell(v)=\log\!\left(\frac{q^2-8(1-c\lambda/k)}{\sigma}\right),
 \qquad
 a(\theta)=
 \frac{2c\,\xi_0(q_c)b(q_c)g_k(q_c)e^{ik\theta}}
      {k(1-c\lambda/k)\mathsf G_{vv}(q_c,\theta)}.
\]
Thus $\Lambda_{g,k}=\chi(v)a(\theta)\ell(v)$.
Since
\[
 \ell'=\frac{2q^2}{q^2-8(1-c\lambda/k)},
 \qquad
 -\ell''=
 \frac{32(1-c\lambda/k)q^2}
      {[q^2-8(1-c\lambda/k)]^2},
\]
the term $-\mathsf G_{vv}a\ell''$ cancels the leading
double pole in the critical angular mode.

Set
\[
 h=\partial_vg-\frac{2k(1-c\lambda/k)}{c}\,a\ell.
\]
We decompose the source in
\eqref{adapted:first-center-source-definition} as
\begin{equation}\label{adapted:first-center-equation}
 F_k=\sum_{j=1}^5F_k^{(j)},
\end{equation}
where the components are defined by
\[
\begin{aligned}
 -\xi_0q^2F_k^{(1)}
 ={}&\bigl(1-\eta(q/q_{0,k}^{I})\bigr)
 \sum_{n\ne0,k}
 \frac{64c\,\xi_0bq^2g_ne^{in\theta}}
      {n[q^2-8(1-c\lambda/n)]^2}\\
 &-\chi(v)\ell(v)\sum_{n\ne0,k}
 \frac{8\xi_0bq^2a_ne^{in\theta}}
      {q^2-8(1-c\lambda/n)},
\end{aligned}
\]
\[
 -\xi_0q^2F_k^{(2)}
 =\bigl(1-\eta(q/q_{0,k}^{I})-\chi(v)\bigr)
 \frac{64c\,\xi_0bq^2g_ke^{ik\theta}}
      {k[q^2-8(1-c\lambda/k)]^2},
\]
\[
\begin{aligned}
 -\xi_0q^2F_k^{(3)}
 ={}&-\frac{32c\chi(v)}{k}
 \partial_v\!\left(
 \frac{\xi_0b(g_k-g_k(q_c))}
      {q^2-8(1-c\lambda/k)}\right)e^{ik\theta}\\
 &+\frac{32c\chi(v)}{k}
 \frac{\partial_v(\xi_0b)(g_k-g_k(q_c))}
      {q^2-8(1-c\lambda/k)}e^{ik\theta},
\end{aligned}
\]
\[
 -\xi_0q^2F_k^{(4)}
 =\frac{32c\chi(v)}{k}
   \frac{\xi_0b\,h_k}{q^2-8(1-c\lambda/k)}e^{ik\theta}
   -8\chi(v)\xi_0b\,a_k\ell(v)e^{ik\theta},
\]
and
\[
\begin{aligned}
 -\xi_0q^2F_k^{(5)}
 ={}&\frac{64c\chi(v)q^2g_k(q_c)e^{ik\theta}}
          {k[q^2-8(1-c\lambda/k)]^2}\\
 &\quad\times
 \left[\xi_0(q)b(q)
 -\frac{\xi_0(q_c)b(q_c)\mathsf G_{vv}(q,\theta)}
        {\mathsf G_{vv}(q_c,\theta)}\right]\\
 &+\operatorname{div}_{v,\theta}
       (\mathsf G\nabla_{v,\theta}\Lambda_{g,k})
       -\chi(v)\mathsf G_{vv}a\ell''.
\end{aligned}
\]
The identity \eqref{adapted:first-center-equation} follows
from the definition of $F_k$, the formulas for
$\ell',\ell''$, and
\[
\begin{aligned}
 &\frac{64c\,\xi_0bq^2(g_k-g_k(q_c))}
       {k[q^2-8(1-c\lambda/k)]^2}\\
 &=-\frac{32c}{k}\partial_v
       \left(\frac{\xi_0b(g_k-g_k(q_c))}
                         {q^2-8(1-c\lambda/k)}\right)
   +\frac{32c}{k}
       \frac{\xi_0b\,\partial_vg_k
                 +\partial_v(\xi_0b)(g_k-g_k(q_c))}
            {q^2-8(1-c\lambda/k)},
\end{aligned}
\]
together with the cancellation
\[
 \frac{32c}{k}
 \frac{\xi_0b(\partial_vg_k-h_k)}
      {q^2-8(1-c\lambda/k)}
 -\frac{8\xi_0bq^2a_k\ell}
       {q^2-8(1-c\lambda/k)}
 =-8\xi_0b\,a_k\ell.
\]
The first four components have zero angular average.
The angular average of the fifth component will be
estimated directly in Step~4.

\smallskip\noindent
\emph{Step 1: General estimates.}
With $W(s)=e^{s/2}(1+e^{2s})^{3/4}$, we have
\begin{equation}\label{adapted:finite-absorption-weight}
 \left\|W(s)^{-1}
 g_k\!\left((|\sigma|+8|\epsilon|/k)^{1/2}e^s\right)
 \right\|_{H^1_k(ds)}
 \le C(|\sigma|+8|\epsilon|/k)k^{-2}M.
\end{equation}
The fixed center cutoff from \eqref{r2:weighted-input}
is understood. Indeed, apply \eqref{r2:weighted-normal}
to \eqref{r2:localized-equation} with parameters
\[
 \left(\frac{\sigma}{|\sigma|+8|\epsilon|/k},
       \frac{8\epsilon/k}{|\sigma|+8|\epsilon|/k}\right).
\]
The four source bounds used for \eqref{r2:weighted-input}
remain valid because
$q^2/W(s)\le|\sigma|+8|\epsilon|/k$ and
\eqref{eq:cylinder-form} is uniform in these parameters.

Since $\sigma\le\sigma+8|\epsilon|/k\le3\sigma/2$,
translation to $v=\log(q/q_c)$ and the trace inequality imply
\begin{equation}\label{adapted:first-center-input}
 \|\partial_vg_k\|_2+k\|g_k\|_2\le C\sigma k^{-2}M,
 \qquad
 |g_k(q_c)|\le C\sigma k^{-5/2}M.
\end{equation}
Here and in the local estimates below, norms use
$|v|<2$ and cylinder measure. Consequently,
\[
 \|a\|_{H^j_\theta}\le Ck^{j-7/2}M,
 \qquad 0\le j\le3.
\]

We also need one radial derivative of $h_k$.
The critical terms are controlled by
\[
 \left|\frac{vq^2}{q^2-8(1-c\lambda/k)}\right|\le C,
 \qquad
 \left\|\frac{g_k-g_k(q_c)}v\right\|_2
 \le2\|\partial_vg_k\|_2,
\]
and
\[
 \left|\frac{\xi_0b}{\mathsf G_{vv}}
       -\frac{\xi_0(q_c)b(q_c)}
              {\mathsf G_{vv}(q_c,\theta)}\right|
 \le C\sigma|v|.
\]
For the other angular indices,
\begin{equation}\label{adapted:first-noncritical-multiplier}
 \left\|\partial_\theta^j
 \sum_{n\ne0,k}
 \frac{8\xi_0bq^2g_ne^{in\theta}}
      {q^2-8(1-c\lambda/n)}\right\|_2
 \le C\sigma\|\partial_\theta^{j+1}g\|_2,
 \qquad j=0,1.
\end{equation}
Expanding \eqref{r2:fixed-source} in cylinder coordinates
and taking one angular derivative before projection onto
mode $k$, it follows from \eqref{r2:angular-g},
\eqref{r2:H-cylinder}, and $\mathsf G-I=O(q^2)$ that
\[
 \|\partial_vh_k\|_2
 \le C\left(k^2\|g_k\|_2+\|\partial_vg_k\|_2
        +\sigma|g_k(q_c)|+\sigma k^{-1}M\right)
 \le C\sigma k^{-1}M.
\]
Since $\|\ell\|_2\le C$, we obtain
\begin{equation}\label{adapted:first-flux-bounds}
 \|h_k\|_2\le C\sigma k^{-2}M,
 \qquad
 \|h_k\|_{H^1_k(|v|<2)}\le C\sigma k^{-1}M.
\end{equation}

We use the physical pairing
\begin{equation}\label{adapted:cylinder-duality}
 \left\langle-\frac{F}{\xi_0q^2},\varphi\right\rangle
 =-\int F\overline\varphi\,dv\,d\theta.
\end{equation}
For nonzero angular modes, integration by parts and
angular Poincar\'e imply
\[
 \left|\int
   (\operatorname{div}_{v,\theta}\boldsymbol U+V)
       \overline{\varphi-\varphi_0}\,dv\,d\theta\right|
 \le C(\|\boldsymbol U\|_2+\|V\|_2)\|\varphi\|_{H^1}.
\]
The corresponding estimate after $L_0=b\partial_\theta$
follows from
\[
 b\partial_\theta(\operatorname{div}\boldsymbol U+V)
 =\operatorname{div}(b\partial_\theta\boldsymbol U)
       +b\partial_\theta V
       -(\partial_vb)\partial_\theta U_v,
\]
where $U_v$ denotes the radial component of $\boldsymbol U$.

\smallskip\noindent
\emph{Step 2: The noncritical component $F_k^{(1)}$
and the critical tail $F_k^{(2)}$.}
The first sum defining $F_k^{(1)}$ satisfies the physical
$L^2$ estimate proved immediately before the lemma.
For its second sum, \eqref{adapted:first-noncritical-multiplier}
applied to $a$, together with $\|\ell\|_2\le C$, implies
\[
 \sum_{j=0}^1k^{-j}
 \left\|\partial_\theta^j\left[
 \chi\ell\sum_{n\ne0,k}
 \frac{8\xi_0bq^2a_ne^{in\theta}}
      {q^2-8(1-c\lambda/n)}\right]\right\|_2
 \le C\sigma
       \bigl(\|a\|_{H^1_\theta}
                    +k^{-1}\|a\|_{H^2_\theta}\bigr).
\]
Thus \eqref{adapted:cylinder-duality} controls this component.

The factor $1-\eta(q/q_{0,k}^{I})-\chi(v)$ vanishes near
$v=0$. Since
\[
 \sup_{|v|\ge1/8}
 \frac{e^{2v}W(v)}
      {|e^{2v}-1-8i\epsilon/(k\sigma)|^2}\le C,
\]
it follows from \eqref{adapted:finite-absorption-weight} that
\[
 \left\|\mathbf1_{\{|v|\ge1/8\}}
       \frac{q^2g_k}
            {k[q^2-8(1-c\lambda/k)]^2}\right\|_2
 \le Ck^{-1}\sigma^{-1}\|W^{-1}g_k\|_2
 \le Ck^{-3}M.
\]
This is a pure $k$th-mode source, so
$\|\varphi_k\|_2\le k^{-1}\|\varphi_k\|_{H^1_k}$
controls its physical dual norm. Applying $L_0$ costs
at most $Ck$. We conclude that
\[
 \sum_{j=1}^2
 \left(\|F_k^{(j)}\|_{H^{-1}}
       +k^{-1}\|L_0F_k^{(j)}\|_{H^{-1}}\right)
 \le C_\kappa k^{-3/2}M.
\]

\smallskip\noindent
\emph{Step 3: The difference quotient $F_k^{(3)}$
and regularized term $F_k^{(4)}$.}
It follows from Hardy's inequality and
\eqref{adapted:first-center-input} that
\[
 \left\|\frac{g_k-g_k(q_c)}
                  {q^2-8(1-c\lambda/k)}\right\|_2
 \le C\sigma^{-1}\|\partial_vg_k\|_2
 \le Ck^{-2}M.
\]
Using $\chi\partial_vu=\partial_v(\chi u)-\chi'u$,
$\|\partial_v(\xi_0b)\|_\infty\le C\sigma$, and
\eqref{adapted:cylinder-duality}, we obtain
\[
 \|F_k^{(3)}\|_{H^{-1}}
 +k^{-1}\|L_0F_k^{(3)}\|_{H^{-1}}
 \le Ck^{-3}M.
\]

For the first term defining $F_k^{(4)}$, the Cauchy estimate
\eqref{eq:cylinder-form} and
\eqref{adapted:first-flux-bounds} imply
\[
 \left\|\frac{\chi(v)\xi_0b\,h_k}
                  {k[q^2-8(1-c\lambda/k)]}\right\|_{H^{-1}_k}
 \le Ck^{-3/2}\sigma^{-1}
                \|h_k\|_{H^1_k(|v|<2)}
 \le Ck^{-5/2}M.
\]
The second term satisfies
\[
 \|\chi\xi_0b\,a_k\ell\|_2\le Ck^{-7/2}M.
\]
Both terms have angular index $k$, and multiplication by
$b$ is bounded on the relevant radial Sobolev spaces.
Therefore
\[
 \|F_k^{(4)}\|_{H^{-1}}
 +k^{-1}\|L_0F_k^{(4)}\|_{H^{-1}}
 \le Ck^{-5/2}M.
\]

\smallskip\noindent
\emph{Step 4: The geometry and cutoff component $F_k^{(5)}$.}
For the first term defining $F_k^{(5)}$, the quotient
\[
 \frac{\xi_0(q)b(q)
       -\xi_0(q_c)b(q_c)\mathsf G_{vv}(q,\theta)
                              /\mathsf G_{vv}(q_c,\theta)}
      {q^2-\sigma}
\]
extends smoothly through $q=q_c$, with bounded derivatives.
The identity
\[
 \frac{q^2(q^2-\sigma)}
      {[q^2-8(1-c\lambda/k)]^2}
 =\frac12\partial_v
 \left[\ell(v)
       -\frac{8i\epsilon/k}{q^2-8(1-c\lambda/k)}\right]
\]
provides a uniformly bounded $L^2$ primitive.
The product rule expresses this term as a radial derivative
of an $L^2$ function plus an $L^2$ function.
Their norms are at most $Ck^{-1}|g_k(q_c)|$,
and one angular derivative costs at most $Ck$.

For the remaining terms, we use
\[
\begin{aligned}
 &\operatorname{div}(\mathsf G\nabla(a\ell))
                    -\mathsf G_{vv}a\ell''\\
 &=\partial_v\!\left(
   \bigl[(\partial_v\mathsf G_{vv}
             +\partial_\theta\mathsf G_{v\theta})a
                 +2\mathsf G_{v\theta}\partial_\theta a\bigr]\ell
   \right)\\
 &\quad+\partial_\theta\!\left(
      [\mathsf G_{\theta\theta}\partial_\theta a
                  -(\partial_v\mathsf G_{v\theta})a]\ell\right)
       -(\partial_v^2\mathsf G_{vv})a\ell.
\end{aligned}
\]
The derivatives of $\chi$ are supported where
$1/8\le|v|\le1/2$, where $\ell,\ell',\ell''$ are uniformly
bounded. Consequently,
\[
\begin{aligned}
 &\sum_{j=0}^1k^{-j}
 \left\|\partial_\theta^j\left[
 \operatorname{div}(\mathsf G\nabla\Lambda_{g,k})
 -\chi\operatorname{div}(\mathsf G\nabla(a\ell))
 \right]\right\|_2\\
 &\qquad\le C\bigl(\|a\|_{H^1_\theta}
                     +k^{-1}\|a\|_{H^2_\theta}\bigr).
\end{aligned}
\]
The displayed identities and \eqref{adapted:cylinder-duality}
control the nonzero angular modes of $F_k^{(5)}$,
also after applying $L_0$.

For its angular average, the other four components have
zero average, so
\[
 [-\xi_0q^2F_k^{(5)}]_0
 =\partial_v\bigl[
    \mathsf G_{vv}\partial_v\Lambda_{g,k}
       +\mathsf G_{v\theta}\partial_\theta\Lambda_{g,k}
 \bigr]_0.
\]
Since $[\mathsf G_{vv}(q_c,\theta)a(\theta)]_0=0$,
the expression inside $\partial_v$ equals
\[
\begin{aligned}
 &\chi\left[
 \frac{2q^2[\mathsf G_{vv}(q,\theta)
                      -\mathsf G_{vv}(q_c,\theta)]a}
      {q^2-8(1-c\lambda/k)}
 +\mathsf G_{v\theta}\partial_\theta a\,\ell
 \right]_0\\
 &\qquad+\chi'
       [(\mathsf G_{vv}(q,\theta)
                      -\mathsf G_{vv}(q_c,\theta))a]_0\ell.
\end{aligned}
\]
Its $L^2_v$ norm is bounded by
$C\sigma\|a\|_{H^1_\theta}$.
Integration by parts in
\eqref{adapted:cylinder-duality} controls the physical
$H^{-1}$ norm of this average; its $L_0$ derivative is zero.
Thus
\[
 \|F_k^{(5)}\|_{H^{-1}}
 +k^{-1}\|L_0F_k^{(5)}\|_{H^{-1}}
 \le C_\kappa k^{-3/2}M.
\]
Summing the five component estimates proves the asserted bound.
The first four components have zero physical mean.
The average of the fifth is a derivative of a compactly
supported function, so $\langle F_k,1\rangle=0$.

\smallskip\noindent
\emph{Step 5: The remaining regimes and boundary values.}
If $\sigma\le0$ or $8|\epsilon|/k>\sigma/2$, then
$\Lambda_{g,k}=0$.
We estimate the critical and noncritical sums directly in
\eqref{adapted:first-center-source-definition}.
Set $q=(|\sigma|+8|\epsilon|/k)^{1/2}e^s$.
Then
\[
 |q^2-8(1-c\lambda/k)|
 \ge C^{-1}(|\sigma|+8|\epsilon|/k)(1+e^{2s}),
 \qquad
 \sup_s\frac{e^{2s}W(s)}{(1+e^{2s})^2}\le C.
\]
It follows from \eqref{adapted:finite-absorption-weight} that
\[
 \left\|\frac{q^2g_k}
                  {k[q^2-8(1-c\lambda/k)]^2}\right\|_2
 \le Ck^{-1}(|\sigma|+8|\epsilon|/k)^{-1}
                    \|W^{-1}g_k\|_2
 \le Ck^{-3}M.
\]
The physical duality and noncritical estimates from Step~2
therefore prove the same bound in these regimes.

For fixed $\lambda\ne k$, Proposition~\ref{global:LAP}
implies local $H^1$ convergence of $g_k$, and also convergence
of $g_k(q_c)$ when $\lambda<k$.
Equation \eqref{adapted:first-flux-bounds} then implies
$L^2$ convergence and weak $H^1$ convergence of $h_k$.
The continuity of the Cauchy form in
Proposition~\ref{compact:normal} identifies the boundary
limit of $F_k^{(4)}$.
The explicit primitives and divergence identities above
identify the limits of the other components.
Their bounds pass to the limits by weak lower semicontinuity.

For the boundary correction used in the statement,
\[
 \log(e^{2v}-1\mp i0)-\log(v\mp i0)
 =\log\!\left(\frac{e^{2v}-1}{v}\right).
\]
The right-hand side is smooth on $\operatorname{supp}\chi$,
with value $\log2$ at $v=0$.
Its product with $\chi a$ has a controlled $H^1$ norm.
The metric estimates and \eqref{eq:cylinder-form} control
the resulting change in $F_k$, also after applying $L_0$,
by the same bounds as above.

\end{proof}

The corrected source estimate also controls differentiation through each
elliptic threshold. We record this consequence for the $H^5$ data used
in the main theorem.
\begin{corollary}\label{r2:threshold-derivative}
For each nonzero integer $m$, let $I_m=(m-\kappa,m+\kappa)$.
For every $f\in H^5\cap Y$ and either boundary sign,
\[
 \psi_{f,\pm}\in W^{1,\infty}(I_m;L^2(\mathbb T^2)),
\]
and
\[
 \operatorname*{ess\,sup}_{\lambda\in I_m}
 \left(\|\psi_{f,\pm}(\lambda)\|_{H^1}
       +\|\partial_\lambda\psi_{f,\pm}(\lambda)\|_2\right)
 \le C_m\|f\|_{H^5}.
\]
\end{corollary}
\begin{proof}
Complex conjugation reduces the proof to $m>0$. Set $M=\|f\|_{H^5}$
and fix
\[
 \operatorname{Re}c=1,\qquad
 0<|\operatorname{Im}c|\le\frac{\kappa}{8(m+\kappa)},\qquad
 g=g_f(ic\lambda),\quad H=(\Delta+1)f,\quad T=T(ic\lambda).
\]
Take $|c-1|$ sufficiently small for Lemma~\ref{aux:logregularity}.
Then $\epsilon=-\lambda\operatorname{Im}c$ satisfies
$0<|\epsilon|\le\kappa/8$ throughout $I_m$.
All derivatives in $\lambda$ below are taken with $c$ fixed.
Use $q_0=q_{0,m}^{I}$ and sum the two local corrections and
sources of Lemma~\ref{adapted:first-center-source}, writing the
result as $\Lambda_g$ and $F$, respectively.

When a correction is nonzero, $\sigma=8(1-\lambda/m)>0$ and
$8|\epsilon|/m\le\sigma/2$. It follows from
\eqref{adapted:first-center-input} and its defining formula that
\[
\begin{aligned}
 \|\Lambda_{g,m}\|_2^2
 &\le C_m\frac{|g_m(\sqrt\sigma)|^2}{\sigma}
   \int_{-1/2}^{1/2}
   \left|\log\!\left(e^{2v}-1-\frac{8i\epsilon}{m\sigma}\right)\right|^2\,dv
 \le C_m\sigma M^2.
\end{aligned}
\]
The integral is uniformly bounded because its integrand is at most
$C(1+|\log|v||^2)$. The corrections vanish in the other parameter
ranges, so $\|\Lambda_g\|_2\le C_mM$.

Let $\Omega_m$ be the complement of the two closed disks
$q\le q_0/2$. Since $g=f+ic\lambda\psi_f(ic\lambda)$, it follows from Proposition~\ref{global:LAP},
Lemma~\ref{aux:logregularity}(iii), and Lemma~\ref{geom:common-field} that
\[
 \|g\|_{H^1}+\|g\|_{W^{1,4}(\Omega_m)}\le C_mM,
 \qquad
 \|\eta V\|_{L^\infty\cap W^{1,4}}
       +\|K_\eta\|_4+\|\zeta\|_4\le C_m.
\]
Here $\eta=\eta(q/q_0)$, and $K_\eta,\zeta$ vanish for $q<q_0$.
Set $R=D_\lambda g-\Lambda_g$. By
\eqref{eq:fixed-source-derivative} and
\eqref{adapted:first-center-source-definition},
\[
 (\Delta+T)R
 =F+\operatorname{div}(H\,\eta V)
     +\operatorname{div}(K_\eta\nabla g)-\zeta Tg.
\]
Lemma~\ref{adapted:first-center-source}, \eqref{cg:dM}, and
\[
 \|K_\eta\nabla g\|_2
 \le\|K_\eta\|_4\|\nabla g\|_{L^4(\Omega_m)},
 \qquad
 \|H\,\eta V\|_2\le\|\eta V\|_\infty\|H\|_2
\]
imply $\|(\Delta+T)R\|_{H^{-1}}\le C_mM$.
For $\operatorname{Re}(ic\lambda)\ne0$, the resolvent equation
and analytic differentiation imply $g,\partial_\lambda g\in H^2$.
Since $\eta V\in L^\infty\cap W^{1,4}$ and $\Lambda_g\in H^1$,
we have $R\in H^1$. Using $\operatorname{mean}g=0$ and
$\zeta=\operatorname{div}(\eta V)$, integration by parts implies
\[
 \operatorname{mean}R=-\operatorname{mean}(\zeta g+\Lambda_g),
 \qquad |\operatorname{mean}R|\le C_mM.
\]
Applying \eqref{r2:full-source-estimate} to
$R-\operatorname{mean}R$, using $T1=0$, proves
$\|R\|_{H^1}\le C_mM$.
The identities
\[
 \partial_\lambda g=R+\Lambda_g-\eta Vg,
 \qquad
 \partial_\lambda\psi_f(ic\lambda)
 =\frac{\partial_\lambda g}{ic\lambda}
       -\frac{g-f}{ic\lambda^2}
\]
therefore imply
$\|\partial_\lambda\psi_f(ic\lambda)\|_2\le C_mM$ uniformly on $I_m$.

Integrating between $\lambda,\mu\in I_m$ and letting $c\to1$
with $\pm\operatorname{Re}(ic\lambda)>0$, it follows from the
strong boundary convergence in Proposition~\ref{global:LAP} that
\[
 \|\psi_{f,\pm}(\lambda)-\psi_{f,\pm}(\mu)\|_2
 \le C_mM|\lambda-\mu|.
\]
This proves the stated $W^{1,\infty}(I_m;L^2)$ conclusion,
including across $\lambda=m$. The $H^1$ bound follows from
Proposition~\ref{global:LAP}.
\end{proof}

\subsubsection{The spectral density near zero}

\begin{lemma}\label{tailD:small-frequency}
After decreasing the fixed constant $\kappa$ if necessary,
let $D_\lambda$ be defined on $I_0=(-\kappa,\kappa)$ by
\eqref{good-der}. Suppose that $f\in H^2\cap Y$ and
\[
 \langle f,g_j\rangle=0,\qquad j=1,2,3.
\]
There exist constants $a,C>0$, independent of $f$ and $\lambda$,
such that, for $0<|\lambda|<\kappa$,
\begin{equation}\label{tailD:F-W14}
 \|\Gamma_f(\lambda)\|_{W^{1,4}}
 \le Ce^{-a/|\lambda|}\|f\|_{H^1},
\end{equation}
and, for almost every such $\lambda$,
\begin{equation}\label{tailD:zero-derivative}
 \|D_\lambda\Gamma_f(\lambda)\|_{H^1}
 +\|\partial_\lambda\Gamma_f(\lambda)\|_2
 \le Ce^{-a/|\lambda|}\|f\|_{H^2}.
\end{equation}
For either boundary sign, we also have
\begin{equation}\label{tailD:individual-W14}
 \|\psi_{f,\pm}(\lambda)\|_{W^{1,4}}
 \le C|\lambda|^{-1}\|f\|_{H^1},
\end{equation}
and, for almost every $0<|\lambda|<\kappa$,
\begin{equation}\label{tailD:Zpolynomial}
 \|D_\lambda\psi_{f,\pm}(\lambda)\|_{H^1}
 +\|\partial_\lambda\psi_{f,\pm}(\lambda)\|_2
 \le C|\lambda|^{-4}\|f\|_{H^2}.
\end{equation}
With $\Gamma_f(0)=0$, the spectral density satisfies
\[
 \Gamma_f\in C(I_0;H^1)\cap W^{1,\infty}(I_0;L^2).
\]
\end{lemma}

\begin{proof}
By complex conjugation, it suffices to establish the estimates
for $0<\lambda<\kappa$. Write
\[
 T_\pm=T_\pm(\lambda),\qquad T_\Delta=T_+-T_-,
 \qquad
 \Gamma_f=\frac{\psi_{f,+}-\psi_{f,-}}{2\pi}.
\]
All spatial norms below are on $\mathbb T^2$.

\smallskip\noindent
\emph{Step 1: Multiplier estimates.}
For $2<q<\infty$, set
\[
 \mathcal B_q
 =\operatorname{div}L^q+L^{2q/(q+2)},\qquad
 \|F\|_{\mathcal B_q}
 =\inf_{F=\operatorname{div}A+B}
 \left(\|A\|_q+\|B\|_{2q/(q+2)}\right),
\]
and write $\mathcal B=\mathcal B_4$.
We claim that
\begin{equation}\label{tailD:Bq}
\begin{aligned}
 \sum_{j=0}^1
 \|L_0^{-j}T_\pm u\|_{\mathcal B_q}
 &\le C_q\|u\|_{H^1},\\
 \sum_{j=0}^1
 \|L_0^{-j}T_\Delta u\|_{\mathcal B_q}
 &\le C_qe^{-a_q/\lambda}\|u\|_{H^1}.
\end{aligned}
\end{equation}
Here $L_0^{-1}$ is the inverse with zero streamline average,
extended by zero on $\ker L_0$.

To verify these estimates, apply the logarithmic decomposition
from Lemma~\ref{aux:logregularity}. The flux terms use
$L^{2q}\cdot L^{2q}\subset L^q$, while the scalar terms use
$L^q\cdot L^2\subset L^{2q/(q+2)}$.
Near the resonant periods
\[
 \tau_n=\frac{2\pi n}{\lambda},\qquad
 w_n=\tau_ne^{-\tau_n/4},\qquad n\ge1,
\]
the corresponding coefficients are bounded by
\[
 C(1+\tau_n)^Nw_n^{1/(2q)}
 \quad\hbox{and}\quad
 C(1+\tau_n)^Nw_n^{1/q}.
\]
For every fixed $N$ and $\alpha>0$,
\[
 \sum_{n\ge1}(1+n/\lambda)^N e^{-\alpha n/\lambda}
 \le C_{N,\alpha}e^{-\alpha/(2\lambda)}.
\]
The jump vanishes where $S<2\pi/\lambda$.
A smooth period cutoff supported in
$S>\pi/(2\lambda)$ therefore leaves the jump unchanged
and preserves these exponential estimates.
Applying $L_0^{-1}$ introduces the factor $S/(2\pi in)$;
its derivatives only increase the polynomial powers of $S$.

Away from the resonant intervals, the individual multipliers
are bounded by a fixed polynomial in $S$. Their contributions
are bounded because $(1+S)^N$ belongs to every finite $L^p$.
The same argument permits multiplication by any fixed radial
weight satisfying \eqref{hg:weight-class}.
It also proves the individual estimates with $T_\pm$ replaced
by $T(ic\lambda)$, uniformly for
$\operatorname{Re}c=1$ and $0<|\operatorname{Im}c|\le\kappa$.

In particular,
\begin{equation}\label{tailD:B}
\begin{aligned}
 \|T_\pm u\|_{\mathcal B}
 +\|L_0^{-1}T_\pm u\|_{\mathcal B}
 &\le C\|u\|_{H^1},\\
 \|T_\Delta u\|_{\mathcal B}
 +\|L_0^{-1}T_\Delta u\|_{\mathcal B}
 &\le Ce^{-a_0/\lambda}\|u\|_{H^1}.
\end{aligned}
\end{equation}
These estimates will also be used with the weight
$\lambda\zeta$, which is independent of $\lambda$ and
satisfies \eqref{hg:weight-class}.

\smallskip\noindent
\emph{Step 2: Bounds without spectral derivatives.}
It follows from \eqref{zeroT:inverse},
\eqref{zeroT:projected-LAP}, and
\eqref{zeroT:exponential-jump} that
\begin{equation}\label{tailD:baseline}
\begin{aligned}
 \|(\Delta+T_\pm)^{-1}\|_{H^{-1}\to H^1}
 &\le C\lambda^{-2},\\
 \|\psi_{f,\pm}\|_{H^1}
 &\le C\lambda^{-1}\|f\|_{H^1},\\
 \|\Gamma_f\|_{H^1}
 &\le Ce^{-a_0/\lambda}\|f\|_{H^1}.
\end{aligned}
\end{equation}
The inverse in this display acts between the spaces of zero
spatial mean. Decrease $a_0$ if necessary so that it is also
admissible in \eqref{tailD:B}.

The individual boundary values satisfy
\[
 \Delta\psi_{f,\pm}
 =-T_\pm\psi_{f,\pm}-L_0^{-1}T_\pm f,
\]
and their difference satisfies
\[
 \Delta\Gamma_f
 =-T_+\Gamma_f
  -\frac1{2\pi}
    \left(L_0^{-1}T_\Delta f+T_\Delta\psi_{f,-}\right).
\]
Since $\mathcal B\hookrightarrow W^{-1,4}$, elliptic
regularity and \eqref{tailD:B} imply
\[
 \|\psi_{f,\pm}\|_{W^{1,4}}
 \le C\left(\|\psi_{f,\pm}\|_{H^1}+\|f\|_{H^1}\right)
 \le C\lambda^{-1}\|f\|_{H^1},
\]
and
\[
\begin{aligned}
 \|\Gamma_f\|_{W^{1,4}}
 &\le C\|\Gamma_f\|_{H^1}
   +Ce^{-a_0/\lambda}
       \left(\|f\|_{H^1}+\|\psi_{f,-}\|_{H^1}\right)\\
 &\le C\lambda^{-1}e^{-a_0/\lambda}\|f\|_{H^1}.
\end{aligned}
\]
After decreasing the exponential rate, these are
\eqref{tailD:individual-W14} and \eqref{tailD:F-W14}.

\smallskip\noindent
\emph{Step 3: Individual spectral derivatives.}
Use the fixed cutoff on $I_0$, with $q_0^2=\kappa/2$.
Write $\eta V$ for the vector field in \eqref{good-der}.
It follows from Lemma~\ref{geom:common-field} that
\[
 \|\eta V\|_\infty+\|\eta V\|_{W^{1,4}}
 +\|K_\eta\|_4+\|\zeta\|_4
 \le C\lambda^{-1},
 \qquad
 \|\eta Vf\|_{H^1}\le C\lambda^{-1}\|f\|_{H^2}.
\]
For $q\le2q_0$ and $n\ne0$,
\[
 |nb(q)-\lambda|
 \ge (1-\kappa/4)|n|-\kappa
 \ge \frac12|n|.
\]
Consequently,
\[
 E_{\eta,+}=E_{\eta,-},\qquad
 \|E_{\eta,\pm}u\|_2\le C\|u\|_2.
\]

Equation \eqref{cg:equation} takes the form
\begin{equation}\label{tailD:individual-common}
\begin{aligned}
 (\Delta+T_\pm)D_\lambda\psi_{f,\pm}
 ={}&
 \frac{i}{\lambda}(T_\pm-I)(\eta Vf)
 -\frac{i}{\lambda^2}(T_\pm-I)f\\
 &+\frac{i}{\lambda}
       \left(E_{\eta,\pm}f+\zeta(T_\pm-I)f\right)\\
 &+\operatorname{div}(K_\eta\nabla\psi_{f,\pm})
   -\zeta T_\pm\psi_{f,\pm}
   -E_{\eta,\pm}\psi_{f,\pm}.
\end{aligned}
\end{equation}
For example,
\[
 \|K_\eta\nabla\psi_{f,\pm}\|_2
 \le\|K_\eta\|_4\|\nabla\psi_{f,\pm}\|_4
 \le C\lambda^{-2}\|f\|_{H^1},
\]
and the weighted version of \eqref{tailD:B} implies
\[
 \|\zeta T_\pm\psi_{f,\pm}\|_{H^{-1}}
 \le C\lambda^{-1}\|\psi_{f,\pm}\|_{H^1}
 \le C\lambda^{-2}\|f\|_{H^1}.
\]
Estimating the other terms in a similar way, we obtain
\[
 \|(\Delta+T_\pm)D_\lambda\psi_{f,\pm}\|_{H^{-1}}
 \le C\lambda^{-2}\|f\|_{H^2}.
\]
Moreover, by integration by parts,
\[
 \operatorname{mean}(D_\lambda\psi_{f,\pm})
 =-\operatorname{mean}(\zeta\psi_{f,\pm}),
 \qquad
 |\operatorname{mean}(D_\lambda\psi_{f,\pm})|
 \le C\lambda^{-2}\|f\|_{H^1}.
\]
Applying \eqref{zeroT:inverse} to the mean-zero part we obtain
\[
 \|D_\lambda\psi_{f,\pm}\|_{H^1}
 \le C\lambda^{-4}\|f\|_{H^2}.
\]
Since
\[
 \partial_\lambda\psi_{f,\pm}
 =D_\lambda\psi_{f,\pm}-\eta V\psi_{f,\pm},
\]
this proves \eqref{tailD:Zpolynomial}.

\smallskip\noindent
\emph{Step 4: The derivative of the spectral density.}
Subtract \eqref{tailD:individual-common} for the two signs.
Since $E_{\eta,+}=E_{\eta,-}$, we obtain
\begin{equation}\label{tailD:jump-common}
\begin{aligned}
 (\Delta+T_+)D_\lambda\Gamma_f
 ={}&
 \frac1{2\pi}D_\lambda\left(\frac{iT_\Delta f}{\lambda}\right)
 +\frac{i}{2\pi\lambda}\zeta T_\Delta f\\
 &+\operatorname{div}(K_\eta\nabla\Gamma_f)
   -\zeta T_+\Gamma_f-E_{\eta,+}\Gamma_f\\
 &-\frac1{2\pi}\zeta T_\Delta\psi_{f,-}
   -\frac1{2\pi}T_\Delta D_\lambda\psi_{f,-}.
\end{aligned}
\end{equation}
The commutation identity for the adapted derivative implies
\begin{equation}\label{tailD:source-jump}
 D_\lambda\left(\frac{iT_\Delta f}{\lambda}\right)
 =\frac{i}{\lambda}T_\Delta(\eta Vf)
  -\frac{i}{\lambda^2}T_\Delta f.
\end{equation}
Indeed,
$\partial_\lambda T_\Delta+[\eta V,T_\Delta]
=E_{\eta,+}-E_{\eta,-}=0$.

It follows from \eqref{tailD:B}, including its version with
weight $\lambda\zeta$, and
\eqref{tailD:F-W14}, \eqref{tailD:Zpolynomial}, and
\eqref{tailD:source-jump} that, for some $a_1>0$,
\[
\begin{aligned}
 \|(\Delta+T_+)D_\lambda\Gamma_f\|_{H^{-1}}
 &\le Ce^{-a_1/\lambda}
 \left(
   \lambda^{-2}\|f\|_{H^2}
   +\lambda^{-1}\|f\|_{H^1}
   +\|D_\lambda\psi_{f,-}\|_{H^1}
 \right)\\
 &\le C\lambda^{-4}e^{-a_1/\lambda}\|f\|_{H^2}.
\end{aligned}
\]
Also,
\[
 \operatorname{mean}(D_\lambda\Gamma_f)
 =-\operatorname{mean}(\zeta\Gamma_f),
 \qquad
 |\operatorname{mean}(D_\lambda\Gamma_f)|
 \le C\lambda^{-1}e^{-a_1/\lambda}\|f\|_{H^1}.
\]
Applying \eqref{zeroT:inverse} and restoring the mean, we obtain
\[
 \|D_\lambda\Gamma_f\|_{H^1}
 \le C\lambda^{-6}e^{-a_1/\lambda}\|f\|_{H^2}.
\]
Together with
\[
 \partial_\lambda\Gamma_f
 =D_\lambda\Gamma_f-\eta V\Gamma_f,\qquad
 \|\eta V\Gamma_f\|_2
 \le C\lambda^{-1}\|\Gamma_f\|_{H^1},
\]
this proves \eqref{tailD:zero-derivative} after choosing
$0<a<a_1$.

Finally, \eqref{tailD:F-W14} and
Proposition~\ref{global:LAP} imply that setting
$\Gamma_f(0)=0$ implies $\Gamma_f\in C(I_0;H^1)$.
For $\varphi\in C_c^\infty(I_0)$ and $\delta>0$,
integration by parts on the two intervals separated by zero we obtain
\[
\begin{aligned}
 \int_{\delta<|\lambda|<\kappa}
       \Gamma_f(\lambda)\varphi'(\lambda)\,d\lambda
 ={}&
 -\int_{\delta<|\lambda|<\kappa}
       \partial_\lambda\Gamma_f(\lambda)\varphi(\lambda)\,d\lambda\\
 &+\Gamma_f(-\delta)\varphi(-\delta)
   -\Gamma_f(\delta)\varphi(\delta).
\end{aligned}
\]
The boundary terms tend to zero in $L^2$ as $\delta\to0$.
It follows from \eqref{tailD:zero-derivative} that the remaining
integrals converge and that
$\Gamma_f\in W^{1,\infty}(I_0;L^2)$.
\end{proof}

\subsubsection{Large spectral parameters}
We next strengthen the fixed-region estimate
\eqref{hg:fixed-exterior-W14} from $W^{1,4}$ to $W^{1,8}$,
uniformly for large $|\lambda|$.
The higher spatial integrability is useful when estimating
two adapted spectral derivatives. Indeed, on a fixed region
$K$ away from the elliptic points,
\[
 \|K_\eta\nabla g\|_{L^4(K)}
 \le \|K_\eta\|_{L^8(K)}\|\nabla g\|_{L^8(K)}
 \le C_K|\lambda|^{-1}\|g\|_{W^{1,8}(K)}.
\]
This controls the gradient flux in
\eqref{eq:fixed-source-derivative} in $L^4$ and is used in
the local elliptic estimate for $D_\lambda g$.
The resulting $W^{1,4}$ bound then controls the flux
\[
 \|K_\eta\nabla D_\lambda g\|_{L^2(K)}
 \le \|K_\eta\|_{L^4(K)}
       \|\nabla D_\lambda g\|_{L^4(K)}
\]
in \eqref{quad:second-weak-equation}; see the proof of
Proposition~\ref{quad:second-regular}.
The following lemma obtains the required uniform
$W^{1,8}$ estimate from the equation
\eqref{r2:fixed-source} and the two transport-derivative
bounds in \eqref{r2:angular-g}.
The region $K$ may meet the separatrix.
\begin{lemma}\label{spectral:exterior-W18}
Let $K$ be a fixed invariant region away from the centers and
$|\lambda|\ge1$. Then
\begin{equation}\label{hg:fixed-exterior-W18}
 \|g\|_{W^{1,8}(K)}
 \le C_K\left(\|H\|_2+\sum_{j=0}^2\|L_0^jg\|_{H^1}\right)
 \le C_KM.
\end{equation}
The bound is also uniform for $g_f(ic\lambda)$ with
$\operatorname{Re}c=1$, $|c-1|$ small, and
$0<|\operatorname{Re}(ic\lambda)|\le1$.
\end{lemma}
\begin{proof}
Choose a fixed invariant cutoff $\xi=1$ near $K$, supported
away from the centers. Applying the logarithmic decomposition
\eqref{w14:ibp} with the Hölder exponents in
Lemma~\ref{aux:logregularity}, we obtain
\begin{equation}\label{hg:exterior-B8}
 \|\xi w_0Th\|_{\mathcal B_8}\le C|\lambda|\|h\|_{H^1}.
\end{equation}
There are $O(|\lambda|)$ resonant indices on regular annuli,
each with bounded cost. Near the separatrix, using fixed-width cutoffs
in the same formula and \eqref{w14:separatrix}, we obtain the bound
\[
 C\|h\|_{H^1}\sum_{\tau_n>\tau_0}(1+\tau_n)^Ne^{-\tau_n/64}
       \le C|\lambda|\|h\|_{H^1}.
\]
The remainder is bounded by
$C\|(1+S)^{N+1}\|_8\|h\|_2$. The polynomial weight $w_0$
and nonzero absorption preserve these estimates, as in
Lemma~\ref{aux:logregularity}. The two-step identity is
\[
 Tg=\frac{i}{\lambda}L_0g+\frac{L_0^2g-TL_0^2g}{\lambda^2},
 \qquad
 \|\xi Tg\|_{\mathcal B_8}
       \le C|\lambda|^{-1}\sum_{j=0}^2\|L_0^jg\|_{H^1};
\]
replace $\lambda$ by $c\lambda$ at nonzero absorption.
Since $L^{8/5}\hookrightarrow W^{-1,8}$ and $H^1\hookrightarrow L^8$,
elliptic regularity applied to
\[
 \Delta(\xi g)=\xi H-\xi Tg
       +2\operatorname{div}(g\nabla\xi)-g\Delta\xi
\]
proves \eqref{hg:fixed-exterior-W18}.
\end{proof}

\begin{lemma}\label{r2:angular-regular-remainder}
Fix $0<\alpha<1$ and take $K$ sufficiently large. On each
interval $I_k$ or $J_k$, use $D_\lambda$ and the fixed spatial
cutoff of Proposition~\ref{prop:spectral-bounds}. For
$|\lambda|\ge K$, we have
\begin{equation}\label{r2:high-regular-g}
 D_\lambda g=\Lambda_g+R_g,\qquad
 \|R_g\|_{H^1}+\|L_0R_g\|_{H^1}
       \le C_\alpha|\lambda|^{-\alpha/2}M.
\end{equation}
On $I_k\cap\{\lambda/k<1\}$, the logarithmic term near
each elliptic point is
\[
 \Lambda_{g,k}=\chi(\log(q/\sqrt\sigma))
 \log(\log(q/\sqrt\sigma)\mp i0\,\operatorname{sgn}k)b_{g,k}(\lambda,\theta),
\qquad \sigma=\sigma_k(\lambda)=8(1-\lambda/k)>0.
\]
$\Lambda_g$ is the sum of these two local terms. It vanishes on $J_k$ and on
$I_k\cap\{\lambda/k\ge1\}$. Its coefficients satisfy
\begin{equation}\label{r2:high-jet-coefficient}
 \|b_{g,k}\|_{H^1_\theta}\le C|k|^{-5/2}M.
\end{equation}
Consequently
\begin{equation}\label{r2:high-regular-density}
 D_\lambda \Gamma_f=\Lambda_f+R_f,\qquad
 \|R_f\|_{H^1}+\|L_0R_f\|_{H^1}
       \le C_\alpha|\lambda|^{-1-\alpha/2}M,
\end{equation}
where
\[
 \Lambda_f=\frac{\Lambda_{g,+}-\Lambda_{g,-}}{2\pi i\lambda}.
\]
Set
$b_{f,k}^{\pm}=\pm b_{g,k}^{\pm}/(2\pi i\lambda)$. Then
\begin{equation}\label{r2:stream-jet-coefficient}
 \|b_{f,k}^{\pm}\|_{L^2(I_k\cap\{\lambda/k<1\};H^1_\theta)}
       \le C|k|^{-7/2}M.
\end{equation}
On a fixed center disk, $L_0R_f$ can be replaced by
$\partial_\theta R_f$. Moreover,
\begin{equation}\label{r2:global-density-tail}
 \|\partial_\lambda \Gamma_f(\lambda)\|_{L^2}
       \le C_\alpha|\lambda|^{-1-\alpha/2}M.
\end{equation}
\end{lemma}

\begin{proof}
We consider $\lambda\ge K$; negative parameters follow by complex
conjugation. On each spectral interval, $q_0$ denotes its fixed
cutoff scale, so $q_0\sim\lambda^{-1/2}$.
All spatial norms below use physical coordinates and measure,
unless $dv\,d\theta$ is specified. A set $q<r$
denotes the union of the two corresponding center disks.

\smallskip\noindent
\emph{Step 1: The logarithmic correction and
\eqref{r2:high-jet-coefficient}.}
For $|\lambda-k|<5\kappa/4$, with the same cutoff radius
$(q_{0,k}^{I})^2=32\kappa/k$, a positive
$\sigma=8(1-\lambda/k)$ satisfies
\[
 0<\sigma<\frac{10\kappa}{k},
 \qquad e\sigma<(q_{0,k}^{I})^2.
\]
Every other critical radius in a center chart satisfies
$q_c^2>4q_0^2$. For
\[
 \lambda\in[k+3\kappa/4,\;k+1-3\kappa/4],
 \qquad q_0=q_{0,k}^{J},
 \qquad q_0^2=\frac{\kappa}{2(1+|k|)},
\]
we have $\operatorname{dist}(\lambda,\mathbb Z)\ge3\kappa/4$.
Since $1-b(q)\le\kappa/[4(1+|k|)]$ for $q\le2q_0$,
splitting into $0<|n|\le2(1+|k|)$ and
$|n|>2(1+|k|)$ implies
\[
 |nb(q)-\lambda|\ge\kappa/4,
 \qquad n\ne0,\quad 0<q\le2q_0.
\]
Thus no critical radius lies in $q\le2q_0$.
In particular, $q_0\le q\le2q_0$ is separated from all critical
circles in both cases.

On $I_k\cap\{\lambda<k\}$, take the boundary correction of
Lemma~\ref{adapted:first-center-source}, whose coefficient is
\[
 b_{g,k}(\lambda,\theta)
 =\frac{16\xi_0(\sqrt\sigma)b(\sqrt\sigma)
          g_{f,\pm,k}(\sqrt\sigma,\lambda)e^{ik\theta}}
        {k\sigma\,\mathsf G_{vv}(\sqrt\sigma,\theta)}.
\]
Since $\lambda=kb(\sqrt\sigma)$, it follows from
\eqref{r2:strong-onshell} and Lemma~\ref{geom:metric} that
\[
 \|b_{g,k}\|_{L^2_\theta}\le Ck^{-7/2}M,
 \qquad
 \|b_{g,k}\|_{H^1_\theta}\le Ck^{-5/2}M.
\]
This proves \eqref{r2:high-jet-coefficient}.
The support condition on $\chi$ implies
$\operatorname{supp}\Lambda_{g,k}\subset\{q^2<e\sigma\}
\subset\{q<q_0\}$. With $s=\log(q/\sqrt\sigma)$,
\[
 dx\,dy=\xi_0(\sqrt\sigma e^s)\sigma e^{2s}\,ds\,d\theta,
 \qquad
 \|\Lambda_g\|_2
 \le C\sqrt\sigma\|b_{g,k}\|_{L^2_\theta}
 \le Ck^{-4}M.
\]
Here we used $\chi(s)\log(s\mp i0)\in L^2(ds)$ and summed
over the two centers. On $J_k$ and on $I_k\cap\{\lambda\ge k\}$,
set $\Lambda_g=0$.

\smallskip\noindent
\emph{Step 2: The $H^1$ bound for $R_g$ in
\eqref{r2:high-regular-g}.}
Set $R_g=D_\lambda g-\Lambda_g$.
By \eqref{eq:fixed-source-derivative} and
\eqref{adapted:first-center-source-definition}, on $I_k$ we have
\[
 (\Delta+T)R_g
 =\operatorname{div}\!\bigl(\eta(q/q_0)VH\bigr)
  +\operatorname{div}(K_\eta\nabla g)-\zeta Tg
  +\sum_{\text{two centers}}F_k.
\]
On $J_k$, we replace the last sum by $-E_\eta g$.
We bound these four source terms step by step.

First, the coefficient bounds for $V$ and
\eqref{r2:H-cylinder} imply
\[
 \|\eta(q/q_0)VH\|_2\le C\lambda^{-1/2}M.
\]
Second, it follows from \eqref{r2:mild-center-weight},
\eqref{r2:H-cylinder}, and Lemma~\ref{geom:metric} that
$\|q^{-\alpha}\nabla g\|_{L^2(q<q_*)}\le C_\alpha M$,
after decreasing the fixed $q_*$ if necessary.
Since $|D_{x,y}(\eta(q/q_0)V)|\le C(\lambda q^2)^{-1}$
for $q\ge q_0$ in these disks,
\begin{equation}\label{r2:high-center-commutator}
 \|D_{x,y}(\eta(q/q_0)V)\nabla g\|_{L^2(q<q_*)}
 \le C\lambda^{-1}q_0^{\alpha-2}
       \|q^{-\alpha}\nabla g\|_{L^2(q<q_*)}
 \le C_\alpha\lambda^{-\alpha/2}M.
\end{equation}
On the complement of $\{q<q_*/2\}$, it follows from
\eqref{hg:fixed-exterior-W18} and Hölder's inequality that
\[
 \|K_\eta\nabla g\|_{L^2(\mathbb T^2\setminus\{q<q_*/2\})}
 \le C\|D_{x,y}(\eta V)\|_{L^4(\mathbb T^2\setminus\{q<q_*/2\})}
       \|\nabla g\|_{L^4(\mathbb T^2\setminus\{q<q_*/2\})}
 \le C\lambda^{-1}M.
\]
Thus $\|K_\eta\nabla g\|_2\le
C_\alpha\lambda^{-\alpha/2}M$ globally.

For the third term, use the exact decomposition
\[
 \zeta=V\eta(q/q_0)+\lambda^{-1}\eta(q/q_0)\zeta_0,
 \qquad \zeta_0=\lambda\operatorname{div}V.
\]
The function $V\eta(q/q_0)$ is uniformly bounded and supported
on $q_0\le q\le2q_0$. On this annulus, separation of the
multipliers implies
\begin{equation}\label{adapted:cutoff-transition}
\begin{aligned}
 &\|Th\|_{L^2(q_0\le q\le2q_0)}
  +\|(\partial_\lambda T)h\|_{L^2(q_0\le q\le2q_0)}\\
 &\qquad\le C\|\partial_\theta h\|_{L^2(q_0\le q\le2q_0)}
 \le Cq_0^{1+\alpha}
       \|q^{-\alpha}\nabla h\|_{L^2(q_0\le q\le2q_0)}.
\end{aligned}
\end{equation}
Indeed, $T_0=\partial_\lambda T_0=0$ and
$|T_n|+|\partial_\lambda T_n|\le C_\kappa|n|$ for $n\ne0$;
apply Parseval's identity and $|\partial_\theta h|\le Cq|\nabla h|$.
It follows from this estimate with $h=g$ and
\eqref{hg:T-improved}, with weight $\zeta_0$, that
\[
 \|\zeta Tg\|_{H^{-1}}
 \le C_\alpha\bigl(q_0^{1+\alpha}+\lambda^{-2}\bigr)M.
\]
Finally, Lemma~\ref{adapted:first-center-source} bounds the
last sum by $Ck^{-3/2}M$. On $J_k$, all multipliers in
$q<2q_0$ are separated, so the proof of
\eqref{adapted:cutoff-transition} instead yields
\[
 \|E_\eta g\|_2
 \le Cq_0^{1+\alpha}
           \|q^{-\alpha}\nabla g\|_{L^2(q<2q_0)}
 \le C_\alpha q_0^{1+\alpha}M.
\]
Consequently,
\[
 \|(\Delta+T)R_g\|_{H^{-1}}
 \le C_\alpha\lambda^{-\alpha/2}M.
\]
Lemma~\ref{r2:full-source-LAP} and
Corollary~\ref{cor:global-high-frequency} control the mean-zero
part of $R_g$. Since $\int\partial_\lambda g=0$, Step~1 implies
\[
 \left|\int R_g\right|
 \le C\bigl(\|\eta(q/q_0)V\|_\infty\|\nabla g\|_2
               +\|\Lambda_g\|_2\bigr)
 \le C\bigl(\lambda^{-1/2}+\lambda^{-4}\bigr)M.
\]
This proves $\|R_g\|_{H^1}\le C_\alpha\lambda^{-\alpha/2}M$.

\smallskip\noindent
\emph{Step 3: The $H^1$ bound for $L_0R_g$ in
\eqref{r2:high-regular-g}.}
Set $u=L_0g$ and write $L_0=U\cdot\nabla$, where
$U=(\sin y,-\sin x)$. Since $[L_0,T]=0$,
\[
 (\Delta+T)u=L_0H+\operatorname{div}(C_U\nabla g),
 \qquad
 C_U=D_{x,y}U+(D_{x,y}U)^t
 =\begin{pmatrix}0&\cos y-\cos x\\
                   \cos y-\cos x&0\end{pmatrix}.
\]
In either center disk,
$|C_U|+|\partial_{v,\theta}C_U|\le Cq^2$.
Moreover, $u=b(q)\partial_\theta g$ implies $u_0(q)=0$.
Thus the nonzero-mode argument in the proof of
Lemma~\ref{lem:weighted-center-energy} applies to $u$;
the additional divergence source contributes at most
$C\|q^{2-\alpha}\nabla g\|_2\le C_\alpha M$.
Using \eqref{r2:H-cylinder} for $L_0H$, we obtain
\begin{equation}\label{r2:angular-mild-weight}
 \|q^{-\alpha}\nabla u\|_{L^2(q<q_*)}\le C_\alpha M.
\end{equation}
The identity
\[
 Tu=\frac{i}{\lambda}L_0u-\frac{i}{\lambda}TL_0u,
 \qquad L_0u=L_0^2g,
\]
together with \eqref{hg:exterior-B} and \eqref{r2:angular-g},
implies, for every admissible weight $w_0$,
\[
 \|\chi_\lambda w_0Tu\|_{\mathcal B}\le CM.
\]
Applying \eqref{hg:B-elliptic} to a fixed cutoff of the
equation for $u$, and using \eqref{hg:fixed-exterior-W18}
for $C_U\nabla g$, proves
\begin{equation}\label{r2:angular-exterior-W14}
 \|u\|_{W^{1,4}(\mathbb T^2\setminus\{q<q_*/2\})}\le CM.
\end{equation}

Since $\operatorname{div}U=0$, we have
\[
 L_0\operatorname{div}F
 =\operatorname{div}(L_0F-(F\cdot\nabla)U),
 \qquad [U,\eta(q/q_0)V]=-\frac{\eta(q/q_0)}{\lambda}U.
\]
Commuting the source equation in Step~2 therefore implies, on $I_k$,
\[
\begin{aligned}
 (\Delta+T)L_0R_g
 ={}&\operatorname{div}\!\left(
       \eta(q/q_0)V L_0H+[U,\eta(q/q_0)V]H\right)\\
 &+\operatorname{div}\!\left(K_\eta\nabla u+
    \{L_0K_\eta-K_\eta(D_{x,y}U)^t-D_{x,y}U K_\eta\}\nabla g\right)\\
 &-\zeta Tu+\sum_{\text{two centers}}L_0F_k
       +\operatorname{div}(C_U\nabla R_g).
\end{aligned}
\]
On $J_k$, the sum is replaced by $-E_\eta u$.
The first line has $H^{-1}$ norm at most $C\lambda^{-1/2}M$.
For the second line, differentiation of the bracket identity implies
\[
\begin{aligned}
 L_0D_{x,y}(\eta V)
 ={}&D_{x,y}[U,\eta V]+D_{x,y}U\,D_{x,y}(\eta V)\\
 &+(\eta V\cdot\nabla)D_{x,y}U
       -D_{x,y}(\eta V)D_{x,y}U.
\end{aligned}
\]
Consequently,
\[
 |L_0D_{x,y}(\eta V)|\le\frac{C}{\lambda q^2}
 \quad(q_0\le q<q_*),\qquad
 \|L_0D_{x,y}(\eta V)\|_{L^4(\mathbb T^2\setminus\{q<q_*/2\})}
 \le C\lambda^{-1}.
\]
It follows from \eqref{r2:high-center-commutator},
\eqref{r2:angular-mild-weight}, and
\eqref{r2:angular-exterior-W14} that the flux in the second
line has $L^2$ norm at most $C_\alpha\lambda^{-\alpha/2}M$.
Since $L_0\zeta=0$, the decomposition of $\zeta$ in Step~2,
\eqref{adapted:cutoff-transition} with $h=u$, and the weighted
bound for $Tu$ imply
\[
 \|\zeta Tu\|_{H^{-1}}
 \le C_\alpha\bigl(q_0^{1+\alpha}+\lambda^{-1}\bigr)M.
\]
Lemma~\ref{adapted:first-center-source} bounds the sum of
$L_0F_k$ by $Ck^{-1/2}M$ in $H^{-1}$; on $J_k$,
$\|E_\eta u\|_2\le C_\alpha q_0^{1+\alpha}M$.
Finally,
\[
 \|\operatorname{div}(C_U\nabla R_g)\|_{H^{-1}}
 \le C\|R_g\|_{H^1}\le C_\alpha\lambda^{-\alpha/2}M.
\]
Thus $\|(\Delta+T)L_0R_g\|_{H^{-1}}
\le C_\alpha\lambda^{-\alpha/2}M$.
Since $\int L_0R_g=0$, Lemma~\ref{r2:full-source-LAP} and
Corollary~\ref{cor:global-high-frequency} prove the second
bound in \eqref{r2:high-regular-g}.

The differentiations in Steps~2--3 are justified first for
$T=T(ic\lambda)$, with
\[
 \operatorname{Re}c=1,\qquad
 0<|\operatorname{Im}c|\le\frac{\kappa}{8\lambda},
\]
using difference quotients along the flow of $L_0$.
The uniform estimates in Lemma~\ref{adapted:first-center-source}
and Corollary~\ref{r2:threshold-derivative} allow passage to
$c\to1$ with $\pm\operatorname{Re}(ic\lambda)>0$.

\smallskip\noindent
\emph{Step 4: The density remainder in
\eqref{r2:high-regular-density} and its angular derivative.}
Applying $D_\lambda$ to
$g_{f,+}-g_{f,-}=2\pi i\lambda\Gamma_f$, we obtain
\[
 D_\lambda\Gamma_f=\Lambda_f+R_f,\qquad
 R_f=\frac{R_{g,+}-R_{g,-}}{2\pi i\lambda}
                         -\frac{\Gamma_f}{\lambda}.
\]
It follows from \eqref{lap:density} and
\eqref{r2:tangential-high} that
$\sum_{j=0}^1\|L_0^j\Gamma_f\|_{H^1}\le C\lambda^{-1}M$.
Hence Steps~2--3 imply
\[
 \sum_{j=0}^1\|L_0^jR_f\|_{H^1}
 \le C_\alpha\bigl(\lambda^{-1-\alpha/2}+\lambda^{-2}\bigr)M
 \le C_\alpha\lambda^{-1-\alpha/2}M.
\]
This proves \eqref{r2:high-regular-density}.
On either fixed center disk, $b^{-1}\in W^{1,\infty}$ and
\[
 \partial_\theta R_f=b^{-1}L_0R_f,\qquad
 \|\partial_\theta R_f\|_{H^1(q<q_*)}
 \le C\|L_0R_f\|_{H^1(q<q_*)}.
\]
This proves the asserted replacement of $L_0$ by
$\partial_\theta$ there.

\smallskip\noindent
\emph{Step 5: The coefficient estimate
\eqref{r2:stream-jet-coefficient}.}
Since $b_{f,k}^{\pm}=\pm b_{g,k}^{\pm}/(2\pi i\lambda)$
and $\lambda\sim k$ on $I_k$, Step~1 implies
\[
 \|b_{f,k}^{\pm}(\lambda)\|_{H^1_\theta}\le Ck^{-7/2}M,
 \qquad
 \|b_{f,k}^{\pm}(\lambda)\|_{L^2_\theta}\le Ck^{-9/2}M.
\]
The interval $I_k\cap\{\lambda<k\}$ has length $\kappa$.
Integrating the first inequality in $\lambda$ proves
\eqref{r2:stream-jet-coefficient}.
The second inequality and the change of variables in Step~1 imply
\[
 \|\Lambda_f\|_2\le C\sqrt\sigma\,k^{-9/2}M
                     \le Ck^{-5}M.
\]

\smallskip\noindent
\emph{Step 6: The ordinary spectral derivative in
\eqref{r2:global-density-tail}.}
By \eqref{good-der},
\[
 \partial_\lambda\Gamma_f
 =R_f+\Lambda_f-\eta(q/q_0)V\Gamma_f,
 \qquad
 \|\eta(q/q_0)V\Gamma_f\|_2
 \le C\lambda^{-1/2}\|\Gamma_f\|_{H^1}
 \le C\lambda^{-3/2}M.
\]
It follows from Steps~4--5 and $0<\alpha<1$ that
\[
 \|\partial_\lambda\Gamma_f\|_2
 \le C_\alpha\bigl(\lambda^{-1-\alpha/2}
                +\lambda^{-5}+\lambda^{-3/2}\bigr)M
 \le C_\alpha\lambda^{-1-\alpha/2}M.
\]
This proves \eqref{r2:global-density-tail}.
\end{proof}

\subsubsection{Completion of the proof}

On bounded $J_k$, it follows from \eqref{cg:good-bound} that (ii) holds with
$\Lambda_f=0$. On $I_k$, the corrected derivative in the proof of
Corollary~\ref{r2:threshold-derivative} satisfies
\[
 D_\lambda\Gamma_f
 =\frac{\Lambda_{g,+}-\Lambda_{g,-}}{2\pi i\lambda}
  +\frac{R_{g,+}-R_{g,-}}{2\pi i\lambda}
  -\frac{\Gamma_f}{\lambda}.
\]
At $q_c=\sqrt{\sigma_k}$, since $b(q_c)=\lambda/k$, the jump coefficient is
\[
 \frac{\pm1}{2\pi i\lambda}
 \frac{16\xi_0(q_c)b(q_c)g_{f,\pm,k}(q_c)e^{ik\theta}}
      {k\sigma_k\mathsf G_{vv}(q_c,\theta)}=b_{f,k}^\pm.
\]
For $\kappa\le|\lambda|\le K$, the angular derivative of the
regular part also satisfies
\begin{equation}\label{r2:compact-regular-density}
 \|R_f\|_{H^1}+\|L_0R_f\|_{H^1}\le C_{\kappa,K}M.
\end{equation}
To verify this across the finitely many thresholds, put
$R_g=D_\lambda g-\Lambda_g$. It follows from Lemma~\ref{aux:logregularity},
\eqref{r2:angular-g}, and
$(\Delta+T)L_0g=L_0H+[\Delta,L_0]g$ that
\[
 \|g\|_{W^{1,4}(q>q_0/2)}
 +\|L_0g\|_{W^{1,4}(q>q_0/2)}\le C_{\kappa,K}M.
\]
Here the indicated region includes the fixed exterior of the
center charts. For $L_0=U\cdot\nabla$, the commutator formulas are
\[
\begin{aligned}
 L_0\operatorname{div}(K_\eta\nabla g)
 &=\operatorname{div}\!\left(K_\eta\nabla L_0g+
       \{L_0K_\eta-K_\eta(DU)^t-DU K_\eta\}\nabla g\right),\\
 L_0(\zeta Tg)&=\zeta T L_0g,
 \qquad [\Delta,L_0]R_g=\operatorname{div}(C_U\nabla R_g).
\end{aligned}
\]
The coefficient bounds and the identity for $L_0D(\eta V)$ in
Lemma~\ref{r2:angular-regular-remainder} are uniform on this
finite collection of intervals. It follows from these bounds,
\eqref{eq:fixed-source-derivative}, \eqref{cg:dM}, and
Lemma~\ref{adapted:first-center-source} that
\[
 \|(\Delta+T)L_0R_g\|_{H^{-1}}
 \le C_{\kappa,K}\bigl(\|f\|_{H^5}+\|R_g\|_{H^1}\bigr)
 \le C_{\kappa,K}M.
\]
On $J_k$, the corrected center source is simply $-E_\eta g$,
and $L_0E_\eta g=E_\eta L_0g$ has the same bound by multiplier
separation. Since $\int L_0R_g=0$, Lemma~\ref{r2:full-source-LAP}
bounds $\|L_0R_g\|_{H^1}$. The calculation is first made at
nonzero absorption and passes to the boundary as in
Corollary~\ref{r2:threshold-derivative}. Taking the boundary
jump and using \eqref{r2:tangential-high} proves
\eqref{r2:compact-regular-density}.

It follows from \eqref{r2:strong-onshell} that
\eqref{spectral:log-coefficient} holds; applying
Lemma~\ref{r2:angular-regular-remainder}, we obtain the uniform
large-parameter bounds. It follows from \eqref{tailD:zero-derivative}
that (iii) holds.

For (i), it follows from Corollary~\ref{r2:threshold-derivative},
\eqref{tailD:zero-derivative}, and \eqref{r2:global-density-tail}
with $\alpha=1/2$ that
$\|\partial_\lambda\Gamma_f\|_{L^2_\lambda L^2_x}\le CM$.
Strong continuity at nonzero thresholds and exponential vanishing
at zero exclude additional distributional terms. Using the resolvent identity
\[
 \Gamma_f=(i\lambda)^{-2}\Gamma_{L^2f},\qquad
 \|L^2f\|_{H^1}\le CM,
\]
and the high-frequency LAP, we obtain
$\|\Gamma_f\|_{H^1}\le C|\lambda|^{-3}M$ at infinity.
Together with \eqref{zeroT:exponential-jump} and
Proposition~\ref{global:LAP} on the remaining compact intervals,
this proves
\[
 \|\Gamma_f\|_{L^1_\lambda H^1_x}
 +\|\Gamma_f\|_{L^2_\lambda L^2_x}\le CM,
\]
and hence \eqref{spectral:integrated}.

\subsection{Pointwise estimates}
\label{subsec:additional-first-spectral}

The vorticity representation in
Subsection~\ref{sec:vorticity-representation} involves the angular
coefficients of the boundary values and their spectral jump.
Lemma~\ref{pointwise:elliptic-inputs} supplies the additional local
elliptic estimates. Corollary~\ref{r2:pointwise-spectral-bounds}
then collects their trace and Hölder consequences, using one
interpolation estimate for spectral differences.
Lemma~\ref{r2:orbit-coefficients}
then converts spatial Hölder bounds into coefficient estimates
with a factor $S(r)^{-1}$, uniformly as $r\to2$.

Let $f\in H^5\cap Y$ satisfy $\langle f,g_j\rangle=0$ for
$1\le j\le3$, and set
\[
 M=\|f\|_{H^5},\qquad
 \Gamma=\Gamma_f=\frac{\psi_{f,+}-\psi_{f,-}}{2\pi},\qquad
 g_{f,\pm}=f+i\lambda\psi_{f,\pm}.
\]
Recall that $r(x,y)=2+\cos x+\cos y$ and $b(r)=2\pi/S(r)$.
For a periodic orbit $X_r(s)$ of $L_0$, we use the normalization
\[
 h_k(r)=\frac1{2\pi}\int_0^{2\pi}h(r,\theta)e^{-ik\theta}\,d\theta
       =\frac1{S(r)}\int_0^{S(r)}h(X_r(s))e^{-ikb(r)s}\,ds.
\]
In a center chart the same coefficients are written $h_k(q)$,
where $q=\sqrt{8(1-b(r))}$. Spectral dependence is displayed
as $h_k(r,\lambda)$ or $h_k(q,\lambda)$.

Fix $0<q_a<q_b<q_*$ so that $b\ge1/2$ on $q<q_b$.
The notation $\{q<a\}$ denotes the union of the two center
disks with this radius. Choose $r_*\in(0,1)$ so that
$\{r<r_*\}\cup\{r>4-r_*\}\subset\{q<q_a\}$, and set
\[
 K_j=\{(x,y)\in\mathbb T^2:
           2^{-j}r_*\le r(x,y)\le4-2^{-j}r_*\},
 \quad 0\le j\le3,\qquad K=K_0.
\]
Decrease $\kappa$, if necessary, so that every center cutoff
in \eqref{good-der} equals one near $K_3$. Thus on these sets
\[
 D_\lambda=\partial_\lambda+V,\qquad
 V=-\frac{S(r)}{\lambda S'(r)}\partial_r,
\]
independently of the interval $I_k$ or $J_k$.
Constants below may depend on these fixed regions and on the
indicated exponents. Spatial norms use physical coordinates
and measure unless another measure is indicated.

The first-derivative bounds already control the corrected remainder
in $H^1$. Away from the centers its logarithmic correction vanishes.
The next lemma applies elliptic regularity to the equation for that
remainder, so that Sobolev embedding can subsequently be used.

\begin{lemma}\label{pointwise:elliptic-inputs}
Under the assumptions above, the following estimates hold.

\smallskip\noindent
\textup{(i) Compact spectral intervals.}
For every compact interval $I\Subset\mathbb R\setminus\{0\}$,
\begin{equation}\label{r2:away-W1q}
 \operatorname*{ess\,sup}_{\lambda\in I}
 \bigl(\|\Gamma\|_{W^{1,8/3}(K)}
       +\|L_0\Gamma\|_{W^{1,8/3}(K)}
       +\|D_\lambda\Gamma\|_{W^{1,8/3}(K)}\bigr)\le C_IM.
\end{equation}
The same estimate holds for each $\psi_{f,\pm}$ in place of
$\Gamma$.

\smallskip\noindent
\textup{(ii) Large spectral parameters.}
There is $\Lambda_0>2$ such that, for $|\lambda|\ge\Lambda_0$,
\begin{equation}\label{hp:G-inputs}
\begin{aligned}
 \|i\lambda\Gamma\|_{H^1}&\le C|\lambda|^{-2}M,
 &\|i\lambda L_0\Gamma\|_{H^1}&\le C|\lambda|^{-1}M,\\
 \|i\lambda\Gamma\|_{W^{1,4}(K)}
 +\|i\lambda L_0\Gamma\|_{W^{1,4}(K)}&\le CM,
\end{aligned}
\end{equation}
and
\begin{equation}\label{pointwise:large-adapted}
 \|D_\lambda(i\lambda\Gamma)\|_{W^{1,8/3}(K)}
 \le C|\lambda|M.
\end{equation}

\smallskip\noindent
\textup{(iii) Small spectral parameters.}
Let $p_\infty$ be defined by \eqref{lap:residue} and set
$f_*=f-p_\infty$. There is $0<\lambda_*<\min\{\kappa,1/4\}$ such that,
for $0<|\lambda|<\lambda_*$,
\begin{equation}\label{pointwise:small-elliptic}
\begin{aligned}
 \|\psi_{f_*,\pm}(\lambda)\|_{W^{1,4}(\mathbb T^2)}&\le CM,\\
 \|\Gamma\|_{W^{1,4}(\mathbb T^2)}
 +\|L_0\Gamma\|_{W^{1,4}(\mathbb T^2)}
 +\|D_\lambda\Gamma\|_{W^{1,3}(\mathbb T^2)}
 &\le Ce^{-c/|\lambda|}M.
\end{aligned}
\end{equation}
Derivative estimates are understood almost everywhere in $\lambda$.
\end{lemma}

\begin{proof}
Choose smooth radial functions $\vartheta_j$, $0\le j\le2$, such that
$\vartheta_j=1$ near $K_j$ and $\vartheta_j=0$ outside $K_{j+1}$.
All derivatives of $\vartheta_j$ are supported in
$K_{j+1}\setminus K_j$. Write
$\mathcal B=\operatorname{div}L^4+L^{4/3}$, as in
\eqref{hg:B-elliptic}; its norms use physical measure.

\smallskip\noindent
\emph{Step 1: Claim (i).}
For either boundary value, set $g=g_{f,\pm}$ and $H=(\Delta+1)f$.
The equations
\[
 (\Delta+T)g=H,\qquad
 (\Delta+T)L_0g=L_0H+[\Delta,L_0]g
\]
and Lemma~\ref{aux:logregularity}\textup{(ii)}, together with
\eqref{r2:angular-g}, imply by local elliptic regularity
\[
 \sup_{\lambda\in I}
 \bigl(\|g\|_{W^{1,4}(K_2)}
       +\|L_0g\|_{W^{1,4}(K_1)}\bigr)\le C_IM.
\]
Here $[\Delta,L_0]g$ is a divergence with smooth coefficients
and flux bounded by $C|\nabla g|$.
Write $D_\lambda g=\Lambda_g+R_g$, using the proof of
Corollary~\ref{r2:threshold-derivative} and
Lemma~\ref{r2:angular-regular-remainder}; on $J_k$ use
\eqref{cg:good-bound}. For $\lambda\notin\mathbb Z$, the finite-$c$
$H^1$ bound for $R_g$ passes weakly to either boundary value.
When $\lambda<k$, replacing $\log(e^{2v}-1\mp i0)$ by
$\log(v\mp i0)$, with $v=\log(q/q_c)$, changes the correction
by a function with $H^1$ norm at most $C_IM$, as proved at the
end of Lemma~\ref{adapted:first-center-source}. Thus
$\|R_g\|_{H^1}\le C_IM$ for almost every $\lambda\in I$.
On $K_3$ the functions $\Lambda_g$ and $E_\eta g$ vanish,
and $T$ preserves the radial variable.
Consequently, on the support of $\vartheta_1$, \eqref{eq:fixed-source-derivative}
reduces to
\[
 \Delta R_g=\operatorname{div}(VH+K_\eta\nabla g)
                   -\zeta Tg-TR_g.
\]
The geometric estimates and Lemma~\ref{aux:logregularity}\textup{(ii)}
with the radial weight $\vartheta_2\zeta$ imply
\[
\begin{aligned}
 \|K_\eta\nabla g\|_{L^{8/3}(K_2)}
 &\le\|K_\eta\|_{L^8(K_2)}\|\nabla g\|_{L^4(K_2)}\le C_IM,\\
 \|VH\|_{L^{8/3}(K_2)}&\le C_I\|H\|_{H^1}\le C_IM,\\
 \|\vartheta_2\zeta Tg\|_{\mathcal B}
 +\|\vartheta_2TR_g\|_{\mathcal B}
 &\le C_I(\|g\|_{H^1}+\|R_g\|_{H^1})\le C_IM.
\end{aligned}
\]
The cutoff terms follow from $\|R_g\|_{H^1}\le C_IM$.
Since $\mathcal B\hookrightarrow W^{-1,8/3}$, elliptic regularity
implies $\|D_\lambda g\|_{W^{1,8/3}(K_1)}\le C_IM$.
Now use
\[
 \Gamma=\frac{g_{f,+}-g_{f,-}}{2\pi i\lambda},\qquad
 D_\lambda\psi_{f,\pm}
 =\frac{D_\lambda g_{f,\pm}-Vf}{i\lambda}
                    -\frac{\psi_{f,\pm}}{\lambda}.
\]
Lemma~\ref{geom:common-field} and $f\in H^5$ imply
$\|Vf\|_{W^{1,8/3}(K_1)}\le C_IM$.
These identities prove claim (i), including the individual boundary values.

\smallskip\noindent
\emph{Step 2: Claim (ii).}
For a general mean-zero datum $h$, Corollary~\ref{cor:global-high-frequency}
and \eqref{r2:tangential-high} imply
\begin{align}
 \|\psi_{h,\pm}\|_{H^1}&\le C|\lambda|^{-1}\|h\|_{H^1},
 &\|\psi_{h,\pm}\|_{W^{1,4}(K_1)}&\le C\|h\|_{H^1},
 \label{hp:base-H1}\\
 \|L_0\psi_{h,\pm}\|_{H^1}&\le C|\lambda|^{-1}\|h\|_{H^2},
 &\|L_0\psi_{h,\pm}\|_{W^{1,4}(K_1)}&\le C\|h\|_{H^2}.
 \label{hp:transport-input}
\end{align}
For the second column, apply \eqref{hg:exterior-B} and
\eqref{hg:B-elliptic} to the equations
\[
\begin{aligned}
 (\Delta+T)\psi_{h,\pm}&=(i\lambda)^{-1}(I-T)h,\\
 (\Delta+T)L_0\psi_{h,\pm}
 &=(i\lambda)^{-1}(I-T)L_0h+[\Delta,L_0]\psi_{h,\pm}.
\end{aligned}
\]
The commutator flux in the second equation is controlled by the
first estimate in $W^{1,4}(K_2)$; the same argument applies on
$K_2$ with a cutoff supported in $K_3$.
The identities
\[
 i\lambda\Gamma_f=\Gamma_{Lf}=(i\lambda)^{-1}\Gamma_{L^2f},
 \qquad \|Lf\|_{H^2}+\|L^2f\|_{H^1}\le CM
\]
then prove \eqref{hp:G-inputs}.

To estimate $D_\lambda(i\lambda\Gamma)$, apply the localized
remainder equation in Step~1. It follows from
\eqref{r2:high-regular-g}, \eqref{hg:fixed-exterior-W18},
\eqref{hg:exterior-B}, and \eqref{hg:T-improved} that
\[
\begin{aligned}
 \|K_\eta\nabla g\|_{L^{8/3}(K_2)}&\le C|\lambda|^{-1}M,\\
 \|D_\lambda g\|_{W^{1,8/3}(K_1)}
 &\le C\bigl(\|R_g\|_{H^1}+|\lambda|^{-1}\|H\|_{H^1}
       +\|K_\eta\nabla g\|_{L^{8/3}(K_2)}\\
 &\hspace{24mm}+\|\vartheta_2\zeta Tg\|_{\mathcal B}
       +\|\vartheta_2TR_g\|_{\mathcal B}\bigr)
 \le C|\lambda|M.
\end{aligned}
\]
Here \eqref{hg:exterior-B} bounds the last term by
$C|\lambda|\|R_g\|_{H^1}$, and \eqref{hg:T-improved} applies
with the radial weight $\lambda\zeta$.
Since $i\lambda\Gamma=(g_{f,+}-g_{f,-})/(2\pi)$,
\eqref{pointwise:large-adapted} follows.

\smallskip\noindent
\emph{Step 3: Claim (iii).}
Choose $0<\lambda_*<\min\{\kappa,1/4\}$ so that \eqref{zp:bounded-full} and
\eqref{zeroT:exponential-jump} hold.
Since $p_\infty\in\ker L$, we have $\Gamma_{f_*}=\Gamma_f$.
The estimate \eqref{tailD:Bq} with $q=4$ implies
\begin{equation}\label{zp:DT-B}
 \|T_\pm h\|_{\mathcal B}
 +\|L_0^{-1}T_\pm h\|_{\mathcal B}\le C\|h\|_{H^1}.
\end{equation}
Applying \eqref{hg:B-elliptic} to
\[
 \Delta\psi_{f_*,\pm}
 =-T_\pm\psi_{f_*,\pm}-L_0^{-1}T_\pm f_*
\]
and using \eqref{zp:bounded-full} proves the first estimate
in \eqref{pointwise:small-elliptic}.
For $L_0\Gamma$, set $T_\Delta=T_+-T_-$ and use
\[
\begin{aligned}
 (\Delta+T_\pm)L_0\psi_{f_*,\pm}
 &=-T_\pm f_*+[\Delta,L_0]\psi_{f_*,\pm},\\
 (\Delta+T_+)L_0\Gamma
 &=[\Delta,L_0]\Gamma
   -\frac1{2\pi}T_\Delta(f_*+L_0\psi_{f_*,-}).
\end{aligned}
\]
The first equation, \eqref{zeroT:inverse}, and \eqref{zp:DT-B}
imply $\|L_0\psi_{f_*,\pm}\|_{H^1}\le C|\lambda|^{-2}M$.
It follows from \eqref{tailD:Bq} and
\eqref{zeroT:exponential-jump} that the source in the second
equation is exponentially small in $H^{-1}$.
After applying \eqref{zeroT:inverse}, elliptic regularity and
\eqref{tailD:F-W14} imply
\[
 \|\Gamma\|_{W^{1,4}}
 +\|L_0\Gamma\|_{H^1}+\|L_0\Gamma\|_{W^{1,4}}
 \le Ce^{-c/|\lambda|}M.
\]
For $D_\lambda\Gamma$, use \eqref{tailD:jump-common},
\eqref{tailD:source-jump}, and \eqref{tailD:zero-derivative}.
The gradient flux satisfies
\[
 \|K_\eta\nabla\Gamma\|_3
 \le\|K_\eta\|_{12}\|\nabla\Gamma\|_4
 \le C|\lambda|^{-1}e^{-c/|\lambda|}M.
\]
By \eqref{tailD:Bq} with $q=3$ and
\eqref{tailD:Zpolynomial}, all other terms, including
$T_+D_\lambda\Gamma$, have exponentially small norm in
$\operatorname{div}L^3+L^{6/5}\hookrightarrow W^{-1,3}$.
Elliptic regularity proves the remaining estimate in
\eqref{pointwise:small-elliptic}. Fixed powers of
$|\lambda|^{-1}$ are absorbed by decreasing $c$.
\end{proof}

The following consequences are stated in the form used in
Section~\ref{sec:main-proof}. The angular trace bounds follow
from \eqref{r2:angular-g}; the Hölder bounds follow from the
preceding lemma and interpolation.

\begin{corollary}\label{r2:pointwise-spectral-bounds}
Fix $0<\beta<1/8$ and $0<\gamma<1/4-\beta$.
With $\lambda_*$ and $\Lambda_0$ as in
Lemma~\ref{pointwise:elliptic-inputs}, the following claims hold.

\smallskip\noindent
\textup{(i) Bounds near the elliptic points.}
For every $\lambda_1>0$, $|\lambda|\ge\lambda_1$,
$k\ne0$, and $0<q\le q_a$,
\begin{equation}\label{r2:tail-trace}
 |g_{f,\pm,k}(q,\lambda)|\le C_{\lambda_1}|k|^{-5/2}M,
 \qquad
 |\Gamma_k(q,\lambda)|\le C_{\lambda_1}|k|^{-5/2}|\lambda|^{-1}M.
\end{equation}
For $|\lambda|\ge\Lambda_0$,
\begin{equation}\label{r2:tail-high-lambda}
 |\Gamma_k(q,\lambda)|\le C|k|^{-1/2}|\lambda|^{-3}M.
\end{equation}
For every integer $N\ge0$ and $0<|\lambda|<\lambda_*$,
\begin{equation}\label{r2:tail-zero}
 \|\Gamma(\lambda)\|_{H^N(q<q_a)}
 \le C_Ne^{-c_N/|\lambda|}M.
\end{equation}

\smallskip\noindent
\textup{(ii) Compact spectral intervals away from zero.}
For every compact interval $I\Subset\mathbb R\setminus\{0\}$,
\begin{equation}\label{r2:away-holder}
 \sup_{\lambda\in I}
 \bigl(\|\Gamma\|_{C^\beta(K)}+\|L_0\Gamma\|_{C^\beta(K)}\bigr)
 +\|\Gamma\|_{C^\gamma(I;C^\beta(K))}\le C_IM.
\end{equation}
The same estimate holds for each $\psi_{f,\pm}$ in place of
$\Gamma$.

\smallskip\noindent
\textup{(iii) Large spectral parameters.}
For $|\lambda|\ge\Lambda_0$,
\begin{equation}\label{hp:G-Holder}
 \|i\lambda\Gamma_f(\lambda)\|_{C^\beta(K)}
 \le C|\lambda|^{-3/2}M,
 \qquad
 \|i\lambda L_0\Gamma_f(\lambda)\|_{C^\beta(K)}
 \le C|\lambda|^{-3/4}M.
\end{equation}
If $|\Lambda|\ge2\Lambda_0$ and
$\lambda,\mu\in[\Lambda-1,\Lambda+1]$, then
\begin{equation}\label{hp:spectral-holder}
 \|i\lambda\Gamma_f(\lambda)-i\mu\Gamma_f(\mu)\|_{C^\beta(K)}
 \le C|\Lambda|\,|\lambda-\mu|^\gamma M.
\end{equation}

\smallskip\noindent
\textup{(iv) Small spectral parameters.}
For $0<|\lambda|<\lambda_*$,
\begin{equation}\label{zp:individual-holder}
 \|\psi_{f_*,\pm}(\lambda)\|_{W^{1,4}(\mathbb T^2)}
 +\|\psi_{f_*,\pm}(\lambda)\|_{C^\beta(\mathbb T^2)}
 \le C_\beta M,
 \qquad f_*=f-p_\infty.
\end{equation}
This estimate also holds for every $0<\beta<1/2$.
For every $0<a<1/3$ and almost every such $\lambda$,
\begin{equation}\label{zp:jump-holder}
 \|\Gamma(\lambda)\|_{C^a(K)}
 +\|L_0\Gamma(\lambda)\|_{C^a(K)}
 +\|D_\lambda\Gamma(\lambda)\|_{C^a(K)}
 \le C_ae^{-c_a/|\lambda|}M.
\end{equation}
For $0<\rho<\lambda_*/2$ and either sign,
\begin{equation}\label{zp:scaled-holder}
 \|\Gamma(\pm\rho\,\cdot)\|_{C^\gamma([1/2,2];C^\beta(K))}
 +\sup_{\rho/2\le|\lambda|\le2\rho}
       \|L_0\Gamma(\lambda)\|_{C^\beta(K)}
 \le Ce^{-c/\rho}M.
\end{equation}
\end{corollary}

\begin{proof}
\emph{Step 1: Claim (i).}
Choose a smooth radial $\chi$ equal to one on $q\le q_a$
and zero on $q\ge q_b$. With $v=\log q$ and
$\|u\|_{H^1_k}^2=\int_{\mathbb R}(|u'|^2+k^2|u|^2)\,dv$,
\eqref{r2:angular-g} and the trace inequality imply
\[
 \sum_{k\ne0}|k|^4\|\chi g_{f,\pm,k}\|_{H^1_k}^2\le CM^2,
 \qquad
 |u(v)|^2\le2\|u\|_2\|u'\|_2
             \le |k|^{-1}\|u\|_{H^1_k}^2.
\]
Together with $g_{f,+}-g_{f,-}=2\pi i\lambda\Gamma_f$,
these estimates prove \eqref{r2:tail-trace}.
For \eqref{r2:tail-high-lambda}, use the same trace inequality and
\[
 \Gamma_f=(i\lambda)^{-2}\Gamma_{L^2f},\qquad
 \|\Gamma_{L^2f}\|_{H^1}
 \le C|\lambda|^{-1}\|L^2f\|_{H^1}\le C|\lambda|^{-1}M.
\]
For small $|\lambda|$, the inequalities $b\ge1/2$ and
$|\lambda|<1/4$ imply $|nb-\lambda|\ge|n|/4$ for $n\ne0$.
Thus $T_+=T_-$ and
\[
 (\Delta+T_\pm)\Gamma=0\qquad(q<q_b).
\]
The smooth periodic flow of $b^{-1}L_0$ and this multiplier
separation imply uniform boundedness of $T_\pm$ on spatial
$H^m$ in these disks. Interior elliptic regularity therefore implies,
for $q_a\le r_1<r_2\le q_b$,
\[
 \|\Gamma\|_{H^{m+2}(q<r_1)}
 \le C_{m,r_1,r_2}\|\Gamma\|_{H^{m+1}(q<r_2)}.
\]
Iteration and \eqref{zeroT:exponential-jump} prove
\eqref{r2:tail-zero} for every $N$.

\smallskip\noindent
\emph{Step 2: Claim (ii) and the interpolation estimate.}
For $p>2$ and $0<\gamma<1-2/p-\beta$, Sobolev embedding and
interpolation on the fixed smooth domain $K$ imply
\begin{equation}\label{pointwise:spectral-interpolation}
\begin{aligned}
 \|u(\lambda)-u(\mu)\|_{C^\beta(K)}
 &\le C\|u(\lambda)-u(\mu)\|_{L^p(K)}^\gamma
        \|u(\lambda)-u(\mu)\|_{W^{1,p}(K)}^{1-\gamma}\\
 &\le C|\lambda-\mu|^\gamma
       \|\partial_\lambda u\|_{L^\infty(I;L^p(K))}^\gamma
       \|u\|_{L^\infty(I;W^{1,p}(K))}^{1-\gamma}.
\end{aligned}
\end{equation}
Indeed, $W^{1-\gamma,p}(K)\hookrightarrow C^\beta(K)$.
For $u=\Gamma$ or $\psi_{f,\pm}$, Corollary~\ref{r2:threshold-derivative}
and \eqref{cg:good-bound} imply local absolute continuity into $L^2$
through every nonzero threshold. Hence the identity
\[
 \partial_\lambda u=D_\lambda u-V\cdot\nabla u
\]
holds distributionally on $I$, without point masses at the thresholds.
By \eqref{r2:away-W1q} and $\|V\|_{L^\infty(I\times K)}\le C_I$,
\[
 \|\partial_\lambda u\|_{L^\infty(I;L^{8/3}(K))}
 +\|u\|_{L^\infty(I;W^{1,8/3}(K))}\le C_IM.
\]
Apply \eqref{pointwise:spectral-interpolation} with $p=8/3$.
The spatial bounds in \eqref{r2:away-holder} follow directly from
$W^{1,8/3}\hookrightarrow C^\beta$; their continuous values follow
from the limiting absorption principle and weak compactness.

\smallskip\noindent
\emph{Step 3: Claim (iii).}
The interpolation inequality
\[
 \|u\|_{W^{1,16/7}(K)}
 \le C\|u\|_{H^1(K)}^{3/4}\|u\|_{W^{1,4}(K)}^{1/4}
\]
and $W^{1,16/7}\hookrightarrow C^\beta$, $\beta<1/8$,
applied to \eqref{hp:G-inputs}, prove \eqref{hp:G-Holder}.
For $u(\lambda)=i\lambda\Gamma_f(\lambda)$ and
$I=[\Lambda-1,\Lambda+1]$, \eqref{hp:G-inputs} and
\eqref{pointwise:large-adapted} imply
\[
 \|u\|_{L^\infty(I;W^{1,8/3}(K))}\le CM,
 \qquad
 \|\partial_\lambda u\|_{L^\infty(I;L^{8/3}(K))}
 \le C|\Lambda|M.
\]
The identity used in Step~2 is valid for this $u$ as well.
Equation \eqref{pointwise:spectral-interpolation} implies
\eqref{hp:spectral-holder}, since $|\Lambda|^\gamma\le|\Lambda|$.

\smallskip\noindent
\emph{Step 4: Claim (iv).}
The embeddings $W^{1,4}\hookrightarrow C^\beta$ for
$\beta<1/2$ and $W^{1,3}\hookrightarrow C^a$ for $a<1/3$
applied to \eqref{pointwise:small-elliptic} prove
\eqref{zp:individual-holder} and \eqref{zp:jump-holder}.
For $u(\nu)=\Gamma(\pm\rho\nu)$, $1/2\le\nu\le2$,
\[
 \partial_\nu u
 =\pm\rho D_\lambda\Gamma(\pm\rho\nu)
  \mp\rho V(\pm\rho\nu)\cdot\nabla\Gamma(\pm\rho\nu),
 \qquad
 \|\rho V(\pm\rho\nu)\|_{L^\infty(K)}\le C.
\]
Consequently, \eqref{pointwise:small-elliptic} implies
\[
 \|u\|_{L^\infty([1/2,2];W^{1,3}(K))}
 +\|\partial_\nu u\|_{L^\infty([1/2,2];L^3(K))}
 \le Ce^{-c/\rho}M.
\]
Apply \eqref{pointwise:spectral-interpolation} with $p=3$;
the condition $\gamma<1/4-\beta$ implies $\gamma<1/3-\beta$.
This proves \eqref{zp:scaled-holder}, including its supremum term
by \eqref{zp:jump-holder}.
\end{proof}

The next lemma applies to the spatial functions controlled
above. It specifies the full range of the orbit coordinate;
the periodic-orbit coefficients are not defined at $r=2$.

\begin{lemma}
\label{r2:orbit-coefficients}
Let $K=\{r_*\le r(x,y)\le4-r_*\}$, with $0<r_*<2$ fixed.
For
\[
 r\in[r_*,2)\cup(2,4-r_*],\qquad k\in\mathbb Z\setminus\{0\},
 \qquad 0<\beta<1,
\]
the following claims hold uniformly in $r$ and $k$.

\smallskip\noindent
\textup{(i)} For $h\in C^\beta(K)$,
\[
 |h_k(r)|\le C_{\beta,K}
       \bigl(S(r)^{-1}+|k|^{-1}\bigr)\|h\|_{C^\beta(K)}.
\]
\textup{(ii)} If $h(p)=0$ at the two hyperbolic fixed points
$p=(0,\pi),(\pi,0)$, then
\[
 |h_k(r)|\le C_{\beta,K}S(r)^{-1}
       \min\{1,|kb(r)|^{-\beta}\}\|h\|_{C^\beta(K)}.
\]
\textup{(iii)} If $u,L_0u\in C^\beta(K)$, with $L_0u$
understood distributionally in the interior of $K$, then
$L_0u(p)=0$ at both hyperbolic fixed points and, for
$|kb(r)|\ge1$,
\[
 |u_k(r)|\le C_{\beta,K}S(r)^{-1}|kb(r)|^{-1-\beta}
                                      \|L_0u\|_{C^\beta(K)}.
\]
\end{lemma}

\begin{proof}
Write $U=(\sin y,-\sin x)$, so $L_0=U\cdot\nabla$.

\smallskip\noindent
\emph{Step 1: Claim (i).}
For $r$ close to $2$, divide one period into the time intervals
$I_j$ in fixed saddle neighborhoods and the remaining arcs,
whose total travel time is uniformly bounded.
In the coordinates $ab=r-2$ of Lemma~\ref{geom:common-field},
set $d=|r-2|$. On each saddle arc,
$ds=j(a)j(d/a)\,da/a$, with $j$ bounded. For fixed $a_0>0$,
\[
 \int_{d/a_0}^{a_0}\left(a+\frac da\right)^\beta\frac{da}{a}
 \le\int_{d/a_0}^{a_0}
       (a^\beta+d^\beta a^{-\beta})\frac{da}{a}
 \le\frac{2a_0^\beta}{\beta}.
\]
Since the coordinate map is regular and $|U(x)|\le C|x-p_j|$
near a saddle, it follows that
\[
 \sum_j\int_{I_j}|X_r(s)-p_j|^\beta\,ds\le C_{\beta,K},
 \qquad
 \int_0^{S(r)}|U(X_r(s))|^\beta\,ds\le C_{\beta,K}.
\]
Set $P_r(s)=h(p_j)$ on $I_j$ and $P_r(s)=0$ on the remaining
arcs. Then
\[
 \int_0^{S(r)}|h(X_r(s))-P_r(s)|\,ds
       \le C_{\beta,K}\|h\|_{C^\beta(K)}.
\]
The periodic step function $P_r$ has a uniformly bounded
number of jumps. Integration on its constant intervals implies
\[
 \left|\frac1{S(r)}\int_0^{S(r)}P_r(s)e^{-ikb(r)s}\,ds\right|
 \le\frac1{2\pi|k|}\sum_s|[P_r](s)|
 \le C_K|k|^{-1}\|h\|_{C^\beta(K)}.
\]
These estimates prove claim (i) near $r=2$.
On the remaining compact range of $r$, $S(r)$ is bounded
above, so $|h_k(r)|\le\|h\|_\infty$ proves the same claim.
The bound for $\int|U(X_r(s))|^\beta ds$ also holds there.

\smallskip\noindent
\emph{Step 2: Claim (ii).}
If $h(p_j)=0$, then $P_r=0$, and Step~1 implies
\[
 |h_k(r)|\le C_{\beta,K}S(r)^{-1}\|h\|_{C^\beta(K)}.
\]
Smoothness of $U$ and Gronwall's inequality imply
\[
 |X_r(s+\tau)-X_r(s)|\le C|\tau|\,|U(X_r(s))|,
 \qquad |\tau|\le1.
\]
Consequently,
\[
 \int_0^{S(r)}|h(X_r(s+\tau))-h(X_r(s))|\,ds
 \le C_{\beta,K}|\tau|^\beta\|h\|_{C^\beta(K)}.
\]
By periodicity,
\[
 (e^{ikb(r)\tau}-1)h_k(r)
 =\frac1{S(r)}\int_0^{S(r)}
 [h(X_r(s+\tau))-h(X_r(s))]e^{-ikb(r)s}\,ds.
\]
For $|kb(r)|\ge\pi$, choose $\tau=\pi/|kb(r)|$; the factor
on the left has absolute value two. For $|kb(r)|<\pi$, use
the first estimate of this step. This proves claim (ii).

\smallskip\noindent
\emph{Step 3: Claim (iii).}
Let $\Phi_t$ be the flow of $U$.
The distributional transport identity and continuity imply
\[
 u(\Phi_t(x))-u(x)=\int_0^tL_0u(\Phi_s(x))\,ds.
\]
This follows in flow coordinates away from fixed points and
extends to them by continuity. Since $\Phi_t(p)=p$, it follows
that $tL_0u(p)=0$, and hence $L_0u(p)=0$.
Integration by parts along each periodic orbit implies
\[
 (L_0u)_k(r)=ikb(r)u_k(r).
\]
Apply claim (ii) to $h=L_0u$ and divide by $|kb(r)|$.
This proves claim (iii).
\end{proof}

\section{The second spectral derivative}\label{sec:second-spectral}
The spectral representation \eqref{lap:spectral-representation} reduces
the $t^{-2}$ bound for the stream function to the following estimate.
The derivative $\partial_\lambda^2\Gamma_f$ contains differentiated
logarithms (and hence Cauchy type singularities), so we estimate its Fourier transform directly.

\begin{proposition}\label{quad:spectral}
Let $f\in H^5\cap Y$ satisfy
$\langle f,g_j\rangle_{L^2}=0$, $1\le j\le3$, and let
$\Gamma_f$ be defined by \eqref{lap:density}. Then
\begin{equation}\label{quad:fourier-second}
 \sup_{t\in\mathbb R}
 \left\|\int_{\mathbb R}e^{-it\lambda}
       \partial_\lambda^2\Gamma_f(\lambda)\,d\lambda\right\|_{L^2}
 \le C\|f\|_{H^5}.
\end{equation}
The derivative and Fourier transform are understood in the sense
of $L^2$-valued distributions.
\end{proposition}

\subsection{Proof of Proposition~\ref{quad:spectral}}
Throughout this subsection, set $M=\|f\|_{H^5}$ and fix
$0<\alpha,\beta<1$. Recall that the first derivative has the form
$D_\lambda\Gamma_f=\Lambda_f+R_f$, where $\Lambda_f$ is the
explicit logarithm in \eqref{spectral:decomposition}.
The proof first establishes integrability of the differentiated
remainder, then bounds the Fourier transform of the logarithm
and more regular terms.

\subsubsection{The differentiated equation and its coefficients}
\label{subsec:second-weak-equation}
On each spectral interval, fix the radial cutoff
$\eta=\eta(q/q_0)$. The adapted derivative and the coefficients
in its commutation with $\Delta+T(ic\lambda)$ are
\[
\begin{aligned}
 D_\lambda
 &=\partial_\lambda+\eta V,
 \qquad
 \zeta=\operatorname{div}(\eta V),\\
 K_\eta
 &=D_x(\eta V)+(D_x(\eta V))^t-\zeta\,\mathrm{Id},\\
 E_\eta(c,\lambda)
 &=\partial_\lambda T(ic\lambda)
   +[\eta V,T(ic\lambda)]
  =(1-\eta)\partial_\lambda T(ic\lambda).
\end{aligned}
\]
Here $\eta V$ acts as a first-order differential operator,
$[D_x(\eta V)]_{ij}=\partial_j(\eta V)_i$, and all spatial
derivatives are taken in the physical coordinates $(x,y)$.
Every derivative in $\lambda$ is taken with $c$ fixed.

Set $H=(\Delta+1)f$ and $g=g_f(ic\lambda)$, where
$\operatorname{Re}c=1$, $0<|c-1|$ is sufficiently small, and
$|\operatorname{Re}(ic\lambda)|\le1$.
Differentiate with $c$ fixed; the boundary formulas follow as
$c\to1$ from either side. All divergences and matrix derivatives are
in the physical coordinates $(x,y)$.
Recall $\zeta=\operatorname{div}(\eta V)$. By direct computation we have the identity
\begin{equation}\label{quad:weak-divergence}
 (D_\lambda+\zeta)\operatorname{div}F
 =\operatorname{div}\{D_\lambda F-D_x(\eta V)F+\zeta F\}.
\end{equation}
Since $H=(\Delta+1)f$ is independent of $\lambda$,
\begin{equation}\label{quad:source-flux}
 (D_\lambda+\zeta)^2H
 =\operatorname{div}\{
 \eta V(\eta V\cdot\nabla H)
       +[\partial_\lambda(\eta V)+\zeta\eta V]H\}.
\end{equation}
Since $\partial_\lambda H=0$ and
$\zeta=\operatorname{div}(\eta V)$, we have
\[
 (D_\lambda+\zeta)H
 =\eta V\cdot\nabla H+\zeta H
 =\operatorname{div}(\eta VH).
\]
Consequently,
\[
 (D_\lambda+\zeta)^2H
 =(D_\lambda+\zeta)\operatorname{div}(\eta VH),
\]
to which \eqref{quad:weak-divergence} applies with $F=\eta VH$.
It follows from \eqref{quad:weak-divergence} and
\eqref{eq:fixed-source-derivative} that
\begin{equation}\label{quad:second-weak-equation}
\begin{aligned}
 (\Delta+T)D_\lambda^2g={}&
 \operatorname{div}\{\eta V(\eta V\cdot\nabla H)
             +[\partial_\lambda(\eta V)+\zeta\eta V]H\}\\
 &+\operatorname{div}\bigl\{
 2K_\eta\nabla D_\lambda g
 +[D_\lambda K_\eta-K_\eta D_x(\eta V)^t
          -D_x(\eta V)K_\eta+\zeta K_\eta]\nabla g\bigr\}\\
 &-(D_\lambda \zeta+\zeta^2)Tg-2\zeta TD_\lambda g
 -(D_\lambda E_\eta+2\zeta E_\eta)g-2E_\eta D_\lambda g.
\end{aligned}
\end{equation}
Here $D_\lambda K_\eta$ denotes componentwise differentiation,
\[
 (D_\lambda K_\eta)_{ij}
 =\partial_\lambda(K_\eta)_{ij}
   +\eta V\cdot\nabla(K_\eta)_{ij},
\]
whereas $D_\lambda E_\eta$ denotes the operator commutator
\[
 (D_\lambda E_\eta)h
 :=D_\lambda(E_\eta h)-E_\eta D_\lambda h
 =(\partial_\lambda E_\eta)h
   +\eta V\cdot\nabla(E_\eta h)
   -E_\eta(\eta V\cdot\nabla h).
\]
In $\partial_\lambda E_\eta$, the operator family is differentiated
with $c$ and its argument held fixed.

The coefficient in the second term on the right-hand side of \eqref{quad:second-weak-equation} satisfies
\begin{equation}\label{quad:weak-coefficients}
\begin{split}
 &|D_\lambda K_\eta-K_\eta D_x(\eta V)^t
          -D_x(\eta V)K_\eta+\zeta K_\eta|\\
 &\qquad\le C\bigl(|D_\lambda D_x(\eta V)|
      +|D_x(\eta V)|^2+|D_\lambda \zeta|+\zeta^2\bigr).
\end{split}
\end{equation}
We note also that the first line of \eqref{quad:second-weak-equation}
contains only $H$ and $\nabla H$.

For the coefficients near the separatrix, it follows from
Lemma~\ref{geom:period-lemma} that
\[
 |r-2|^j\left|\partial_r^j\frac{S}{S'}\right|
 \le C|r-2|(1+|\log|r-2||),\qquad 0\le j\le2.
\]
On a regular separatrix arc, in the coordinates of
Lemma~\ref{geom:common-field},
$\lambda V=-(S/S')\partial_r-\sigma\partial_\sigma$.
Since the chart and its derivatives are bounded, the product rule implies
\[
 |D_x(\lambda V)|+|(\lambda V)\cdot\nabla D_x(\lambda V)|
 \le C\left(1+\left|\left(\frac S{S'}\right)'\right|^2
       +\left|\frac S{S'}\left(\frac S{S'}\right)''\right|\right)
 \le C(1+|\log|r-2||)^2.
\]

By symmetry it suffices to work near the saddle $(\pi,0)$, where
\[
 a=\sqrt{2}\sin\frac{x-\pi+y}{2},\qquad
 b=\sqrt{2}\sin\frac{x-\pi-y}{2},\qquad r-2=ab.
\]
In the sector $0<|b|\le |a|\le a_0$, the forward travel time from
$(a,b)$ along its streamline to the section $|a|=\sqrt{2}$ is
\[
 t(a,b)=\int_{|a|}^{\sqrt{2}}j(u)j(ab/u)\,\frac{du}{u},
 \qquad j(u)=(1-u^2/2)^{-1/2},
\]
as in \eqref{geom:field-travel-time}.
On $0<|b|\le |a|\le a_0$,
expand the integrand as
\[
 j(u)j(ab/u)\frac1u
 =\frac{j(u)}u+
 \sum_{n\ge1}c_n(ab)^{2n}\frac{j(u)}{u^{2n+1}},
 \qquad j(z)=\sum_{n\ge0}c_nz^{2n}.
\]
The series and its first three $b$ derivatives converge uniformly
for $|ab/u|\le a_0$. For $1\le i+j\le3$, differentiation of
the integral, including its lower endpoint, therefore implies
\[
 |a|^i|\partial_a^it|\le C_i\quad(i\ge1),\qquad
 |a|^i|\partial_a^i\partial_b^jt|
 \le C_{i,j}|b|^{\max(2-j,0)}\quad(j\ge1),
 \qquad |t|\le C(1+|\log|a||).
\]
The first derivatives of the numerators in
\eqref{geom:field-components} are
\[
 \begin{aligned}
 \partial_a(J\lambda V^a)
 &=t+at_a-b\left(\frac S{S'}\right)'t_b-\frac S{S'}t_{ab},\\
 \partial_b(J\lambda V^a)
 &=a\left[1-\left(\frac S{S'}\right)'\right]t_b-\frac S{S'}t_{bb},\\
 \partial_a(J\lambda V^b)
 &=b\left[\left(\frac S{S'}\right)'-1\right]t_a+\frac S{S'}t_{aa},\\
 \partial_b(J\lambda V^b)
 &=a\left(\frac S{S'}\right)'t_a+\frac S{S'}t_{ab}-t-bt_b.
 \end{aligned}
\]
The quotient and its derivatives are evaluated at $r=2+ab$.
After applying $a\partial_a$ or $b\partial_b$, the terms with
second quotient derivatives or three $a$ derivatives of $t$ satisfy
\[
 \begin{split}
 &\left|ab^2\left(\frac S{S'}\right)''t_a\right|
 +\left|a\frac S{S'}t_{aaa}\right|
 +\left|ab\left(\frac S{S'}\right)'t_{aa}\right|
 \le C\frac{|b|}{|a|}(1+|\log|ab||),\\
 &\left|a^2b\left(\frac S{S'}\right)''t_a\right|
 \le C(1+|\log|ab||).
 \end{split}
\]
The other terms obey the same bounds, and $J,J^{-1}$ have
bounded derivatives through order two. Since $|b/a|\le1$, it follows that
\[
 \left|\frac{\lambda V^a}{a}\right|
 +\left|\frac{\lambda V^b}{b}\right|
 +|D_{a,b}(\lambda V)|
 +|a\partial_aD_{a,b}(\lambda V)|
 +|b\partial_bD_{a,b}(\lambda V)|
 \le C(1+|\log|a||+|\log|b||)^2.
\]
The same bounds hold on $0<|a|\le|b|\le a_0$ after exchanging
$a,b$. Since the physical coordinate change is smooth, it follows that
\begin{equation}\label{quad:saddle-coefficients}
 \|D_{x,y}V\|_{L^p}+\|V\cdot\nabla D_{x,y}V\|_{L^p}
 +\|\partial_\lambda D_{x,y}V\|_{L^p}\le C_{I,p},
 \qquad 1\le p<\infty,\quad \lambda\in I\Subset\mathbb R\setminus\{0\},
\end{equation}
on every fixed region away from the centers. Here
$\partial_\lambda V=-V/\lambda$. The derivatives of $\eta(q/q_0)$ are supported in
$q_0\le q\le2q_0$, where $V$ is smooth. The product rule
therefore extends \eqref{quad:saddle-coefficients} to
$\eta(q/q_0)V$, with constants allowed to depend on $q_0$.

To justify the distributional derivatives across the separatrix,
first consider a regular separatrix chart. On $r-2=\pm\varepsilon$,
\eqref{geom:field-period-bounds} and
$\lambda V=-(S/S')\partial_r-\sigma\partial_\sigma$ imply
\[
 |V\cdot n|\le C_I\varepsilon(1+|\log\varepsilon|),
 \qquad
 |D_xV|\le C_I(1+|\log\varepsilon|).
\]
These curves have uniformly bounded length on the chart, so
\[
 \int_{r-2=\pm\varepsilon}|(V\cdot n)D_xV|\,dS
 \le C_I\varepsilon(1+|\log\varepsilon|)^2
 \longrightarrow0.
\]
Here $n$ is the physical unit normal.

Near a saddle, integrate over
$\{\varepsilon<|a|<a_0,\ \varepsilon<|b|<a_0\}$.
Using
\[
 |V^a/a|+|V^b/b|+|D_xV|
 \le C_I(1+|\log|a||+|\log|b||)^2
\]
and the bounded coordinate Jacobian, we obtain
\[
\begin{aligned}
 &\int_{\substack{a=\pm\varepsilon\\
                  \varepsilon<|b|<a_0}}
       |(V\cdot n)D_xV|\,dS\\
 &\qquad\le C_I\varepsilon
       \int_{\varepsilon<|b|<a_0}
       (1+|\log\varepsilon|+|\log|b||)^4\,db\\
 &\qquad\le C_I\varepsilon(1+|\log\varepsilon|)^4
       \longrightarrow0.
\end{aligned}
\]
The same estimate holds on $b=\pm\varepsilon$ after exchanging
$a$ and $b$.

Consequently, integration by parts and $\varepsilon\to0+$
imply, componentwise, for every smooth test function $\varphi$
supported away from the centers,
\[
 \int (V\cdot\nabla D_xV)\varphi\,dx\,dy
 =-\int D_xV\bigl(V\cdot\nabla\varphi
                  +(\operatorname{div}V)\varphi\bigr)\,dx\,dy.
\]
Thus the classical expression for $V\cdot\nabla D_xV$
off the separatrix equals its distributional expression.
The corresponding identification for $D_xV$ follows from
Lemma~\ref{geom:common-field}.

In the region where $\eta=1$, set $\tau=S$ and recall
$w(\tau)=\tau/|S'(r(\tau))|$ from \eqref{eq:weight}.
Since $D_\lambda\tau=-\tau/\lambda$,
\[
 \zeta=-\frac1\lambda\left(1+\frac{\tau w'}w\right),\qquad
 D_\lambda\zeta=\frac1{\lambda^2}
 \left[1+\frac{\tau w'}w+
       \tau\partial_\tau\left(\frac{\tau w'}w\right)\right].
\]
For a polynomially growing weight $P$, the integration by parts
in Lemma~\ref{aux:logregularity} represents $PTu$ by the flux
\[
 -\frac{\tau_n}{c}\eta_n\log(\tau-\frac{\tau_n}{c})P(\tau)
                      u_ne^{in\theta}X_\tau
\]
and scalar terms containing $P$ and $P'$. The coefficients are bounded by a convergent sum
$C\sum_{n:\,\tau_n>0}(1+\tau_n)^Ne^{-\tau_n/C}$, uniformly on
compact $I\subset\mathbb R\setminus\{0\}$.
Consequently, $\zeta Tu$ and $(D_\lambda\zeta+\zeta^2)Tu$
satisfy the weighted
$\operatorname{div}L^p+L^{2p/(p+2)}$ bounds in
Lemma~\ref{aux:logregularity}. The matching traces across the
separatrix exclude an additional interface distribution.

\subsubsection{Preliminary bounds for the second adapted derivative}
\label{subsec:second-preliminary}
We first estimate the second adapted derivative of the spectral density away from the
integer thresholds; these bounds also provide control of the spectral density in the exterior of the
center disks near the integer thresholds.

\begin{lemma}\label{quad:second-away}
 Let $I\Subset\mathbb R\setminus\mathbb Z$, and choose $q_0>0$
so that
\[
 nb(q)-\lambda\ne0
 \qquad
 (n\ne0,\quad q\le2q_0,\quad \lambda\in I).
\]
Let $\eta$ be a smooth radial cutoff equal to zero for
$q\le q_0$ and equal to one outside the center disks
$q<2q_0$. Recall that
\[
 V=-\frac{S(r)}{\lambda S'(r)}\partial_r,
 \qquad
 D_\lambda=\partial_\lambda+\eta V
 =\partial_\lambda
   -\eta\,\frac{S(r)}{\lambda S'(r)}\partial_r.
\]
Thus $D_\lambda=\partial_\lambda$ on $q\le q_0$, where
the denominators are separated from zero, and
$D_\lambda=\partial_\lambda+V$ outside $q<2q_0$. Then
\[
 \sup_{\lambda\in I}
   (\|D_\lambda g\|_{H^1}+\|D_\lambda^2g\|_{H^1})\le C_I M.
\]
Away from these disks, $g,D_\lambda g\in W^{1,p}$ locally for
every finite $p$, with bounds $C_{I,p}M$.
\end{lemma}
\begin{proof}
\emph{Step 1: Local regularity and the differentiated estimates.}
Fix $p>2$ and open sets $U,U'$ whose closures avoid the disks
$q\le q_0$, with $\overline U\subset U'$.
For $\Delta u=\operatorname{div}F+r$, the local elliptic estimate is
\[
 \|u\|_{W^{1,p}(U)}\le C\left(\|u\|_{H^1(U')}
       +\|F\|_{L^p(U')}+\|r\|_{L^{2p/(p+2)}(U')}\right).
\]
It follows from Lemma~\ref{aux:logregularity},
$(\Delta+T)g=H$, and this estimate (on a slightly larger open set
containing $\overline{U'}$) that $\|g\|_{W^{1,2p}(U')}\le C_{I,p}M$.
Then
\[
 \|\eta VH+K_\eta\nabla g\|_{L^p(U')}
 \le \|\eta V\|_\infty\|H\|_p
     +\|K_\eta\|_{2p}\|\nabla g\|_{2p}\le C_{I,p}M.
\]
It follows from \eqref{spectral:log-regularity},
with $\rho=\zeta$, that $\zeta Tg$ is bounded in
$\operatorname{div}L^p+L^{2p/(p+2)}$. Near the elliptic centers
$E_\eta$ is smooth.
It follows from \eqref{eq:fixed-source-derivative},
Lemma~\ref{r2:full-source-LAP}, and the local elliptic estimate that
\[
 \|D_\lambda g\|_{H^1}
       +\|D_\lambda g\|_{W^{1,p}(U)}\le C_{I,p}M.
\]
In particular, taking $p=4$ and using \eqref{quad:saddle-coefficients}, we get
\[
 \begin{split}
 &\|2K_\eta\nabla D_\lambda g\|_2
 +\|[D_\lambda K_\eta-K_\eta D_x(\eta V)^t
       -D_x(\eta V)K_\eta+\zeta K_\eta]\nabla g\|_2\\
 &\qquad\le C_I\bigl(\|\nabla D_\lambda g\|_4
                                  +\|\nabla g\|_4\bigr)\le C_I M.
 \end{split}
\]
The vector field in \eqref{quad:source-flux} satisfies
\[
\begin{aligned}
 &\|\eta V(\eta V\cdot\nabla H)
       +[\partial_\lambda(\eta V)+\zeta\eta V]H\|_2\\
 &\qquad\le C_I\bigl(\|\nabla H\|_2+\|H\|_2+\|H\|_4\bigr)
 \le C_I\|H\|_{H^1},
\end{aligned}
\]
because $\eta V,\partial_\lambda(\eta V)\in L^\infty$,
$\zeta\in L^4$, and $H^1\hookrightarrow L^4$.
For the scalar terms in \eqref{quad:second-weak-equation},
Lemma~\ref{aux:logregularity}\textup{(ii)} with
$\rho=\zeta$ and $\rho=D_\lambda\zeta+\zeta^2$, implies
\[
 \|\zeta TD_\lambda g\|_{H^{-1}}
 +\|(D_\lambda\zeta+\zeta^2)Tg\|_{H^{-1}}
 \le C_I(\|D_\lambda g\|_{H^1}+\|g\|_{H^1})
 \le C_I M.
\]
Moreover,
\[
 E_\eta=(1-\eta)\partial_\lambda T,\qquad
 D_\lambda E_\eta
 =(1-\eta)^2\partial_\lambda^2T
 -\eta\left(V\eta+\frac{1-\eta}{\lambda}\right)
                                  \partial_\lambda T.
\]
These operators are supported in $q\le2q_0$, where the
denominators are uniformly separated from zero. Hence
\[
 \|E_\eta D_\lambda g\|_2
 +\|(D_\lambda E_\eta+2\zeta E_\eta)g\|_2
 \le C_I(\|D_\lambda g\|_2+\|g\|_2)
 \le C_I M.
\]
Substitution in \eqref{quad:second-weak-equation}, together
with the gradient-flux estimate, therefore implies
\[
 \|(\Delta+T)D_\lambda^2g\|_{H^{-1}}\le C_I M.
\]
Since $T1=0$, Lemma~\ref{r2:full-source-LAP} applies to the
mean-zero part and implies the a priori estimate
\[
 \|D_\lambda^2g-\operatorname{mean}(D_\lambda^2g)\|_{H^1}
 \le C_I M.
\]
Step~2 bounds the mean using \eqref{quad:restore-mean}
and justifies differentiation and passing to the boundary values.

\smallskip\noindent
\emph{Step 2: Spatial means and boundary differentiation.}
Fix $\operatorname{Re}c=1$, $c\ne1$, and set
$g_c=g_f(ic\lambda)$, with $c$ held fixed under differentiation.
The formula
\[
 g_c=f+ic\lambda\,\Delta^{-1}(ic\lambda-L)^{-1}f
\]
and holomorphy of the resolvent imply
$\partial_\lambda^jg_c\in H^2$ for $0\le j\le2$.
Since $\eta V\in L^\infty\cap W^{1,4}$,
$\partial_\lambda(\eta V)=-\eta V/\lambda$, and
$H^2\hookrightarrow W^{1,4}$, the product rule implies
\[
 D_\lambda g_c,\ \partial_\lambda D_\lambda g_c\in H^1,
 \qquad
 D_\lambda^2g_c
 =\partial_\lambda D_\lambda g_c
       +\eta V\cdot\nabla D_\lambda g_c\in L^2.
\]
The distributional identities
\eqref{quad:weak-divergence}--\eqref{quad:second-weak-equation},
with the coefficient bounds \eqref{quad:saddle-coefficients},
therefore apply. By Step~1 their right-hand side belongs to
$H^{-1}$. Since $T(ic\lambda)$ is bounded on $L^2$,
\[
 \Delta D_\lambda^2g_c
 =(\Delta+T(ic\lambda))D_\lambda^2g_c
       -T(ic\lambda)D_\lambda^2g_c\in H^{-1}.
\]
Elliptic regularity implies $D_\lambda^2g_c\in H^1$.
These statements do not require bounds uniform in $c$.

We now complete the uniform estimates by controlling the means.
Since $\int_{\mathbb T^2}g_c=0$ and
$\int\eta V\cdot\nabla u=-\int\zeta u$, differentiation and
integration by parts imply
\begin{equation}\label{quad:restore-mean}
\begin{aligned}
 \int D_\lambda g_c
 &=-\int\zeta g_c,\\
 \int D_\lambda^2g_c
 &=-\int\{2\zeta D_\lambda g_c
                 +(D_\lambda\zeta+\zeta^2)g_c\}.
\end{aligned}
\end{equation}
In particular,
\[
\begin{aligned}
 \left|\int D_\lambda g_c\right|
 &\le \|\zeta\|_2\|g_c\|_2,\\
 \left|\int D_\lambda^2g_c\right|
 &\le 2\|\zeta\|_2\|D_\lambda g_c\|_2
       +\|D_\lambda\zeta+\zeta^2\|_2\|g_c\|_2.
\end{aligned}
\]
The coefficient norms are bounded uniformly for $\lambda\in I$.
Applying Lemma~\ref{r2:full-source-LAP} to the mean-zero
parts, first for $D_\lambda g_c$ and then for
$D_\lambda^2g_c$, completes the estimates
\[
 \sup_{\lambda\in I}
 \bigl(\|D_\lambda g_c\|_{H^1}
            +\|D_\lambda^2g_c\|_{H^1}\bigr)\le C_I M,
\]
uniformly as $c\to1$ from either boundary side.

By Proposition~\ref{global:LAP},
$g_c\to g_{f,\pm}$ in $L^2(I;H^1)$.
After passing to a subsequence, write
\[
 D_\lambda^jg_c\rightharpoonup^\ast h_j
       \quad\text{in }L^\infty(I;H^1),\qquad j=1,2.
\]
For $\varphi\in C_c^\infty(I^\circ\times\mathbb T^2)$,
integration by parts implies
\[
 \int_I\!\int D_\lambda^jg_c\,\varphi
 =-\int_I\!\int D_\lambda^{j-1}g_c
       \bigl(\partial_\lambda\varphi
              +\eta V\cdot\nabla\varphi+\zeta\varphi\bigr),
 \qquad j=1,2.
\]
The expression in parentheses belongs to $L^1(I;L^2)$.
Passing to the limit successively for $j=1,2$ identifies
$h_1=D_\lambda g_{f,\pm}$ and
$h_2=D_\lambda^2g_{f,\pm}$ in distributions. Consequently,
\[
 \operatorname*{ess\,sup}_{\lambda\in I}
 \bigl(\|D_\lambda g_{f,\pm}\|_{H^1}
       +\|D_\lambda^2g_{f,\pm}\|_{H^1}\bigr)
 \le C_I M.
\]

\end{proof}

To estimate spectral derivatives uniformly near the nonzero integers,
we first record the bound
\begin{equation}\label{quad:uniform-full-inverse}
 \|u\|_{H^1}\le C\left(
 \|(\Delta+T(\epsilon+i\lambda))u\|_{H^{-1}}
       +\left|\int_{\mathbb T^2}u\right|\right),
 \qquad |\lambda|\ge\tfrac12,\quad 0\le|\epsilon|\le1.
\end{equation}
This holds for every $u\in H^1(\mathbb T^2)$; at $\epsilon=0$
either boundary value of $T$ is allowed, including at integer thresholds.
Indeed, $T$ annihilates constants, so it suffices to consider
$u-|\mathbb T^2|^{-1}\int u$. For $|\lambda|\ge\Lambda_0$,
with $\Lambda_0$ sufficiently large, Corollary~\ref{cor:global-high-frequency}
implies
\[
 (\Delta+T)^{-1}=(I+K)^{-1}\Delta^{-1},\qquad
 \|(I+K)^{-1}\|_{H^1\to H^1}\le2
\]
on the spaces of zero spatial mean. On
$\tfrac12\le|\lambda|\le\Lambda_0$, the required bound is
\eqref{r2:full-source-estimate} in Lemma~\ref{r2:full-source-LAP}.
Adding back the mean proves \eqref{quad:uniform-full-inverse}.

The next lemma bounds two adapted derivatives near $\lambda=k$,
allowing the factor $(k\delta^2)^{-1}$ for each derivative.
It also proves uniform bounds for differentiation along
$\lambda=k(1\mp\delta^2/8)$ on explicitly specified annuli.
These bounds will be useful in the weighted estimates near the centers below.

\begin{lemma}\label{quad:coarse-second}
Let $f\in H^5\cap Y$ satisfy $\langle f,g_j\rangle=0$ for
$1\le j\le3$, and set $M=\|f\|_{H^5}$.
Fix $\eta\in C^\infty([0,\infty))$ with $0\le\eta\le1$,
$\eta=0$ on $[0,1]$ and $\eta=1$ on $[2,\infty)$, and fix
\[
 0<c_0\le e^{-5}\min\{q_*/2,1/4\},\qquad
 0<c_1\le\tfrac12e^{-6}c_0.
\]
For an integer $k\ge1$, $0<\delta<c_1k^{-1/2}$, and either choice
\[
 \lambda_*=k(1-\delta^2/8)\quad\text{or}\quad
 \lambda_*=k(1+\delta^2/8),
\]
consider $|\lambda-\lambda_*|\le k\delta^2/64$.
Define $\chi_\delta(q)=\eta(8q/\delta)$. Set
\[
 D_\delta=\partial_\lambda+\chi_\delta V,
 \qquad V=-\frac{S(r)}{\lambda S'(r)}\partial_r,
 \qquad g=g_{f,\pm}(\lambda)=f+i\lambda\psi_{f,\pm}(\lambda).
\]
The parameter $\delta$ is held fixed when $D_\delta$ is applied.
Then, for either boundary value and every $\lambda$ in the
above,
\begin{equation}\label{quad:coarse-bound}
 \sum_{\ell=0}^2\sum_{j=0}^2
 (k\delta^2)^j\|L_0^\ell D_\delta^jg\|_{H^1(\mathbb T^2)}
 \le CM.
\end{equation}

Set $q_k=c_0k^{-1/2}$. In either chart centered at $(0,0)$ or
$(\pi,\pi)$, use $v=\log(q/\delta)$ and $w=\log\delta$.
Along the curve $\lambda=k(1-\delta^2/8)$, set
\[
 D_w=q\partial_q-\frac{k\delta^2}{4}\partial_\lambda.
\]
Thus $D_w=\partial_w$ with $v$ fixed. On the annulus
$e^{-4}q_k<q<e^4q_k$,
\begin{equation}\label{h5c:collar-input}
 \sum_{\ell=0}^2\sum_{j=0}^2
 \|\nabla_{v,\theta}\partial_\theta^\ell D_w^jg\|_
 {L^2(dv\,d\theta;\,e^{-4}q_k<q<e^4q_k)}\le CM.
\end{equation}
Along the curve $\lambda=k(1+\delta^2/8)$, the same bound
holds with $D_w=q\partial_q+k\delta^2\partial_\lambda/4$.
The constant may depend on $c_0,c_1,\eta$ and the fixed coordinate charts,
but is independent of $k,\delta$, the boundary sign, and the
choice of center.
\end{lemma}

\begin{proof}
We first prove the estimates for spectral parameters off the
imaginary axis. Fix $c\in\mathbb C$ satisfying
\[
 \operatorname{Re}c=1,\qquad 0<|c-1|\le2\delta^2,
\]
and consider the family $\lambda\mapsto g_f(ic\lambda)$ on the
spectral neighborhood specified in the lemma. Since
\[
 ic\lambda=\epsilon(\lambda)+i\lambda,
 \qquad \epsilon(\lambda)=-\lambda\operatorname{Im}c,
\]
the condition $c\ne1$ ensures that $\epsilon(\lambda)\ne0$.
Throughout Steps~1--2, spectral derivatives are taken with
$c$ fixed; thus
\[
 \partial_\lambda g_f(ic\lambda)
       =ic\,\partial_z g_f(ic\lambda).
\]
This is the differentiation convention in
\eqref{eq:fixed-source-derivative}. In particular, since
$Vb=b/\lambda$, it preserves the identity
\[
 (\partial_\lambda+V)
 \left(\frac{nb}{nb-c\lambda}\right)=0.
\]

On the stated neighborhood, $1/2\le\lambda\le3k/2$, and hence
\[
 0<|\epsilon(\lambda)|
   \le2\lambda\delta^2
   \le3k\delta^2<1.
\]
Therefore \eqref{quad:uniform-full-inverse} applies.
We obtain estimates uniform in $c$ in this range.
Step~3 justifies passage to the boundary values as $c\to1$
and separately treats differentiation along the curves
$\lambda=k(1\mp\delta^2/8)$.

\smallskip\noindent
\emph{Step 1: The estimate \eqref{quad:coarse-bound}.}
The identities $b(q)=1-q^2/8$ and
$|\lambda-k|\in[7k\delta^2/64,9k\delta^2/64]$ imply
\[
 \begin{aligned}
 |kb(q)-c\lambda|&\ge |kb(q)-\lambda|
                    \ge\frac{13}{128}k\delta^2,\qquad q\le\delta/4,\\
 |nb(q)-c\lambda|&\ge\tfrac12|n-k|,
                    \qquad n\ne0,k,\quad q\le\delta/4.
 \end{aligned}
\]
Thus no critical circle lies in the support of $1-\chi_\delta$.
On the lower neighborhood, the $k$th critical circle satisfies
$7\delta^2/8\le q^2\le9\delta^2/8$; on the upper neighborhood
there is no $k$th critical circle.

On the fixed center disks, direct differentiation of the cutoff
field implies
\[
 \begin{aligned}
 |D_x(\chi_\delta V)|+|\operatorname{div}(\chi_\delta V)|
     &\le C(k\delta^2)^{-1},\\
 |D_\delta D_x(\chi_\delta V)|
 +|D_\delta\operatorname{div}(\chi_\delta V)|
 +|D_x(\chi_\delta V)|^2
     &\le C(k\delta^2)^{-2}.
 \end{aligned}
\]
The corresponding coefficient bounds hold after at most two
$L_0$ derivatives; here
$[L_0,\chi_\delta V]=-(\chi_\delta/\lambda)L_0$.
For the separated multipliers, the preceding denominator bounds imply
\[
 \sup_{n\ne0}\sup_{q<\delta/4}
 \frac{\delta^2}{n^2}
 \left|\frac{d^j}{d\lambda^j}
           \frac{nb(q)}{nb(q)-c\lambda}\right|
 \le C(k\delta^2)^{-j},\qquad j=1,2.
\]
Since
\[
 \|u_n\|_{L^2(q<\delta/4;\,dx\,dy)}
 \le C\frac{\delta}{|n|}\|\nabla u_n\|_{L^2(q<\delta/4;\,dx\,dy)},
\]
applying this inequality to the solution and test function bounds
the multiplier errors in $H^1\to H^{-1}$ by
$C(k\delta^2)^{-j}$. Derivatives of $\chi_\delta$ are supported
in $\delta/8<q<\delta/4$, where the same bounds apply.
For the scalar terms containing $T$, set
$\zeta=\operatorname{div}(\chi_\delta V)$. In the center coordinates,
\[
 \zeta=-\frac{4}{\lambda\xi_0(q)q}
          \partial_q\bigl(\chi_\delta(q)\xi_0(q)b(q)\bigr).
\]
Direct differentiation, with $\partial_v=q\partial_q$, implies
\[
 \|\zeta\|_{W^{1,\infty}_v}\le C(k\delta^2)^{-1},\qquad
 \|D_\delta\zeta+\zeta^2\|_{W^{1,\infty}_v}
       \le C(k\delta^2)^{-2}.
\]
These norms are taken on each fixed center disk. The coefficients
are radial, so multiplication by them preserves each cylinder
space $H^1_n$ with the same bounds. The scalar terms are
therefore controlled by \eqref{eq:cylinder-form}; $T$ annihilates
the zeroth mode.

On the complement of fixed smaller center disks, the first and
second coefficient $L^p$ norms are bounded by $C_p/k$ and
$C_p/k^2$, respectively, by \eqref{quad:saddle-coefficients}.
For the transport commutators, write $L_0=U\cdot\nabla$,
where $U=(-r_y,r_x)$ and $\operatorname{div}U=0$. The identities
\[
 [L_0,\operatorname{div}]F=-\operatorname{div}((D_xU)F),
 \qquad
 [\Delta,L_0]u=\operatorname{div}
       ((D_xU+(D_xU)^t)\nabla u)
\]
have smooth coefficients. The source flux in
\eqref{quad:source-flux} contains $H=(\Delta+1)f$ and $\nabla H$;
two transport commutations require at most
$\|L_0^2\nabla H\|_2\le C\|H\|_{H^3}\le CM$.

We now apply the differentiated equations. First,
\eqref{r2:angular-g}, Lemma~\ref{aux:logregularity}, and the local
elliptic estimate in Lemma~\ref{quad:second-away} imply
\[
 \sum_{\ell=0}^2\|L_0^\ell g\|_{W^{1,p}(q>q_*/2)}\le C_pkM
 \qquad(2<p<\infty).
\]
The right-hand sides of the local equations are estimated on
$q>q_*/4$. Using this estimate in
\eqref{eq:fixed-source-derivative}, with $\eta$ replaced by
$\chi_\delta$, and applying \eqref{quad:uniform-full-inverse}
proves the global $H^1$ bound for $D_\delta g$. Applying the same
local elliptic estimate to this first differentiated equation
then implies
\[
 \sum_{\ell=0}^2\|L_0^\ell D_\delta g\|_{W^{1,p}(q>q_*/2)}
       \le C_pk(k\delta^2)^{-1}M.
\]
The second exterior gradient flux is consequently bounded in
$L^2$ by
\[
 C\bigl(k^{-1}k(k\delta^2)^{-1}+k^{-2}k\bigr)M
       \le C(k\delta^2)^{-2}M.
\]
Apply \eqref{quad:second-weak-equation} with $\eta$ replaced
throughout by $\chi_\delta(q)=\eta(8q/\delta)$. Thus the derivative
is $D_\delta=\partial_\lambda+\chi_\delta V$, and its coefficients are
\[
\begin{aligned}
 \zeta_\delta&=\operatorname{div}(\chi_\delta V),\\
 K_\delta&=D_{x,y}(\chi_\delta V)
          +\bigl(D_{x,y}(\chi_\delta V)\bigr)^t
          -\zeta_\delta\,\mathrm{Id},\\
 E_\delta&=[D_\delta,T(ic\lambda)]
          =(1-\chi_\delta)\partial_\lambda T(ic\lambda).
\end{aligned}
\]
Here $c$ and $\delta$ are fixed when differentiating in $\lambda$.
In particular, $\partial_\lambda\chi_\delta=0$, and all spatial
derivatives of $\chi_\delta$ are supported in
$\delta/8<q<\delta/4$. The preceding estimates for the terms on
the right-hand side imply
\[
 \|(\Delta+T(ic\lambda))D_\delta^2g\|_{H^{-1}}
 \le C(k\delta^2)^{-2}M.
\]
Applying \eqref{quad:uniform-full-inverse} and using
\eqref{quad:restore-mean} with $D_\delta,\zeta_\delta$ yields
$\|D_\delta^2g\|_{H^1}\le C(k\delta^2)^{-2}M$.
At each of these two stages, commuting $L_0^\ell$ yields
\[
 \|L_0^\ell D_\delta^jg\|_{H^1}
 \le C(k\delta^2)^{-j}M
       +C\sum_{a<\ell}\|L_0^aD_\delta^jg\|_{H^1}.
\]
For $\ell\ge1$ the spatial mean is zero because $L_0$ is
divergence-free. For $\ell=0$, restore the means using
\eqref{quad:restore-mean} with
$\zeta=\operatorname{div}(\chi_\delta V)$.
This proves \eqref{quad:coarse-bound}.

\smallskip\noindent
\emph{Step 2: The estimate \eqref{h5c:collar-input}.}
The choices of $c_0,c_1$ ensure that the larger annulus
$e^{-5}q_k<q<e^5q_k$ lies in the center chart and
$\delta<\tfrac12e^{-6}q_k$. On this annulus,
\[
 |kb(q)-\lambda|\ge\frac{kq^2}{16},\qquad
 |nb(q)-\lambda|\ge\tfrac12|n-k|\quad(n\ne0,k).
\]
Consequently, for the derivative orders needed below,
\[
 \sup_{n\ne0}\left|\partial_v^m\left(
 q^2\frac{d^j}{d\lambda^j}
             \frac{nb(q)}{nb(q)-c\lambda}\right)\right|
 \le C_{m,j},\qquad 0\le m\le3,\quad 0\le j\le2.
\]
Set $u_j=(k\delta^2)^j\partial_\lambda^jg$, $0\le j\le2$;
the prefactor is fixed when these spectral derivatives are taken.
The cylinder equations (in $v, \theta$) are
\[
 \begin{aligned}
 \mathcal A_\lambda u_j
 &=-\mathbf1_{\{j=0\}}\xi_0q^2H
   +\xi_0q^2\sum_{a=1}^j\binom ja(k\delta^2)^a
                    (\partial_\lambda^aT)u_{j-a},\\
 \mathcal A_\lambda
 &=-\operatorname{div}_{v,\theta}(\mathsf G\nabla_{v,\theta})
                              -\xi_0q^2T(ic\lambda).
 \end{aligned}
\]
The multiplier bounds above and
$\|\xi_0q^2H\|_{H^3(dv\,d\theta)}\le CM$ are uniform
on the larger annulus.

Let $U,U'$ be full annuli with
$\overline U\subset U'\subset\{e^{-5}q_k<q<e^5q_k\}$,
whose boundary radii have fixed ratios to $q_k$.
Expanding $D_\delta$ in \eqref{quad:coarse-bound}, and then
applying interior elliptic regularity to the displayed equations,
we obtain
\[
 \begin{aligned}
 \|\nabla u_1\|_{L^2(U)}
 &\le C\bigl(M+\|\nabla g\|_{H^1(U')}\bigr),\\
 \|\nabla u_2\|_{L^2(U)}
 &\le C\bigl(M+\|\nabla g\|_{H^2(U')}
                         +\|\nabla u_1\|_{H^1(U')}\bigr),\\
 \|\nabla g\|_{H^4(U)}&\le CM,\\
 \|\nabla u_j\|_{H^{4-j}(U)}
 &\le C\left(\|\nabla u_j\|_{L^2(U')}
          +\sum_{a<j}\|\nabla u_a\|_{H^{2-j}(U')}\right),
                      \qquad j=1,2.
 \end{aligned}
\]
Here and below, $\nabla=\nabla_{v,\theta}$ and the norms on
annuli use $dv\,d\theta$. Subtract a scalar constant before
localization; $\mathcal A_\lambda$ and its spectral derivatives
annihilate constants. Applying these estimates successively on
annuli with boundary radii $e^{\pm(4+a/6)}q_k$, $0\le a\le6$,
implies
\[
 \|\nabla g\|_{H^4}
 +k\delta^2\|\nabla\partial_\lambda g\|_{H^3}
 +(k\delta^2)^2\|\nabla\partial_\lambda^2g\|_{H^2}\le CM
 \quad(e^{-4}q_k<q<e^4q_k).
\]
For the lower curve, direct differentiation at fixed $v$ implies
\[
 D_w^2g=(q\partial_q)^2g
 -\frac{k\delta^2}{2}q\partial_q\partial_\lambda g
 +\frac{k^2\delta^4}{16}\partial_\lambda^2g
 -\frac{k\delta^2}{2}\partial_\lambda g.
\]
The last two displays imply \eqref{h5c:collar-input} for a
fixed-$c$ family. For the upper curve, the two terms linear
in $k\delta^2$ have the opposite signs.

\smallskip\noindent
\emph{Step 3: Boundary values and differentiation along the curves.}
The difference-quotient argument in Lemma~\ref{quad:second-away}
applies with $q_0=\delta/8$, since Step~1 proves its
nonresonance condition on $q\le\delta/4$. All estimates are uniform
as $c\to1$ from either side. This proves
\eqref{quad:coarse-bound} for $g_{f,\pm}$.
For the annulus estimate along the lower curve, set
\[
 z=\epsilon+ik(1-\delta^2/8),\qquad
 \epsilon=\kappa k\delta^2,\quad 0<|\kappa|\le1,
 \qquad s=2\epsilon-ik\delta^2/4.
\]
Hold $\kappa$ and $v$ fixed when differentiating with $w$.
Write $\partial_\lambda^{\,c}$ for the derivative along
$z=ic\lambda$ with $c$ fixed. Since $g_f$ is holomorphic for $\operatorname{Re}z\ne0$,
differentiation along $z=ic\lambda$ with $c$ fixed satisfies
$\partial_\lambda^{\,c}g_f=ic\,\partial_zg_f$.
Along the path parametrized by $w$, with $\kappa$ and $v$ fixed,
we have
\[
 \partial_w q=q,\qquad \partial_w z=s,\qquad \partial_ws=2s.
\]
It therefore follows from the chain rule that
\begin{equation}\label{quad:absorption-paths}
 \begin{aligned}
 \partial_wg_f(z)&=\partial_vg_f(z)+s\partial_zg_f(z),\\
 \partial_w^2g_f(z)
 &=\partial_v^2g_f(z)+2s\partial_v\partial_zg_f(z)
                 +s^2\partial_z^2g_f(z)+2s\partial_zg_f(z),\\
 \partial_zg_f(z)&=\frac{\lambda}{z}\partial_\lambda^{\,c}g_f(z),
 \qquad
 \partial_z^2g_f(z)=\left(\frac{\lambda}{z}\right)^2
                         (\partial_\lambda^{\,c})^2g_f(z).
 \end{aligned}
\end{equation}
At each point of this path, $c=z/(i\lambda)$ satisfies
$|c-1|\le2\delta^2$. Moreover,
$|\lambda/z|\le1$, $|s|\le(9/4)k\delta^2$, and
$\partial_ws=2s$. It follows from Step~2 and
\eqref{quad:absorption-paths} that \eqref{h5c:collar-input}
holds uniformly for these regularized families. On each compact
$\delta$-subinterval of $(0,c_1k^{-1/2})$, let $\kappa\to0$ from either
side in distributions in $w=\log\delta$. The uniform Sobolev bounds
identify the first two $w$ derivatives of $g_{f,\pm}$ and preserve
the estimates by weak lower semicontinuity. For the upper curve,
use $z=\epsilon+ik(1+\delta^2/8)$ and
$s=2\epsilon+ik\delta^2/4$; the same calculation applies.
\end{proof}

The exponential bound below removes the endpoint at $\lambda=0$
in the spectral integrations by parts.
\begin{lemma}\label{quad:zero-second}
Set $q_0=q_{0,0}^{I}=\sqrt{\kappa/2}$. As in
\eqref{good-der}, use $D_\lambda=\partial_\lambda+\eta(q/q_0)V$,
where the radial cutoff is zero for $q\le q_0$ in each center
chart and equals one outside the two disks $q<2q_0$.
Since $b(q)=1-q^2/8$, for $q\le2q_0$ and $0<|\lambda|<\kappa$ we have
\[
 |nb(q)-\lambda|
 \ge |n|(1-\kappa/4)-\kappa
 \ge |n|/2,\qquad n\ne0.
\]
Then, for $0<|\lambda|<\kappa$ and $j=0,1,2$,
\[
 \|D_\lambda^j\Gamma_f(\lambda)\|_{H^1}
                    \le Ce^{-c/|\lambda|}M.
\]
\end{lemma}
\begin{proof}
\emph{Step 1: Bounds for the individual boundary values.}
Set $T_\Delta=T_+-T_-$, and recall
$g_{f,\pm}(\lambda)=f+i\lambda\psi_{f,\pm}(\lambda)$.
For $j=0,1,2$, set $J_j=D_\lambda^j(g_{f,+}-g_{f,-})$.
It follows from \eqref{tailD:baseline} and \eqref{tailD:Bq} that $\|g_{f,\pm}\|_{H^1\cap W^{1,p}}\le C_pM$ for finite $p>2$.
The first and second weak-equation coefficients have finite
$L^p$ norms $O(|\lambda|^{-1})$ and $O(|\lambda|^{-2})$.
Apply \eqref{zeroT:inverse} and then the local elliptic estimate
in Lemma~\ref{quad:second-away} to the first equation. This proves
$\|D_\lambda g_{f,\pm}\|_{H^1\cap W^{1,p}}\le C_p|\lambda|^{-3}M$.
It then follows from \eqref{quad:second-weak-equation} that
\[
 \|D_\lambda^2g_{f,\pm}\|_{H^1}
 \le C|\lambda|^{-2}\left(|\lambda|^{-2}M
    +|\lambda|^{-1}\|D_\lambda g_{f,\pm}\|_{W^{1,4}}\right)
 \le C|\lambda|^{-6}M.
\]
The means obey the same bounds by \eqref{quad:restore-mean}.

\smallskip\noindent
\emph{Step 2: Exponential decay of the jump.}
Let $F_{2,\pm}$ denote the right-hand side of
\eqref{quad:second-weak-equation}. Subtracting the two equations, we obtain
\[
 (\Delta+T_+)J_2=F_{2,+}-F_{2,-}-T_\Delta D_\lambda^2g_{f,-}.
\]
Since $E_{\eta,+}=E_{\eta,-}$ on the nonresonant cutoff disk,
\[
 \begin{split}
 \|F_{2,+}-F_{2,-}\|_{H^{-1}}
 \le C|\lambda|^{-N}\bigl(&\|J_0\|_{W^{1,4}}+\|J_1\|_{W^{1,4}}\\
 &+e^{-c/|\lambda|}
       (\|g_{f,-}\|_{H^1}+\|D_\lambda g_{f,-}\|_{H^1})\bigr)
 \end{split}
\]
for a fixed $N$. The coefficient fluxes use Hölder's inequality;
the factors containing $T_\Delta$ use \eqref{tailD:B}, including
its polynomial period weights. It follows from
\eqref{tailD:F-W14}, \eqref{tailD:zero-derivative}, and the equation
\eqref{tailD:jump-common}, with the local elliptic estimate and
\eqref{tailD:Bq}, that
$\|J_0\|_{W^{1,4}}+\|J_1\|_{W^{1,4}}\le Ce^{-c/|\lambda|}M$.
It also follows from \eqref{tailD:B} and the individual bound that
$\|T_\Delta D_\lambda^2g_{f,-}\|_{H^{-1}}
 \le C|\lambda|^{-6}e^{-c/|\lambda|}M$.
Apply \eqref{zeroT:inverse} to $J_2$ and use
\eqref{quad:restore-mean} for its mean. Decreasing $c$ absorbs
all polynomial factors. Finally,
\[
 J_0=2\pi i\lambda\Gamma_f,\qquad
 J_1=2\pi i(\lambda D_\lambda\Gamma_f+\Gamma_f),\qquad
 J_2=2\pi i(\lambda D_\lambda^2\Gamma_f+2D_\lambda\Gamma_f).
\]
It follows from these identities that the stated bounds hold for $j=0,1,2$.
\end{proof}

\subsubsection{Weighted estimates near an elliptic point}
The change of variables $v=\log(q/\delta)$ fixes the singular circle at
$v=0$. We first record the weighted inverse estimate used to differentiate
with respect to $w=\log\delta$ while keeping $v$ fixed.

\begin{lemma}\label{spectral:weighted-normal-family}
Let $k\ge k_0\ge1$, $|a|+|\eta|=1$, and
\[
 W(v)=e^{\mu_-v}(1+e^{2v})^{(\mu_+-\mu_-)/2},
 \qquad |\mu_-|<k_0,\qquad |\mu_+|<\sqrt{k_0^2+8}.
\]
The operator in Proposition~\ref{compact:normal} satisfies for any $u\in H^1_k$,
\[
 \|u\|_{H^1_k}\le C\|W^{-1}\mathcal N_{k,a,\eta}Wu\|_{H^{-1}_k}.
\]
The constant may depend on $k_0,\mu_-,\mu_+$, but is uniform in
$k,a,\eta$, including both boundary
values at $(a,\eta)=(1,0)$.
\end{lemma}
\begin{proof}
\emph{Step 1. The conjugated form.}
Set
\[
 \chi(v)=\frac{e^{2v}}{1+e^{2v}},\qquad
 p(v)=\frac{W'(v)}{W(v)}=(1-\chi(v))\mu_-+\chi(v)\mu_+,\qquad
 V_{a,\eta}(v)=\frac{8e^{2v}}{e^{2v}-a-i\eta}.
\]
Since $W^{-1}\partial_v W=\partial_v+p$, the differential expression is
\[
 W^{-1}\mathcal N_{k,a,\eta}Wu
 =-u''-2pu'+(k^2+V_{a,\eta}-p'-p^2)u .
\]
For $u,\varphi\in H^1_k(\mathbb R)$ its sesquilinear form is
\begin{equation}\label{spectral:wn-explicit-form}
\begin{aligned}
 \langle W^{-1}\mathcal N_{k,a,\eta}Wu,\varphi\rangle
  ={}&\int_{\mathbb R}
   \big[(u'+pu)(\overline{\varphi'}-p\bar\varphi)
                       +k^2u\bar\varphi\big]\,dv
      +\langle V_{a,\eta},u\bar\varphi\rangle .
\end{aligned}
\end{equation}
Here the last pairing is $\int_{\mathbb R}V_{a,\eta}u\bar\varphi\,dv$
when $(a,\eta)\ne(1,0)$. For the denominator $e^{2v}-1\mp i0$ it is
\[
 \lim_{\varepsilon\to0^+}
 \int_{|e^{2v}-1|>\varepsilon}
       \frac{8e^{2v}u(v)\overline{\varphi(v)}}{e^{2v}-1}\,dv
 \ \pm\ 4\pi i\,u(0)\overline{\varphi(0)} .
\]
Boundedness of the form follows from \eqref{eq:cylinder-form}.
For $\varphi=u$, the cross term $p(u\overline{u'}-u'\bar u)$ is purely
imaginary. Thus, after replacing $V_{a,\eta}$ by $8\chi$, the real part
of \eqref{spectral:wn-explicit-form} is
\[
 \|u'\|_2^2+\int_{\mathbb R}(k^2+8\chi-p^2)|u|^2\,dv .
\]
By convexity, $p^2\le(1-\chi)\mu_-^2+\chi\mu_+^2$. It follows that
\begin{equation}\label{spectral:wn-coercivity}
 k^2+8\chi-p^2
 \ge (1-\chi)(k^2-\mu_-^2)+\chi(k^2+8-\mu_+^2)
 \ge c k^2 ,
\end{equation}
where one may take
$c=\min\{1,1-\mu_-^2/k_0^2,1+(8-\mu_+^2)/k_0^2\}>0$.
The reference form with potential $8\chi$ is therefore coercive.
By \eqref{eq:model-log-identity} and the compactness argument in
Proposition~\ref{compact:normal}, multiplication by $V_{a,\eta}-8\chi$
is compact from $H^1_k$ to $H^{-1}_k$. Hence the operator in
\eqref{spectral:wn-explicit-form} is Fredholm of index zero.

\smallskip\noindent
\emph{Step 2. Injectivity.}
Suppose $u\in H^1_k(\mathbb R)$ and
$W^{-1}\mathcal N_{k,a,\eta}Wu=0$ in $H^{-1}_k(\mathbb R)$.
Then $U=Wu$ belongs to $H^1_{\mathrm{loc}}(\mathbb R)$ and satisfies
$\mathcal N_{k,a,\eta}U=0$ distributionally. To apply
Proposition~\ref{compact:normal}, we prove that $U\in H^1_k(\mathbb R)$.

The explicit potential satisfies $V_{a,\eta}=O(e^{2v})$ as $v\to-\infty$
and $V_{a,\eta}=8+O(e^{-2v})$ as $v\to+\infty$, uniformly in the
parameters. By \eqref{spectral:wn-coercivity}, there is $R>1$ such that
\begin{equation}\label{spectral:wn-tail}
 k^2+\operatorname{Re}V_{a,\eta}-p^2\ge ck^2/2
 \qquad (|v|\ge R).
\end{equation}
Choose a smooth $\zeta$ equal to zero on $|v|\le R$ and to one on
$|v|\ge R+1$. For $N>\max_{|v|\le R+1}W(v)$, set $m_N=\min\{1,N/W\}$.
Then $|m_N'|\le |p|m_N$ almost everywhere, and $m_N=1$ on the support
of $\zeta'$. Moreover,
\[
 m_NU=\min\{W,N\}u,\qquad m_NU'=\min\{W,N\}(u'+pu)
\]
belong to $L^2(\mathbb R)$. Testing the equation for $U$ against
$(\zeta m_N)^2 U$ and taking real parts yields
\[
 \|(\zeta m_N U)'\|_2^2+
 \int_{\mathbb R}(k^2+\operatorname{Re}V_{a,\eta})|\zeta m_N U|^2\,dv
 =\int_{\mathbb R}|(\zeta m_N)'|^2|U|^2\,dv .
\]
The test is justified by first inserting a compact cutoff and then
using $m_NU,m_NU'\in L^2$. Since $\zeta'm_N'=0$ almost everywhere,
the right-hand side is at most
\[
 \int_{R<|v|<R+1}|\zeta'U|^2\,dv+
 \int_{\mathbb R}p^2|\zeta m_N U|^2\,dv .
\]
It follows from \eqref{spectral:wn-tail} that
\[
 \|(\zeta m_N U)'\|_2^2+\tfrac c2 k^2\|\zeta m_N U\|_2^2
 \le \int_{R<|v|<R+1}|\zeta'U|^2\,dv .
\]
The right-hand side is independent of $N$. On any fixed bounded
interval $m_N=1$ for all sufficiently large $N$. Exhausting $\mathbb R$
by such intervals proves $\zeta U\in H^1_k$, hence $U\in H^1_k$.
Proposition~\ref{compact:normal} now implies $U=0$, and therefore $u=0$.

\smallskip\noindent
\emph{Step 3. Uniformity in $k,a,\eta$.}
For fixed $k$, bijectivity follows from Steps~1 and~2.
On each compact parameter set $\{|a|+|\eta|=1,\ \pm\eta\ge0\}$,
with its corresponding boundary value at $(1,0)$, the operator is
continuous in norm by \eqref{eq:model-log-identity} and
Proposition~\ref{compact:normal}. Its inverse is consequently bounded
uniformly for $k$ in any bounded range of angular indices.
For large $k$, set $\varphi=u$ in \eqref{spectral:wn-explicit-form}.
It follows from \eqref{eq:cylinder-form} that
\[
 \operatorname{Re}\langle W^{-1}\mathcal N_{k,a,\eta}Wu,u\rangle
 \ge \bigl(1-k^{-2}\|p\|_\infty^2-Ck^{-1/2}\bigr)\|u\|_{H^1_k}^2 .
\]
The coefficient is at least $1/2$ for all sufficiently large $k$.
Duality proves the desired estimate for these $k$ and completes the proof.
\end{proof}

To estimate two spectral derivatives near $\lambda=k$, we need
weighted bounds for $g$, $D_wg$, and $D_w^2g$ near each elliptic
center. The next lemma first controls their spatial gradients,
with up to two angular derivatives, in the norm weighted by
$q^{-\beta}$ in \eqref{h5c:mild}. It then establishes the stronger
bounds \eqref{h5c:critical}--\eqref{h5c:trace} for the $k$th angular
Fourier coefficient.

For $\lambda=k(1-\delta^2/8)$, the critical circle is $q=\delta$,
or equivalently $v=0$. The bounds on $g_k(0,w)$ and its first two
$w$ derivatives control the coefficients of the logarithmic
correction and the singular terms produced by its spectral
differentiation in Lemma~\ref{quad:critical-subtraction}.
\begin{lemma}\label{quad:weighted-center}
Let $k\ge1$, $q_k=c_0k^{-1/2}$ and $0<\delta<c_1k^{-1/2}$, where
$\sqrt\kappa\ll c_1<e^{-6}c_0\ll1$ and $e^5c_0<q_*$. Recall that $M=\|f\|_{H^5(\mathbb{T}^2)}$. Set
\[
 \lambda=k(1-\delta^2/8),\qquad
 v=\log(q/\delta),\qquad w=\log\delta,\qquad
 D_w=q\partial_q-\frac{k\delta^2}{4}\partial_\lambda.
\]
In $(v,w)$ coordinates, $D_w=\partial_w$ with $v$ fixed. Let
$g=g_{f,\pm}(\lambda)$ and $H=(\Delta+1)f$. For $0<\beta<1$,
\begin{equation}\label{h5c:mild}
 \sum_{j=0}^2\sum_{a=0}^2
 \|q^{-\beta}\nabla_{v,\theta}\partial_\theta^aD_w^jg
                  \|_{L^2(dv\,d\theta;q<q_k)}
 \le C_\beta k^{\beta/2}M.
\end{equation}
Choose a smooth radial function $\chi_k$ equal to one for $q<q_k/2$
and zero for $q>q_k$, with $|\partial_v^m\chi_k|\le C_m$, and set
\[
 W(v)=e^{v/2}(1+e^{2v})^{3/4+\beta/2}.
\]
Writing $g_k(v,w)$ for the $k$th angular Fourier coefficient of $g$
in these coordinates, we have
\begin{align}
 \|W^{-1}\chi_k\partial_w^jg_k\|_{H^1_k(\mathbb R_v)}
 &\le C_\beta k^{-2+\beta/2}\delta^{2+\beta}M,
                 &&0\le j\le2,\label{h5c:critical}\\
 |\partial_w^jg_k(0,w)|
 &\le C_\beta k^{-5/2+\beta/2}\delta^{2+\beta}M,
                 &&0\le j\le2.\label{h5c:trace}
\end{align}
The same estimates hold for $\lambda=k(1+\delta^2/8)$ with
$D_w=q\partial_q+k\delta^2\partial_\lambda/4$.
\end{lemma}
\begin{proof}
For $\lambda=k(1-\delta^2/8)$ the critical potential is
$8e^{2v}/(e^{2v}-1\mp i0)$, independent of $w$. Thus the
change of variables has the same purpose as the renormalization
in~\cite{JiaVortex}.
All norms in $(v,\theta)$ below use $dv\,d\theta$.

\emph{Step 1: the source and the commuted equation.}
It follows from the zero orbit average of $f$ and
$X(q,\theta)=q(\cos\theta,-\sin\theta)+O(q^3)$ that
$f(0)=\Delta f(0)=0$: these are the constant and quadratic
coefficients of the angular average. Hence $H(0)=0$. For
$0<\beta<1$,
\begin{equation}\label{h5c:source}
 \sum_{m=0}^3\ \sum_{Z_1,\ldots,Z_m}
 \|q^{-\beta}Z_1\cdots Z_mH\|_{L^2(dv\,d\theta;q<q_*)}
 \le C_\beta\|H\|_{H^3}\le C_\beta M,
\end{equation}
where each $Z_i$ is $\partial_v$ or $\partial_\theta$.
Indeed, a product of $m\ge1$ cylinder derivatives is bounded by
$C\sum_{\ell=1}^m q^\ell|\nabla^\ell H|$. Since cylinder measure
is comparable to $q^{-2}\,dx\,dy$, the squared weights with respect to physical measure
are $q^{2\ell-2\beta-2}$. For $\ell\ge2$ these are bounded;
for $\ell=0,1$ use $H(0)=0$ and $\|\nabla H\|_\infty\le CM$.

Set $h_j=D_w^jg$, $0\le j\le2$. The exact
center equation is
\begin{equation}\label{h5c:full-center-equation}
 -\operatorname{div}_{v,\theta}(\mathsf G\nabla_{v,\theta}g)
 +\xi_0(q)b(q)\sum_{n\ne0}\frac{8q^2}{q^2-8(1-\lambda/n)\mp i0\,\operatorname{sgn}n}\Pi_ng
 =-\xi_0q^2H.
\end{equation}
Here $\Pi_n$ denotes the projection onto the $n$th angular Fourier
mode, with $v,w$ fixed:
\[
 (\Pi_ng)(v,\theta,w)
 :=\frac{e^{in\theta}}{2\pi}
   \int_0^{2\pi}g(v,\vartheta,w)e^{-in\vartheta}\,d\vartheta
 =g_n(v,w)e^{in\theta}.
\]
For $n\ne0$, set $V_n=8q^2/(q^2-8(1-\lambda/n)\mp i0\,\operatorname{sgn}n)$,
where the boundary sign agrees with that of $g$. Then
\[
 \begin{aligned}
& D_wV_k=0,\qquad
 |\partial_v^bD_w^jV_n|\le C_{b,j}\frac{|n|q^2}{|n-k|}
 \quad(n\ne0,k),\\
& |D_w^j\partial_{v,\theta}^{\gamma}(\mathsf G-I,\xi_0b-1)|\le C_{j,\gamma}q^2.
 \end{aligned}
\]
For the above estimates, $0\le b,j\le2$ and $|\gamma|\le2$ suffice.
Since $\xi_0b$ and $V_n$ are radial, the commuted equation is
\begin{align}\label{h5c:commuted-center-equation}
 &-\operatorname{div}(\mathsf G\nabla\partial_\theta^ah_j)
       +\xi_0(q)b(q)\sum_{n\ne0}V_n\Pi_n\partial_\theta^ah_j\notag\\
 &=-\partial_\theta^a\partial_v^j(\xi_0q^2H)
 +\sum_{\substack{0\le i\le j,\ 0\le b\le a\\(i,b)\ne(j,a)}}
 \binom{j}{i}\binom{a}{b}\operatorname{div}
 ((D_w^{j-i}\partial_\theta^{a-b}\mathsf G)
                      \nabla\partial_\theta^bh_i)\notag\\
 &\quad-\sum_{i<j}\binom{j}{i}\sum_{n\ne0}
             D_w^{j-i}(\xi_0bV_n)\Pi_n\partial_\theta^ah_i,
 \qquad j,a\le2.
\end{align}

\emph{Step 2: the bound \eqref{h5c:mild}.}
We first prove the unweighted estimate and then insert the weight
$q^{-\beta}$. This keeps the choice of $c_0$ independent of $\beta$.

\smallskip\noindent
\emph{(a) The angular average.}
For an angular function $u$, write
$[u]_0=(2\pi)^{-1}\int_0^{2\pi}u\,d\theta$.
Since the potentials are radial, averaging
\eqref{h5c:full-center-equation} we obtain
\[
 \partial_v[\mathsf G_{vv}g_v+\mathsf G_{v\theta}g_\theta]_0
       =[\xi_0q^2H]_0.
\]
The flux belongs to $L^2$ near $v=-\infty$ and
$[\xi_0q^2H]_0=O(q^3)$. Hence the integration constant is zero:
\begin{equation}\label{h5c:mean-flux}
 [\mathsf G_{vv}g_v+\mathsf G_{v\theta}g_\theta]_0(v,w)
   =\int_{-\infty}^v[\xi_0q^2H]_0(u,w)\,du.
\end{equation}
Since the source is independent of the spectral parameter,
$D_w^j(\xi_0q^2H)=\partial_v^j(\xi_0q^2H)$. It follows that
\[
\begin{aligned}
 \partial_v[h_j]_0
 ={}&\int_{-\infty}^v[\partial_v^j(\xi_0q^2H)]_0(u,w)\,du\\
 &-[(\mathsf G_{vv}-1)\partial_vh_j
                      +\mathsf G_{v\theta}\partial_\theta h_j]_0\\
 &-\sum_{i=0}^{j-1}\binom ji
   [(D_w^{j-i}\mathsf G_{vv})\partial_vh_i
       +(D_w^{j-i}\mathsf G_{v\theta})\partial_\theta h_i]_0.
\end{aligned}
\]
For $0\le\beta<1$ and $R\le3q_k$, Young's inequality implies
\[
 \left\|q^{-\beta}\int_{-\infty}^v q(u)^2F(u)\,du\right\|_{L^2(q<R)}
 \le\frac{R^{2-\beta}}2\|F\|_{L^2(q<R)};
\]
indeed, the expression inside the norm is
$q(v)^{2-\beta}\int_0^\infty e^{-2t}F(v-t)\,dt$.
It follows from \eqref{h5c:source} and the coefficient bounds
in Step~1 that
\[
 \|q^{-\beta}\partial_v[h_j]_0\|_{L^2(q<R)}
 \le CR^{2-\beta}
       \bigg(M+\sum_{i=0}^j\|\nabla h_i\|_{L^2(q<R)}\bigg).
\]
For $a\ge1$, $[\partial_\theta^ah_j]_0=0$.

\smallskip\noindent
\emph{(b) The nonzero angular modes.}
Choose smooth radial cutoffs $\chi_0,\chi_1$: $\chi_0=1$ for
$q\le q_k$ and $\chi_0=0$ for $q\ge2q_k$; $\chi_1=1$ for
$q\le2q_k$ and $\chi_1=0$ for $q\ge3q_k$.
We can arrange $|\partial_v^m\chi_\ell|\le C_m$ for
$\ell=0,1$ and $0\le m\le3$.
Set
\begin{equation}\label{h5c:step2b-reference-operator}
 \mathcal P=-\partial_v^2-\partial_\theta^2+V_k\Pi_k.
\end{equation}
For $0\le\beta<1$, Lemma~\ref{spectral:weighted-normal-family},
with $\mu_-=\mu_+=\beta$, applies to the $k$th angular mode.
For the regularized critical mode, the translation
$v\mapsto v-\tfrac12\log(1+8|\epsilon|/(k\delta^2))$
reduces the parameters to those of Lemma~\ref{spectral:weighted-normal-family}.
For $n\ne0,k$, we use
\begin{equation}\label{h5c:step2b-free-coercivity}
 \begin{aligned}
 q^{-\beta}(-\partial_v^2+n^2)q^\beta
   &=-(\partial_v+\beta)^2+n^2,\\
 \operatorname{Re}\left\langle
   \bigl[-(\partial_v+\beta)^2+n^2\bigr]\varphi,\varphi
   \right\rangle
   &=\|\varphi'\|_2^2+(n^2-\beta^2)\|\varphi\|_2^2\\
   &\ge(1-\beta^2)
          \bigl(\|\varphi'\|_2^2+n^2\|\varphi\|_2^2\bigr).
 \end{aligned}
\end{equation}
Summing over the nonzero modes, and using the
Poincar\'e inequality, we obtain, for $[U]_0=0$,
\begin{equation}\label{h5c:step2b-weighted-inverse}
 \|q^{-\beta}\nabla U\|_{L^2(\mathbb R_v\times\mathbb T_\theta)}
 \le C_\beta
 \|q^{-\beta}\mathcal PU\|_{H^{-1}(\mathbb R_v\times\mathbb T_\theta)}.
\end{equation}
Here the angular Poincar\'e inequality also implies the
equivalence of $\|q^{-\beta}\nabla U\|_2$ and
$\|q^{-\beta}U\|_{H^1}$.

Moving all perturbations in \eqref{h5c:commuted-center-equation}
to the right, we obtain
\begin{equation}\label{h5c:step2b-rearranged-equation}
\begin{aligned}
 \mathcal P\partial_\theta^ah_j
 ={}&-\partial_\theta^a\partial_v^j(\xi_0q^2H)\\
 &+\sum_{i=0}^j\sum_{b=0}^a\binom ji\binom ab
  \operatorname{div}\!\left(
   D_w^{j-i}\partial_\theta^{a-b}(\mathsf G-I)
                          \nabla\partial_\theta^bh_i\right)\\
 &-\sum_{i=0}^j\binom ji\bigg\{
   D_w^{j-i}(\xi_0b-1)V_k\Pi_k\partial_\theta^ah_i
   +\sum_{n\ne0,k}D_w^{j-i}(\xi_0bV_n)
                           \Pi_n\partial_\theta^ah_i\bigg\},
 \qquad 0\le j,a\le2.
\end{aligned}
\end{equation}
Since $\Pi_0U=[U]_0$ and $\mathcal P$ commutes with $\Pi_0$, set
\begin{equation}\label{h5c:step2b-localized-function}
 u=\chi_0(I-\Pi_0)\partial_\theta^ah_j,
 \qquad \Pi_0u=0.
\end{equation}
The product rule implies
\begin{equation}\label{h5c:step2b-localized-equation}
 \begin{aligned}
 \mathcal Pu
 ={}&\chi_0(I-\Pi_0)\mathcal P\partial_\theta^ah_j\\
 &-2\chi_0'(I-\Pi_0)\partial_v\partial_\theta^ah_j
   -\chi_0''(I-\Pi_0)\partial_\theta^ah_j,
 \end{aligned}
\end{equation}
where $\chi_0'=\partial_v\chi_0$ and
$\chi_0''=\partial_v^2\chi_0$.
It follows from \eqref{h5c:step2b-weighted-inverse} and
$\chi_0=1$ on $q<q_k$ that
\begin{equation}\label{h5c:step2b-localized-bound}
 \begin{aligned}
 &\|q^{-\beta}\nabla(I-\Pi_0)\partial_\theta^ah_j
                         \|_{L^2(q<q_k)}\\
 &\qquad\le
   \|q^{-\beta}\nabla u\|_{L^2(\mathbb R_v\times\mathbb T_\theta)}
 \le C_\beta
   \|q^{-\beta}\mathcal Pu\|_{H^{-1}(\mathbb R_v\times\mathbb T_\theta)}.
 \end{aligned}
\end{equation}
All dual norms below are taken against
$H^1(\mathbb R_v\times\mathbb T_\theta)$, with measure
$dv\,d\theta$. Since $q^{-\beta}\mathcal Pu$ has zero angular
average, its dual norm may be computed using test functions
with zero angular average.

Substitute \eqref{h5c:step2b-rearranged-equation} into the first
term of \eqref{h5c:step2b-localized-equation}. The fixed source is
$-\chi_0(I-\Pi_0)\partial_\theta^a\partial_v^j(\xi_0q^2H)$.
The perturbations are $\chi_0(I-\Pi_0)$ applied to the metric
and potential sums in \eqref{h5c:step2b-rearranged-equation}.
The last two terms of \eqref{h5c:step2b-localized-equation}
are the cutoff contributions.

\emph{The source.}
We claim that
\begin{equation}\label{h5c:step2b-source-bound}
 \bigl\|q^{-\beta}\chi_0(I-\Pi_0)
       \partial_\theta^a\partial_v^j(\xi_0q^2H)\bigr\|_{H^{-1}}
 \le C_\beta M,\qquad 0\le j,a\le2.
\end{equation}
For $j+a\le3$, this follows from \eqref{h5c:source}.
For $j=a=2$, integration by parts implies
\begin{equation}\label{h5c:step2b-source-ibp}
\begin{aligned}
 &\int q^{-\beta}\chi_0
       \partial_v^2\partial_\theta^2(\xi_0q^2H)\overline\varphi\\
 &\quad=-\int q^{-\beta}\partial_v\partial_\theta^2(\xi_0q^2H)
       \left(\chi_0\partial_v\overline\varphi
                  +(\chi_0'-\beta\chi_0)\overline\varphi\right).
\end{aligned}
\end{equation}
It follows from \eqref{h5c:source} and
\eqref{h5c:step2b-source-ibp} that
\eqref{h5c:step2b-source-bound} holds, using only three
derivatives of $H$. The source estimate is also valid at
$\beta=0$, by the proof of \eqref{h5c:source}.

\emph{The perturbations.}
For the matrix coefficients
$A=D_w^{j-i}\partial_\theta^{a-b}(\mathsf G-I)$ in
\eqref{h5c:step2b-rearranged-equation}, and the radial
coefficients $\rho=D_w^m(\xi_0b-1)$, Step~1 implies
\begin{equation}\label{h5c:step2b-coefficient-bound}
 |A|+|\rho|+|\partial_v\rho|\le Cq^2,
 \qquad q<3q_k.
\end{equation}
For any function $U$, integration by parts and
\eqref{eq:cylinder-form} imply
\begin{equation}\label{h5c:step2b-perturbation-bound}
\begin{aligned}
 &\|q^{-\beta}\chi_0\operatorname{div}(A\nabla U)\|_{H^{-1}}
   +\|q^{-\beta}\chi_0\rho V_k\Pi_kU\|_{H^{-1}}\\
 &\qquad\le Cq_k^{2-\beta}
                         \|\nabla U\|_{L^2(q<3q_k)}.
\end{aligned}
\end{equation}
For the second term of \eqref{h5c:step2b-perturbation-bound},
use \eqref{eq:cylinder-form} with the coefficient
$q^{-\beta}\chi_0\rho$ and the input $\chi_1\Pi_kU$.
The required coefficient bound follows from
\eqref{h5c:step2b-coefficient-bound}:
\begin{equation}\label{h5c:step2b-critical-multiplier-bound}
 \|q^{-\beta}\chi_0\rho\|_{W^{1,\infty}(\mathbb R_v)}
       \le Cq_k^{2-\beta}.
\end{equation}

For the noncritical potentials, Step~1 implies
\begin{equation}\label{h5c:step2b-noncritical-coefficient-bound}
 \sup_{n\ne0,k}\frac{|D_w^m(\xi_0bV_n)|}{n^2}
       \le Cq^2,\qquad 0\le m\le2.
\end{equation}
It follows from \eqref{h5c:step2b-noncritical-coefficient-bound}
and Parseval's identity that
\begin{equation}\label{h5c:step2b-noncritical-bound}
 \left\|q^{-\beta}\chi_0
       \sum_{n\ne0,k}D_w^m(\xi_0bV_n)\Pi_nU\right\|_{H^{-1}}
 \le Cq_k^{2-\beta}\|\nabla U\|_{L^2(q<2q_k)}.
\end{equation}
Apply \eqref{h5c:step2b-perturbation-bound} and
\eqref{h5c:step2b-noncritical-bound} to the functions
$\partial_\theta^bh_i$ and $\partial_\theta^ah_i$
occurring in \eqref{h5c:step2b-rearranged-equation}.
The projection $I-\Pi_0$ does not increase the dual norms.

\emph{The cutoff terms.}
The commutator identity is
\begin{equation}\label{h5c:step2b-cutoff-identity}
 [\mathcal P,\chi_0]U
       =-2\chi_0'\partial_vU-\chi_0''U.
\end{equation}
Its coefficients are supported in $q_k<q<2q_k$.
It follows from \eqref{h5c:step2b-cutoff-identity},
the angular Poincar\'e inequality, and
\eqref{h5c:collar-input} that
\begin{equation}\label{h5c:step2b-cutoff-bound}
 \bigl\|q^{-\beta}[\mathcal P,\chi_0]
       (I-\Pi_0)\partial_\theta^ah_j\bigr\|_{H^{-1}}
       \le Cq_k^{-\beta}M.
\end{equation}
The annulus estimate \eqref{h5c:collar-input} also controls
the gradient norms on $q_k<q<3q_k$ occurring in
\eqref{h5c:step2b-perturbation-bound} and
\eqref{h5c:step2b-noncritical-bound}.

\smallskip\noindent
\emph{(c) Summation and the weight.}
Combining (a) and (b), and summing over $0\le j,a\le2$, we obtain
\[
\begin{aligned}
 &\sum_{j,a=0}^2
   \|q^{-\beta}\nabla\partial_\theta^ah_j\|_{L^2(q<q_k)}\\
 &\qquad\le C_\beta q_k^{-\beta}M
   +C_\beta q_k^{2-\beta}\sum_{j,a=0}^2
                 \|\nabla\partial_\theta^ah_j\|_{L^2(q<q_k)}.
\end{aligned}
\]
At $\beta=0$, we choose the fixed constant $c_0$ sufficiently small such that
$C_0q_k^2\le C_0c_0^2<1/2$. Absorbing the last term proves
$\sum_{j,a=0}^2\|\nabla\partial_\theta^ah_j\|_{L^2(q<q_k)}\le CM$.
Plug this bound into the same inequality for $0<\beta<1$.
Since $q_k<1$ and
$q_k^{-\beta}=c_0^{-\beta}k^{\beta/2}$, we then obtain
\eqref{h5c:mild}.

\emph{Step 3: the bounds \eqref{h5c:critical} and
\eqref{h5c:trace}.}
For a scalar function $F$, write
\[
 [F]_k(v,w)=\frac1{2\pi}\int_0^{2\pi}
                 F(v,\theta,w)e^{-ik\theta}\,d\theta.
\]
Thus $h_{j,k}=[h_j]_k=\partial_w^jg_k$.
Throughout this step, $H^{-1}_k(\mathbb R_v)$ denotes the dual
of $H^1_k(\mathbb R_v)$, where
\[
 \|u\|_{H^1_k}^2
 =\|\partial_vu\|_{L^2(dv)}^2+k^2\|u\|_{L^2(dv)}^2.
\]
All derivatives in $v$ are taken with $w$ fixed.

\emph{(a) The equation for the localized Fourier coefficient.}
Recall that
\[
 V_k(v)=\frac{8e^{2v}}{e^{2v}-1\mp i0},
 \qquad
 \mathcal N_{k,1,\pm0}=-\partial_v^2+k^2+V_k.
\]
Take the $k$th Fourier coefficient of
\eqref{h5c:step2b-rearranged-equation} with $a=0$, and then
multiply by $\chi_k$. Since every potential coefficient is
radial, all terms indexed by $n\ne0,k$ vanish under this
projection. We obtain
\begin{equation}\label{h5c:step3-localized-equation}
\begin{aligned}
 \mathcal N_{k,1,\pm0}(\chi_kh_{j,k})
 ={}&-\chi_k[\partial_v^j(\xi_0q^2H)]_k\\
 &+\chi_k\sum_{i=0}^j\binom ji
   \bigl[\operatorname{div}
       \bigl(D_w^{j-i}(\mathsf G-I)\nabla h_i\bigr)\bigr]_k\\
 &-\sum_{i=0}^j\binom ji
       D_w^{j-i}(\xi_0b-1)V_k\chi_kh_{i,k}\\
 &-2\chi_k'\partial_vh_{j,k}-\chi_k''h_{j,k},
 \qquad 0\le j\le2,
\end{aligned}
\end{equation}
where $\chi_k'=\partial_v\chi_k$ and
$\chi_k''=\partial_v^2\chi_k$.
We estimate the four lines on the right of
\eqref{h5c:step3-localized-equation}, respectively, as the fixed
source, metric terms, critical potential terms, and cutoff terms.

\emph{(b) The fixed source and metric terms.}
Since $q=\delta e^v$, the weight satisfies
\begin{equation}\label{h5c:step3-weight}
 W^{-1}q^2\le\delta^{2+\beta}q^{-\beta},
 \qquad
 |(\log W)'|\le2+\beta.
\end{equation}
The identity $[F]_k=-k^{-2}[\partial_\theta^2F]_k$,
integration by parts, and \eqref{h5c:step3-weight} imply
\begin{equation}\label{h5c:step3-projected-dual}
\begin{aligned}
 &\|W^{-1}\chi_k[q^2F]_k\|_{H^{-1}_k}
 +\|W^{-1}\chi_k
          [\operatorname{div}(q^2\boldsymbol F)]_k\|_{H^{-1}_k}\\
 &\quad\le C_\beta k^{-2}\delta^{2+\beta}
 \left(
  \|q^{-\beta}\partial_\theta^2F\|_{L^2(dv\,d\theta;q<q_k)}
  +\|q^{-\beta}\partial_\theta^2\boldsymbol F
                         \|_{L^2(dv\,d\theta;q<q_k)}
 \right).
\end{aligned}
\end{equation}
Here $F$ is scalar and $\boldsymbol F$ is a cylinder vector field.
In the divergence term, the derivatives falling on
$W^{-1}\chi_k$ are controlled by
\eqref{h5c:step3-weight} and the bounds on $\chi_k'$.

For the fixed source with $j=0,1$, apply
\eqref{h5c:step3-projected-dual} with
$F=q^{-2}\partial_v^j(\xi_0q^2H)$ and use
\eqref{h5c:source}. For $j=2$, integrate once in $v$:
\begin{equation}\label{h5c:step3-source-ibp}
\begin{aligned}
 &\left\langle
 W^{-1}\chi_k[\partial_v^2(\xi_0q^2H)]_k,\varphi
 \right\rangle\\
 &\qquad=
 k^{-2}\int_{\mathbb R}
 [\partial_v\partial_\theta^2(\xi_0q^2H)]_k
       \partial_v(W^{-1}\chi_k\overline\varphi)\,dv.
\end{aligned}
\end{equation}
Only three derivatives of $H$ occur in
\eqref{h5c:step3-source-ibp}.

For each metric term, apply
\eqref{h5c:step3-projected-dual} with
\[
 \boldsymbol F
 =q^{-2}D_w^{j-i}(\mathsf G-I)\nabla h_i.
\]
It follows from \eqref{h5c:step2b-coefficient-bound} that
\[
 \bigl|\partial_\theta^2
       \bigl(q^{-2}D_w^{j-i}(\mathsf G-I)\nabla h_i\bigr)\bigr|
 \le C\sum_{a=0}^2|\nabla\partial_\theta^ah_i|.
\]
Consequently, \eqref{h5c:source}, \eqref{h5c:mild},
\eqref{h5c:step3-projected-dual}, and
\eqref{h5c:step3-source-ibp} imply
\begin{equation}\label{h5c:step3-source-metric-bound}
\begin{aligned}
 &\|W^{-1}\chi_k[\partial_v^j(\xi_0q^2H)]_k\|_{H^{-1}_k}\\
 &\quad+\sum_{i=0}^j\binom ji
 \left\|W^{-1}\chi_k
   \bigl[\operatorname{div}
       \bigl(D_w^{j-i}(\mathsf G-I)\nabla h_i\bigr)\bigr]_k
 \right\|_{H^{-1}_k}\\
 &\qquad\le C_\beta k^{-2+\beta/2}\delta^{2+\beta}M,
 \qquad 0\le j\le2.
\end{aligned}
\end{equation}

\emph{(c) The critical potential terms.}
Use the cutoff $\chi_0$ from Step~2(b), which equals one for
$q\le q_k$ and vanishes for $q\ge2q_k$. It follows from
\eqref{h5c:step2b-coefficient-bound} and
\eqref{h5c:step3-weight} that
\begin{equation}\label{h5c:step3-potential-multiplier}
 \left\|\chi_0q^\beta W^{-1}
                 D_w^m(\xi_0b-1)\right\|_{W^{1,\infty}(\mathbb R_v)}
 \le C_\beta\delta^{2+\beta},
 \qquad 0\le m\le2.
\end{equation}
Moreover, angular projection, the bounds on $\chi_k'$,
and \eqref{h5c:mild} with $a=2$ imply
\begin{equation}\label{h5c:step3-mode-energy}
\begin{aligned}
 \|q^{-\beta}\chi_kh_{i,k}\|_{H^1_k}
 &\le C_\beta k^{-2}
       \|q^{-\beta}\nabla\partial_\theta^2h_i
                          \|_{L^2(dv\,d\theta;q<q_k)}\\
 &\le C_\beta k^{-2+\beta/2}M,
 \qquad 0\le i\le2.
\end{aligned}
\end{equation}
Indeed, the terms containing derivatives of $q^{-\beta}\chi_k$
are controlled by $k|h_{i,k}|$ because $k\ge1$.

Since $\chi_0=1$ on the support of $\chi_k$, we have
\[
 W^{-1}\chi_kD_w^m(\xi_0b-1)V_kh_{i,k}
 =
 \bigl(\chi_0q^\beta W^{-1}D_w^m(\xi_0b-1)\bigr)
 V_k\bigl(q^{-\beta}\chi_kh_{i,k}\bigr).
\]
Applying \eqref{eq:cylinder-form} to this factorization, and
using \eqref{h5c:step3-potential-multiplier} and
\eqref{h5c:step3-mode-energy}, we obtain
\begin{equation}\label{h5c:step3-potential-bound}
 \|W^{-1}\chi_kD_w^m(\xi_0b-1)V_kh_{i,k}\|_{H^{-1}_k}
 \le C_\beta k^{-5/2+\beta/2}\delta^{2+\beta}M,
 \qquad 0\le m,i\le2.
\end{equation}
Thus every critical potential term in
\eqref{h5c:step3-localized-equation} is controlled by the
already proved estimate \eqref{h5c:mild}.

\emph{(d) The cutoff terms.}
The derivatives of $\chi_k$ are supported in
$q_k/2<q<q_k$. On this annulus,
\eqref{h5c:step3-weight} and \eqref{h5c:collar-input} imply
\begin{equation}\label{h5c:step3-annulus}
\begin{aligned}
 \sup_{q_k/2<q<q_k}W^{-1}
 &\le C_\beta\delta^{2+\beta}q_k^{-2-\beta},\\
 \|h_{j,k}\|_{L^2(dv;q_k/2<q<q_k)}
 &\le Ck^{-3}
   \|\nabla\partial_\theta^2h_j
                  \|_{L^2(dv\,d\theta;q_k/2<q<q_k)}
 \le Ck^{-3}M.
\end{aligned}
\end{equation}
We use the identity
\begin{equation}\label{h5c:step3-cutoff-identity}
 -2\chi_k'\partial_vh_{j,k}-\chi_k''h_{j,k}
 =-2\partial_v(\chi_k'h_{j,k})+\chi_k''h_{j,k}.
\end{equation}
Integration by parts in \eqref{h5c:step3-cutoff-identity},
followed by \eqref{h5c:step3-weight} and
\eqref{h5c:step3-annulus}, implies
\begin{equation}\label{h5c:step3-cutoff-bound}
\begin{aligned}
 &\|W^{-1}(-2\chi_k'\partial_vh_{j,k}
                         -\chi_k''h_{j,k})\|_{H^{-1}_k}\\
 &\qquad\le C_\beta\delta^{2+\beta}q_k^{-2-\beta}k^{-3}M
 \le C_\beta k^{-2+\beta/2}\delta^{2+\beta}M.
\end{aligned}
\end{equation}

\emph{(e) The weighted bound and the trace.}
First, \eqref{h5c:mild} with any fixed exponent in $(1/2,1)$
and the angular Poincar\'e inequality imply
$W^{-1}\chi_kh_{j,k}\in H^1_k(\mathbb R_v)$ for each fixed
$\delta>0$. Indeed, $W^{-1}$ is comparable to $e^{-v/2}$
as $v\to-\infty$, and $\chi_k$ vanishes for sufficiently
large $v$.

We now apply Lemma~\ref{spectral:weighted-normal-family} with
$k_0=1$, $\mu_-=1/2$, and $\mu_+=2+\beta<3$.
It follows from \eqref{h5c:step3-localized-equation},
\eqref{h5c:step3-source-metric-bound},
\eqref{h5c:step3-potential-bound}, and
\eqref{h5c:step3-cutoff-bound} that
\begin{equation}\label{h5c:step3-inverse-bound}
\begin{aligned}
 \|W^{-1}\chi_kh_{j,k}\|_{H^1_k}
 &\le C_\beta
   \|W^{-1}\mathcal N_{k,1,\pm0}
                       (\chi_kh_{j,k})\|_{H^{-1}_k}\\
 &\le C_\beta k^{-2+\beta/2}\delta^{2+\beta}M.
\end{aligned}
\end{equation}
Since $h_{j,k}=\partial_w^jg_k$, this proves
\eqref{h5c:critical}.

We use the inequality
\begin{equation}\label{h5c:step3-trace-inequality}
 |u(0)|\le Ck^{-1/2}\|u\|_{H^1_k(\mathbb R_v)}.
\end{equation}
At $v=0$, we have $q=\delta<q_k/2$, so $\chi_k=1$;
also $W(0)=2^{3/4+\beta/2}$.
Applying \eqref{h5c:step3-trace-inequality} to
$u=W^{-1}\chi_kh_{j,k}$ and using
\eqref{h5c:step3-inverse-bound}, we obtain
\[
 |\partial_w^jg_k(0,w)|
 \le C_\beta k^{-5/2+\beta/2}\delta^{2+\beta}M,
 \qquad 0\le j\le2,
\]
which proves \eqref{h5c:trace}.

\emph{Step 4: $\lambda>k$ and the resolvent limits.}
For $\lambda=k(1+\delta^2/8)$, use $D_w=q\partial_q+k\delta^2\partial_\lambda/4$.
The potential $8e^{2v}/(e^{2v}+1)$ corresponds to $(a,\eta)=(-1,0)$,
so \eqref{h5c:mild}--\eqref{h5c:trace} remain valid.
For $\epsilon\ne0$, set $\sigma=\epsilon/(k\delta^2)$ and fix $\sigma$. Then the critical potential is exactly
$8e^{2v}/(e^{2v}-1-8i\sigma)$, independent of $w$.
It follows from \eqref{quad:absorption-paths} that the estimates \eqref{h5c:collar-input} on $e^{-4}q_k<q<e^4q_k$ hold for these derivatives. For bounded $\sigma$, the translation
$v\mapsto v-\tfrac12\log(1+8|\sigma|)$ normalizes the parameters to
$|a|+|\eta|=1$; the weights before and after this bounded
translation are equivalent. The weighted estimates are therefore
uniform and pass to the boundary as $\sigma\to0$ on either side.

\end{proof}

\subsubsection{Subtraction at a critical circle}\label{subsec:critical-subtraction}

The bounds of Lemma~\ref{quad:weighted-center} control derivatives
with respect to $w$ at fixed $v$. Ordinary spectral differentiation
also differentiates $v=\log(q/\delta)$ and therefore produces
singular terms at $q=\delta$. We identify these terms from the
radial equation and subtract them before estimating the second
spectral derivative.

Fix $k\ge1$, $k-5\kappa/4<\lambda<k$, and $0<\beta<1$. Set
\[
 \sigma=8(1-\lambda/k)=\delta^2,\qquad
 v=\log(q/\delta),\qquad w=\log\delta,
 \qquad g=g_{f,\pm}(\lambda).
\]
Then $\delta^2<10\kappa/k<c_1^2/k$, so
\eqref{h5c:mild}--\eqref{h5c:trace} apply. The critical circle is
$v=0$, whereas the elliptic center corresponds to $v\to-\infty$.
The coefficients $\mathsf G_{ij}(v,\theta)$ below are evaluated at
$q=\delta e^v$, and $g_k(0,w)$ is the $k$th angular Fourier
coefficient of $g$ on $q=\delta$. All $w$ derivatives are taken
at fixed $v$, and $\nabla=\nabla_{v,\theta}$.

In \eqref{h5c:full-center-equation}, the critical potential satisfies
\[
 V_k(v)=\frac{8e^{2v}}{e^{2v}-1\mp i0},\qquad
 V_k(v)-\frac4{v\mp i0}\in C^\infty((-1,1)).
\]
The distributional identities
\begin{equation}\label{h5c:singular-primitives}
 \partial_v^2\{v\log(v\mp i0)\}=\frac1{v\mp i0},\qquad
 \partial_v^2\{v^2\log(v\mp i0)\}=2\log(v\mp i0)+3
\end{equation}
suggest the singular profiles $a_1v\log(v\mp i0)$ and
$a_2v^2\log(v\mp i0)$. The first cancels the pole in the equation
for the remainder; the second cancels the logarithm left by the
first. Accordingly, define
\begin{equation}\label{h5c:profile-coefficients}
 \begin{aligned}
 a_1(w,\theta)
 &=\frac{4\xi_0(\delta)b(\delta)g_k(0,w)e^{ik\theta}}
         {\mathsf G_{vv}(0,\theta)},\\
 2\mathsf G_{vv}(0,\theta)a_2(w,\theta)
 &=4\xi_0(\delta)b(\delta)\Pi_ka_1
   -\bigl(\partial_v\mathsf G_{vv}(0,\theta)
          +\partial_\theta\mathsf G_{v\theta}(0,\theta)\bigr)a_1\\
 &\quad-2\mathsf G_{v\theta}(0,\theta)\partial_\theta a_1.
 \end{aligned}
\end{equation}
Use the cutoff $\chi\in C_c^\infty((-1/2,1/2))$ from
Proposition~\ref{prop:spectral-bounds}, with $\chi=1$ on
$[-1/8,1/8]$, and set
\begin{equation}\label{h5c:regularized-density}
 \mathcal R=g-\chi(v)
       \{a_1v\log(v\mp i0)+a_2v^2\log(v\mp i0)\}.
\end{equation}
Step~1 below verifies the two cancellations. Both profiles vanish
at $v=0$, so $\mathcal R(0,\theta,w)=g(0,\theta,w)$; no vanishing
of $g$ on the critical circle is assumed.

At fixed $q$, the ordinary derivatives are
\begin{equation}\label{h5c:sigma-derivatives}
 \partial_\sigma=\frac{\partial_w-\partial_v}{2\sigma},\qquad
 \partial_\sigma^2=
 \frac{(\partial_w-\partial_v)^2-2(\partial_w-\partial_v)}{4\sigma^2}.
\end{equation}
Applying \eqref{h5c:sigma-derivatives} to the subtracted profiles
in \eqref{h5c:regularized-density} produces a Cauchy term and a
logarithm. Define their coefficients and their sum by
\begin{equation}\label{h5c:second-singular-profile}
 \begin{aligned}
 \mathcal C_2&=\frac{a_1}{4\sigma^2},\qquad
 \mathcal L_2=\frac{a_2+a_1-\partial_wa_1}{2\sigma^2},\\
 \mathcal S_2&=\chi(v)\left\{
       \frac{\mathcal C_2}{v\mp i0}
                +\mathcal L_2\log(v\mp i0)\right\}.
 \end{aligned}
\end{equation}
After subtracting $\mathcal S_2$, the remaining terms from
differentiating $g-\mathcal R$ are in
$H^1((-1,1)\times\mathbb T)$; derivatives of $\chi$ are supported
away from $v=0$. The next lemma bounds the mixed derivatives of
$\mathcal R$, the coefficients $\mathcal C_2,\mathcal L_2$, and the
localized source for $g_{\sigma\sigma}-\mathcal S_2$.
Since $q^2<e\sigma<32\kappa/k$ on the support of $\chi(v)$,
the cutoff adapted field vanishes on the support of $\mathcal S_2$.

\begin{lemma}\label{quad:critical-subtraction}
The remainder and the two coefficients satisfy
\begin{equation}\label{h5c:full-mixed}
 \sum_{j=0}^2\|\nabla \partial_w^j\mathcal R\|_{H^{2-j}}
 \le C\delta^\beta k^{\beta/2}M,
\end{equation}
\begin{equation}\label{h5c:critical-mixed}
 \|\mathcal R_k\|_{H^3_v}+\|(\partial_w\mathcal R)_k\|_{H^2_v}
                  +\|(\partial_w^2\mathcal R)_k\|_{H^1_v}
 \le Ck^{\beta/2}\delta^{2+\beta}M,
\end{equation}
\begin{equation}\label{h5c:second-coeff}
 \|\mathcal C_2\|_{H^2_\theta}+\|\mathcal L_2\|_{H^2_\theta}
 \le Ck^{-1/2+\beta/2}\delta^{-2+\beta}M,
\end{equation}
where the norms involving $v$ are on $(-1,1)$ with measure
$dv\,d\theta$. Let $\chi_k$ be the radial cutoff in
\eqref{h5c:critical}, equal to one for $q<q_k/2$ and zero for
$q>q_k$, where $q_k=c_0k^{-1/2}$. Then
\begin{equation}\label{h5c:final-local}
 \left\|\chi_k(\Delta+T_\pm(\lambda))
                 (g_{\sigma\sigma}-\mathcal S_2)\right\|_{H^{-1}(\mathbb T^2)}
 \le C\left(k^{\beta/2}\delta^{-2+\beta}+k^{3/2}\right)M.
\end{equation}
Here the localized source is extended by zero outside the selected
center chart; its definition uses only $q<q_k$. For
$\lambda=k(1+\delta^2/8)$, \eqref{h5c:final-local} holds with $\mathcal S_2=0$.
The analogous estimates at negative thresholds hold with $k$ replaced
by $|k|$ and the boundary signs conjugated.
\end{lemma}
\begin{proof}
\noindent\emph{Step 1. Cancellation of the pole and logarithm.}
It follows from \eqref{h5c:trace},
$\mathsf G_{v\theta}(0,\theta)=O(\delta^2)$, and $k\delta^2\le c_1^2$ that
\begin{equation}\label{h5c:jet-coeff}
 \|\partial_w^ja_i\|_{H^m_\theta}
 \le C_mk^{m-5/2+\beta/2}\delta^{2+\beta}M,
 \qquad i=1,2,\quad j\le2,\quad m\le3.
\end{equation}
The formula for $a_1$ implies this bound at every fixed angular
order; $a_2$ in $H^3_\theta$ uses $a_1$ in $H^4_\theta$ and
$k\delta^2\le c_1^2$.

For $j=0,1,2$, we use the one-dimensional Hardy estimates
\begin{equation}\label{h5c:hardy-quotients}
 \left\|\partial_v^j\frac{u(v)-\sum_{m=0}^j
                     u^{(m)}(0)v^m/m!}{v}\right\|_2
 \le C_j\|u^{(j+1)}\|_2.
\end{equation}
Here the norms are on $(-1,1)$ and $u\in H^{j+1}(-1,1)$.
The differentiated quotient equals
\[
 \int_0^1t^ju^{(j+1)}(tv)\,dt.
\]
Thus \eqref{h5c:hardy-quotients} follows from Minkowski's inequality
and $\int_0^1t^{j-1/2}\,dt<\infty$.
Set
\[
 \mathcal A_\sigma=-\operatorname{div}(\mathsf G\nabla)
                         +\xi_0b\sum_{n\ne0}V_n\Pi_n,\qquad
 \mathcal F=-\xi_0q^2H-\mathcal A_\sigma(g-\mathcal R)
                    -\frac{4\xi_0b\,g_k(0,w)e^{ik\theta}}{v\mp i0}.
\]
The equation for the remainder is exactly
\begin{equation}\label{h5c:regularized-core-equation}
 \begin{split}
 -\operatorname{div}(\mathsf G\nabla \mathcal R)
 +4\xi_0b\frac{\mathcal R_k(v)-\mathcal R_k(0,w)}v e^{ik\theta}
 +\xi_0b\left(V_k-\frac4{v\mp i0}\right)\mathcal R_ke^{ik\theta}
 +\xi_0b\sum_{n\ne0,k}V_n\mathcal R_ne^{in\theta}=\mathcal F.
 \end{split}
\end{equation}
The function $V_k-4/(v\mp i0)$ is smooth on $|v|<1$.
The pole coefficient of $\mathcal A_\sigma(g-\mathcal R)$ is $-\mathsf G_{vv}(0,\theta)a_1$;
its logarithmic coefficient is
\[
 \begin{aligned}
 &-2\mathsf G_{vv}(0,\theta)a_2+4\xi_0(\delta)b(\delta)\Pi_ka_1\\
 &\qquad-\bigl(\partial_v\mathsf G_{vv}(0,\theta)
              +\partial_\theta\mathsf G_{v\theta}(0,\theta)\bigr)a_1
       -2\mathsf G_{v\theta}(0,\theta)\partial_\theta a_1=0.
 \end{aligned}
\]
Thus $\mathcal F+\xi_0q^2H$ is a smooth-coefficient sum of
$1$, $v\log(v\mp i0)$ and $v^2\log(v\mp i0)$; coefficient
quotients such as $(\mathsf G_{vv}-\mathsf G_{vv}(0,\theta))/v$ have derivatives $O(\delta^2)$.
It follows from \eqref{h5c:source} and \eqref{h5c:jet-coeff} that
\begin{align}
 \|\mathcal F\|_{H^1}+\|\partial_w\mathcal F\|_2
 &\le C\delta^{2+\beta}M
   +C\sum_{i=1}^2\sum_{j=0}^1\|\partial_w^ja_i\|_{H^3_\theta}
 \le Ck^{\beta/2}\delta^\beta M,
 \label{h5c:regularized-core-source}\\
 \|\mathcal F_k\|_{H^1_v}+\|\partial_w\mathcal F_k\|_{L^2_v}
 &\le C\delta^{2+\beta}M
   +C\sum_{i=1}^2\sum_{j=0}^1\|\partial_w^ja_i\|_{H^2_\theta}
 \le Ck^{\beta/2}\delta^{2+\beta}M.
 \label{h5c:regularized-projected-source}
\end{align}
The projected estimate uses only $m\le2$; the full $H^1$
estimate uses $m=3$ and $k^{1/2}\delta^2\le c_1^2$.

It follows from \eqref{h5c:mild}, \eqref{h5c:critical}, and
the explicit subtractions that
\begin{equation}\label{h5c:jet-initial-energy}
 \begin{split}
 \sum_{j=0}^2\sum_{a=0}^2
       \|\nabla\partial_\theta^a\partial_w^j\mathcal R\|_2
       &\le Ck^{\beta/2}\delta^\beta M,\\
 \sum_{j=0}^2\|(\partial_w^j\mathcal R)_k\|_{H^1_v}
       &\le Ck^{-2+\beta/2}\delta^{2+\beta}M.
 \end{split}
\end{equation}
Both subtracted profiles belong to $H^1$, and the first line uses
at most three angular derivatives of $a_i$. The same estimates hold with $(-1,1)$ replaced by $(-2,2)$.

\medskip
\noindent\emph{Step 2. Mixed derivatives of the remainder.}
Solving \eqref{h5c:regularized-core-equation} for $\mathcal R_{vv}$ and
using the $j=0$ Hardy estimate, we obtain
\[
 \begin{split}
 \|\mathcal R_{vv}\|_2
 \le C\biggl(&\|\nabla\partial_\theta \mathcal R\|_2+\|\nabla \mathcal R\|_2
 +\|\mathcal R_k'\|_2+\|\mathcal F\|_2
 +\delta^2\|\mathcal R_{\theta\theta}\|_2\biggr)
 \le Ck^{\beta/2}\delta^\beta M.
 \end{split}
\]
Here $|V_n|\le C\delta^2n^2$ for $n\ne0,k$.
Applying $\partial_\theta$ to the same radial identity, we obtain
\[
 \|\mathcal R_{vv\theta}\|_2
 \le C\left(\|\nabla\partial_\theta^2\mathcal R\|_2
       +\|\nabla \mathcal R\|_{H^1}+k\|\mathcal R_k'\|_2
       +\|\partial_\theta\mathcal F\|_2\right)
 \le Ck^{\beta/2}\delta^\beta M.
\]
Thus all third derivatives containing an angular derivative
are bounded. Projecting the original radial identity onto $k$, we obtain
\begin{equation}\label{h5c:regularized-projected-equation}
 \mathcal R_k''=k^2\mathcal R_k-[\operatorname{div}((\mathsf G-I)\nabla \mathcal R)]_k
  +4\xi_0b\frac{\mathcal R_k(v)-\mathcal R_k(0,w)}v
  +\xi_0b\left(V_k-\frac4{v\mp i0}\right)\mathcal R_k-\mathcal F_k.
\end{equation}
It follows from \eqref{h5c:jet-initial-energy},
\eqref{h5c:regularized-projected-source}, and
$\|\mathsf G-I\|_{C^2}\le C\delta^2$ that
\[
 \|\mathcal R_k''\|_2\le
 C\left(k^2\|\mathcal R_k\|_2+\|\mathcal R_k\|_{H^1_v}
     +\delta^2\|\nabla \mathcal R\|_{H^1}+\|\mathcal F_k\|_2\right)
 \le Ck^{\beta/2}\delta^{2+\beta}M.
\]
Differentiate the radial equation in \eqref{h5c:regularized-core-equation}:
\[
 \begin{aligned}
 &\partial_v\{\mathsf G_{vv}\mathcal R_{vv}+2\mathsf G_{v\theta}\mathcal R_{v\theta}
  +\mathsf G_{\theta\theta}\mathcal R_{\theta\theta}
  +(\partial_v\mathsf G_{vv}+\partial_\theta\mathsf G_{v\theta})\mathcal R_v
  +(\partial_v\mathsf G_{v\theta}+\partial_\theta\mathsf G_{\theta\theta})\mathcal R_\theta\}\\
 &\quad=\partial_v\biggl\{4\xi_0b\frac{\mathcal R_k(v)-\mathcal R_k(0,w)}v e^{ik\theta}
  +\xi_0b\left(V_k-\frac4{v\mp i0}\right)\mathcal R_ke^{ik\theta}
  +\xi_0b\sum_{n\ne0,k}V_n\mathcal R_ne^{in\theta}-\mathcal F\biggr\}.
 \end{aligned}
\]
The right side is in $L^2$: use the $j=1$ Hardy estimate,
the bound for $\mathcal R_k''$, and $|\partial_v^bV_n|\le C\delta^2n^2$
for $b\le1$. The terms $\mathcal R_{v\theta\theta},\mathcal R_{\theta\theta\theta}$ are
bounded by \eqref{h5c:jet-initial-energy}, and the displayed
estimate bounds $\mathcal R_{vv\theta}$. Since $\mathsf G_{vv}\ge1/2$, it follows that
$\|\nabla \mathcal R\|_{H^2}\le Ck^{\beta/2}\delta^\beta M$.

At fixed $v$, $\partial_w$ commutes with the trace-subtracted
quotient. Its equation is
\begin{equation}\label{h5c:regularized-w-equation}
 \begin{aligned}
 &\mathcal A_\sigma \partial_w\mathcal R-\frac{4\xi_0b(\partial_w\mathcal R)_k(0,w)e^{ik\theta}}{v\mp i0}
 =\partial_w\mathcal F
 +\operatorname{div}((\partial_w\mathsf G)\nabla \mathcal R)\\
 &\qquad-4(\partial_w(\xi_0b))\frac{\mathcal R_k(v)-\mathcal R_k(0,w)}v e^{ik\theta}
 -(\partial_w(\xi_0b))\left(V_k-\frac4{v\mp i0}\right)\mathcal R_ke^{ik\theta}\\
 &\qquad-\sum_{n\ne0,k}\partial_w(\xi_0bV_n)\mathcal R_ne^{in\theta}.
 \end{aligned}
\end{equation}
Its right-hand side has $L^2$ norm at most
\[
 \|\partial_w\mathcal F\|_2+
 C\delta^2\left(\|\nabla \mathcal R\|_{H^1}+\|\mathcal R_k\|_{H^1_v}
                  +\|\mathcal R_{\theta\theta}\|_2\right)
 \le Ck^{\beta/2}\delta^\beta M.
\]
It follows from the radial identity, \eqref{h5c:jet-initial-energy},
and \eqref{h5c:hardy-quotients} that
$\|\nabla \partial_w\mathcal R\|_{H^1}\le Ck^{\beta/2}\delta^\beta M$. The energy of $\partial_w^2\mathcal R$ is in
\eqref{h5c:jet-initial-energy}. This proves \eqref{h5c:full-mixed}.
Finally, differentiating \eqref{h5c:regularized-projected-equation}
and projecting \eqref{h5c:regularized-w-equation}, we obtain
\[
 \begin{split}
 \|\mathcal R_k'''\|_2
 &\le C\left(k^2\|\mathcal R_k'\|_2+\|\mathcal R_k\|_{H^2_v}
       +\delta^2\|\nabla \mathcal R\|_{H^2}+\|\mathcal F_k\|_{H^1_v}\right),\\
 \|(\partial_w\mathcal R)_k''\|_2
 &\le C\left(k^2\|(\partial_w\mathcal R)_k\|_2+\|(\partial_w\mathcal R)_k\|_{H^1_v}
       +\delta^2\|\nabla \partial_w\mathcal R\|_{H^1}
       +\delta^2\|\nabla \mathcal R\|_{H^1}+\|\partial_w\mathcal F_k\|_2\right).
 \end{split}
\]
The bound \eqref{h5c:critical-mixed} therefore follows from
\eqref{h5c:full-mixed} and \eqref{h5c:regularized-projected-source}.

\medskip
\noindent\emph{Step 3. The corrected source on $|v|<1$.}
It follows from \eqref{h5c:singular-primitives} and
\eqref{h5c:sigma-derivatives} that
\[
 \partial_\sigma^2(g-\mathcal R)-\mathcal S_2
       \in H^1((-1,1)\times\mathbb T).
\]
The coefficient estimate \eqref{h5c:second-coeff}
follows from \eqref{h5c:profile-coefficients},
\eqref{h5c:second-singular-profile}, and \eqref{h5c:jet-coeff}.
Set $u=g_{\sigma\sigma}-\mathcal S_2$.
It follows from \eqref{h5c:full-mixed},
\eqref{h5c:sigma-derivatives}, and \eqref{h5c:jet-coeff} that,
on $|v|<1$,
$\|\nabla u\|_2\le Ck^{\beta/2}\delta^{-4+\beta}M$.
It follows from \eqref{h5c:critical} and \eqref{h5c:critical-mixed} that
\[
 \|\mathcal R_k'\|_2\le Ck^{-2+\beta/2}\delta^{2+\beta}M,
 \quad \|\mathcal R_k''\|_2\le Ck^{-1+\beta/2}\delta^{2+\beta}M,
 \quad \|u_k\|_{H^1_k}\le Ck^{\beta/2}\delta^{-2+\beta}M.
\]
The middle estimate follows by interpolation; the last also uses
\eqref{h5c:jet-coeff} and \eqref{h5c:sigma-derivatives}.
Since $\xi_0b$ is radial, the critical equation is exactly
\[
 (\mathcal A_\sigma u)_k
 =(-\partial_v^2+k^2+\xi_0bV_k)u_k
       -[\operatorname{div}((\mathsf G-I)\nabla u)]_k.
\]
It follows from \eqref{eq:cylinder-form} and
$\mathsf G-I=O(\delta^2)$ that its local dual norm is bounded by $C(\|u_k\|_{H^1_k}+\delta^2\|\nabla u\|_2)$.
We localize the source using $\chi\operatorname{div}F
=\operatorname{div}(\chi F)-\chi'F_v$.
For $n\ne0,k$, ordinary differentiation implies
\[
 (\mathcal A_\sigma u)_n
 =-2\xi_0b(\partial_\sigma V_n)g_{\sigma,n}
       -\xi_0b(\partial_\sigma^2V_n)g_n-(\mathcal A_\sigma \mathcal S_2)_n.
\]
The cubic and double pole coefficients of $\mathcal A_\sigma \mathcal S_2$ are
\[
 \begin{aligned}
 &-2\mathsf G_{vv}(0,\theta)\mathcal C_2=-2\xi_0(\delta)b(\delta)g_k(0,w)e^{ik\theta}/\sigma^2,\\
 &\mathsf G_{vv}(0,\theta)\mathcal L_2-\partial_v\mathsf G_{vv}(0,\theta)\mathcal C_2
       +2\mathsf G_{v\theta}(0,\theta)\partial_\theta \mathcal C_2\\
 &\qquad +(\partial_\theta\mathsf G_{v\theta}(0,\theta))\mathcal C_2
       +4\xi_0(\delta)b(\delta)\Pi_k\mathcal C_2\\
 &\qquad=\frac{8\xi_0(\delta)b(\delta)\Pi_ka_1+2\mathsf G_{vv}(0,\theta)a_1-2\partial_w(\mathsf G_{vv}(0,\theta)a_1)}{4\sigma^2}.
 \end{aligned}
\]
Both coefficients have only angular mode $k$. This follows from the
identity for $a_2$ in \eqref{h5c:profile-coefficients},
$\partial_w\mathsf G_{vv}(0,\theta)=\partial_v\mathsf G_{vv}(0,\theta)$, and $\mathsf G_{vv}(0,\theta)a_1=4\xi_0(\delta)b(\delta)g_k(0,w)e^{ik\theta}$.
All remaining terms have at most a simple Cauchy or logarithmic
Cauchy singularity, with primitives $\log(v\mp i0)$ and
$\log^2(v\mp i0)/2$. Their dual norms are bounded by
$C(\|\mathcal C_2\|_{H^2_\theta}+\|\mathcal L_2\|_{H^2_\theta})$.
For the angular average,
\[
 \begin{aligned}
 &\mathsf G_{vv}(0,\theta)\mathcal L_2-\partial_v\mathsf G_{vv}(0,\theta)\mathcal C_2
                  +\mathsf G_{v\theta}(0,\theta)\partial_\theta \mathcal C_2\\
 &\quad=\frac{1}{4\sigma^2}\bigl\{
       4\xi_0(\delta)b(\delta)\Pi_ka_1+2\mathsf G_{vv}(0,\theta)a_1\\
 &\hspace{38mm}-2\partial_w(\mathsf G_{vv}(0,\theta)a_1)
       -\partial_\theta(\mathsf G_{v\theta}(0,\theta)a_1)\bigr\}.
 \end{aligned}
\]
Its average and $[\mathsf G_{vv}(0,\theta)\mathcal C_2]_0$ vanish. Taylor expansion
$\mathsf G_{vv}=\mathsf G_{vv}(0,\theta)+v\partial_v\mathsf G_{vv}(0,\theta)+O(\delta^2v^2)$, $\mathsf G_{v\theta}=\mathsf G_{v\theta}(0,\theta)+O(\delta^2v)$ then implies
\[
 \|[(\mathsf G\nabla \mathcal S_2)_v]_0\|_{L^2_v}
 \le C\delta^2(\|\mathcal C_2\|_{H^1_\theta}+\|\mathcal L_2\|_{H^1_\theta}).
\]
Thus the average source is the derivative of an $L^2_v$ flux.
Use the smooth function $\eta$ from the definition of $D_\lambda$,
with $\eta=0$ on $[0,1]$ and $\eta=1$ on $[2,\infty)$.
The preceding estimates, also on $|v|<2$, imply
\begin{equation}\label{h5c:second-residual}
 \begin{split}
 &\left\|\frac{1-\eta(|v|)}{\xi_0q^2}
   \left[\mathcal A_\sigma u+\xi_0b\sum_{n\ne0,k}
     \{2(\partial_\sigma V_n)g_{\sigma,n}
                +(\partial_\sigma^2V_n)g_n\}e^{in\theta}\right]
       \right\|_{H^{-1}(\mathbb T^2)}\\
 &\hspace{40mm}\le Ck^{\beta/2}\delta^{-2+\beta}M.
 \end{split}
\end{equation}
Indeed, $\mathcal A_\sigma=-\xi_0q^2(\Delta+T_\pm(\lambda))$, so the physical
pairing is exactly the cylinder pairing against the test function.
The mean of $\mathcal S_2$ is bounded by pairing it against the smooth density
$\xi_0\sigma e^{2v}$ and satisfies the same bound.

\medskip
\noindent\emph{Step 4. The first ordinary derivative away from $v=0$.}
For $\gamma=\beta,2+\beta$, set
\[
 W_\gamma(v)=e^{v/2}(1+e^{2v})^{\gamma/2-1/4}.
\]
The first ordinary derivative satisfies
\begin{equation}\label{h5c:first-ordinary}
 \chi_kg_{\sigma,k}=\ell_k\chi(v)\log(v\mp i0)+r_k,
 \qquad
 |\ell_k|+\|W_\beta^{-1}r_k\|_{H^1_k}
 \le Ck^{\beta/2}\delta^\beta M.
\end{equation}
Indeed, $\ell_k=- (a_1)_k/(2\sigma)$ and
$g_\sigma=-(k/8)g_\lambda$. Use the fixed threshold cutoff
$\eta=\eta(q/\sqrt{32\kappa/k})$ from
Proposition~\ref{prop:spectral-bounds}, so
$D_\lambda=\partial_\lambda+\eta V$. The identity
\[
 g_\sigma=-\frac k8(\Lambda_g+R_g-\eta V g),
 \qquad -\frac k8\Lambda_g
       =-\frac{a_1}{2\sigma}\chi(v)\log(v\mp i0),
\]
uses the same logarithmic cutoff as the first subtraction.
On $q<\sqrt{32\kappa/k}$ the field $\eta V$ vanishes; on its
remaining support in $q<q_k$, one has $q\sim k^{-1/2}$ and
no critical circle.
On this annulus,
$\eta V=-4\eta b(\lambda q^2)^{-1}\partial_v$ has bounded
cylinder coefficients and first derivatives. Solving
\eqref{h5c:full-center-equation} for $g_{vv}$ and using
\eqref{h5c:collar-input}, we obtain
\[
 \|\eta Vg\|_{H^1(q<q_k)}
 \le C\|\nabla_{v,\theta}g\|_{H^1(q\sim k^{-1/2};\,dv\,d\theta)}
 \le CM.
\]
It follows from \eqref{r2:high-regular-g} that
\[
 \left\|g_\sigma+\frac{a_1}{2\sigma}
                       \chi(v)\log(v\mp i0)\right\|_{H^1(q<q_k)}
 \le CkM.
\]
It follows from
\eqref{r2:high-jet-coefficient} that $\|g_\sigma\|_{L^2(q<q_k;dx\,dy)}\le CkM$.
Use $1-\eta(|v|)=1$ on $|v|\le1$, supported in $|v|<2$.
It follows from \eqref{h5c:critical}, \eqref{h5c:critical-mixed}, and
$\partial_\sigma=(2\sigma)^{-1}(\partial_w-\partial_v)$ that
\[
 \|(1-\eta(|v|))r_k\|_{H^1_k}\le Ck^{\beta/2}\delta^\beta M.
\]
For $\eta(|v|)\chi_kg_{\sigma,k}$ the logarithm vanishes and
this bound on $|v|<2$ controls its cutoff commutators. On the two tails,
\[
 \begin{aligned}
 W_\beta^{-1}q^2&\le
 \begin{cases}C\delta^\beta q_k^{2-\beta},&v>1,\ q<q_k,\\
 C\delta^2,&v<-1,\end{cases}\\
 W_\beta^{-1}|\partial_\sigma V_k|W_{2+\beta}
 &\le C\sigma^{-1},\qquad |v|>1.
 \end{aligned}
\]
These inequalities bound the projected metric fluxes by
$Ck^{\beta/2}\delta^\beta M$ and, with \eqref{h5c:critical},
the remaining multiplier source by $Ck^{-2+\beta/2}\delta^\beta M$.
On $q\sim q_k$, project the $L_0R_g$ bound in
\eqref{r2:high-regular-g} and use the nonresonant equation and
$g_\sigma=-(k/8)(D_\lambda g-\eta Vg)$ to obtain
\[
 \|(D_\lambda g)_k\|_{H^1_k}\le Ck^{-1}M,
 \qquad \|\partial_vg_k\|_{H^1_k}\le Ck^{-1}M,
 \qquad \|g_{\sigma,k}\|_{H^1_k}\le CM.
\]
The commutator
$-2\partial_v(\chi_k'g_{\sigma,k})+\chi_k''g_{\sigma,k}$
therefore has weighted dual norm at most
$C\delta^\beta k^{\beta/2}M$, since
$W_\beta^{-1}\le C\delta^\beta k^{\beta/2}$ on this annulus.
It follows from Lemma~\ref{spectral:weighted-normal-family}
with exponents $1/2,\beta$ that \eqref{h5c:first-ordinary} holds.

\medskip
\noindent\emph{Step 5. The remaining source terms.}
For $|v|>1$, direct differentiation of $V_k$ implies
\[
\begin{array}{c|cc}
 &v>1&v<-1\\ \hline
 |\partial_\sigma V_k|&Cq^{-2}&Cq^2\delta^{-4}\\
 |\partial_\sigma^2V_k|&Cq^{-4}&Cq^2\delta^{-6}.
\end{array}
\]
For $|v|>1$ these inequalities imply
\[
 \begin{split}
 \|\mathbf1_{|v|>1}(\partial_\sigma V_k)\chi_kg_{\sigma,k}\|_{H^{-1}_k}
 &\le C\delta^{-2}\|W_\beta^{-1}r_k\|_{H^1_k},\\
 \|\mathbf1_{|v|>1}(\partial_\sigma^2V_k)\chi_kg_k\|_{H^{-1}_k}
 &\le C\delta^{-4}\|W_{2+\beta}^{-1}\chi_kg_k\|_{H^1_k}.
 \end{split}
\]
Indeed, the respective weighted multipliers are bounded by
$C\delta^{-2}e^{-(2-\beta)v}$ and
$C\delta^{-4}e^{-(2-\beta)v}$ for $v>1$, and by
$C\delta^{-2}e^{5v/2}$ and $C\delta^{-4}e^{5v/2}$ for $v<-1$.
It follows from \eqref{h5c:critical} and
\eqref{h5c:first-ordinary} that \eqref{h5c:second-residual}
holds on the critical tails as well. For $n\ne0,k$,
\[
 |\partial_\sigma V_n|\le Cq^2\frac{k|n|}{|n-k|^2},\qquad
 |\partial_\sigma^2V_n|\le Cq^2\frac{k^2|n|}{|n-k|^3}.
\]
If $\phi$ has unit physical energy, the angular Poincar\'e inequality on
$q<q_k$ implies
\begin{align*}
 \left|\sum_{n\ne0,k}\int
       (\partial_\sigma V_n)g_{\sigma,n}\overline{\phi_n}
                                          \,dv\,d\theta\right|
 &\le Cq_k\sup_{n\ne0,k}\frac{k}{|n-k|^2}
                         \|g_\sigma\|_{L^2(dxdy)}
 \le Ck^{3/2}M,\\
 \left|\sum_{n\ne0,k}\int
       (\partial_\sigma^2V_n)g_n\overline{\phi_n}
                                          \,dv\,d\theta\right|
 &\le Cq_k^2\sup_{n\ne0,k}
       \frac{k^2}{|n|^3|n-k|^3}
       \|\nabla\partial_\theta^2g\|_2\le CM.
\end{align*}
On the transition
annulus $q\sim q_k$, all denominators satisfy
$|nb-\lambda|\ge c(1+|n-k|)$, and
\[
 \sup_{n\ne0}\frac{|\partial_\sigma V_n|}{|n|}\le C,
 \qquad
 \sup_{n\ne0}\frac{|\partial_\sigma^2V_n|}{|n|}\le Ck.
\]
On the annulus $A_k=\{q_k/2<q<q_k\}$, the logarithmic
correction vanishes, and
\[
 g_\sigma=-\frac{k}{8}(R_g-\eta Vg),\qquad
 \eta V=-\frac{4\eta b(q)}{\lambda q^2}\partial_v,
 \qquad
 \left|\frac{4\eta b(q)}{\lambda q^2}\right|\le C.
\]
It follows from \eqref{r2:high-regular-g},
\eqref{h5c:collar-input}, and the angular Poincar\'e inequality
that
\[
 \|(I-\Pi_0)g_\sigma\|_{L^2(A_k;dv\,d\theta)}
 +k\|(I-\Pi_0)g\|_{L^2(A_k;dv\,d\theta)}
 \le CkM.
\]
For $\|\phi\|_{H^1(\mathbb T^2)}\le1$, we also have
$\|\partial_\theta\phi\|_{L^2(A_k;dv\,d\theta)}\le C$.
Parseval's identity and the displayed multiplier bounds
therefore imply
\begin{equation}\label{h5c:annular-second-source}
 \left\|
 \frac{\mathbf1_{A_k}b(q)}{q^2}
 \sum_{n\ne0}
 \bigl\{2(\partial_\sigma V_n)g_{\sigma,n}
          +(\partial_\sigma^2V_n)g_n\bigr\}e^{in\theta}
 \right\|_{H^{-1}(\mathbb T^2)}
 \le CkM.
\end{equation}
The factor $q^{-2}$ cancels against the physical measure
$dx\,dy=\xi_0q^2\,dv\,d\theta$ in the dual pairing.
The same estimate holds with an additional bounded radial cutoff,
in particular with $\chi_k$.

Combining \eqref{h5c:annular-second-source} with
\eqref{h5c:second-residual}, the critical-mode estimates on
$|v|>1$, and the two estimates for $n\ne0,k$ in this step,
we obtain
\[
 \|\chi_k(\Delta+T_\pm(\lambda))
           (g_{\sigma\sigma}-\mathcal S_2)\|_{H^{-1}}
 \le C\bigl(k^{\beta/2}\delta^{-2+\beta}
                  +k^{3/2}+k\bigr)M.
\]
Since $k\ge1$, this proves \eqref{h5c:final-local}.

For $\lambda>k$, set
$\sigma=8(\lambda/k-1)>0$, $\delta=\sqrt{\sigma}$, and
$v=\log(q/\delta)$. Differentiation at fixed $q$ implies
\[
 V_k=\frac{8q^2}{q^2+\sigma},\qquad
 \partial_\sigma V_k=-\frac{8q^2}{(q^2+\sigma)^2},\qquad
 \partial_\sigma^2V_k=\frac{16q^2}{(q^2+\sigma)^3}.
\]
With the weights $W_\gamma$ defined in Step~4, these formulas
imply, for every $v\in\mathbb R$,
\[
\begin{aligned}
 W_\beta^{-1}|\partial_\sigma V_k|W_{2+\beta}
 &=\frac8{\sigma}\frac{e^{2v}}{1+e^{2v}}
 \le 8\sigma^{-1},\\
 |\partial_\sigma V_k|W_\beta&\le C_\beta\sigma^{-1},\\
 |\partial_\sigma^2V_k|W_{2+\beta}&\le C_\beta\sigma^{-2}.
\end{aligned}
\]
Thus the weighted inverse argument in Step~4 applies with
$\mathcal N_{k,-1,0}$, using
Lemma~\ref{spectral:weighted-normal-family}, and implies
\[
 \|W_\beta^{-1}\chi_k g_{\sigma,k}\|_{H^1_k}
 \le C_\beta k^{\beta/2}\delta^\beta M.
\]
There is no logarithmic term in this estimate.
Together with \eqref{h5c:critical}, which also holds for
$\lambda>k$, it follows that
\[
\begin{aligned}
 &\|2(\partial_\sigma V_k)\chi_k g_{\sigma,k}
       +(\partial_\sigma^2V_k)\chi_k g_k\|_{H^{-1}_k}\\
 &\qquad\le C_\beta\left(
 \delta^{-2}\|W_\beta^{-1}\chi_k g_{\sigma,k}\|_{H^1_k}
 +\delta^{-4}\|W_{2+\beta}^{-1}\chi_k g_k\|_{H^1_k}
 \right)\\
 &\qquad\le C_\beta k^{\beta/2}\delta^{-2+\beta}M.
\end{aligned}
\]
This estimate holds on the entire $v$-axis, including $|v|\le1$.
The physical dual pairing again cancels the factor $q^{-2}$.
The estimates for $n\ne0,k$ and for $A_k$ remain valid, so
\[
 \|\chi_k(\Delta+T_\pm(\lambda))g_{\sigma\sigma}\|_{H^{-1}}
 \le C_\beta\bigl(k^{\beta/2}\delta^{-2+\beta}
                        +k^{3/2}\bigr)M.
\]
This is \eqref{h5c:final-local} with $\mathcal S_2=0$.
For negative $k$, complex conjugation and
$(k,\lambda)\mapsto(-k,-\lambda)$ imply the same estimate
with $k$ replaced by $|k|$.
\end{proof}

Differentiating $\psi_f=(g-f)/z$, $z=i\lambda$, including $1/z$, we obtain
\begin{equation}\label{h5c:stream-conversion}
 \psi_{f,\lambda\lambda}
 =\frac{64}{k^2z}g_{\sigma\sigma}
       +\frac{16i}{kz^2}g_\sigma-\frac{2(g-f)}{z^3}.
\end{equation}
It follows from \eqref{h5c:first-ordinary} and
\eqref{h5c:stream-conversion} that the singular profile for
$\psi_{f,\lambda\lambda}$ is
\begin{equation}\label{h5c:stream-singular-profile}
 \chi(v)\left\{\frac{64}{k^2z}
       \left(\frac{\mathcal C_2}{v\mp i0}+\mathcal L_2\log(v\mp i0)\right)
              -\frac{8ia_1}{kz^2\sigma}\log(v\mp i0)\right\}.
\end{equation}

\subsubsection{The integrable remainder}
We now subtract the derivative of the logarithmic term in
Proposition~\ref{prop:spectral-bounds}. The remaining second derivative
is integrable in the spectral parameter, including across each
integer threshold.

We slightly enlarge the intervals of Lemma~\ref{r2:angular-regular-remainder} to
\[
 (k-5\kappa/4,k+5\kappa/4),\qquad
 (k+3\kappa/4,k+1-3\kappa/4),
 \qquad k\in\mathbb Z,
\]
and use their respective spatial cutoffs from
Proposition~\ref{prop:spectral-bounds}. Choose a subordinate smooth
partition $\sum_\nu\chi_\nu=1$ with bounded overlap and uniformly bounded
first two derivatives. On the interval indexed by $\nu$, set
\[
 V_\nu=\eta(q/q_\nu)V=\lambda^{-1}V_{0,\nu},\qquad
\mathcal{D}_{\nu}=\partial_\lambda+V_\nu,\qquad
\mathcal{D}_{\nu}\Gamma_f=\Lambda_{f,\nu}+R_{f,\nu}.
\]
Denote the corresponding integer by $k_\nu$, and write $k=k_\nu$
when estimating one interval. Here $q_\nu=q_{0,k}^{I}$ or
$q_{0,k}^{J}$, according to the interval, and $V_{0,\nu}$ is
independent of $\lambda$. For the estimates as $|k|\to\infty$,
write $K=1+|k|$, so $q_\nu\sim K^{-1/2}$.
The transition region $q_\nu<q<2q_\nu$ contains no critical circle.
Indeed, on the enlarged $I_k$ the resonant circle satisfies
$q^2=\sigma_k\le10\kappa/|k|$, whereas
$(q_{0,k}^{I})^2=32\kappa/|k|$; the other circles remain outside
$q\le2q_{0,k}^{I}$. On the enlarged $J_k$, the distance from
$\lambda$ to either adjacent integer is at least $3\kappa/4$, while
$\max\{|k|,|k+1|\}q^2/8\le\kappa/4$ for $q\le2q_{0,k}^{J}$.
The logarithmic cutoff has support $|\log(q/\sqrt{\sigma_k})|<1/2$,
so
\[
 \operatorname{supp}\Lambda_{f,\nu}\subset\{q<q_\nu\},
 \qquad V_\nu\Lambda_{f,\nu}=0.
\]

\begin{proposition}\label{quad:second-regular}
Under the hypotheses of Proposition~\ref{quad:spectral},
\[
 \sum_\nu\int_{\operatorname{supp}\chi_\nu}
          \|\mathcal{D}_{\nu} R_{f,\nu}(\lambda)\|_2\,d\lambda\le CM.
\]
Moreover,
\[
 \sum_\nu\sum_{n\in\operatorname{supp}\chi_\nu
                         \cap(\mathbb Z\setminus\{0\})}
 \bigl(\|R_{f,\nu}(n-)\|_2+\|R_{f,\nu}(n+)\|_2\bigr)
 \le CM.
\]
\end{proposition}
\begin{proof}
Fix $0<\alpha,\beta<1$. We first estimate the intervals with $K$
large; all norms on subsets of a center chart use physical measure.

\medskip\noindent\textit{Step 1: the gradient and source terms.}
It follows from \eqref{r2:high-regular-g},
\eqref{r2:mild-center-weight}, and \eqref{r2:H-cylinder} that,
for either boundary value,
\[
\mathcal{D}_{\nu} g=\Lambda_g+R_g,\qquad
 \|R_g\|_{H^1}\le CK^{-\alpha/2}M,\qquad
 \|q^{-\alpha}\nabla g\|_{L^2(q<q_*)}\le C_\alpha M.
\]
Here $\eta=\eta(q/q_\nu)$, $V_\nu=\eta V$, and
$\zeta=\operatorname{div}V_\nu$. It follows from
$V_\nu=-4\eta b(q)(\lambda q)^{-1}\partial_q$,
Lemma~\ref{geom:metric}, and \eqref{quad:weak-coefficients}
that, for $q_\nu\le q<q_*$,
\[
\begin{aligned}
 |K_\eta|&\le CK^{-1}q^{-2},\\
 \bigl|\mathcal D_\nu K_\eta
       -K_\eta(D_xV_\nu)^t-(D_xV_\nu)K_\eta+\zeta K_\eta\bigr|
 &\le CK^{-2}(q^{-4}+q^{-2}).
\end{aligned}
\]
Both matrices vanish for $q\le q_\nu$. Since
$\operatorname{supp}\Lambda_g\subset\{q<q_\nu\}$ and
$\mathcal D_\nu g=\Lambda_g+R_g$, we have
$K_\eta\nabla\mathcal D_\nu g=K_\eta\nabla R_g$.
Using $q_\nu\sim K^{-1/2}$ and the preceding bounds on
$R_g$ and $q^{-\alpha}\nabla g$, we obtain
\[
\begin{aligned}
 \|K_\eta\nabla\mathcal D_\nu g\|_{L^2(q<q_*)}
 &\le CK^{-1}q_\nu^{-2}\|R_g\|_{H^1}
 \le C_\alpha K^{-\alpha/2}M,
\end{aligned}
\]
and
\[
\begin{aligned}
 &\bigl\|
 [\mathcal D_\nu K_\eta-K_\eta(D_xV_\nu)^t
       -(D_xV_\nu)K_\eta+\zeta K_\eta]\nabla g
 \bigr\|_{L^2(q<q_*)}\\
 &\qquad\le CK^{-2}
       (q_\nu^{\alpha-4}+q_\nu^{\alpha-2})
       \|q^{-\alpha}\nabla g\|_{L^2(q<q_*)}\\
 &\qquad\le C_\alpha K^{-\alpha/2}M.
\end{aligned}
\]
Thus the vector field inside the second divergence in
\eqref{quad:second-weak-equation} has $L^2(q<q_*)$ norm
at most $C_\alpha K^{-\alpha/2}M$.

On $q>q_*/4$, the respective coefficient $L^4$ norms are
$CK^{-1}$ and $CK^{-2}$. Apply the local elliptic estimate in
Lemma~\ref{quad:second-away} on $q>q_*/4$, using its right-hand
side on $q>q_*/8$. It follows from \eqref{hg:exterior-B},
\eqref{eq:fixed-source-derivative}, and
\eqref{hg:fixed-exterior-W18} that
\[
 \|\mathcal{D}_{\nu} g\|_{W^{1,4}(q>q_*/4)}
 \le C\bigl(K\|R_g\|_{H^1}+K^{-1}M\bigr)
 \le CK^{1-\alpha/2}M,
 \qquad
 \|g\|_{W^{1,8}(q>q_*/4)}\le CM.
\]
Hölder's inequality therefore bounds the gradient flux in
\eqref{quad:second-weak-equation} by $CK^{-\alpha/2}M$ globally.

For \eqref{quad:source-flux}, use
$H(0)=0$, $\|H\|_{H^3}\le CM$, and
\[
 |V_\nu|\le\frac{C}{Kq},\qquad |H|\le CqM,
 \qquad |\nabla H|\le CM
 \quad(q_\nu<q<q_*).
\]
Consequently,
\[
 \|V_\nu(V_\nu\cdot\nabla H)\|_{L^2(q<q_*)}
 \le CK^{-2}M\left(\int_{q_\nu}^{q_*}q^{-3}\,dq\right)^{1/2}
 \le CK^{-3/2}M.
\]
Since $\partial_\lambda V_\nu=-V_\nu/\lambda$ and
$\zeta=\operatorname{div}V_\nu$ is $O(K^{-1})$ outside the
transition region and $O(1)$ inside it,
\[
 \| (\partial_\lambda V_\nu+\zeta V_\nu)H
       \|_{L^2(q<q_*)}\le CK^{-3/2}M.
\]
Here the transition region has area $O(K^{-1})$ and
$|V_\nu H|\le CK^{-1}M$. The coefficient bounds on $q>q_*/4$
bound the remaining source flux by $CK^{-2}M$.

\medskip\noindent\textit{Step 2: the scalar and cutoff terms.}
In each center chart, set $V_0=\lambda V$ and
$\zeta_0=\operatorname{div}V_0$. The identities
\[
 V_0=-\frac{4b(q)}q\partial_q,\qquad
 \zeta_0=-\frac4{\xi_0(q)q}\partial_q(\xi_0(q)b(q))
\]
show that $\zeta_0$ and
$-\zeta_0+V_0\cdot\nabla\zeta_0+\zeta_0^2$ extend smoothly to
$q=0$. Outside $q_\nu<q<2q_\nu$,
\[
 \zeta=\lambda^{-1}\eta\zeta_0,
 \qquad
\mathcal{D}_{\nu}\zeta+\zeta^2
 =\lambda^{-2}\eta
       (-\zeta_0+V_0\cdot\nabla\zeta_0+\zeta_0^2).
\]
On the transition region the differences from these expressions
are bounded uniformly. There the separated multipliers satisfy
\[
 \sup_{n\ne0}\frac{|\partial_\lambda^jT_n|}{|n|}\le C_j,
 \qquad j=0,1,2,
\]
and hence
\[
 \|\partial_\lambda^jTu\|_{L^2(q_\nu<q<2q_\nu)}
 \le C\|\partial_\theta u\|_{L^2(q_\nu<q<2q_\nu)}
 \le Cq_\nu\|\nabla u\|_2.
\]
Thus the transition contributions to $\zeta T$ and
$(\mathcal{D}_{\nu}\zeta+\zeta^2)T$ have norm $CK^{-1/2}$ from $H^1$ to
$H^{-1}$.

It follows from Corollary~\ref{cor:global-high-frequency} and
$T=\Delta K$ that
\[
 \|Th\|_{H^{-1}}\le CK^{-1/2}\|h\|_{H^1}.
\]
For a fixed smooth radial function $a$ in a center chart and
$\varphi\in H_0^1(q<q_*)$, multiplication by $\eta a$ satisfies
\[
 \|\eta a\varphi\|_{H^1}
 \le Cq_\nu^{-1}\|\varphi\|_{H^1}
 \le CK^{1/2}\|\varphi\|_{H^1}.
\]
Consequently, for $j=1,2$,
\[
 \begin{aligned}
 |\langle \lambda^{-j}\eta aTh,\varphi\rangle|
 &\le C|\lambda|^{-j}K^{-1/2}
          \|h\|_{H^1}\|\eta a\varphi\|_{H^1}\\
 &\le CK^{-j}\|h\|_{H^1}\|\varphi\|_{H^1}.
 \end{aligned}
\]
Here localized $H^{-1}$ norms are dual to $H_0^1$ on the
indicated region. Applying this estimate with
$a=\zeta_0$ and
$a=-\zeta_0+V_0\cdot\nabla\zeta_0+\zeta_0^2$, respectively, we obtain
\[
 \begin{aligned}
 \|\lambda^{-1}\eta\zeta_0Th\|_{H^{-1}(q<q_*)}
 &\le CK^{-1}\|h\|_{H^1},\\
 \|\lambda^{-2}\eta
   (-\zeta_0+V_0\cdot\nabla\zeta_0+\zeta_0^2)Th
                                      \|_{H^{-1}(q<q_*)}
 &\le CK^{-2}\|h\|_{H^1}.
 \end{aligned}
\]
These are the two contributions outside the transition annulus.
Since $T$ acts separately on each streamline,
$\operatorname{supp}(T\Lambda_g)\subset\{q<q_\nu\}$.
The function $\zeta$ vanishes there, and hence
\[
 \zeta T\mathcal D_\nu g=\zeta TR_g.
\]

On $q>q_*/4$, we have $\eta=1$ for sufficiently large $K$.
The radial weights
$\zeta_0$ and $-\zeta_0+V_0\cdot\nabla\zeta_0+\zeta_0^2$
satisfy \eqref{hg:weight-class}. Indeed, with $\tau=S(r)$,
\[
 V_0\cdot\nabla\tau=-\tau,\qquad
 -\zeta_0+V_0\cdot\nabla\zeta_0+\zeta_0^2
 =-\zeta_0-\tau\partial_\tau\zeta_0+\zeta_0^2,
\]
and the required polynomial bounds follow from
Lemma~\ref{geom:period-lemma}.
It follows from \eqref{hg:exterior-B},
\eqref{hg:T-improved}, and
$\mathcal B\hookrightarrow H^{-1}$ that
\[
 \begin{aligned}
 \|\zeta TR_g\|_{H^{-1}(q>q_*/4)}
 &\le C|\lambda|^{-1}K\|R_g\|_{H^1}
 \le CK^{-\alpha/2}M,\\
 \|(\mathcal D_\nu\zeta+\zeta^2)Tg\|_{H^{-1}(q>q_*/4)}
 &\le C|\lambda|^{-2}K^{-1}
             \sum_{j=0}^2\|L_0^jg\|_{H^1}
 \le CK^{-3}M,
 \end{aligned}
\]
where the last inequalities follow from
\eqref{r2:high-regular-g} and \eqref{r2:angular-g}.

For the remaining cutoff terms in
\eqref{quad:second-weak-equation}, recall that
$E_\eta=(1-\eta)\partial_\lambda T$.
Since $\partial_\lambda\eta=0$,
$\partial_\lambda V=-V/\lambda$, and
$[V,T]=-\partial_\lambda T$, we have
\[
 \begin{aligned}
 \mathcal D_\nu E_\eta
 &=(1-\eta)^2\partial_\lambda^2T
   -\eta\left(V\cdot\nabla\eta+\frac{1-\eta}{\lambda}\right)
                              \partial_\lambda T,\\
 (\mathcal D_\nu E_\eta+2\zeta E_\eta)g
 &=(1-\eta)^2(\partial_\lambda^2T)g\\
 &\quad+
 \left[2\zeta(1-\eta)-\eta V\cdot\nabla\eta
                   -\frac{\eta(1-\eta)}{\lambda}\right]
                              (\partial_\lambda T)g.
 \end{aligned}
\]
Set $A_\nu=\{q_\nu<q<2q_\nu\}$. On this annulus,
\[
 |V\cdot\nabla\eta|+|\zeta|\le C,
 \qquad \mathcal D_\nu g=R_g.
\]
The preceding separated multiplier estimates therefore imply
\[
 \begin{aligned}
 &\|2E_\eta\mathcal D_\nu g\|_{L^2(A_\nu)}
 +\|(\mathcal D_\nu E_\eta+2\zeta E_\eta)g\|_{L^2(A_\nu)}\\
 &\qquad\le Cq_\nu
       \bigl(\|R_g\|_{H^1}+\|g\|_{H^1}\bigr)
 \le CK^{-1/2}M.
 \end{aligned}
\]
This also bounds their $H^{-1}$ contributions.

On a threshold interval, $\eta=0$ in $q<q_\nu$, so the
remaining cutoff terms there are
\[
 -2(\partial_\lambda T)\partial_\lambda g
 -(\partial_\lambda^2T)g
 =(\Delta+T)\partial_\lambda^2g.
\]
They are estimated together with the singular correction:
for $\lambda=k(1\pm\sigma/8)$,
\[
 (\Delta+T)\left(\mathcal D_\nu^2g
                       -\frac{64}{k^2}\mathcal S_2\right)
 =\frac{64}{k^2}(\Delta+T)
                       (g_{\sigma\sigma}-\mathcal S_2)
 \qquad(q<q_\nu).
\]
Thus \eqref{h5c:final-local} applies after multiplication by
$64/k^2$, with $\mathcal S_2=0$ on the nonresonant side.
On $J_k$, the multipliers are separated throughout
$q<2q_\nu$, and the preceding $CK^{-1/2}M$ estimate applies
on that entire region.

\medskip\noindent
\textit{Step 3: inversion and identification of the remainder.}
We first consider $k\ge1$ and $\lambda<k$, and set
$\sigma=8(1-\lambda/k)=\delta^2$.
For either boundary value, define
\[
 u=\mathcal D_\nu^2g
       -\frac{64}{k^2}\sum_{\text{two center charts}}\mathcal S_2,
\]
where each profile from \eqref{h5c:second-singular-profile}
is extended by zero outside its chart.
It follows from \eqref{h5c:sigma-derivatives},
\eqref{h5c:regularized-density}, and \eqref{h5c:full-mixed}
that $u\in H^1$ near $q=\delta$.
On the complement, the local argument in
Lemma~\ref{quad:second-away}, applied to
\eqref{quad:second-weak-equation}, proves the same regularity.
Thus $u\in H^1(\mathbb T^2)$.

We first specify how \eqref{h5c:final-local} enters the global
source estimate. On $q_\nu<q<q_k$, all angular multipliers are
separated from their poles and
$\partial_\lambda g=R_g-\eta Vg$. Since $V$ and $\eta$ are radial,
\eqref{r2:angular-g} and \eqref{r2:high-regular-g} imply
\[
 \begin{aligned}
 \|\partial_\theta\partial_\lambda g\|_{L^2(q_\nu<q<q_k)}
 &\le Cq_k\|R_g\|_{H^1}
   +CK^{-1/2}\|\nabla\partial_\theta g\|_{L^2(q<q_k)}
 \le CK^{-1/2}M,\\
 \|\partial_\theta g\|_{L^2(q_\nu<q<q_k)}
 &\le Cq_k\|g\|_{H^1}\le CK^{-1/2}M.
 \end{aligned}
\]
Consequently,
\[
 \|2(\partial_\lambda T)\partial_\lambda g
             +(\partial_\lambda^2T)g\|_{L^2(q_\nu<q<q_k)}
 \le CK^{-1/2}M.
\]
In $q<q_\nu$, the cutoff terms in
\eqref{quad:second-weak-equation} equal
$-2(\partial_\lambda T)\partial_\lambda g
-(\partial_\lambda^2T)g$.
Since $\chi_k=1$ on the support of $\mathcal S_2$, we may
therefore use \eqref{h5c:final-local} with its original cutoff
$\chi_k$; the additional annular terms satisfy the preceding
$L^2$ bound. Together with Steps~1--2, this proves
\[
 \|(\Delta+T)u\|_{H^{-1}}
 \le CM\left(
 K^{-\alpha/2}+K^{-1/2}
 +K^{-2+\beta/2}\sigma^{-1+\beta/2}\right).
\]

To estimate the mean, use \eqref{quad:restore-mean},
$\zeta\Lambda_g=0$, and the coefficient bounds
\[
 \|\zeta\|_2\le CK^{-1/2},\qquad
 \|\mathcal D_\nu\zeta+\zeta^2\|_{4/3}\le CK^{-3/4}.
\]
They imply
\[
 |\langle\mathcal D_\nu^2g,1\rangle|
 \le C\left(K^{-1/2}\|R_g\|_2
                    +K^{-3/4}\|g\|_4\right)
 \le C\left(K^{-(1+\alpha)/2}+K^{-3/4}\right)M.
\]
For each subtracted profile, the physical measure is
$\delta^2\xi_0(\delta e^v)e^{2v}\,dv\,d\theta$.
The identity
$(v\mp i0)^{-1}=\operatorname{pv}(1/v)\pm i\pi\delta_0$
and \eqref{h5c:second-coeff} imply
\[
 \begin{aligned}
 |\langle\mathcal S_2,1\rangle|
 &\le C\delta^2
       \left(\|\mathcal C_2\|_{L^2_\theta}
                         +\|\mathcal L_2\|_{L^2_\theta}\right)\\
 &\le Ck^{-1/2+\beta/2}\delta^\beta M.
 \end{aligned}
\]
Since $\delta\sqrt{k}\le c_1$, it follows that
$|\int u|\le CK^{-1/2}M$.
Applying \eqref{quad:uniform-full-inverse}, we obtain
\begin{equation}\label{quad:subtracted-g-estimate}
 \|u\|_{H^1}
 \le CM\left(
 K^{-\alpha/2}+K^{-1/2}
 +K^{-2+\beta/2}\sigma^{-1+\beta/2}\right).
\end{equation}

We next compare this subtraction with
$\mathcal D_\nu\Lambda_g=\partial_\lambda\Lambda_g$.
In each center chart,
\[
 \Lambda_g=\frac{4a_1}{k\sigma}\chi(v)\log(v\mp i0),
 \qquad
 \partial_\lambda w=-\frac4{k\sigma},\qquad
 \partial_\lambda v=\frac4{k\sigma}.
\]
Differentiating this expression and using
\eqref{h5c:second-singular-profile}, we obtain the exact
cancellation
\begin{equation}\label{quad:log-subtraction-comparison}
 \frac{64}{k^2}\mathcal S_2-\partial_\lambda\Lambda_g
 =\frac{16}{k^2\sigma^2}
   \left[(2a_2-\partial_wa_1)\chi(v)-a_1\chi'(v)\right]
                                      \log(v\mp i0).
\end{equation}
In particular, no Cauchy term remains.
The physical measure and \eqref{h5c:jet-coeff} imply
\[
 \left\|\frac{64}{k^2}\mathcal S_2
                   -\partial_\lambda\Lambda_g\right\|_2
 \le Ck^{-9/2+\beta/2}\sigma^{-1/2+\beta/2}M
 \le CK^{-2+\beta/2}\sigma^{-1+\beta/2}M.
\]
Since
\[
 \mathcal D_\nu R_g
 =u+\sum_{\text{two center charts}}
       \left(\frac{64}{k^2}\mathcal S_2
                            -\partial_\lambda\Lambda_g\right),
\]
the same right-hand side as in
\eqref{quad:subtracted-g-estimate} bounds
$\|\mathcal D_\nu R_g\|_2$.

Finally, the boundary-value identities imply
\[
 \begin{aligned}
 R_{f,\nu}
 &=\frac{R_{g,+}-R_{g,-}}{2\pi i\lambda}
                                      -\frac{\Gamma_f}{\lambda},\\
 \mathcal D_\nu R_{f,\nu}
 &=\frac{\mathcal D_\nu R_{g,+}
                  -\mathcal D_\nu R_{g,-}}{2\pi i\lambda}
       -\frac{2R_{f,\nu}}{\lambda}
       -\frac{\Lambda_{f,\nu}}{\lambda}.
 \end{aligned}
\]
By \eqref{r2:high-regular-density},
\[
 |\lambda|^{-1}\|R_{f,\nu}\|_2
 \le CK^{-2-\alpha/2}M.
\]
Moreover, \eqref{r2:high-jet-coefficient} and
$\Lambda_{f,\nu}=(\Lambda_{g,+}-\Lambda_{g,-})/(2\pi i\lambda)$
imply
\[
 |\lambda|^{-1}\|\Lambda_{f,\nu}\|_2
 \le CK^{-2}\sqrt{\sigma}
       \sum_\pm\|b_{g,k}^{\pm}\|_{L^2_\theta}
 \le CK^{-5}M.
\]
It follows that
\begin{equation}\label{quad:second-regular-bound}
 \|\mathcal D_\nu R_{f,\nu}\|_2
 \le CM\left(
 K^{-1-\alpha/2}
 +K^{-3+\beta/2}\sigma^{-1+\beta/2}
 +K^{-3/2}\right).
\end{equation}

For $\lambda>k$, set $\sigma=8(\lambda/k-1)$.
The same argument applies with
$\mathcal S_2=\Lambda_g=0$, using the nonresonant conclusion
of \eqref{h5c:final-local}.
Thus \eqref{quad:second-regular-bound} holds on both sides
of a threshold with $\sigma=8|1-\lambda/k|$.
On $J_k$, the middle term is absent.
Negative thresholds follow by complex conjugation, with
$k$ replaced by $|k|$.
Bounded intervals separated from zero are treated by
Lemmas~\ref{quad:second-away} and
\ref{quad:critical-subtraction}; near zero, we can apply
Lemma~\ref{quad:zero-second}.

\medskip\noindent\textit{Step 4: integration and threshold traces.}
On a resonant interval, $|d\lambda|=(|k|/8)d\sigma$ and
$0<\sigma<C/K$, so
\[
 \int K^{-3+\beta/2}\sigma^{-1+\beta/2}\,d\lambda
 \le CK^{-2+\beta/2}\int_0^{C/K}\sigma^{-1+\beta/2}\,d\sigma
 \le C_\beta K^{-2}.
\]
The other terms in \eqref{quad:second-regular-bound} are summable
on intervals of bounded length. We conclude that
\begin{equation}\label{quad:regular-sum}
 \sum_\nu\int_{\operatorname{supp}\chi_\nu}
         \|\mathcal{D}_{\nu} R_{f,\nu}(\lambda)\|_2\,d\lambda\le CM.
\end{equation}
Moreover, \eqref{r2:high-regular-density} and
$\|V_\nu\|_\infty\le CK^{-1/2}$ imply
\[
 \|V_\nu R_{f,\nu}\|_2\le CK^{-3/2-\alpha/2}M.
\]
Thus $\partial_\lambda R_{f,\nu}=\mathcal{D}_{\nu} R_{f,\nu}-V_\nu R_{f,\nu}$
is integrable on either side of each nonzero integer. The
one-sided $L^2$ limits exist and, by
\eqref{r2:high-regular-density} and
\eqref{r2:compact-regular-density}, satisfy the stated summable
trace bound.
\end{proof}

\subsubsection{The logarithmic term and inverse spectral moments}
Proposition~\ref{quad:second-regular} bounds
$\mathcal D_\nu R_{f,\nu}$ in $L^1_\lambda L^2_x$.
The identity \eqref{quad:identity} also contains
$\partial_\lambda\Lambda_{f,\nu}$, the derivative of the logarithmic
correction in \eqref{spectral:decomposition}. This derivative has a
Cauchy singularity. The next lemma bounds its Fourier transform
uniformly in time, as required in \eqref{quad:fourier-second}.
Equivalently, the Fourier transform of the logarithmic correction
is $O(|t|^{-1})$ in $L^2_x$.

\begin{lemma}\label{quad:log-Fourier}
Fix $k>0$, $\gamma\ge0$ and $\chi\in C_c^\infty((-2,2))$.
In one center chart $0<q<q_*$, set
\[
 \sigma=8(1-\lambda/k),\qquad \rho=q^2,\qquad
 v=\frac12\log(\rho/\sigma).
\]
Choose $\sigma_0>0$ so that $\chi(v)$ is supported in $q<q_*$
whenever $0<\sigma<\sigma_0$. Let
$a\in W^{1,1}_{\mathrm{loc}}((0,\sigma_0);L^2(\mathbb T_\theta))$
vanish near $\sigma_0$ and satisfy, for almost every $\sigma$,
\[
 \|a(\sigma)\|_{L^2_\theta}
 +\|\sigma\partial_\sigma a(\sigma)\|_{L^2_\theta}
 \le M\sigma^\gamma.
\]
For either boundary sign, define
\[
 \Lambda_a(\lambda,q,\theta)
   =a(\sigma,\theta)\chi(v)\log(v\mp i0)
       \quad\text{when }0<\sigma<\sigma_0,
\]
and extend it by zero outside this spectral interval and the
chosen center chart. Then
\begin{equation}\label{quad:log-derivative-Fourier}
 \sup_{t\in\mathbb R}
 \left\|\int_{\mathbb R}e^{-it\lambda}
             \partial_\lambda\Lambda_a(\lambda)\,d\lambda\right\|_{L^2_x}
 \le C_{\gamma,\chi}M\sigma_0^{\gamma+1/2}.
\end{equation}
Here the derivative and Fourier transform are understood as
$L^2_x$-valued distributions. Consequently, for $t\ne0$,
\begin{equation}\label{quad:log-decay}
 \left\|\int_{\mathbb R}e^{-it\lambda}
                    \Lambda_a(\lambda)\,d\lambda\right\|_{L^2_x}
 \le \frac{C_{\gamma,\chi}M\sigma_0^{\gamma+1/2}}{|t|}.
\end{equation}
The constant may depend on the fixed chart, but is independent
of $k$, $\sigma_0$ and $M$.
\end{lemma}
\begin{proof}
Since $\partial_\sigma v=-1/(2\sigma)$, distributional
differentiation implies
\begin{equation}\label{quad:log-sigma-derivative}
\begin{aligned}
 \partial_\sigma\{a\chi(v)\log(v\mp i0)\}
 ={}&-\frac{a(\sigma,\theta)\chi(v)}
              {\sigma[\log(\rho/\sigma)\mp i0]}\\
 &+\left\{(\partial_\sigma a)\chi(v)
               -\frac{a}{2\sigma}\chi'(v)\right\}\log(v\mp i0).
\end{aligned}
\end{equation}
By Lemma~\ref{geom:metric},
$dx\,dy=\frac12\xi_0(\sqrt\rho)\,d\rho\,d\theta$, with $\xi_0$
bounded above and below. The change of variables
$q=\sqrt\sigma e^v$ therefore implies
\begin{equation}\label{quad:log-regular-term}
\begin{aligned}
 \left\|\left\{(\partial_\sigma a)\chi(v)
       -\frac{a}{2\sigma}\chi'(v)\right\}\log(v\mp i0)\right\|_2
 &\le C_\chi M\sigma^{\gamma-1/2},\\
 \|a(\sigma)\chi(v)\log(v\mp i0)\|_2
 &\le C_\chi M\sigma^{\gamma+1/2}.
\end{aligned}
\end{equation}

For the first term in \eqref{quad:log-sigma-derivative}, consider
\[
 (\mathcal K_\pm h)(\rho,\theta)
 =\int_0^\infty
 \frac{\chi(\frac12\log(\rho/\sigma))h(\sigma,\theta)}
      {\sigma[\log(\rho/\sigma)\mp i0]}\,d\sigma.
\]
Under the unitary map $(Uh)(u,\theta)=e^{u/2}h(e^u,\theta)$,
$\mathcal K_\pm$ becomes convolution in $u$ with
$m(r)\operatorname{p.v.}(1/r)\pm i\pi m(0)\delta_0$, where
$m(r)=e^{r/2}\chi(r/2)$. Since
$\widehat{\operatorname{p.v.}(1/r)}=-i\pi\operatorname{sgn}$,
Plancherel's theorem implies
\begin{equation}\label{quad:log-Cauchy-bound}
 \|\mathcal K_\pm h\|_{L^2(d\rho\,d\theta)}
 \le\bigl(\tfrac12\|\widehat m\|_1+\pi|m(0)|\bigr)
                   \|h\|_{L^2(d\sigma\,d\theta)}
 \le C_\chi\|h\|_{L^2(d\sigma\,d\theta)}.
\end{equation}

The second bound in \eqref{quad:log-regular-term} tends to zero
as $\sigma\to0+$, and $a$ vanishes near $\sigma_0$.
Changing variables from $\lambda$ to $\sigma$ thus implies
\[
 \int_{\mathbb R}e^{-it\lambda}\partial_\lambda\Lambda_a\,d\lambda
 =-e^{-itk}\int_0^{\sigma_0}e^{itk\sigma/8}
       \partial_\sigma\{a\chi(v)\log(v\mp i0)\}\,d\sigma.
\]
This identity is justified first on $\varepsilon<\sigma<\sigma_0$;
\eqref{quad:log-regular-term} and \eqref{quad:log-Cauchy-bound}
justify passage to the limit in $L^2_x$ as $\varepsilon\to0+$.
Apply \eqref{quad:log-Cauchy-bound} to
$h=e^{itk\sigma/8}a\mathbf1_{(0,\sigma_0)}$ and integrate the first
bound in \eqref{quad:log-regular-term}. It follows that
\[
 \left\|\int e^{-it\lambda}\partial_\lambda\Lambda_a\,d\lambda\right\|_2
 \le C_\chi\|a\|_{L^2((0,\sigma_0)\times\mathbb T)}
       +C_\chi M\int_0^{\sigma_0}\sigma^{\gamma-1/2}\,d\sigma
 \le C_{\gamma,\chi}M\sigma_0^{\gamma+1/2}.
\]
Finally, \eqref{quad:log-decay} follows from
$\widehat{\partial_\lambda\Lambda_a}(t)=it\widehat{\Lambda_a}(t)$.
\end{proof}

To treat the term $\lambda^{-2}\Gamma_f$ we need the following result. Recall that $QL$ is
skew-adjoint on $Y\cap E^\perp$ for the inner product
$\langle(I+\Delta^{-1})\cdot,\cdot\rangle$; its inverse acts
on the subspace \eqref{zero:continuous-space}.

\begin{lemma}\label{r2:inverse-moments}
Suppose $f\in H^2\cap Y$ and $\langle f,g_j\rangle=0$ for
$j=1,2,3$. For every integer $N\ge1$,
\[
 Qf\in\operatorname{Dom}((QL)^{-N}),\qquad
 \|(QL)^{-N}Qf\|_{L^2}\le C_N\|f\|_{H^2}.
\]
\end{lemma}
\begin{proof}
By Proposition~\ref{zero:constraints}, $Qf$ is energy-orthogonal
to $\ker(QL)$. Let $\mu_{Qf}$ be its spectral measure for
$-iQL$. It follows from \eqref{lap:measure-density} that,
for $\lambda\ne0$,
\begin{equation}\label{r2:measure-density}
\begin{aligned}
 \frac{d\mu_{Qf}}{d\lambda}
 &=\langle\Delta Q\Gamma_f,(I+\Delta^{-1})Qf\rangle\\
 &=\langle Q\Gamma_f,(\Delta+1)Qf\rangle
  =\langle\Gamma_f,(\Delta+1)f\rangle.
\end{aligned}
\end{equation}
The last equality follows
from $(\Delta+1)Pf=0$ and $(\Delta+1)f\perp E$.
The exponential jump bound in Corollary~\ref{global:projected}
implies, near zero,
\[
 0\le\frac{d\mu_{Qf}}{d\lambda}
 \le Ce^{-c/|\lambda|}\|f\|_{H^2}^2.
\]
Moreover, $\mu_{Qf}(\{0\})=0$ and the measure is absolutely
continuous away from zero. Consequently,
\[
 \int_{\mathbb R}|\lambda|^{-2N}\,d\mu_{Qf}(\lambda)
 \le C_N\|f\|_{H^2}^2,
\]
where the integral over $|\lambda|\ge1$ is bounded by the total
mass. The claim follows from the spectral theorem and equivalence
of the energy and $L^2$ norms.
\end{proof}

\subsubsection{Completion of the proof}
\begin{proof}[Proof of Proposition~\ref{quad:spectral}]
\emph{Step 1: The ordinary second derivative and threshold terms.}
We use the partition $\{\chi_\nu\}$ and decompositions
$\mathcal{D}_{\nu}\Gamma_f=\Lambda_{f,\nu}+R_{f,\nu}$ specified before
Proposition~\ref{quad:second-regular}. Since
$V_\nu\Lambda_{f,\nu}=0$ thanks to the separation of support, we obtain
\begin{equation}\label{quad:identity}
 \partial_\lambda^2\Gamma_f
 =\partial_\lambda\Lambda_{f,\nu}
       +\mathcal{D}_{\nu} R_{f,\nu}-2V_\nu R_{f,\nu}
       +(V_\nu^2-\partial_\lambda V_\nu)\Gamma_f.
\end{equation}
For a nonzero integer $n$, write
$[R_{f,\nu}]_n=R_{f,\nu}(n+)-R_{f,\nu}(n-)$.
Continuity of $\Gamma_f$ and $V_\nu\Gamma_f$ at $n$, and
$\Lambda_{f,\nu}\to0$ in $L^2$, imply
\begin{equation}\label{quad:partition-identity}
\begin{aligned}
 \partial_\lambda^2(\chi_\nu\Gamma_f)
 ={}&\chi_\nu\bigl\{\partial_\lambda\Lambda_{f,\nu}
       +\mathcal{D}_{\nu} R_{f,\nu}-2V_\nu R_{f,\nu}
       +(V_\nu^2-\partial_\lambda V_\nu)\Gamma_f\bigr\}\\
 &+2\chi_\nu'\partial_\lambda\Gamma_f+\chi_\nu''\Gamma_f
       +\sum_{n\ne0}\chi_\nu(n)[R_{f,\nu}]_n\delta_n.
\end{aligned}
\end{equation}

\smallskip\noindent
\emph{Step 2: The integrable remainder and the logarithms.}
It follows from Proposition~\ref{quad:second-regular} that
\[
 \sum_\nu\int\|\chi_\nu\mathcal{D}_{\nu} R_{f,\nu}\|_2\,d\lambda
 +\sum_{\nu,n\ne0}\|\chi_\nu(n)[R_{f,\nu}]_n\|_2
 \le CM.
\]
On an interval indexed by $k_\nu$ with $|k_\nu|$ large,
\[
 \|V_\nu R_{f,\nu}\|_2
 \le C(1+|k_\nu|)^{-1/2}\|R_{f,\nu}\|_{H^1}
 \le C(1+|k_\nu|)^{-3/2-\alpha/2}M.
\]
These bounds are summable. The corresponding integral on bounded
intervals follows from Proposition~\ref{prop:spectral-bounds}
and, near zero, Lemma~\ref{quad:zero-second}.

For $b_{f,k}^{\pm}$, it follows from
\eqref{h5c:trace} and the coefficient formula preceding
Proposition~\ref{prop:spectral-bounds} that
\[
 \|b_{f,k}^{\pm}(\sigma)\|_{L^2_\theta}
 +\|\sigma\partial_\sigma b_{f,k}^{\pm}(\sigma)\|_{L^2_\theta}
 \le C|k|^{-5/2+\beta/2}\sigma^{\beta/2}M,
 \qquad 0<\sigma<C/|k|.
\]
The same estimate holds after multiplication by $\chi_\nu$, since
$\sigma|d\lambda/d\sigma|\le C$ on its support.
Apply Lemma~\ref{quad:log-Fourier} with $\gamma=\beta/2$ and
$\sigma_0=C/|k|$. We get
\[
 \sup_t\left\|\int e^{-it\lambda}
       \partial_\lambda(\chi_\nu\Lambda_{f,\nu})\,d\lambda
       \right\|_2
 \le C|k|^{-5/2+\beta/2}|k|^{-\beta/2-1/2}M
 =C|k|^{-3}M.
\]

\smallskip\noindent
\emph{Step 3: The spatial derivatives on bounded spectral intervals.}
For $a\in C_c^\infty(\mathbb R\setminus\{0\})$, set
\[
 \psi_a(t)=\int e^{-it\lambda}a(\lambda)\Gamma_f(\lambda)\,d\lambda,
 \qquad \eta_a(t)=e^{-tQL}a(-iQL)Qf.
\]
It follows from \eqref{lap:Stone-formula} and
\eqref{lap:P-resolvent} that, for every $t\in\mathbb R$,
\[
 \Delta\psi_a=\eta_a+PL(QL)^{-1}\eta_a.
\]
Here $e^{-tQL}$ is the unitary group on $Y\cap E^\perp$
for the inner product
$\langle(I+\Delta^{-1})\cdot,\cdot\rangle_{L^2}$,
as established in Proposition~\ref{zero:constraints}.
The moment conditions imply $Qf\in\mathcal Z$, where
$\mathcal Z$ is defined in \eqref{zero:continuous-space};
hence $\eta_a(t)\in\mathcal Z$ for every $t$.
The map $PL:L^2\to E$ is bounded. For a fixed
$\chi\in C_c^\infty(\mathbb R)$, Lemma~\ref{r2:inverse-moments}
and cutoff approximation at zero therefore imply
\[
 \sup_{t\in\mathbb R}\|\psi_{\chi/\lambda^2}(t)\|_{H^2}
 \le C_\chi\bigl(\|(QL)^{-2}Qf\|_2+
                       \|(QL)^{-3}Qf\|_2\bigr)
 \le C_\chi M.
\]
Write $V_\nu=\lambda^{-1}V_{0,\nu}$, where $V_{0,\nu}$ is
independent of $\lambda$. For each of the finitely many bounded
spectral intervals, \eqref{quad:saddle-coefficients} and the fixed
center cutoff imply
$V_{0,\nu}\in L^\infty\cap W^{1,4}$. Consequently,
\[
 \|(V_{0,\nu}\cdot\nabla)^2u\|_2
 \le \|V_{0,\nu}\|_\infty^2\|\nabla^2u\|_2
 +\|V_{0,\nu}\cdot\nabla V_{0,\nu}\|_4\|\nabla u\|_4
 \le C_\nu\|u\|_{H^2}.
\]
Applying this estimate and
$\|V_{0,\nu}u\|_2\le C_\nu\|u\|_{H^1}$
to $\psi_{\chi_\nu/\lambda^2}$ bounds the Fourier transforms
of $\chi_\nu V_\nu^2\Gamma_f$ and
$-\chi_\nu(\partial_\lambda V_\nu)\Gamma_f$
uniformly in $L^2_x$, for each of these finitely many
bounded spectral intervals.

\smallskip\noindent
\emph{Step 4: Summation of the spatial derivatives at large frequency.}
Recall that $V_{0,\nu}=\eta(q/q_\nu)\lambda V$, where
$\eta(q/q_\nu)=0$ for $q\le q_\nu$ and
$q_\nu\sim_\kappa(1+|k_\nu|)^{-1/2}$.
In each center chart, $V_{0,\nu}$ therefore vanishes for
$q\le q_\nu$, while on $q_\nu<q<q_*$,
\[
 |V_{0,\nu}|\le Cq^{-1},\qquad
 |V_{0,\nu}\cdot\nabla V_{0,\nu}|\le Cq^{-3}.
\]
Integration against
$q\,dq$, together with \eqref{quad:saddle-coefficients} on the
fixed exterior and the Sobolev embedding $H^2\hookrightarrow W^{1,4}$,
implies
\begin{equation}\label{quad:vector-norms}
\begin{aligned}
 \|V_{0,\nu}\|_\infty&\le C(1+|k_\nu|)^{1/2},\\
 \|V_{0,\nu}\cdot\nabla V_{0,\nu}\|_4
   &\le C(1+|k_\nu|)^{5/4},\\
 \|(V_{0,\nu}\cdot\nabla)^2u\|_2
   &\le C(1+|k_\nu|)^{5/4}\|u\|_{H^2}.
\end{aligned}
\end{equation}
Since $|\lambda|\ge c(1+|k_\nu|)$ on
$\operatorname{supp}\chi_\nu$, the spectral theorem and
\eqref{zero:energy-coercivity} imply
\[
 \|(QL)^{-1}\chi_\nu(-iQL)Qf\|_2
 \le C(1+|k_\nu|)^{-1}
          \|\chi_\nu(-iQL)Qf\|_2.
\]
Here $(QL)^{-1}$ is the inverse with values in
$\mathcal Z$, defined in \eqref{zero:continuous-space}.
The spectral identity in Step~3
therefore implies
\[
 \sup_t\|\psi_{\chi_\nu/\lambda^2}(t)\|_{H^2}
 \le C(1+|k_\nu|)^{-2}\|\chi_\nu(-iQL)Qf\|_2.
\]
It follows from bounded overlap and the spectral theorem that
\[
 \sum_\nu\|\chi_\nu(-iQL)Qf\|_2^2\le C\|f\|_2^2.
\]
Hence
\begin{align*}
 \sum_{\nu:\,|k_\nu|\text{ large}}\sup_t
 \left\|\int e^{-it\lambda}\chi_\nu V_\nu^2\Gamma_f\,d\lambda
                                                    \right\|_2
 &\le C\sum_\nu(1+|k_\nu|)^{-3/4}
                           \|\chi_\nu(-iQL)Qf\|_2\\
 &\le C\left(\sum_\nu(1+|k_\nu|)^{-3/2}\right)^{1/2}\|f\|_2.
\end{align*}
Since $-\partial_\lambda V_\nu=\lambda^{-2}V_{0,\nu}$,
\[
\begin{aligned}
 &\sup_{t\in\mathbb R}
 \left\|\int e^{-it\lambda}\chi_\nu(\lambda)
       (-\partial_\lambda V_\nu)\Gamma_f(\lambda)\,d\lambda
 \right\|_2\\
 &\qquad=\sup_{t\in\mathbb R}
       \|V_{0,\nu}\psi_{\chi_\nu/\lambda^2}(t)\|_2
 \le C(1+|k_\nu|)^{-3/2}\|\chi_\nu(-iQL)Qf\|_2.
\end{aligned}
\]
These bounds are summable by Cauchy--Schwarz, since
\[
 \sum_\nu(1+|k_\nu|)^{-3}<\infty,\qquad
 \sum_\nu\|\chi_\nu(-iQL)Qf\|_2^2\le C\|f\|_2^2.
\]

Bounded overlap of $\chi_\nu$ and $\Gamma_f\in L^1_\lambda L^2_x$ imply
$\sum_\nu\chi_\nu\Gamma_f=\Gamma_f$ in $L^1_\lambda L^2_x$.
It follows from Steps~1--4, after summing
\eqref{quad:partition-identity} in distributions, that
\[
\Big \|\widehat{\partial_\lambda^2\Gamma_f}\Big\|_{L^\infty_tL^2_x}
 \le \sum_\nu
 \|\widehat{\partial_\lambda^2(\chi_\nu\Gamma_f)}
                                      \|_{L^\infty_tL^2_x}
 \le CM.
\]
Finally,
$\widehat{\partial_\lambda^2\Gamma_f}(t)
=-t^2\widehat{\Gamma_f}(t)$ in distributions, and the
right-hand side is continuous with values in $L^2_x$.
The bound therefore holds for every $t\in\mathbb R$,
which proves \eqref{quad:fourier-second}.
\end{proof}

\section{Cubic bounds near elliptic centers}\label{sec:cubic-centers}
In this section we establish the cubic vanishing at the critical streamlines and additional bounds for the spectral density function needed for the pointwise vorticity bound \eqref{Intm9}. We note that the bound \eqref{h5c:trace} has the depletion factor $\delta^{2+\beta}$, with $0<\beta<1$, but it is not strong enough for our purposes.
 We first treat the case $|k|\ge2$ in
Subsection~\ref{subsec:cubic-higher-modes}. The first modes have the
 exponent $\sqrt{1+8}=3$ which leads to a logarithmic divergence using standard perturbation arguments, and therefore require a separate argument; we refer to~\cite{JiaVortex} and the earlier work~\cite{Bed2} for an analysis of this issue in the simpler case of monotonic vortices.

We study the region near the positive elliptic point; the other center can be
treated by the cell symmetries. We use
$q=\sqrt{8(1-b)}$ and the area element
$dx\,dy=\xi_0(q)q\,dq\,d\theta$ from Lemma~\ref{geom:metric}.
Below $v=\log(q/\delta)$, $w=\log\delta$, and in all cylinder
norms we use $dv\,d\theta$.

\subsection{Higher angular modes}
\label{subsec:cubic-higher-modes}

For $|k|\ge2$, the inequality $\sqrt{k^2+8}>3$ allows the
stronger weights in Lemma~\ref{spectral:weighted-normal-family}
used below.

\begin{lemma}[The higher angular modes]\label{cubic:profiles}
Let $f\in H^5\cap Y$ satisfy $\langle f,g_j\rangle=0$ for
$j=1,2,3$, and set $M=\|f\|_{H^5}$. Suppose $|k|\ge2$,
$q_k=c_0|k|^{-1/2}$, and $0<\delta<c_1|k|^{-1/2}$, where
$c_0,c_1$ are fixed sufficiently small constants satisfying
$c_1<e^{-6}c_0$ and $e^5c_0<q_*$. Choose a smooth radial cutoff
$0\le\chi_k\le1$ such that
\[
 \chi_k=1\quad(q\le q_k/2),\qquad
 \chi_k=0\quad(q\ge q_k),\qquad
 |(q\partial_q)^m\chi_k|\le C_m\quad(0\le m\le3).
\]
Fix $\vartheta\in C^\infty([0,\infty))$ with $0\le\vartheta\le1$,
$\vartheta=1$ on $[0,1/2]$, and $\vartheta=0$ on $[3/4,\infty)$.
Write $g_{\pm,k}=g_{f,\pm,k}$ for the $k$th angular Fourier
coefficient of $g_{f,\pm}$. For $0<q<q_k$ and
$\lambda=k(1-\delta^2/8)$, define
\[
 \begin{aligned}
 G_k(q,\delta)
 &=\chi_k(q)\vartheta\!\left(\frac{\delta\sqrt{|k|}}{c_1}\right)
       [g_{+,k}-g_{-,k}](q,\lambda),\\
 U_k(w,v)&=\delta^{-3}G_k(\delta e^v,\delta),\qquad
 v=\log(q/\delta),\quad w=\log\delta.
 \end{aligned}
\]
Define $G_k=0$ for $q\ge q_k$ or $\delta\ge c_1|k|^{-1/2}$.
For an interval $I\subset\mathbb R$, set
\[
 \|u\|_{H^1_k(I)}^2
 =\int_I\bigl(|\partial_vu|^2+k^2|u|^2\bigr)\,dv.
\]
Then
\begin{equation}\label{cubic:profile-core}
 \|U_k(w,\cdot)\|_{H^1_k((-2,2))}
 +\|\partial_wU_k(w,\cdot)\|_{H^1_k((-2,2))}
 \le C|k|^{-3/2}M.
\end{equation}
For $0<q<q_k$, we have
\begin{equation}\label{cubic:profile-outer}
 |G_k(q,\delta)|
 \le C|k|^{-2}q^{13/4}\delta^{-1/4}M,
 \qquad 0<\delta<q/2,
\end{equation}
and
\begin{equation}\label{cubic:profile-inner}
 \begin{aligned}
 |G_k(q,\delta)|
 &\le C|k|^{-2}q^{5/4}\delta^{7/4}M,\\
 |\delta\partial_\delta G_k(q,\delta)|
 &\le C|k|^{-1}q^{5/4}\delta^{7/4}M,
 \end{aligned}
 \qquad \delta>2q.
\end{equation}
The estimates \eqref{cubic:profile-core}--\eqref{cubic:profile-inner}
also hold when $\lambda=k(1+\delta^2/8)$ is used in the definition
of $G_k$.

For the original boundary values, without either cutoff, we have
\begin{equation}\label{cubic:higher-trace}
 \left|(\delta\partial_\delta)^j
       g_{\pm,k}(\delta,kb(\delta))\right|
 \le C|k|^{-2}\delta^3M,
 \qquad j=0,1,\quad 0<\delta<c_1|k|^{-1/2}.
\end{equation}
In \eqref{cubic:higher-trace}, the derivative acts on the composite
function $\delta\mapsto g_{\pm,k}(\delta,kb(\delta))$.
The constant $C$ is independent of $k$, $\delta$, and $f$.
\end{lemma}

\medskip

\begin{proof}
We take $k\ge2$ without loss of generality.
For $\lambda=k(1+\delta^2/8)$, replace
$D_w=\partial_v-k\delta^2\partial_\lambda/4$ by
$\partial_v+k\delta^2\partial_\lambda/4$.
Throughout the proof, $h_j=D_w^jg$, $j=0,1$, and $V_n$ is the
coefficient in \eqref{h5c:full-center-equation}.

\emph{Step 1: the source on a small annulus.}
Let $\ell$ be the linear Taylor polynomial of $H$ at the center.
For $0<\alpha<1$ and $0<a<q_*/2$, we have
\begin{equation}\label{cubic:source-annulus}
 \sum_{m=0}^3\ \sum_{Z_1,\ldots,Z_m}
 \|Z_1\cdots Z_m(H-\ell)\|_{L^2(dv\,d\theta;a/2<q<2a)}
 \le C_\alpha a^{1+\alpha}M.
\end{equation}
Here the inner sum runs over all ordered choices
$Z_1,\ldots,Z_m\in\{\partial_v,\partial_\theta\}$, where
$\partial_v=q\partial_q$; for $m=0$, the product denotes
the identity operator.
On $a/2<q<2a$, Taylor's theorem and
$\|H\|_{C^{1,\alpha}}\le C_\alpha M$ imply
\[
 |H-\ell|+a|\nabla(H-\ell)|
 \le C_\alpha a^{1+\alpha}M.
\]
Since $dv\,d\theta\sim a^{-2}dx\,dy$ on this annulus,
the product rule implies
\[
 \sum_{m=0}^3\sum_{Z_1,\ldots,Z_m}
 \|Z_1\cdots Z_m(H-\ell)\|_{L^2(dv\,d\theta;\,a/2<q<2a)}
 \le C\sum_{j=0}^3a^{j-1}
       \|\nabla^j(H-\ell)\|_{L^2(dx\,dy;\,a/2<q<2a)}.
\]
The terms with $j=0,1$ are bounded by
$C_\alpha a^{1+\alpha}M$ by the preceding pointwise estimate.
For $p\ge2/(1-\alpha)$, Hölder's inequality and
$\|H\|_{H^3}\le CM$ imply
\[
 \begin{aligned}
 a\|\nabla^2H\|_{L^2(a/2<q<2a)}
 &\le Ca^{2-2/p}\|\nabla^2H\|_{L^p}
 \le C_\alpha a^{1+\alpha}M,\\
 a^2\|\nabla^3H\|_{L^2(a/2<q<2a)}
 &\le Ca^2M\le Ca^{1+\alpha}M.
 \end{aligned}
\]
These norms use physical measure, and
$\nabla^j(H-\ell)=\nabla^jH$ for $j\ge2$.
This proves \eqref{cubic:source-annulus}.
Finally, $H(0)=0$ implies $\ell(x,y)=\nabla H(0)\cdot(x,y)$,
whose cylinder derivatives of order at most three have
$L^2(dv\,d\theta)$ norm at most $CaM$ on the same annulus.

\emph{Step 2: the estimate \eqref{cubic:dyadic-energy}.}
Fix either boundary value $g=g_{f,\pm}(\lambda)$.
Recall that
\[
 D_w=q\partial_q-\frac{k\delta^2}{4}\partial_\lambda
     =\partial_w\big|_v.
\]
In this step we prove
\begin{equation}\label{cubic:dyadic-energy}
 \sum_{j=0}^1\sum_{a=0}^2
 \|\nabla_{v,\theta}\partial_\theta^aD_w^jg\|_
                       {L^2(q/2<q'<2q)}
 \le Ck^{1/2}qM,\qquad 0<q<q_k/2.
\end{equation}
All norms in this step use cylinder measure.
The starting estimate is \eqref{h5c:mild} with $\beta=1/2$.
Together with \eqref{h5c:collar-input}, it implies
\begin{equation}\label{cubic:preliminary-annular}
 \sum_{j=0}^1\sum_{a=0}^2
 \|\nabla_{v,\theta}\partial_\theta^aD_w^jg\|_
                {L^2(e^{-2}\rho<q'<e^2\rho)}
 \le CM\left(\frac{\rho}{q_k}\right)^{1/2},
 \qquad 0<\rho<q_k.
\end{equation}
We use this bound to estimate the perturbation terms before
applying the inverse operators.

\smallskip\noindent
\emph{The nonzero angular modes.}
In \eqref{h5c:commuted-center-equation}, keep
$-\partial_v^2-\partial_\theta^2+V_k\Pi_k$ on the left.
For $j=0,1$ and $a=0,1,2$, the fixed source satisfies
\[
 \|\partial_\theta^a\partial_v^j(\xi_0q^2H)\|_
                   {L^2(\rho/e<q<e\rho)}
 \le C\rho^3M
\]
by Step~1.
The other sources contain the metric coefficients
$D_w^{j-i}\partial_\theta^{a-b}(\mathsf G-I)$,
the critical coefficients $D_w^{j-i}(\xi_0b-1)$,
and the noncritical coefficients
$D_w^{j-i}(\xi_0bV_n)$.
Their localized $H^{-1}$ norms are bounded by
\[
 C\rho^2
 \sum_{i=0}^1\sum_{b=0}^2
 \|\nabla\partial_\theta^bD_w^ig\|_
                     {L^2(e^{-2}\rho<q<e^2\rho)}
 \le CMq_k^{-1/2}\rho^{5/2}.
\]
Here localization means multiplication by a smooth cutoff
supported in $\rho/e<q<e\rho$, with uniformly bounded
$v$ derivatives, and $H^{-1}$ is the dual of $H^1$ on
$\mathbb R_v\times\mathbb T_\theta$.
The first inequality follows from
\eqref{h5c:step2b-coefficient-bound},
\eqref{eq:cylinder-form}, and
\eqref{h5c:step2b-noncritical-coefficient-bound};
the second follows from \eqref{cubic:preliminary-annular}.

Apply the inverse operators to
$\chi_k(I-\Pi_0)\partial_\theta^aD_w^jg$.
The terms containing $\chi_k'$ or $\chi_k''$ are supported
in $q_k/2<q<q_k$ and have dual norm at most $CM$, by
\eqref{h5c:step2b-cutoff-identity} and
\eqref{h5c:collar-input}.
For the modes $n=\pm1$, the inverse kernel is
$e^{-|v-v'|}/2$.
For the other free modes and the critical mode,
Lemma~\ref{spectral:weighted-normal-family}, with
$\mu_-=\mu_+=\pm3/2$, implies the stronger factor
$e^{-3|v-v'|/2}$.
Since
\[
 e^{-|v-v'|}
 =\min\!\left(\frac q\rho,\frac\rho q\right),
 \qquad q=\delta e^v,\quad \rho=\delta e^{v'},
\]
localizing the sources in fixed-length $v'$ intervals and
summing the inverse estimates proves, for $q<q_k/4$,
\begin{equation}\label{cubic:strip-inequality}
 \begin{aligned}
 &\sum_{j=0}^1\sum_{a=0}^2
 \|\nabla(I-\Pi_0)\partial_\theta^aD_w^jg\|_
                                      {L^2(q/2<q'<2q)}\\
 &\qquad\le CM\frac q{q_k}
 +CM\int_0^{q_k}
       \min\!\left(\frac q\rho,\frac\rho q\right)
       \left(\rho^3+q_k^{-1/2}\rho^{5/2}\right)
                                      \frac{d\rho}{\rho}.
 \end{aligned}
\end{equation}
The first term is the contribution from the derivatives
of $\chi_k$. For $p>1$,
\[
 \begin{aligned}
 \int_0^{q_k}\min\!\left(\frac q\rho,\frac\rho q\right)
                  \rho^p\,\frac{d\rho}{\rho}
 &=\frac1q\int_0^q\rho^p\,d\rho
       +q\int_q^{q_k}\rho^{p-2}\,d\rho\\
 &\le C_pq q_k^{p-1}.
 \end{aligned}
\]
Taking $p=3$ and $p=5/2$ in
\eqref{cubic:strip-inequality}, we obtain
\[
 \sum_{j=0}^1\sum_{a=0}^2
 \|\nabla(I-\Pi_0)\partial_\theta^aD_w^jg\|_
                                      {L^2(q/2<q'<2q)}
 \le CM\left(\frac q{q_k}+qq_k^2+qq_k\right)
 \le CM\frac q{q_k}.
\]

\smallskip\noindent
\emph{The angular average.}
Differentiating \eqref{h5c:mean-flux}, we obtain
\[
 \begin{aligned}
 \partial_v[D_w^jg]_0(v,w)
 ={}&\int_{-\infty}^v
       [\partial_v^j(\xi_0q^2H)]_0(v',w)\,dv'\\
 &-\sum_{i=0}^j\binom ji
   \left[(D_w^{j-i}(\mathsf G-I)\nabla D_w^ig)_v\right]_0,
 \qquad j=0,1.
 \end{aligned}
\]
In the integral, $q=\delta e^{v'}$.
Since $H(0)=0$ and $\|H\|_{C^1}\le CM$, the integral is
bounded by $Cq^3M$.
The metric terms are controlled by
\eqref{cubic:preliminary-annular}. Hence
\[
 \sum_{j=0}^1
 \|\partial_v[D_w^jg]_0\|_{L^2(q/2<q'<2q)}
 \le CM\left(q^3+q_k^{-1/2}q^{5/2}\right)
 \le CM\frac q{q_k}.
\]
For $a\ge1$, $[\partial_\theta^aD_w^jg]_0=0$.
Combining the two angular estimates proves
\eqref{cubic:dyadic-energy} for $q<q_k/4$, since
$q/q_k=c_0^{-1}k^{1/2}q$.
For $q_k/4\le q<q_k/2$, the desired conclusion follows
directly from \eqref{h5c:collar-input}.

\emph{Step 3: the weighted estimate and conclusions
\eqref{cubic:profile-core}, \eqref{cubic:profile-outer},
and \eqref{cubic:higher-trace}.}
Throughout this step, $g_{\pm,k}$ is evaluated at
$q=\delta e^v$ and $\lambda=k(1-\delta^2/8)$.
Since $D_w$ commutes with angular projection,
\[
 [D_w^jg_{f,\pm}]_k=D_w^jg_{\pm,k},
 \qquad
 [F]_k=\frac1{2\pi}\int_0^{2\pi}F e^{-ik\theta}\,d\theta.
\]
The norm $H^{-1}_k$ is the dual of $H^1_k(\mathbb R)$.

For $j=0,1$, the localized equation
\eqref{h5c:step3-localized-equation} is
\begin{equation}\label{eq:localized-J0}
\begin{split}
 (-\partial_v^2+k^2+V_k)(\chi_kD_w^jg_{\pm,k})
 ={}&-\chi_k[\partial_v^j(\xi_0q^2H)]_k\\
 &+\chi_k\sum_{i=0}^j\binom ji
 \left[\operatorname{div}
   \left(D_w^{j-i}(\mathsf G-I)\nabla D_w^ig_{f,\pm}\right)\right]_k\\
 &-\sum_{i=0}^j\binom ji
 D_w^{j-i}(\xi_0b-1)V_k\chi_kD_w^ig_{\pm,k}\\
 &-2\chi_k'\partial_vD_w^jg_{\pm,k}
                         -\chi_k''D_w^jg_{\pm,k},
\end{split}
\end{equation}
where $\chi_k'=\partial_v\chi_k$ and
$\chi_k''=\partial_v^2\chi_k$.

We first estimate the source terms on an annulus of radius $q$.
The local dual norms below are taken after multiplication by
a smooth cutoff supported in $q/2<q'<2q$, with uniformly
bounded $v$ derivatives.
Angular projection and \eqref{cubic:dyadic-energy} imply
\[
 \|D_w^ig_{\pm,k}\|_{H^1_k(q/4<q'<4q)}
 \le Ck^{-2}\sum_{a=0}^2
       \|\nabla\partial_\theta^aD_w^ig_{f,\pm}\|_
                                  {L^2(q/4<q'<4q)}
 \le Ck^{-3/2}qM.
\]
Near $q=q_k$, use \eqref{h5c:collar-input} in the last
inequality. Since
$[F]_k=-k^{-2}[\partial_\theta^2F]_k$ and
$\|F\|_{H^{-1}_k}\le k^{-1}\|F\|_2$, Step~1,
\eqref{cubic:dyadic-energy}, and \eqref{eq:cylinder-form}
imply, for $0\le i\le j\le1$,
\[
\begin{aligned}
 \|[\partial_v^j(\xi_0q^2H)]_k\|_{H^{-1}_k}
 &\le Ck^{-3}q^3M,\\
 \left\|
 \left[\operatorname{div}
   \left(D_w^{j-i}(\mathsf G-I)\nabla D_w^ig_{f,\pm}\right)\right]_k
                                      \right\|_{H^{-1}_k}
 &\le Ck^{-2}q^2
       \sum_{a=0}^2\|\nabla\partial_\theta^aD_w^ig_{f,\pm}\|_2
 \le Ck^{-3/2}q^3M,\\
 \|D_w^{j-i}(\xi_0b-1)V_kD_w^ig_{\pm,k}\|_{H^{-1}_k}
 &\le Cq^2\|D_w^ig_{\pm,k}\|_{H^1_k(q/4<q'<4q)}
 \le Ck^{-3/2}q^3M.
\end{aligned}
\]
Thus the critical potential error is controlled by the
already proved annular estimate.

Apply Lemma~\ref{spectral:weighted-normal-family} with weight
$e^{5v/4}(1+e^{2v})$, which is allowed since
\[
 \frac54<k,\qquad \frac{13}4<\sqrt{k^2+8},\qquad k\ge2.
\]
The logarithmic derivative of this weight is bounded, and
\[
 \left\|q^3\frac{e^{-5v/4}}{1+e^{2v}}\right\|_{L^2(dv)}
 =\delta^3
   \left(\int_{\mathbb R}
      \frac{e^{7v/2}}{(1+e^{2v})^2}\,dv\right)^{1/2}
 \le C\delta^3.
\]
Consequently, the first three lines of \eqref{eq:localized-J0},
after multiplication by $e^{-5v/4}/(1+e^{2v})$, have total
$H^{-1}_k$ norm at most $Ck^{-3/2}\delta^3M$.

For the last line, \eqref{h5c:collar-input} implies
\[
 \|D_w^jg_{\pm,k}\|_{L^2(q_k/2<q<q_k)}
 \le Ck^{-3}M,\qquad
 \sup_{q_k/2<q<q_k}\frac{e^{-5v/4}}{1+e^{2v}}
 \le C\delta^{13/4}k^{13/8}.
\]
Using
\[
 -2\chi_k'\partial_vD_w^jg_{\pm,k}
                         -\chi_k''D_w^jg_{\pm,k}
 =-2\partial_v(\chi_k'D_w^jg_{\pm,k})
                         +\chi_k''D_w^jg_{\pm,k},
\]
integration by parts bounds its weighted dual norm by
\[
 C\delta^{13/4}k^{-11/8}M
 =C\delta^3k^{-3/2}(\delta\sqrt{k})^{1/4}M
 \le C\delta^3k^{-3/2}M.
\]
The weighted inverse estimate therefore proves
\begin{equation}\label{cubic:weighted-higher}
 \left\|\frac{e^{-5v/4}}{1+e^{2v}}
                 \chi_kD_w^jg_{\pm,k}\right\|_{H^1_k(\mathbb R)}
 \le Ck^{-3/2}\delta^3M,\qquad j=0,1.
\end{equation}

The trace inequality
$\|u\|_\infty\le Ck^{-1/2}\|u\|_{H^1_k}$ implies
\begin{equation}\label{cubic:weighted-pointwise}
 |\chi_kD_w^jg_{\pm,k}|
 \le Ck^{-2}\delta^3 e^{5v/4}(1+e^{2v})M,
 \qquad j=0,1.
\end{equation}
At $v=0$, $\chi_k=1$, and
$(\delta\partial_\delta)^j g_{\pm,k}(\delta,kb(\delta))
=D_w^jg_{\pm,k}|_{v=0}$.
This proves \eqref{cubic:higher-trace}.

On $|v|<2$, the condition $c_1<e^{-6}c_0$ implies
$\chi_k=1$. The definition of $U_k$, the product rule, and
the bounded logarithmic derivative of $\vartheta$ imply
\[
 \|U_k\|_{H^1_k((-2,2))}
 +\|\partial_wU_k\|_{H^1_k((-2,2))}
 \le C\delta^{-3}\sum_\pm\sum_{j=0}^1
       \|D_w^jg_{\pm,k}\|_{H^1_k((-2,2))}
 \le Ck^{-3/2}M.
\]
This proves \eqref{cubic:profile-core}.
Finally, \eqref{cubic:weighted-pointwise} and the definition
of $G_k$ imply
\[
 |G_k(q,\delta)|
 \le Ck^{-2}M
 \begin{cases}
 q^{13/4}\delta^{-1/4},&\delta<q/2,\\
 q^{5/4}\delta^{7/4},&\delta>2q.
 \end{cases}
\]
Thus \eqref{cubic:profile-outer} and the amplitude bound in
\eqref{cubic:profile-inner} are established.

\medskip\noindent
\emph{Step 4: the derivative bound in
\eqref{cubic:profile-inner}.}
Fix $v<-\log2$ and set $q=\delta e^v$.
All local norms in this step are taken over
$|v'-v|<1/4$; there
$q'=\delta e^{v'}<2\delta/3$ and $\chi_k=1$.
It follows from \eqref{cubic:weighted-higher} that
\[
 \|\partial_vg_{\pm,k}\|_2+k\|g_{\pm,k}\|_2
 \le Ck^{-3/2}q^{5/4}\delta^{7/4}M.
\]

Taking the $k$th coefficient of
\eqref{h5c:full-center-equation}, we obtain
\[
\begin{aligned}
 &[\mathsf G_{vv}\partial_vg_{f,\pm}
       +\mathsf G_{v\theta}\partial_\theta g_{f,\pm}]_k\\
 &\qquad=\partial_vg_{\pm,k}
   +[(\mathsf G_{vv}-1)\partial_vg_{f,\pm}
       +\mathsf G_{v\theta}\partial_\theta g_{f,\pm}]_k,\\
 &\partial_v[\mathsf G_{vv}\partial_vg_{f,\pm}
       +\mathsf G_{v\theta}\partial_\theta g_{f,\pm}]_k\\
 &\qquad=k^2g_{\pm,k}+\xi_0bV_kg_{\pm,k}
                      +[\xi_0q^2H]_k\\
 &\qquad\quad
   -ik[\mathsf G_{v\theta}\partial_vg_{f,\pm}
       +(\mathsf G_{\theta\theta}-1)
                              \partial_\theta g_{f,\pm}]_k.
\end{aligned}
\]
Here $|V_k|\le C$. Two angular integrations by parts and
\eqref{cubic:dyadic-energy} bound the metric terms in these
two identities by $Ck^{-3/2}q^3M$ and
$Ck^{-1/2}q^3M$, respectively.
Since $q^3\le q^{5/4}\delta^{7/4}$, it follows that
\[
\begin{aligned}
 &\|[\mathsf G_{vv}\partial_vg_{f,\pm}
       +\mathsf G_{v\theta}\partial_\theta g_{f,\pm}]_k\|_2
 \le Ck^{-3/2}q^{5/4}\delta^{7/4}M,\\
 &\|\partial_v[\mathsf G_{vv}\partial_vg_{f,\pm}
       +\mathsf G_{v\theta}\partial_\theta g_{f,\pm}]_k\|_2
 \le Ck^{-1/2}q^{5/4}\delta^{7/4}M.
\end{aligned}
\]
The one-dimensional trace inequality therefore implies
\[
 \left|[\mathsf G_{vv}\partial_vg_{f,\pm}
       +\mathsf G_{v\theta}\partial_\theta g_{f,\pm}]_k(v)\right|
 \le Ck^{-1}q^{5/4}\delta^{7/4}M.
\]

To recover $\partial_vg_{\pm,k}(v)$, we estimate the metric
correction point-wisely.
On the chosen interval,
$\sup_{n\ne0}|V_n|/n^2\le C$.
Solving \eqref{h5c:full-center-equation} for
$\partial_v^2g_{f,\pm}$ and using
\eqref{cubic:dyadic-energy}, we obtain
\[
 \|\nabla g_{f,\pm}\|_{L^2_{v'}H^2_\theta}
 +\|\partial_{v'}\nabla g_{f,\pm}\|_{L^2_{v'}H^1_\theta}
 \le Ck^{1/2}qM.
\]
The trace inequality then implies
\[
\begin{aligned}
 \|\nabla g_{f,\pm}(v)\|_{H^{3/2}_\theta}^2
 &\le C\|\nabla g_{f,\pm}\|_{L^2_{v'}H^2_\theta}
 \left(\|\partial_{v'}\nabla g_{f,\pm}\|_{L^2_{v'}H^1_\theta}
             +\|\nabla g_{f,\pm}\|_{L^2_{v'}H^1_\theta}\right)\\
 &\le Ckq^2M^2.
\end{aligned}
\]
Since the angular derivatives of $\mathsf G-I$ are $O(q^2)$,
\[
\begin{aligned}
 &\left|[(\mathsf G_{vv}-1)\partial_vg_{f,\pm}
       +\mathsf G_{v\theta}\partial_\theta g_{f,\pm}]_k(v)\right|\\
 &\qquad\le Ck^{-3/2}q^2
                    \|\nabla g_{f,\pm}(v)\|_{H^{3/2}_\theta}
 \le Ck^{-1}q^3M.
\end{aligned}
\]
Consequently,
\[
 |\partial_vg_{\pm,k}(v)|
 \le Ck^{-1}q^{5/4}\delta^{7/4}M.
\]
Together with \eqref{cubic:weighted-pointwise}, this proves
\[
 |(D_w-\partial_v)g_{\pm,k}|
 \le Ck^{-1}q^{5/4}\delta^{7/4}M.
\]

Finally, $\delta\partial_\delta|_q=D_w-\partial_v$.
Differentiating the definition of $G_k$ at fixed $q$ gives
\[
\begin{aligned}
 \delta\partial_\delta G_k
 ={}&\chi_k(q)\vartheta\!\left(\frac{\delta\sqrt{k}}{c_1}\right)
       (D_w-\partial_v)(g_{+,k}-g_{-,k})\\
 &+\chi_k(q)
   \left[\delta\partial_\delta
       \vartheta\!\left(\frac{\delta\sqrt{k}}{c_1}\right)\right]
                                      (g_{+,k}-g_{-,k}).
\end{aligned}
\]
Both cutoff coefficients are uniformly bounded. Hence
\[
 |\delta\partial_\delta G_k(q,\delta)|
 \le Ck^{-1}q^{5/4}\delta^{7/4}M,
 \qquad \delta>2q,
\]
which completes \eqref{cubic:profile-inner}.

For $\lambda=k(1+\delta^2/8)$, we use
$D_w=q\partial_q+k\delta^2\partial_\lambda/4$ and
$V_k=8e^{2v}/(1+e^{2v})$ and the same estimates apply. The Lemma is then proved.
\end{proof}

\subsection{The critical trace for the first angular modes}
\label{subsec:cubic-first-trace}

For $k=\pm1$, the normal exponent equals three. We obtain the
cubic trace bound from an exact homogeneous test function.
In this subsection and Subsection~\ref{subsec:cubic-first-jump},
the constants are uniform for $0<\delta<q_0/8$, where $q_0$
is fixed and sufficiently small such that the cutoff in
\eqref{h5c:critical}, with $k=1$, equals one on $q<q_0$.

\begin{lemma}\label{cubic:first-mode}
For $j=0,1,2$,
\begin{equation}\label{cubic:first-trace}
 \left|(\delta\partial_\delta)^j
     g_{\pm,1}(\delta,b(\delta))\right|
       \le C\delta^3M.
\end{equation}
The same estimate holds with $1,b(\delta)$ replaced by
$-1,-b(\delta)$, and at the other elliptic point.
\end{lemma}

\begin{proof}
We treat the first angular mode at one center; either boundary
value is allowed. Set
\[
 \lambda=b(\delta)=1-\delta^2/8,\qquad
 v=\log(q/\delta),\qquad w=\log\delta,\qquad
 D_w=q\partial_q-\frac{\delta^2}{4}\partial_\lambda.
\]
Thus $D_w=\partial_w$ with $v$ fixed, and $q=\delta$ is $v=0$.
Write $g=g_{f,\pm}(\lambda)$ and $H=(\Delta+1)f$ in these
coordinates, and define the differential operator
\[
 P=\mathsf G_{vv}\partial_v+\mathsf G_{v\theta}\partial_\theta.
\]
All cylinder norms use $dv\,d\theta$; norms of individual
Fourier coefficients use $dv$. Fourier subscripts denote angular
coefficients, and $(0,w)$ denotes evaluation at $q=\delta$.

\emph{Step 1: proof assuming the estimates below.}
On the chosen cell, define
\[
 \phi=(L_0-i\lambda)(x-iy)
      =i\bigl(\sin x-i\sin y-\lambda(x-iy)\bigr).
\]
The first coefficient of $\phi$ satisfies
\[
 \phi_1(v,w)=i[b(\delta e^v)-b(\delta)](x-iy)_1(\delta e^v),
 \qquad \phi_1(0,w)=0.
\]
It follows from \eqref{h5c:full-center-equation} and the definition
of $P$ that
\begin{equation}\label{cubic:green-flux-equation}
 \partial_v(Pg)
 =-\partial_\theta\!\left(
    \mathsf G_{v\theta}\partial_vg
       +\mathsf G_{\theta\theta}\partial_\theta g\right)
   +\xi_0(q)q^2\bigl(H-T_\pm(\lambda)g\bigr).
\end{equation}
The identity $(\Delta+T_\pm(\lambda))\phi=0$, verified in
Step~2 below, implies
\[
 \partial_v(P\phi)
 =-\partial_\theta\!\left(
    \mathsf G_{v\theta}\partial_v\phi
       +\mathsf G_{\theta\theta}\partial_\theta\phi\right)
   -\xi_0(q)q^2T_\pm(\lambda)\phi.
\]

For $q<\delta$, we have $b(q)>b(\delta)=\lambda$.
The angular multipliers
\[
 \frac{nb(q)}{nb(q)-\lambda},\qquad n\in\mathbb Z,
\]
are therefore real and nonsingular. Consequently,
\[
 \int_{\mathbb T}
 \left((T_\pm(\lambda)g)\overline\phi
       -g\,\overline{T_\pm(\lambda)\phi}\right)\,d\theta=0.
\]
Recall that
\[
 \Delta g+T_\pm(\lambda)g=H,
 \qquad
 \Delta\phi+T_\pm(\lambda)\phi=0.
\]
Subtracting the equation for $\overline\phi$ multiplied by $g$
from the equation for $g$ multiplied by $\overline\phi$,
and applying Green's identity, we obtain
\begin{equation}\label{cubic:green-before-limit}
 \begin{aligned}
 \frac1{2\pi}\int_{q'<\delta e^v}H\overline\phi\,dx\,dy
 &=\frac1{2\pi}\int_{q'<\delta e^v}
       \left((\Delta g)\overline\phi
                     -g\,\overline{\Delta\phi}\right)\,dx\,dy\\
 &=\frac1{2\pi}\int_{\mathbb T}
       \left((Pg)\overline\phi-g\,\overline{P\phi}\right)
                                      (v,\theta,w)\,d\theta\\
 &=\sum_{n\equiv1\!\!\pmod4}
       \left((Pg)_n\overline{\phi_n}
                 -g_n\overline{(P\phi)_n}\right)(v,w).
 \end{aligned}
\end{equation}
Here the boundary term uses
$\partial_{\boldsymbol n}h\,dS=(Ph)\,d\theta$
on $q'=\delta e^v$, where $\boldsymbol n$ is the outward unit
normal. The last equality follows from Parseval's identity
and the fact that $\phi$ and $P\phi$ contain only angular
modes $n\equiv1\pmod4$. Indeed, by the rotation identity established in the proof of
Lemma~\ref{lem:geometry},
\[
 X(q,\theta+\pi/2)=RX(q,\theta),
 \qquad R(x,y)=(y,-x).
\]
Since $\phi\circ R=i\phi$ and
$\mathsf G(q,\theta+\pi/2)=\mathsf G(q,\theta)$,
the definition of $P$ implies
\[
 \phi(q,\theta+\pi/2)=i\phi(q,\theta),
 \qquad
 (P\phi)(q,\theta+\pi/2)=i(P\phi)(q,\theta).
\]
Taking angular Fourier coefficients, we obtain
\[
 (e^{in\pi/2}-i)\phi_n
 =(e^{in\pi/2}-i)(P\phi)_n=0.
\]
Thus both coefficients vanish unless $n\equiv1\pmod4$.

The first coefficient of $Pg$ need not have a finite trace
at $v=0$. More precisely, \eqref{cubic:first-flux-log},
proved in Step~4 below, states that
\[
 (Pg)_1(v,w)
 =4\xi_0(\delta)b(\delta)g_1(0,w)\log(v\mp i0)+r_0(v,w),
 \qquad r_0(\cdot,w)\in H^1(-1/2,0).
\]
On the other hand,
\[
 \phi_1(v,w)
 =\frac{i\delta^2}{8}(1-e^{2v})
                     (x-iy)_1(\delta e^v),
 \qquad
 |\phi_1(v,w)|\le C\delta^3|v|,
 \quad -1/2<v<0.
\]
Since $H^1(-1/2,0)\hookrightarrow L^\infty(-1/2,0)$
and $v\log(v\mp i0)\to0$, it follows that
\[
 \lim_{v\to0^-}
       (Pg)_1(v,w)\overline{\phi_1(v,w)}=0.
\]
The limits are justified in Step~3.
Taking $v\to0^-$ in \eqref{cubic:green-before-limit}
and separating the term with $n=1$, we obtain
\begin{equation}\label{cubic:exact-green}
 \begin{aligned}
 g_1(0,w)\overline{(P\phi)_1(0,w)}
 ={}&-\frac1{2\pi}\int_{q<\delta}H\overline\phi\,dx\,dy\\
 &+\sum_{\substack{n\equiv1\!\!\pmod4\\n\ne1}}
   \left(\overline{\phi_n}(Pg)_n
              -g_n\overline{(P\phi)_n}\right)(0,w).
 \end{aligned}
\end{equation}
Thus the boundary term involving $(Pg)_1$ vanishes, whereas
the coefficient multiplying $g_1$ satisfies
\[
 (P\phi)_1(0,w)=-\frac{i\delta^3}{4}+O(\delta^5),
\]
as proved in Step~2. In particular, this coefficient is nonzero
for sufficiently small $\delta$. Step~4 also justifies
\eqref{cubic:exact-green} after up to two $w$ derivatives.

Assume for the moment that, for $0\le j\le2$,
\begin{equation}\label{cubic:exact-test-traces}
 \begin{aligned}
 \left|\partial_w^j\left[\delta^{-3}(P\phi)_1(0,w)+i/4\right]\right|
 &\le C\delta^2,\\
 \bigg(\sum_{\substack{n\equiv1\!\!\pmod4\\n\ne1}}
   \left[|\partial_w^j\phi_n(0,w)|^2
        +|\partial_w^j(P\phi)_n(0,w)|^2\right]\bigg)^{1/2}
 &\le C\delta^3,
 \end{aligned}
\end{equation}
\begin{equation}\label{cubic:first-volume}
 \left|\partial_w^j\int_{q<\delta}H\overline\phi\,dx\,dy\right|
 \le C\delta^6M,
\end{equation}
and
\begin{equation}\label{cubic:noncritical-traces}
 \bigg(\sum_{\substack{n\equiv1\!\!\pmod4\\n\ne1}}
   \left[|\partial_w^jg_n(0,w)|^2
        +|\partial_w^j(Pg)_n(0,w)|^2\right]\bigg)^{1/2}
 \le C\delta^3M.
\end{equation}
It follows from \eqref{cubic:exact-green}, the product rule,
and Cauchy--Schwarz that
\[
 \left|\partial_w^j\left[g_1(0,w)
                   \overline{(P\phi)_1(0,w)}\right]\right|
 \le C\delta^6M,\qquad 0\le j\le2.
\]
The first line of \eqref{cubic:exact-test-traces} implies
\[
 \left|\partial_w^j
       \left(\frac1{\overline{(P\phi)_1(0,w)}}\right)\right|
 \le C\delta^{-3},\qquad 0\le j\le2.
\]
Applying the product rule once more, we obtain
\[
 |\partial_w^jg_1(0,w)|\le C\delta^3M,\qquad 0\le j\le2.
\]
Since $g_1(0,w)=g_{f,\pm,1}(\delta,b(\delta))$, this proves
\eqref{cubic:first-trace}. It remains
to establish the estimates and identity used above.

\emph{Step 2: the test function and the volume estimate.}
The identities
\[
 \Delta(x-iy)=0,\qquad L_0(x-iy)=i(\sin x-i\sin y),\qquad
 \Delta L_0(x-iy)=-L_0(x-iy)
\]
imply
\[
 T(i\lambda\pm0)\phi=L_0(x-iy),\qquad
 (\Delta+T(i\lambda\pm0))\phi=0.
\]
Indeed, $\phi_n=i(nb(q)-\lambda)(x-iy)_n$, so the first
identity follows by cancellation in the angular multiplier
$nb(q)/(nb(q)-\lambda\mp i0)$, in the distributional sense.

By Lemma~\ref{geom:metric},
\[
 (x-iy)_1=q+O(q^3),\qquad
 b(q)-b(\delta)=\frac{\delta^2-q^2}{8},\qquad
 \mathsf G-I=O(q^2).
\]
By Lemma~\ref{geom:metric},
\[
 \left\|(q\partial_q)^m
       (x-iy-qe^{i\theta})\right\|_{C^2_\theta}
 \le Cq^3,\qquad 0\le m\le3.
\]
Since $\phi=(b(q)\partial_\theta-i b(\delta))(x-iy)$,
we have the exact decomposition
\[
 \phi
 =\frac{i}{8}q(\delta^2-q^2)e^{i\theta}
  +\bigl(b(q)\partial_\theta-i b(\delta)\bigr)
                       (x-iy-qe^{i\theta}).
\]
Differentiation in $w$ with $v$ fixed acts as
$q\partial_q+\delta\partial_\delta$ on expressions written
in $(q,\delta,\theta)$. Hence the product rule implies
\[
 \sum_{a+b\le1}
 \left\|D_w^j\partial_v^a\partial_\theta^b\phi
                                      \right\|_{L^\infty_\theta}
 \le C(\delta^2q+q^3),\qquad 0\le j\le2.
\]
In particular,
\[
 |D_w^j\phi|\le C\delta^2q,
 \qquad 0<q<\delta,\quad 0\le j\le2,
\]
and, at $v=0$, the metric bounds imply
\begin{equation}\label{cubic:test-full-traces}
 \|\partial_w^j\phi(0,w)\|_{L^2_\theta}
 +\|\partial_w^j(P\phi)(0,w)\|_{L^2_\theta}
 \le C\delta^3,\qquad 0\le j\le2.
\end{equation}

To identify the leading first coefficient, use
$\phi_1=i[b(q)-b(\delta)](x-iy)_1(q)$ and
$q b'(q)=-q^2/4$. Since $b(q)-b(\delta)=0$ at $q=\delta$,
\[
 \begin{aligned}
 (P\phi)_1(0,w)
 ={}&-\frac{i\delta^2}{4}(x-iy)_1(\delta)\\
 &+\left[
   (\mathsf G_{vv}-1)\partial_v\phi
       +\mathsf G_{v\theta}\partial_\theta\phi
   \right]_1(0,w).
 \end{aligned}
\]
The preceding gradient estimate and
$D_w^j(\mathsf G-I)=O(q^2)$ imply
\[
 \left|
 \partial_w^j\left[
   (\mathsf G_{vv}-1)\partial_v\phi
       +\mathsf G_{v\theta}\partial_\theta\phi
 \right]_1(0,w)
 \right|\le C\delta^5,\qquad 0\le j\le2.
\]
Moreover, the coordinate expansion implies
\[
 \left|\partial_w^j
   \left[\delta^{-1}(x-iy)_1(\delta)-1\right]\right|
 \le C\delta^2,\qquad 0\le j\le2.
\]
It follows from these two estimates and the product rule that
\[
 \left|\partial_w^j
   \left[\delta^{-3}(P\phi)_1(0,w)+\frac{i}{4}\right]\right|
 \le C\delta^2,\qquad 0\le j\le2.
\]
This proves the first line of \eqref{cubic:exact-test-traces}.
The second line follows from \eqref{cubic:test-full-traces}
and Parseval's identity, restricting the sum to
$n\equiv1\pmod4$, $n\ne1$.

The function $H$ vanishes at the elliptic center, as proved in
Lemma~\ref{quad:weighted-center}.
The linear Taylor term of $H$ and
\eqref{cubic:source-annulus} imply
\[
 \sum_{m=0}^2\int_{q<\delta}q|\partial_v^mH|\,dx\,dy
 \le C\sum_{a\le\delta}a^3\sum_{m=0}^2
       \|\partial_v^mH\|_{L^2(q\sim a;dv\,d\theta)}
 \le CM\sum_{a\le\delta}a^4\le C\delta^4M,
\]
where the sums are over dyadic radii. In the fixed domain
$v<0$, the area measure is $\xi_0(q)q^2\,dv\,d\theta$, so
\[
 \partial_w^j\int_{q<\delta}H\overline\phi\,dx\,dy
 =\sum_{a+b+c=j}\frac{j!}{a!b!c!}
  \int_{v<0}(\partial_v^aH)\overline{D_w^b\phi}
                    \partial_v^c(\xi_0q^2)\,dv\,d\theta.
\]
Since $|\partial_v^c(\xi_0q^2)|\le Cq^2$, the preceding bound
and $|D_w^b\phi|\le C\delta^2q$ prove
\eqref{cubic:first-volume}.

\emph{Step 3: the traces
\eqref{cubic:noncritical-traces}.}
By the rotation identity in the proof of
Lemma~\ref{lem:geometry},
\[
 \mathsf G(q,\theta+\pi/2)=\mathsf G(q,\theta),
\]
and the area density $\xi_0(q)q$ is radial. Moreover, $T$
acts diagonally on angular Fourier coefficients, and the
coefficients of $D_w$ are independent of $\theta$.
Consequently, the orthogonal projection
\[
 u(q,\theta)\longmapsto
 \frac14\sum_{m=0}^3i^{-m}u(q,\theta+m\pi/2)
 =\sum_{n\equiv1\!\!\pmod4}u_n(q)e^{in\theta}
\]
commutes with $\Delta,T,D_w$, and $P$.

Apply this projection to $g$ and $H$, and with a slight abuse of notation denote
the resulting functions still by $g$ and $H$. Their equation is
unchanged. The projection preserves $g_1$ and all coefficients
appearing in \eqref{cubic:noncritical-traces}. It also preserves
$\int_{q<\delta}H\overline\phi\,dx\,dy$, by angular
orthogonality, since $\phi$ contains only modes
$n\equiv1\pmod4$ and the area density is radial.

For functions with this angular support, write
\[
 u_{\ne1}=(I-\Pi_1)u
 =\sum_{\substack{n\equiv1\!\!\pmod4\\n\ne1}}
                         u_n e^{in\theta},
 \qquad \Pi_1u=u_1e^{i\theta}.
\]
Every index in this sum satisfies $|n|\ge3$, so Parseval's
identity implies
\[
 \|u_{\ne1}\|_2\le\frac13\|\partial_\theta u_{\ne1}\|_2.
\]

Set $h_j=D_w^jg$ for $0\le j\le2$. Since $D_w$ commutes
with $\partial_v$ and $\partial_\theta$, the product rule
implies
\[
 \begin{aligned}
 D_w^j(Pg)
 ={}&Ph_j+\sum_{\ell=0}^{j-1}\binom j\ell
 \left[
  (D_w^{j-\ell}\mathsf G_{vv})\partial_vh_\ell
  +(D_w^{j-\ell}\mathsf G_{v\theta})
                                  \partial_\theta h_\ell
 \right],
 \end{aligned}
\]
where the sum is empty for $j=0$.

Unless another domain is specified, $\|\cdot\|_2$ in this
step denotes the norm on
\[
 L^2\!\left(
 (-\infty,\log(q_0/\delta))\times\mathbb T;\,dv\,d\theta
 \right),
\]
with $w$ fixed. For an individual angular coefficient,
the norm uses only $dv$ on the same interval.

\smallskip\noindent
\emph{(a) The equations for the modes $n\equiv1\pmod4$, $n\ne1$.}
It follows from \eqref{h5c:full-center-equation} that
\begin{equation}\label{cubic:conormal-system}
 \begin{aligned}
 g_v&=\mathsf G_{vv}^{-1}
                  (Pg-\mathsf G_{v\theta}g_\theta),\\
 \partial_v(Pg)
 &=-\partial_\theta\!\left(
       \mathsf G_{v\theta}\mathsf G_{vv}^{-1}Pg
       +(\mathsf G_{\theta\theta}
         -\mathsf G_{v\theta}^2/\mathsf G_{vv})g_\theta\right)
       -\xi_0q^2Tg+\xi_0q^2H.
 \end{aligned}
\end{equation}
Apply $I-\Pi_1$ to the equation for $h_j=D_w^jg$.
In the resulting equation for $(h_j)_{\ne1}$, define
the two terms explicitly by
\[
 \begin{aligned}
 \mathcal F_j
 &=D_w^j(Pg)-(h_j)_v\\
 &=\sum_{\ell=0}^j\binom j\ell
   \left[\partial_v^{j-\ell}(\mathsf G_{vv}-1)(h_\ell)_v
       +\partial_v^{j-\ell}\mathsf G_{v\theta}(h_\ell)_\theta\right],\\
 \mathcal H_j
 &=\partial_\theta\sum_{\ell=0}^j\binom j\ell
   \left[\partial_v^{j-\ell}\mathsf G_{v\theta}(h_\ell)_v
       +\partial_v^{j-\ell}(\mathsf G_{\theta\theta}-1)
                                      (h_\ell)_\theta\right]\\
 &\quad+D_w^j(\xi_0q^2Tg)-D_w^j(\xi_0q^2H).
 \end{aligned}
\]
For $n\equiv1\pmod4$, $n\ne1$, we have
\begin{equation}\label{cubic:noncritical-equation}
 -\partial_v^2(h_j)_n+n^2(h_j)_n
 =\partial_v(\mathcal F_j)_n+(\mathcal H_j)_n,
 \qquad (D_w^j(Pg))_n=(h_j)_n'+(\mathcal F_j)_n.
\end{equation}
Moreover, for $0\le a,b\le2$,
\begin{equation}\label{cubic:noncritical-multiplier}
 \left|D_w^a\partial_v^b
      \left(q^2\frac{nb(q)}{nb(q)-\lambda}\right)\right|
 \le C_{a,b}q^2,
 \qquad n\equiv1\pmod4,\quad n\ne1.
\end{equation}

\smallskip\noindent
\emph{(b) The weight and the source estimates.}
Choose the fixed radius $q_0$ so that
\[
 4c_1<q_0<c_0/2,
\]
where $0<\delta<c_1$. In particular, the cutoff in
\eqref{h5c:critical}, with $k=1$, equals one throughout
$q<q_0$. Set
\[
 \xi_\delta(q)=q^{5/2}(q+\delta)^{1/2},
 \qquad
 \partial_v\log\xi_\delta
 =\frac52+\frac{q}{2(q+\delta)}\in[5/2,3].
\]
Choose a smooth radial cutoff $\chi$ satisfying
\[
 \chi=1\quad(q\le q_0/2),\qquad
 \chi=0\quad(q\ge3q_0/4),\qquad
 |\partial_v^m\chi|\le C_m,\quad 0\le m\le3.
\]

Define
\[
 h_{j,\mathrm h}
 =\sum_{\substack{n\equiv1\!\!\pmod4\\|n|\ge5}}
                       (h_j)_n e^{in\theta},
\]
and
\[
 \begin{aligned}
 \mathcal B_j
 &=\sup_{v<\log(q_0/\delta)}
   \xi_\delta(\delta e^v)^{-1}
   \left(
      |(h_j)_{-3}(v,w)|
      +|(D_w^j(Pg))_{-3}(v,w)|
   \right),\\
 \mathcal E_j
 &=\sum_{a=0}^1
   \left\|\xi_\delta^{-1}
     \bigl(\partial_v\partial_\theta^ah_{j,\mathrm h},
           \partial_\theta^{a+1}h_{j,\mathrm h}\bigr)
   \right\|_{L^2(q<q_0;dv\,d\theta)}.
 \end{aligned}
\]
The separate norm for $n=-3$ avoids integrating the square
of a constant in $v$ when a term proportional to $q^3$ is
divided by $\xi_\delta\sim q^3$ on $q>\delta$.
For $|n|\ge5$, the weighted energy has the strict bound
\[
 n^2-|\partial_v\log\xi_\delta|^2
 \ge \frac{16}{25}n^2.
\]
The desired conclusion is
\begin{equation}\label{cubic:noncritical-weighted}
 \mathcal B_j+\mathcal E_j\le C_{q_0}M,
 \qquad 0\le j\le2.
\end{equation}
Here $C_{q_0}$ is independent of $\delta$. In the perturbation
estimates below, the constants multiplying $q_0^2$ are
independent of both $q_0$ and $\delta$.

\smallskip\noindent
\emph{The outer annulus.}
It follows from \eqref{h5c:mild} and the angular
Poincar\'e inequality that
\[
 \sum_{j=0}^2\sum_{a=0}^2
 \|\partial_\theta^ah_j\|_
 {H^1(q_0/2<q<q_0;dv\,d\theta)}
 \le CM.
\]
Differentiating \eqref{cubic:conormal-system} and using
\eqref{cubic:noncritical-multiplier}, we obtain
\[
 \sum_{j=0}^2
 \|(I-\Pi_1)D_w^j(Pg)\|_
 {H^1_vL^2_\theta(q_0/2<q<q_0)}
 \le C_{q_0}M.
\]
The corresponding $v$ interval has length $\log2$.
The embedding $H^1_vL^2_\theta\hookrightarrow L^\infty_vL^2_\theta$
therefore bounds the traces of $(I-\Pi_1)h_j$ and
$(I-\Pi_1)D_w^j(Pg)$ on this annulus.
Since $\xi_\delta(q)\sim q_0^3$ there, its contributions
to $\mathcal B_j,\mathcal E_j$, and the weighted terms
containing derivatives of $\chi$ are bounded by $C_{q_0}M$.

\smallskip\noindent
\emph{The metric terms involving the first mode.}
For $\beta=1/2$, the weight in \eqref{h5c:critical} is
\[
 W(v)=e^{v/2}(1+e^{2v}),
 \qquad |\partial_v\log W|\le5/2.
\]
Consequently, \eqref{h5c:critical} implies
\[
 \left\|W^{-1}
   \bigl((h_j)_1,\partial_v(h_j)_1\bigr)
 \right\|_{L^2(q<q_0;dv)}
 \le C\delta^{5/2}M,\qquad 0\le j\le2.
\]
The multiplier relating this estimate to $\xi_\delta$ satisfies
\[
 \frac{\delta^{5/2}q^2W(\log(q/\delta))}{\xi_\delta(q)}
 =\frac{q^2+\delta^2}{(q+\delta)^{1/2}}
 \le(q+\delta)^{3/2}
 \le Cq_0^{3/2}.
\]
It follows that
\begin{equation}\label{cubic:resonant-metric-source}
 \left\|\xi_\delta^{-1}q^2
   \bigl((h_j)_1,\partial_v(h_j)_1\bigr)
 \right\|_{L^2(q<q_0;dv)}
 \le Cq_0^{3/2}M,\qquad 0\le j\le2.
\end{equation}

\smallskip\noindent
\emph{The fixed source.}
Let $\ell$ be the linear Taylor polynomial of the projected
function $H$ at the center. Its rotation symmetry implies
\[
 \ell(x,y)=\partial_xH(0,0)(x-iy).
\]
It follows from Lemma~\ref{geom:metric} that
\[
 \sum_{r+s\le3}
 \|(I-\Pi_1)\partial_v^r\partial_\theta^s\ell
                                      \|_{L^2_\theta}
 \le Cq^3M.
\]
Moreover, \eqref{cubic:source-annulus}, with $\alpha=1/2$,
implies
\begin{equation}\label{cubic:source-remainder}
 \sum_{r+s\le3}
 \|\partial_v^r\partial_\theta^s(H-\ell)\|_
 {L^2(a/2<q<2a;dv\,d\theta)}
 \le Ca^{3/2}M,\qquad 0<a\le q_0.
\end{equation}
Since $\xi_0q^2H$ is independent of the spectral parameter,
\[
 D_w^j(\xi_0q^2H)=\partial_v^j(\xi_0q^2H).
\]
Thus, for $0\le j\le2$ and $0\le m\le1$, the product rule
and the preceding estimates imply
\[
 \begin{aligned}
 &\|(I-\Pi_1)\partial_\theta^mD_w^j(\xi_0q^2H)\|_
       {L^2(a/2<q<2a;dv\,d\theta)}\\
 &\qquad\le C(a^5+a^{7/2})M
 \le Ca^{7/2}M.
 \end{aligned}
\]
On each such annulus, $\xi_\delta(q)\sim\xi_\delta(a)$.
Summing over $a=q_0,q_0/2,q_0/4,\ldots$, we have
\[
 \sum_a\frac{a^7}{\xi_\delta(a)^2}
 =\sum_a\frac{a^2}{a+\delta}
 \le\sum_a a
 \le2q_0.
\]
Consequently,
\begin{equation}\label{cubic:fixed-source-weighted}
 \sum_{j=0}^2\sum_{m=0}^1
 \left\|\xi_\delta^{-1}(I-\Pi_1)
       \partial_\theta^mD_w^j(\xi_0q^2H)
 \right\|_{L^2(q<q_0;dv\,d\theta)}
 \le Cq_0^{1/2}M.
\end{equation}
\smallskip\noindent
\emph{(c) The modes $|n|\ge5$.}
After conjugation by $\xi_\delta$ and integration
by parts, with $u=\xi_\delta^{-1}\chi h_{j,\mathrm h}$, we obtain
\[
\begin{aligned}
 \operatorname{Re}\langle\xi_\delta^{-1}
          (-\partial_v^2-\partial_\theta^2)(\xi_\delta u),u\rangle
 &=\int\left(|u_v|^2+|u_\theta|^2
              -|\partial_v\log\xi_\delta|^2|u|^2\right)\\
 &\ge\|u_v\|_2^2+\frac{16}{25}\|u_\theta\|_2^2
 \ge c\|u\|_{H^1}^2.
\end{aligned}
\]
Here $|n|\ge5$. The metric and noncritical potential errors have
relative form bound $Cq_0^2$. Moreover,
\[
 \|\xi_\delta^{-1}q^2((h_j)_{-3},(h_j)_{-3}')\|_2
       \le Cq_0^2\mathcal B_j+Cq_0^4\mathcal E_j+CM.
\]
Use $(h_j)_{-3}'=(D_w^j(Pg))_{-3}-(\mathcal F_j)_{-3}$, absorb its
$O(q_0^2)$ derivative term, and integrate $q^4\,dv$ to obtain $q_0^4/4$.
Testing also after one angular derivative, and using
\eqref{cubic:resonant-metric-source}--\eqref{cubic:source-remainder},
it follows that
\begin{equation}\label{cubic:high-closure}
 \mathcal E_j\le CM+Cq_0^2(\mathcal B_j+\mathcal E_j)
             +Cq_0^2\sum_{\ell<j}(\mathcal B_\ell+\mathcal E_\ell).
\end{equation}
Only three derivatives of $H$ occur, since $j\le2$ and there is
at most one angular commutation.

\smallskip\noindent
\emph{(d) The mode $n=-3$.}
Set $v_0=\log(q_0/\delta)$ and use the Dirichlet Green kernel on $(-\infty,v_0)$:
\[
 K(v,u)=\frac{e^{-3|v-u|}-e^{3(v+u-2v_0)}}6.
\]
For $v<v_0$, the exact localized Green formula is
\[
\begin{split}
 \chi(h_j)_{-3}={}&-\int_{-\infty}^{v_0}K_u(v,u)
             [\chi \mathcal F_j-2\chi'h_j]_{-3}(u)\,du\\
 &+\int_{-\infty}^{v_0}K(v,u)[\chi \mathcal H_j-\chi'\mathcal F_j+\chi''h_j]_{-3}(u)\,du.
\end{split}
\]
Here primes on $\chi$ denote cylinder derivatives. For $v,u<v_0$, the Green kernel satisfies the distributional
identity
\[
 -\partial_v\partial_uK(v,u)
 =-\delta_u(v)
   +\frac32\left(e^{-3|v-u|}+e^{3(v+u-2v_0)}\right),
\]
where $\delta_u$ denotes the Dirac mass at $u$.
Since $\chi=1$ and $\chi'=0$ on $q<q_0/2$,
differentiating the localized Green formula produces
$-(\mathcal F_j)_{-3}(v)$. This term cancels against
$(\mathcal F_j)_{-3}(v)$ in
\[
 (D_w^j(Pg))_{-3}
 =\partial_v(h_j)_{-3}+(\mathcal F_j)_{-3}.
\]
Moreover,
\[
 e^{3(v+u-2v_0)}\le e^{-3|v-u|}.
\]
It follows that, for $q=\delta e^v<q_0/2$,
\begin{equation}\label{cubic:low-kernel-bound}
\begin{aligned}
 |(h_j)_{-3}(v)|+|(D_w^j(Pg))_{-3}(v)|
 &\le C\int_{-\infty}^{v_0}e^{-3|v-u|}
       \left(|(\chi\mathcal F_j)_{-3}(u)|
             +|(\chi\mathcal H_j)_{-3}(u)|\right)\,du\\
 &\quad+C_{q_0}M\xi_\delta(q).
\end{aligned}
\end{equation}
Here the last term bounds the contributions containing
$\chi'$ or $\chi''$, using the annulus estimates in part~(b)
and $q^3\le\xi_\delta(q)$.

We first estimate the terms involving modes other than $1$.
Since
\[
 \frac{\xi_\delta(a)}{\xi_\delta(q)}
 \le
 \begin{cases}
 (a/q)^{5/2},&0<a\le q,\\
 (a/q)^3,&q\le a<q_0,
 \end{cases}
\]
the change of variables $a=\delta e^u$ implies, for $r=1,2$,
\begin{equation}\label{cubic:complementary-kernel-moments}
\begin{aligned}
 &\frac1{\xi_\delta(q)^r}
   \int_0^{q_0}\min(q/a,a/q)^{3r}
          a^{2r}\xi_\delta(a)^r\,\frac{da}{a}\\
 &\qquad\le
 q^{2r}\int_0^q(a/q)^{15r/2}\,\frac{da}{a}
       +\int_q^{q_0}a^{2r}\,\frac{da}{a}=
 \frac{2q^{2r}}{15r}
 +\frac{q_0^{2r}-q^{2r}}{2r}
 \le C_rq_0^{2r}.
\end{aligned}
\end{equation}
The case $r=1$ bounds the contributions controlled by
$\mathcal B_\ell$; the case $r=2$, followed by
Cauchy--Schwarz, bounds those controlled by $\mathcal E_\ell$.
For terms containing $\partial_v(h_\ell)_{-3}$, use
\[
 \partial_v(h_\ell)_{-3}
 =(D_w^\ell(Pg))_{-3}-(\mathcal F_\ell)_{-3},
\]
as in part~(c), and absorb the coefficient of size $Cq_0^2$
multiplying this derivative. It follows from
\eqref{cubic:noncritical-multiplier} and the bounds on the
metric coefficients that the total contribution of these
modes to the weighted supremum in
\eqref{cubic:low-kernel-bound} is at most
\[
 Cq_0^2\sum_{\ell=0}^j(\mathcal B_\ell+\mathcal E_\ell).
\]

For the terms containing the first angular mode, recall from
\eqref{h5c:critical}, with $\beta=1/2$, that
\[
 \left\|
 \frac{e^{-v/2}}{1+e^{2v}}
 \bigl((h_\ell)_1,\partial_v(h_\ell)_1\bigr)
 \right\|_{L^2(v<v_0)}
 \le C\delta^{5/2}M,\qquad 0\le\ell\le2.
\]
Their metric coefficients are $O(q^2)$. The squared kernel
factor in the corresponding Cauchy--Schwarz estimate satisfies
\[
\begin{aligned}
 &\frac{\delta^5}{\xi_\delta(q)^2}
   \int_{-\infty}^{v_0}
      e^{-6|v-u|}\delta^4e^{5u}(1+e^{2u})^2\,du\\
 &\quad=
 \frac{
 q^{-6}\displaystyle\int_0^q a^{10}(\delta^2+a^2)^2\,da
 +q^6\displaystyle\int_q^{q_0}a^{-2}(\delta^2+a^2)^2\,da}
 {q^5(q+\delta)}\\
 &\quad\le
 \frac{(\delta^2+q^2)^2}{11(q+\delta)}
 +\frac{q}{q+\delta}
       \left(\frac{\delta^4}{q}
                    +2\delta^2q_0+\frac{q_0^3}{3}\right)
 \le Cq_0^3.
\end{aligned}
\]
Thus the first-mode contribution is bounded by
$Cq_0^{3/2}M$.

For the fixed source, \eqref{cubic:source-remainder} bounds
its $L^2$ norm on $a/2<q<2a$ by $Ca^{7/2}M$ after
multiplication by $\xi_0q^2$ and the indicated derivatives.
Cauchy--Schwarz on each annulus therefore reduces its
contribution to
\[
\begin{aligned}
 \sum_{a\le q_0}\min(q/a,a/q)^3a^{7/2}
 &=q^{-3}\sum_{a\le q}a^{13/2}
   +q^3\sum_{q<a\le q_0}a^{1/2}\\
 &\le Cq^{7/2}+Cq^3q_0^{1/2}
 \le Cq_0^{1/2}\xi_\delta(q),
\end{aligned}
\]
where $a=q_0\,2^{-m}$, $m=0,1,\ldots$.

Combining these bounds with \eqref{cubic:low-kernel-bound}
and the annulus estimates in part~(b), we obtain
\begin{equation}\label{cubic:low-closure}
 \mathcal B_j
 \le C_{q_0}M
     +Cq_0^2(\mathcal B_j+\mathcal E_j)
     +Cq_0^2\sum_{\ell<j}(\mathcal B_\ell+\mathcal E_\ell).
\end{equation}
The constant $C$ multiplying $q_0^2$ is independent of
$q_0$ and $\delta$, whereas $C_{q_0}$ is independent of
$\delta$.

\smallskip\noindent
\emph{(e) Closing the estimates and taking traces.}
Since the finiteness of $\mathcal B_j+\mathcal E_j$ has not
yet been established, we first replace $\xi_\delta$ by
$\xi_{\delta,N}=\xi_\delta+N^{-1}$, so that the corresponding
weighted norms are finite:
\[
 \xi_{\delta,N}(q)=\xi_\delta(q)+N^{-1},\qquad N\ge1,
\]
and denote the corresponding quantities by
$\mathcal B_{j,N}$ and $\mathcal E_{j,N}$.
It follows from \eqref{h5c:mild},
\eqref{cubic:conormal-system}, and
\eqref{cubic:noncritical-multiplier} that
\[
 (h_j)_{-3},\ (D_w^j(Pg))_{-3}
       \in H^1((-\infty,v_0)),
 \qquad
 \nabla\partial_\theta^a h_{j,\mathrm h}\in L^2,
 \quad a=0,1.
\]
In particular, $\mathcal B_{j,N}+\mathcal E_{j,N}<\infty$.
For the flux, the first-mode potential disappears upon
projection onto mode $-3$; all remaining terms in its radial
derivative are bounded by the cited estimates.

The estimates in parts~(c)--(d) hold uniformly in $N$.
Indeed,
\[
 0\le\partial_v\log\xi_{\delta,N}\le3,
 \qquad \xi_{\delta,N}\ge\xi_\delta,
\]
so the coercivity and source estimates remain valid.
For the kernel estimate, the corresponding calculation is
\[
\begin{aligned}
 &\frac1{\xi_{\delta,N}(q)^r}
   \int_0^{q_0}\min(q/a,a/q)^{3r}
           a^{2r}\xi_{\delta,N}(a)^r\,\frac{da}{a}\\
 &\qquad\le
 q^{2r}\int_0^q(a/q)^{5r}\,\frac{da}{a}
       +\int_q^{q_0}a^{2r}\,\frac{da}{a}\\
 &\qquad=
 \frac{q^{2r}}{5r}
 +\frac{q_0^{2r}-q^{2r}}{2r}
 \le C_rq_0^{2r},\qquad r=1,2.
\end{aligned}
\]
Thus \eqref{cubic:high-closure} and
\eqref{cubic:low-closure}, with this weight, imply
\[
 \mathcal B_{j,N}+\mathcal E_{j,N}
 \le C_{q_0}M
   +Cq_0^2\sum_{\ell=0}^j
          (\mathcal B_{\ell,N}+\mathcal E_{\ell,N}).
\]
Choose $q_0$ so that $Cq_0^2\le1/2$. By induction
over $j=0,1,2$ we get
\[
 \mathcal B_{j,N}+\mathcal E_{j,N}\le C_{q_0}M,
\]
uniformly in $N$. Since $\xi_{\delta,N}$ decreases to
$\xi_\delta$, passage to the limit proves
\eqref{cubic:noncritical-weighted}.

We next establish the traces at $v=0$. On $|v|<1$,
$\xi_\delta(\delta e^v)\sim\delta^3$. Therefore
\eqref{cubic:noncritical-weighted} implies
\[
 \|h_{j,\mathrm h}\|_{H^1_vL^2_\theta}
 +\|\partial_v\partial_\theta h_{j,\mathrm h}\|_2
 +\|\partial_\theta^2h_{j,\mathrm h}\|_2
 \le C_{q_0}\delta^3M.
\]
All norms in the remainder of this step are over
$(-1,1)\times\mathbb T$, unless otherwise indicated.
The first identity in \eqref{cubic:conormal-system},
the bound on $\mathcal B_j$, and absorption of the
$O(\delta^2)$ metric coefficient also imply
\[
 \|(h_j)_{-3}\|_{H^1(-1,1)}
 \le C_{q_0}\delta^3M.
\]

To estimate the flux, differentiate the divergence identity:
\[
\begin{aligned}
 \partial_vD_w^j(Pg)
 ={}&-\partial_\theta D_w^j
       \bigl(\mathsf G_{v\theta}\partial_vg
                    +\mathsf G_{\theta\theta}\partial_\theta g\bigr)-D_w^j(\xi_0q^2Tg)+D_w^j(\xi_0q^2H).
\end{aligned}
\]
Project this identity onto the modes
$n\equiv1\pmod4$ with $|n|\ge5$.
The critical potential contributes zero to this projection.
The remaining potential terms satisfy
\[
 \|[D_w^j(\xi_0q^2Tg)]_{\mathrm h}\|_2
 \le C\delta^2\sum_{\ell=0}^j\|h_{\ell,\mathrm h}\|_2
\]
by \eqref{cubic:noncritical-multiplier}.
The metric terms containing the first mode satisfy
\[
 \delta^2\sum_{\ell=0}^j
       \|(h_\ell)_1\|_{H^1(-1,1)}
 \le C\delta^{9/2}M
\]
by \eqref{h5c:critical} with $\beta=1/2$.
Moreover, the source calculation in part~(b), based on
\eqref{cubic:source-remainder}, implies
\[
 \|[D_w^j(\xi_0q^2H)]_{\mathrm h}\|_2
 \le C\delta^{7/2}M.
\]
It follows from these estimates and the preceding bounds
on the other angular modes that
\[
 \|(D_w^j(Pg))_{\mathrm h}\|_{H^1_vL^2_\theta}
 \le C_{q_0}\delta^3M.
\]

Finally, apply
\[
 \|U(0,\cdot)\|_{L^2_\theta}
 \le C\|U\|_{H^1((-1,1);L^2_\theta)}
\]
to $h_{j,\mathrm h}$ and $(D_w^j(Pg))_{\mathrm h}$.
For mode $-3$, we use
\[
 |(h_j)_{-3}(0,w)|+|(D_w^j(Pg))_{-3}(0,w)|
 \le \xi_\delta(\delta)\mathcal B_j
 =\sqrt2\,\delta^3\mathcal B_j.
\]
It follows from Parseval's identity that
\eqref{cubic:noncritical-traces} holds.

\emph{Step 4: justification of Green's identity.}
We can first apply Green's identity on $q<\delta e^v$, with $v<0$.
All multipliers of $T$ in the chosen rotation class are real on
this domain, so their contributions cancel. It remains to justify
the limit of the boundary integral as $v\to0^-$.

The limits exist by Step~3.
\eqref{cubic:conormal-system} implies
\begin{equation}\label{cubic:first-flux-log}
 (D_w^j(Pg))_1(v,w)
 =4\partial_w^j\!\left[\xi_0(\delta)b(\delta)g_1(0,w)\right]
       \log(v\mp i0)+r_j(v,w),
 \qquad r_j(\cdot,w)\in H^1(-1/2,0).
\end{equation}
To verify this, note
$8e^{2v}/(e^{2v}-1)=4/v+O(1)$. After subtracting the derivative
of the above logarithm, the singular term is
\[
 \frac4v\left\{
 D_w^j[\xi_0(q)b(q)g_1(v,w)]
 -\partial_w^j[\xi_0(\delta)b(\delta)g_1(0,w)]\right\}.
\]
It belongs to $L^2(-1/2,0)$ by \eqref{h5c:mild} and
\[
 \left\|\frac{F(v)-F(0)}{v}\right\|_{L^2(-1/2,0)}
 \le2\|F'\|_{L^2(-1/2,0)}.
\]
The remaining angular and source terms are in $L^2$ by
\eqref{h5c:mild} and \eqref{h5c:source}. This proves
\eqref{cubic:first-flux-log}.
On the other hand,
\[
 \phi_1(v,w)=i[b(\delta e^v)-b(\delta)](x-iy)_1(\delta e^v),
 \qquad
 |\partial_w^j\phi_1(v,w)|\le C\delta^3|v|,
 \quad -1/2<v<0,\quad j\le2.
\]
Since $v\log(v\mp i0)\to0$, the products containing
$\overline{\phi_1}(Pg)_1$, including their first two $w$ derivatives,
vanish in the boundary limit. Thus \eqref{cubic:exact-green} holds through two $w$ derivatives.
All ingredients used in Step~1 have now been established.

\end{proof}

\subsection{The spectral jump for the first angular modes}
\label{subsec:cubic-first-jump}

Lemma~\ref{cubic:first-mode} controls the first angular
coefficient on the critical circle $q=\delta$. We now obtain
bounds for the spectral jump throughout the center neighborhood.
For $\delta\ll q\ll1$, the leading radial operator in the first
angular mode is $-(q\partial_q)^2+9$. Since
\[
 \bigl[-(q\partial_q)^2+9\bigr]q^3=0,\qquad
 \bigl[-(q\partial_q)^2+9\bigr](q^3\log q)=-6q^3,
\]
the cubic term in the fixed source $\xi_0q^2H_1$ may produce
a $q^3\log q$ contribution to each boundary value.
Subtracting the equations
$(\Delta+T_\pm)g_\pm=H$ cancels this source exactly:
\[
 (\Delta+T_+)(g_+-g_-)=-(T_+-T_-)g_-.
\]
For $\lambda=b(\delta)$, within the chosen center chart and
rotation class, the right-hand side is supported on
$q=\delta$ and depends only on $g_{-,1}(\delta,b(\delta))$.
Its coefficient, together with its first two normalized
spectral derivatives, is controlled by
\eqref{cubic:first-trace}. After localization, the additional
cutoff terms are supported in $q_0/2\le q<q_0$; the precise
equation is \eqref{cubic:first-jump-source}.
\begin{lemma}\label{cubic:first-profiles}
Let $\chi$ be supported in $q<q_0$ and equal to one on
$q<q_0/2$. For either $\lambda=1-\delta^2/8$ or
$\lambda=1+\delta^2/8$, set
\[
 G_1(q,\delta)=\chi(q)
           [g_{+,1}(q,\lambda)-g_{-,1}(q,\lambda)],
 \qquad U_1(w,v)=\delta^{-3}G_1(\delta e^v,\delta).
\]
Then
\begin{align}
 \|U_1\|_{H^1_v(|v|<2)}+
       \|\partial_wU_1\|_{H^1_v(|v|<2)}
       &\le CM,\label{cubic:first-profile-core}\\
 |G_1(q,\delta)|&\le CMq^3,
             &&\delta<q/2,\label{cubic:first-profile-outer}\\
 \sum_{j=0}^2|(\delta\partial_\delta)^jG_1(q,\delta)|
       &\le CMq\delta^2,
             &&\delta>2q.\label{cubic:first-profile-inner}
\end{align}
The derivatives in the last line hold $q$ fixed.
The analogous assertions hold at the negative threshold and at
the other center.
\end{lemma}

\begin{proof}
We first consider $\lambda=b(\delta)=1-\delta^2/8$ and set
\[
 v=\log(q/\delta),\qquad w=\log\delta,\qquad
 D_w=q\partial_q-\frac{\delta^2}{4}\partial_\lambda.
\]
Thus $D_w=\partial_w$ with $v$ fixed. Apply the projection
$h\mapsto\frac14\sum_{m=0}^3i^{-m}h(q,\theta+m\pi/2)$ to
$g_\pm=g_{f,\pm}(\lambda)$. As in the proof of
Lemma~\ref{cubic:first-mode}, this projection commutes with the
equation and preserves the first coefficient. In this proof,
$g_\pm$ denote the projected functions. Set
\[
 u=\chi(g_+-g_-),\qquad
 u_{\ne1}=(I-\Pi_1)u,\qquad \Pi_1u=u_1e^{i\theta}.
\]
All modes of $u_{\ne1}$ satisfy $n\equiv1\pmod4$, $n\ne1$,
and hence $|n|\ge3$. Extend $u$ by zero outside the support of
$\chi$. Cylinder norms use $dv\,d\theta$, and norms of individual
coefficients use $dv$. We choose $q_0$ sufficiently small that
$q_0<c_0$ and $q_0/8<c_1$, so Lemma~\ref{quad:weighted-center}
applies throughout the parameter range.

\emph{Step 1: the decisive cancellation and the estimates to be proved.}
Subtracting the two boundary equations cancels the fixed source.
After multiplication by $-\xi_0q^2$ and localization, the exact
equation is
\begin{equation}\label{cubic:first-jump-source}
 \begin{aligned}
 \relax[-\operatorname{div}(\mathsf G\nabla)-\xi_0q^2T_+]u
 ={}&-8\pi i\,\xi_0(\delta)b(\delta)
       g_{-,1}(\delta,b(\delta))\delta_0(v)e^{i\theta}
       +\mathcal F_\chi,\\
 \mathcal F_\chi={}&-2\mathsf G\nabla\chi\cdot\nabla(g_+-g_-)
       -\operatorname{div}(\mathsf G\nabla\chi)(g_+-g_-).
 \end{aligned}
\end{equation}
Here $T_+=T(i\lambda+0)$, all derivatives in the divergence
terms use $(v,\theta)$, and
$\operatorname{supp}\mathcal F_\chi\subset\{q_0/2\le q<q_0\}$.
Indeed, the first multiplier is
$T_{+,1}=-8b/(q^2-\delta^2-i0)$, and
\[
 8q^2\left(\frac1{q^2-\delta^2-i0}
              -\frac1{q^2-\delta^2+i0}\right)
 =8\pi i\delta_0(v).
\]

The Green kernel of the first-mode normal operator satisfies
\begin{equation}\label{cubic:first-green-kernel}
 \begin{gathered}
 N_1=-\partial_v^2+1+\frac{8e^{2v}}{e^{2v}-1-i0},
 \qquad p(v)=e^v(e^{2v}-1),\\
 N_{1,v}K_1(v,t)=\delta_t(v),\qquad t\in\mathbb R,\\
 \lim_{v\to-\infty}K_1(v,t)
 =\lim_{v\to+\infty}K_1(v,t)=0,\\
 K_1(v,t)=p(\min(v,t))z(\max(v,t)),\qquad
 K_1(v,0)=\tfrac12p(v)\mathbf1_{\{v<0\}}.
 \end{gathered}
\end{equation}
Here $N_{1,v}$ acts on the variable $v$, $\delta_t$ is the
Dirac mass at $v=t$, and the equation holds in the
distributional sense.
In Step~2 we will define $z$ and verify this formula. By
\eqref{cubic:first-trace} and the product rule,
\[
 \left|\partial_w^j\bigl[
 \xi_0(\delta)b(\delta)g_{-,1}(\delta,b(\delta))\bigr]\right|
 \le C\delta^3M,\qquad 0\le j\le2.
\]
By \eqref{cubic:first-green-kernel}, the contribution of the
$j$th $w$ derivative of the Dirac source in
\eqref{cubic:first-jump-source} is
\[
\begin{aligned}
 &-8\pi i\,\partial_w^j
   \bigl[\xi_0(\delta)b(\delta)
          g_{-,1}(\delta,b(\delta))\bigr]K_1(v,0)\\
 &\qquad=-4\pi i\,\partial_w^j
   \bigl[\xi_0(\delta)b(\delta)
          g_{-,1}(\delta,b(\delta))\bigr]
          p(v)\mathbf1_{\{v<0\}}.
\end{aligned}
\]
Since $p(0)=0$, its first distributional $v$ derivative
contains no Dirac mass. Moreover,
\[
 |p(v)|+|p'(v)|\le Ce^v,\qquad v<0.
\]
It therefore follows from \eqref{cubic:first-trace} that
this contribution and its first $v$ derivative vanish for
$q>\delta$ and satisfy
\[
 C\delta^3e^vM=Cq\delta^2M,\qquad 0<q<\delta,
 \quad 0\le j\le2.
\]
These are the required bounds for the Dirac source.

It remains to estimate the metric and potential perturbations
after applying the Green operators. In the weighted norms
defined in \eqref{cubic:first-jump-weighted-norm}, we will prove
\[
 \mathcal E_j
 \le CM+Cq_0^2\sum_{\ell=0}^j\mathcal E_\ell,
 \qquad 0\le j\le2,
\]
with $C$ multiplying $q_0^2$ independent of $q_0$ and $\delta$.
For sufficiently small $q_0$, moving
$Cq_0^2\mathcal E_j$ to the left and proceeding by induction
over $j$ proves the desired weighted bounds.

For this purpose define
\[
 W_\delta(q)=
 \begin{cases}q\delta^2,&q\le\delta,\\q^3,&q\ge\delta,\end{cases}
 \qquad
 V_\delta(q)=
 \begin{cases}\delta^{1/2}q^{5/2},&q\le\delta,\\q^3,&q\ge\delta.
 \end{cases}
\]
The stronger left weight $V_\delta$ is used for the modes
$|n|\ge3$. For $0\le j\le2$, set
\begin{equation}\label{cubic:first-jump-weighted-norm}
 \begin{aligned}
 \mathcal E_j={}&\sup_{\sigma\in\mathbb R}
 \frac{\|D_w^ju_1\|_{H^1(\sigma-1,\sigma+1)}}
      {W_\delta(\delta e^\sigma)}+\sum_{a=0}^2\sup_{\sigma\in\mathbb R}
 \frac{\|\partial_\theta^aD_w^ju_{\ne1}\|_
 {H^1((\sigma-1,\sigma+1)\times\mathbb T)}}
      {V_\delta(\delta e^\sigma)}.
 \end{aligned}
\end{equation}
We will prove in Steps~2--3 that
\begin{equation}\label{cubic:first-jump-weighted-bound}
 \mathcal E_j\le CM,\qquad 0\le j\le2.
\end{equation}
We first deduce all three conclusions from this estimate.

\smallskip\noindent
\emph{(a) The profile on $|v|<2$.}
On these intervals $W_\delta(\delta e^v)\sim\delta^3$ and
\[
 U_1=\delta^{-3}u_1,\qquad
 \partial_wU_1=\delta^{-3}(D_wu_1-3u_1).
\]
Covering $(-2,2)$ by a fixed number of intervals of length two,
we obtain \eqref{cubic:first-profile-core} from
\eqref{cubic:first-jump-weighted-bound} for $j=0,1$.

\smallskip\noindent
\emph{(b) The bound for $\delta<q/2$.}
The one-dimensional embedding $H^1(-1,1)\hookrightarrow L^\infty(-1,1)$
and \eqref{cubic:first-jump-weighted-bound} imply
\[
 |G_1(q,\delta)|=|u_1(\log(q/\delta),w)|
 \le CMW_\delta(q)=CMq^3.
\]
This proves \eqref{cubic:first-profile-outer}.

\smallskip\noindent
\emph{(c) The fixed-$q$ derivatives for $q<\delta/2$.}
Fix $v_*< -\log2$ and use intervals of radii
$\log2/4$ and $\log2/8$ about $v_*$. They lie in $v<0$ and
$q<\delta<q_0/8$, where $\chi=1$ and
\eqref{cubic:first-jump-source} is homogeneous. In particular,
\begin{equation}\label{cubic:first-jump-radial-equation}
 \begin{aligned}
 \mathsf G_{vv}u_{vv}={}&-2\mathsf G_{v\theta}u_{v\theta}
       -\mathsf G_{\theta\theta}u_{\theta\theta}
       -(\partial_v\mathsf G_{vv}
                         +\partial_\theta\mathsf G_{v\theta})u_v\\
 &-(\partial_v\mathsf G_{v\theta}
                         +\partial_\theta\mathsf G_{\theta\theta})u_\theta
       -\xi_0q^2T_+u.
 \end{aligned}
\end{equation}
Here $\mathsf G_{vv}$ is bounded below, and the potential and
its $\partial_v$ and $D_w$ derivatives are bounded angular multipliers,
uniformly on the larger interval. Solve
\eqref{cubic:first-jump-radial-equation} for $u_{vv}$;
then differentiate it in $\theta$ and $v$, and apply $D_w$.
Using $[D_w,\partial_v]=[D_w,\partial_\theta]=0$, we obtain
on the smaller interval
\begin{equation}\label{cubic:first-jump-radial-bound}
 \begin{aligned}
 &\|u_1\|_{H^3_v}+\|(D_wu)_1\|_{H^2_v}
                         +\|(D_w^2u)_1\|_{H^1_v}\\
 &\qquad\le C\sum_{j=0}^2\sum_{a=0}^2
       \|\partial_\theta^aD_w^ju\|_{H^1}
 \le CMW_\delta(\delta e^{v_*}).
 \end{aligned}
\end{equation}
The norms on the right are on the larger interval times
$\mathbb T$; the last inequality uses $V_\delta\le W_\delta$.
At fixed $q$, $\delta\partial_\delta=D_w-\partial_v$, so
\[
 (D_w-\partial_v)^2u_1
 =D_w^2u_1-2\partial_vD_wu_1+\partial_v^2u_1.
\]
Applying the same one-dimensional embedding to
\eqref{cubic:first-jump-radial-bound} proves
\eqref{cubic:first-profile-inner}. It remains to establish
\eqref{cubic:first-jump-weighted-bound}.

\emph{Step 2: the normal inverse in the two weights.}
For $v>0$, set
$z(v)=p(v)\int_v^\infty p(t)^{-2}\,dt$ and continue it across
zero with the boundary prescription in $N_1$. Direct
differentiation shows that
\[
 N_1p=N_1z=0,\qquad pz'-p'z=-1,\qquad z(0)=\tfrac12.
\]
Since $p$ decays at $-\infty$ and $z$ decays at $+\infty$,
these identities prove \eqref{cubic:first-green-kernel}.
For $|v|,|t|>1$ and $v\ne t$, the explicit formulas for
$p$ and $z$ imply
\begin{equation}\label{cubic:first-green-decay}
 \sum_{a,b=0}^1
       |\partial_v^a\partial_t^bK_1(v,t)|
 \le C
 \begin{cases}
 e^{-|v-t|},&v,t<-1,\\
 e^{-3|v-t|},&v,t>1,\\
 e^{v-3t},&v<-1,\ t>1,\\
 e^{t-3v},&t<-1,\ v>1.
 \end{cases}
\end{equation}
The mixed derivative in this estimate is evaluated away
from the diagonal $v=t$.

For $n\equiv1\pmod4$, $n\ne1$, the inverse of
$-\partial_v^2+n^2$ has kernel
\[
 K_n(v,t)=\frac{e^{-|n||v-t|}}{2|n|},
 \qquad
 |\partial_v^a\partial_t^bK_n(v,t)|
 =\frac12|n|^{a+b-1}e^{-|n||v-t|},
 \quad a,b\in\{0,1\},\quad v\ne t.
\]
Let $I,J\subset\mathbb R$ be disjoint intervals and set
$d=\operatorname{dist}(I,J)$. It follows from this formula
and Cauchy--Schwarz that
\[
 \left\|\mathbf1_I\partial_v^a
       (-\partial_v^2+n^2)^{-1}
       \partial_v^b\mathbf1_J\right\|_{L^2\to L^2}
 \le C|n|^{a+b-2}e^{-|n|d}
 \le Ce^{-3d}.
\]

When the intervals intersect, we use the weak equation
instead of pointwise bounds on the mixed derivative.
Indeed, if $u_n\in H^1(\mathbb R)$ satisfies
\[
 -u_n''+n^2u_n=\partial_vF_n^v+inF_n^\theta,
 \qquad F_n^v,F_n^\theta\in L^2(\mathbb R),
\]
then testing against $\overline{u_n}$ implies
\[
\begin{aligned}
 \|u_n'\|_2^2+n^2\|u_n\|_2^2
 &=\operatorname{Re}\int_{\mathbb R}
       \left(-F_n^v\overline{u_n'}
                    +inF_n^\theta\overline{u_n}\right)\,dv,\\
 \bigl(\|u_n'\|_2^2+n^2\|u_n\|_2^2\bigr)^{1/2}
 &\le
 \bigl(\|F_n^v\|_2^2+\|F_n^\theta\|_2^2\bigr)^{1/2}.
\end{aligned}
\]
This estimates divergence sources uniformly in $n$ and
accounts for the diagonal term in
\[
 \partial_v\partial_tK_n(v,t)
 =\delta_t(v)-\frac{|n|}{2}e^{-|n||v-t|}.
\]

Set $v_0=\log(q_0/\delta)$, and denote the bound in
\eqref{cubic:first-green-decay} by $\mathcal K_1(v,t)$.
The required integral estimates are
\begin{equation}\label{cubic:first-endpoint-kernels}
 \begin{aligned}
 \int_{\substack{t<v_0\\|t|>1}}\mathcal K_1(v,t)\delta^2e^{2t}
                   W_\delta(\delta e^t)\,dt
 &\le Cq_0^2W_\delta(\delta e^v),\qquad |v|>1,\\
 \int_{-\infty}^{v_0}e^{-3|v-t|}\delta^2e^{2t}
                   V_\delta(\delta e^t)\,dt
 &\le Cq_0^2V_\delta(\delta e^v),\\
 \int_{-\infty}^{v_0}e^{-3|v-t|}\delta^2e^{2t}
                   W_\delta(\delta e^t)\,dt
 &\le Cq_0^2V_\delta(\delta e^v).
 \end{aligned}
\end{equation}
The second line follows from
\eqref{cubic:complementary-kernel-moments} with $r=1$.
For the first and third, set $x=e^v$ and $R=q_0/\delta$.
When $x\le1$, the respective ratios of the integrals to their
output weights are bounded by
\[
 C\delta^2\left(\frac{R^2}{2}-\frac{x^2}{4}\right),\qquad
 C\delta^2x^{1/2}
       \left(\frac16+\log(1/x)+\frac{R^2-1}{2}\right)
 \le Cq_0^2.
\]
When $1\le x\le R$, the parts with $e^t<1$ are bounded by
$C\delta^2/(4x^6)$ and $C\delta^2/(6x^6)$; the common remaining
bound is
\[
 C\delta^2\left(\frac{x^2-x^{-6}}8+\frac{R^2-x^2}2\right)
 \le Cq_0^2.
\]
For $v\ge v_0$, each response decays at least as
$e^{-3(v-v_0)}$, while the weights increase as $e^{3(v-v_0)}$.
The modes $n\equiv1\pmod4$, $n\ne1$, are estimated with the
weight $V_\delta$, whereas the first mode is estimated with
$W_\delta$. Their contribution to the first-mode equation
comes from the metric terms, whose coefficients are $O(q^2)$.
Since $V_\delta\le W_\delta$, the first line of
\eqref{cubic:first-endpoint-kernels} implies, for $|v|>1$,
\[
\begin{aligned}
 &\int_{\substack{t<v_0\\|t|>1}}
   \mathcal K_1(v,t)\delta^2e^{2t}
                    V_\delta(\delta e^t)\,dt\le
 \int_{\substack{t<v_0\\|t|>1}}
   \mathcal K_1(v,t)\delta^2e^{2t}
                    W_\delta(\delta e^t)\,dt
 \le Cq_0^2W_\delta(\delta e^v).
\end{aligned}
\]
Thus the metric coupling from these modes to the first mode
satisfies the required bound with input weight $V_\delta$,
output weight $W_\delta$, and factor $Cq_0^2$.

We now estimate the contributions for which either the source
variable $t$ or the output variable $v$ is close to zero.
Choose $\rho\in C_c^\infty((-3,3))$ with $\rho=1$ on
$[-2,2]$. Let $F$ be a metric or potential error term acting
on $D_w^ju$, as in \eqref{cubic:first-jump-commuted}, and define
$U^{(0)},U^{(1)}\in H^1(\mathbb R\times\mathbb T)$ by
\[
\begin{aligned}
 \Big[N_1\Pi_1+(-\partial_v^2-\partial_\theta^2)(I-\Pi_1)\Big]U^{(0)}
 &=\rho F,\\
 \Big[N_1\Pi_1+(-\partial_v^2-\partial_\theta^2)(I-\Pi_1)\Big]U^{(1)}
 &=(1-\rho)F.
\end{aligned}
\]
These equations are restricted to the rotation class
$n\equiv1\pmod4$.

On $|t|<4$, both weights are comparable to $\delta^3$,
while the error coefficients and their required cylinder
derivatives are bounded by $C\delta^2$.
Integration by parts for the metric terms and
\eqref{eq:cylinder-form} for the critical potential imply
\[
 \sum_{a=0}^2
 \|\partial_\theta^a(\rho F)\|_{H^{-1}(\mathbb R\times\mathbb T)}
 \le C\delta^2
       \sum_{a=0}^2
       \|\partial_\theta^aD_w^ju\|_{H^1((-4,4)\times\mathbb T)}
 \le C\delta^5\mathcal E_j.
\]
It follows from Proposition~\ref{compact:normal} and the
preceding energy estimate for the other angular modes that
\[
 \sum_{a=0}^2
 \|\partial_\theta^aU^{(0)}\|_{H^1(\mathbb R\times\mathbb T)}
 \le C\delta^5\mathcal E_j.
\]

Outside the support of $\rho F$, the explicit Green kernels
also imply, for $|v|>4$,
\[
\begin{aligned}
 \|\Pi_1U^{(0)}\|_{H^1((v-1,v+1)\times\mathbb T)}
 &\le C\delta^5\mathcal E_j
 \begin{cases}
 e^v,&v<-4,\\
 e^{-3v},&v>4,
 \end{cases}\\
 \sum_{a=0}^2
 \|\partial_\theta^a(I-\Pi_1)U^{(0)}
                   \|_{H^1((v-1,v+1)\times\mathbb T)}
 &\le C\delta^5\mathcal E_j
 \begin{cases}
 e^{3v},&v<-4,\\
 e^{-3v},&v>4.
 \end{cases}
\end{aligned}
\]
For $v<0$, the normalization weights are
\[
 W_\delta(\delta e^v)=\delta^3e^v,\qquad
 V_\delta(\delta e^v)=\delta^3e^{5v/2};
\]
for $v>0$, both equal $\delta^3e^{3v}$.
Thus these estimates, together with the local $H^1$ bound,
bound the contribution of $U^{(0)}$ to the weighted norms
by $C\delta^2\mathcal E_j\le Cq_0^2\mathcal E_j$.

For $U^{(1)}$, consider first $|v|<1$.
The source is supported in $|t|>2$, and the kernel formula
implies
\[
 \sum_{b=0}^1
 \|\partial_t^bK_1(\cdot,t)\|_{H^1(-1,1)}
 \le C
 \begin{cases}
 e^t,&t<-2,\\
 e^{-3t},&t>2.
 \end{cases}
\]
The kernels of the modes $|n|\ge3$ satisfy the corresponding
operator estimates, with the stronger factor $e^{3t}$ for
$t<-2$. Applying Cauchy--Schwarz on intervals of length one
and using \eqref{cubic:first-jump-weighted-norm}, we obtain
\[
\begin{aligned}
 \sum_{a=0}^2
 \|\partial_\theta^aU^{(1)}\|_{H^1((-1,1)\times\mathbb T)}
 &\le C\delta^5\mathcal E_j
 \left(
 \int_{-\infty}^{-2}e^{4t}\,dt
       +\int_2^{v_0}e^{2t}\,dt\right)\\
 &\le C\delta^3q_0^2\mathcal E_j.
\end{aligned}
\]
The integrals bound the corresponding sums over those
intervals, since the exponential weights vary by bounded
factors on each interval. Division by the output weights,
which are comparable to $\delta^3$ on $|v|<1$, proves the
required $Cq_0^2\mathcal E_j$ bound.

For $|v|>1$, the contribution of $U^{(1)}$ is controlled by
\eqref{cubic:first-endpoint-kernels}.
For an output interval crossing $v=1$ or $v=-1$, split its
$H^1$ integrals into the portions inside and outside $(-1,1)$
and apply the respective estimates. This controls every
interval in \eqref{cubic:first-jump-weighted-norm}.
The logarithmic singularity of $z'$ at zero is handled by
the inverse estimate and the $H^1_v$ kernel bounds above.

\emph{Step 3: the commuted equation and absorption of the errors.}
Moving the coefficient errors in
\eqref{cubic:first-jump-source} to the right, and applying
$D_w^j$, we obtain
\begin{equation}\label{cubic:first-jump-commuted}
 \begin{aligned}
 &[N_1\Pi_1+(-\partial_v^2-\partial_\theta^2)(I-\Pi_1)]D_w^ju\\
 &=-8\pi i\,\partial_w^j
   [\xi_0(\delta)b(\delta)g_{-,1}(\delta,b(\delta))]
                                    \delta_0(v)e^{i\theta}
       +D_w^j\mathcal F_\chi\\
 &\quad+\sum_{\ell=0}^j\binom j\ell
       \operatorname{div}\bigl(D_w^{j-\ell}(\mathsf G-I)
                                            \nabla D_w^\ell u\bigr)\\
 &\quad-\sum_{\ell=0}^j\binom j\ell
       D_w^{j-\ell}(\xi_0b-1)
       \frac{8e^{2v}}{e^{2v}-1-i0}(D_w^\ell u)_1e^{i\theta}\\
 &\quad+\sum_{\ell=0}^j\binom j\ell
       \sum_{\substack{n\equiv1\!\!\pmod4\\n\ne1}}
       D_w^{j-\ell}\left(\xi_0q^2\frac{nb}{nb-\lambda}\right)
                                   (D_w^\ell u)_ne^{in\theta}.
 \end{aligned}
\end{equation}
Multiply each
coefficient error by a fixed smooth radial function equal to one
on a neighborhood of $\operatorname{supp}\chi$ and zero for
$q\ge q_0$; this does
not change \eqref{cubic:first-jump-commuted} for the localized
function. By Lemma~\ref{geom:metric} and
\eqref{cubic:noncritical-multiplier}, these extended coefficients
satisfy
\begin{equation}\label{cubic:first-jump-coefficients}
 \begin{aligned}
 &|D_w^j\partial_\theta^a(\mathsf G-I)|
       +|D_w^j(\xi_0b-1)|
       +|\partial_vD_w^j(\xi_0b-1)|\\
 &\qquad+\sup_{\substack{n\equiv1\!\!\pmod4\\n\ne1}}
    \left|D_w^j\left(\xi_0q^2\frac{nb}{nb-\lambda}\right)\right|
 \le Cq^2,\qquad 0\le j,a\le2.
 \end{aligned}
\end{equation}
Radial derivatives of the last multiplier satisfy the same
bound. Applying $\partial_\theta^a$, $0\le a\le2$, to the
metric sum in \eqref{cubic:first-jump-commuted} replaces each
term by
\[
 \sum_{b=0}^a\binom ab\operatorname{div}
 \bigl(D_w^{j-\ell}\partial_\theta^{a-b}(\mathsf G-I)
                       \nabla\partial_\theta^bD_w^\ell u\bigr).
\]
The potential coefficients are radial. Integration by parts for
the metric terms, \eqref{eq:cylinder-form} for the singular
potential, and Parseval's identity for the other potentials
therefore reduce their bounds to Step~2. The total contribution
to \eqref{cubic:first-jump-weighted-norm} is at most
\begin{equation}\label{cubic:first-jump-perturbation}
 Cq_0^2\sum_{\ell=0}^j\mathcal E_\ell.
\end{equation}

The first source in \eqref{cubic:first-jump-commuted} was
estimated in Step~1. For the second source, the product rule,
the angular Poincar\'e inequality, and \eqref{h5c:mild} imply
\begin{equation}\label{cubic:first-jump-cutoff}
 \sum_{j=0}^2\sum_{a=0}^2
 \|\partial_\theta^aD_w^j\mathcal F_\chi\|_{H^{-1}}
 \le C_{q_0}M.
\end{equation}
Every term in the above is supported in
$v_0-\log2\le v\le v_0$. By
\eqref{cubic:first-green-decay}, its first-mode response on
intervals to the left of this set is bounded by
\[
 C_{q_0}M
 \begin{cases}e^{v-3v_0},&v<0,\\
 e^{3(v-v_0)},&0<v<v_0-2,
 \end{cases}
 \le C_{q_0}MW_\delta(\delta e^v).
\]
For the modes $n\equiv1\pmod4$, $n\ne1$, the free Green
kernels and \eqref{cubic:first-jump-cutoff} imply
\[
 \sum_{a=0}^2
 \left\|\partial_\theta^a
 (-\partial_v^2-\partial_\theta^2)^{-1}
 (I-\Pi_1)D_w^j\mathcal F_\chi
 \right\|_{H^1((v-1,v+1)\times\mathbb T)}
 \le C_{q_0}M e^{3(v-v_0)},
 \qquad v\le v_0-2.
\]
Since $e^{v_0}=q_0/\delta$,
\[
 \frac{e^{3(v-v_0)}}{V_\delta(\delta e^v)}
 =q_0^{-3}
 \begin{cases}
 e^{v/2},&v<0,\\
 1,&v\ge0,
 \end{cases}
 \le q_0^{-3}.
\]
Thus these contributions satisfy the required weighted bound,
with the factor $q_0^{-3}$ included in $C_{q_0}$.

For the remaining intervals, whose centers satisfy
$v\ge v_0-2$, the condition $\delta<q_0/8$ implies
$v_0-2>0$, and hence
\[
 W_\delta(\delta e^v)=V_\delta(\delta e^v)
 =\delta^3e^{3v}\ge e^{-6}q_0^3.
\]
It follows from Proposition~\ref{compact:normal}, the free-mode
energy estimate, and \eqref{cubic:first-jump-cutoff} that the
contributions of $D_w^j\mathcal F_\chi$ have $H^1$ norms
bounded by $C_{q_0}M$ on these intervals, also after the
angular derivatives appearing in
\eqref{cubic:first-jump-weighted-norm}.
Division by the displayed lower bound for the weights
therefore proves the same $C_{q_0}M$ bound for the normalized
norms.

For sufficiently small fixed $q_0$, the perturbation bound
\eqref{cubic:first-jump-perturbation} is less than $1/2$.
The Neumann series therefore constructs solutions of
\eqref{cubic:first-jump-commuted} in the weighted space
\eqref{cubic:first-jump-weighted-norm}, successively for
$j=0,1,2$.
These solutions belong to $H^1$: the weights decay
exponentially as $v\to-\infty$, and for $v>v_0$ the Green
formulas express their angular coefficients as multiples of
$z(v)$ in mode $1$ and $e^{-|n|v}$ in the other modes.

For the difference $h$ of two $H^1$ solutions with the same
lower-order terms, Proposition~\ref{compact:normal}, the
free-mode inverse estimate,
\eqref{cubic:first-jump-coefficients}, and
\eqref{eq:cylinder-form} imply
\[
\begin{aligned}
 \|h\|_{H^1}
 &\le C\|
 [N_1\Pi_1+(-\partial_v^2-\partial_\theta^2)(I-\Pi_1)]h
 \|_{H^{-1}}\\
 &\le Cq_0^2\|h\|_{H^1}.
\end{aligned}
\]
Thus $h=0$. Since $D_w^ju\in H^1$ by \eqref{h5c:mild},
induction identifies the constructed solutions with $D_w^ju$.

The weighted norms of $D_w^ju$ are consequently finite, and the
preceding estimates imply
\begin{equation}\label{cubic:first-jump-closure}
 \mathcal E_j
 \le CM+Cq_0^2\mathcal E_j
       +Cq_0^2\sum_{\ell<j}\mathcal E_\ell,
 \qquad 0\le j\le2.
\end{equation}
Moving $Cq_0^2\mathcal E_j$ to the left and proceeding by
induction proves \eqref{cubic:first-jump-weighted-bound}.
The conclusions of the lemma then follow from Step~1.

\emph{Step 4: $\lambda>1$ and the other cases.}
For $\lambda=1+\delta^2/8$, set
\[
 D_w=q\partial_q+\frac{\delta^2}{4}\partial_\lambda,
 \qquad N_1=-\partial_v^2+1+\frac{8e^{2v}}{e^{2v}+1},
 \qquad p(v)=e^v(e^{2v}+1).
\]
Now $z(v)=p(v)\int_v^\infty p(t)^{-2}\,dt$ is defined for
all $v\in\mathbb R$. The same Wronskian calculation proves
the kernel bounds \eqref{cubic:first-green-decay}.
The two local multipliers agree, so
\eqref{cubic:first-jump-source} contains only
$\mathcal F_\chi$. Steps~2--3 apply without the point source, and the identities in Step~1 prove the three desired
conclusions.

A smooth factor depending on $\delta$ as in the statement has
bounded $(\delta\partial_\delta)^j$ derivatives for $j\le2$. By complex conjugation
and the cell symmetries we obtain the same conclusions at the negative
threshold and at the other center.
\end{proof}

\section{Proof of the main theorem}\label{sec:main-proof}
In this section we give the proof of our main Theorem~\ref{th1:Int}. Throughout this section, $f=\omega_0$ and $M=\|f\|_{H^5}$.
\subsection{Inviscid damping}

\begin{theorem}\label{r2:global-damping}\label{r2:compact-damping}
Under the hypotheses of Theorem~\ref{th1:Int}, let
$\psi_c=\Delta^{-1}(\omega-p_\infty)$, where $p_\infty\in E$
is the stationary part in \eqref{lap:spectral-representation}.
Then, for $t\ge0$,
\begin{align}
 \|\psi_c(t)\|_2&\le C(1+t)^{-2}M,
                                      \label{quad:stream-decay}\\
 \|\psi_c(t)\|_{H^1}&\le C(1+t)^{-1}M.
                                      \label{r2:global-velocity}
\end{align}
\end{theorem}
\begin{proof}
It follows from Proposition~\ref{qual:theorem} and
\eqref{quad:fourier-second} that
\[
 \psi_c(t)=\int_{\mathbb R}e^{-it\lambda}\Gamma_f(\lambda)\,d\lambda,
 \qquad
 -t^2\psi_c(t)=\int_{\mathbb R}e^{-it\lambda}
                 \partial_\lambda^2\Gamma_f(\lambda)\,d\lambda.
\]
Moreover, it follows from
\eqref{lap:lift} that
\[
 \sup_{t\ge0}\|\Delta\psi_c(t)\|_2
 =\sup_{t\ge0}\|\omega(t)-p_\infty\|_2
 \le C\|f\|_{H^1}\le CM.
\]
It sufficies to consider $t\ge1$. Using integration by parts we obtain that
\[
 \|\nabla\psi_c(t)\|_2^2
 =-\operatorname{Re}\langle\psi_c(t),\Delta\psi_c(t)\rangle_{L^2}
 \le \|\psi_c(t)\|_2\|\Delta\psi_c(t)\|_2
 \le C(1+t)^{-2}M^2,
\]
which proves \eqref{r2:global-velocity}.
\end{proof}

\subsection{The vorticity representation}
\label{sec:vorticity-representation}

Recall that
\[
 g_f(z)=f+z\psi_f(z),\qquad
 \Gamma_f=\frac{\psi_{f,+}-\psi_{f,-}}{2\pi}.
\]
On a regular streamline, $L_0=b(r)\partial_\theta$.
It follows from \eqref{rayleigh1} that
\[
 \bigl[(\Delta+1)\psi_f(z)\bigr]_k(r)
 =\frac{g_{f,k}(r,z)}{z-ikb(r)},
 \qquad \operatorname{Re}z\ne0.
\]
The angular coefficient is taken after applying the physical
operator $\Delta+1$. This operator annihilates the stationary
resolvent term $-p_\infty/z$, which was identified in the proof of
Proposition~\ref{qual:theorem}, since $(\Delta+1)p_\infty=0$.

We have the identity
\[
 \lim_{\varepsilon\to0^+}
 \frac1{\pm\varepsilon+i(\lambda-kb(r))}
 =-i\operatorname{PV}\frac1{\lambda-kb(r)}
       \pm\pi\delta(\lambda-kb(r)).
\]
Using
$g_{f,+}-g_{f,-}=2\pi i\lambda\Gamma_f$, we obtain
\begin{equation}\label{sd:jump}
\begin{aligned}
 \Delta\Gamma_f
 &+\sum_{k\ne0}\operatorname{PV}
       \frac{kb(r)\Gamma_{f,k}(r,\lambda)}{kb(r)-\lambda}
                         e^{ik\theta}\\
 &=\frac12\sum_{k\ne0}\delta(\lambda-kb(r))
   \bigl[g_{f,+,k}(r,kb(r))
        +g_{f,-,k}(r,kb(r))\bigr]e^{ik\theta}.
\end{aligned}
\end{equation}
The zeroth coefficient of $\Delta\Gamma_f$ vanishes because
$f\in Y$ and $(\Delta\psi_f(z))_0=f_0/z=0$.

For $t>0$,
\[
 \operatorname{PV}\int_{\mathbb R}
       \frac{e^{-it\lambda}}{\lambda-kb(r)}\,d\lambda
 =-i\pi e^{-itkb(r)}.
\]
For $t>0$, the Dirac source in \eqref{sd:jump} and the pole term contribute
a term proportional to $e^{-itkb(r)}$. Its coefficient is
\begin{equation}\label{sd:amplitude}
 \frac{g_{f,+,k}(r,kb(r))+g_{f,-,k}(r,kb(r))}{2}
 -i\pi kb(r)\Gamma_{f,k}(r,kb(r))
 =g_{f,-,k}(r,kb(r)).
\end{equation}
The first term on the left comes from the Dirac source.
The second follows from the Fourier transform of
$\operatorname{PV}(\lambda-kb(r))^{-1}$.
The equality follows from \eqref{lap:density} and
\[
 g_{f,+,k}(r,kb(r))-g_{f,-,k}(r,kb(r))
 =2\pi i\,kb(r)\Gamma_{f,k}(r,kb(r)).
\]

In either center chart, this coefficient is
$g_{f,-,k}(q,kb(q))$, where $b(q)=1-q^2/8$.
With $M=\|f\|_{H^5}$, it follows from
\eqref{spectral:critical-amplitude} that
\[
 |g_{f,-,k}(q,kb(q))|
 \le Cq^2|k|^{-5/2}M,
 \qquad k\ne0,\quad 0<q<q_0/4.
\]
The estimates \eqref{cubic:first-trace} and
\eqref{cubic:higher-trace}, respectively, imply the stronger
cubic bounds
\[
\begin{aligned}
 |g_{f,-,k}(q,kb(q))|
 &\le Cq^3M,
 &&|k|=1,\quad 0<q<q_0/8,\\
 |g_{f,-,k}(q,kb(q))|
 &\le C|k|^{-2}q^3M,
 &&|k|\ge2,\quad 0<q<c_1|k|^{-1/2}.
\end{aligned}
\]
These estimates also hold for the other boundary value
$g_{f,+,k}$.

\subsection{The large angular and spectral parameters}

We use the following elementary principal-value bound. Suppose
$H\in C^\gamma(\mathbb R)$ is supported on an interval of length
$CL$, $L\ge1$, and
$\|H\|_\infty\le A$, $[H]_{C^\gamma}\le B$, then
\begin{equation}\label{r2:small-amplitude-PV}
 \sup_{a,t}\left|\operatorname{PV}\int
       e^{-it\lambda}\frac{H(\lambda)}{\lambda-a}\,d\lambda\right|
 \le C_\gamma A\bigl[1+\log(2+L(B/A)^{1/\gamma})\bigr].
\end{equation}

\begin{proposition}\label{r2:high-output-tail}
Let $f\in H^5\cap Y$ satisfy the moment conditions in
Theorem~\ref{th1:Int}, and $M=\|f\|_{H^5}$. For every fixed
$c_*>0$, the angular coefficients of the vorticity satisfy
\begin{equation}\label{cubic:large-angular}
 \sum_{|k|q^2\ge c_*}|(\omega-p_\infty)_k(t,q)|
 \le C_{c_*}Mq^3,\qquad t\ge0.
\end{equation}
Fix $c_{\rm spec}>0$ and
$0<c_{\rm low}\le c_{\rm spec}/2$.
For the modes $0<|k|q^2<c_{\rm low}$, the contribution from
$|\lambda|q^2\ge c_{\rm spec}$ satisfies
\begin{equation}\label{cubic:far-spectral-integral}
 \sum_{0<|k|q^2<c_{\rm low}}
 \int_{|\lambda|q^2\ge c_{\rm spec}}
 \left|
   \frac{kb(q)\Gamma_{f,k}(q,\lambda)}
        {\lambda-kb(q)}
 \right|\,d\lambda
 \le C_{c_{\rm spec}}Mq^3.
\end{equation}
\end{proposition}
\begin{remark}
    Since $0<b(q)\le1$ in the center chart,
$|\lambda-kb(q)|\ge|\lambda|/2$ on this integration region,
so no principal value is needed.
In particular, inserting $e^{-it\lambda}$ and taking the
absolute value of each integral before summing in $k$
preserves the bound \eqref{cubic:far-spectral-integral}, uniformly for $t\ge0$.
\end{remark}
\begin{proof}
The third commuted resolvent equation is
\begin{align*}
 (\Delta+T)L_0^3\psi_f
 & =-(L_0-z)^{-1}L_0^3f+3[\Delta,L_0]L_0^2\psi_f\\
 &\quad-3[[\Delta,L_0],L_0]L_0\psi_f
       +[[[\Delta,L_0],L_0],L_0]\psi_f.
\end{align*}
The commutators are second-order divergence-form operators with
smooth coefficients. It follows from Lemma~\ref{r2:full-source-LAP},
\eqref{r2:tangential-high}, and
$(L_0-z)^{-1}=(T-I)/z$ that
$\|L_0^3\psi_f\|_{H^1}\le C|\lambda|^{-1}\|f\|_{H^4}$.
Since $L_0=b(q)\partial_\theta$ and $b$ is smooth and positive in
the center disk, it follows from the one-dimensional trace inequality that
\begin{equation}\label{cubic:tail-traces}
 |g_{f,\pm,k}(q,\lambda)|\le C|k|^{-7/2}M,
 \qquad
 |\Gamma_{f,k}(q,\lambda)|\le C|k|^{-7/2}|\lambda|^{-1}M
 \quad (|\lambda|\ge\lambda_0).
\end{equation}
For the far spectral region we use
$\Gamma_f=(i\lambda)^{-2}\Gamma_{L^2f}$. Here
$\|L^2f\|_{H^2}\le C\|f\|_{H^4}$, so it follows from one angular commutation,
applied to this $H^2$ datum, that
\begin{equation}\label{cubic:tail-far}
 |\Gamma_{f,k}(q,\lambda)|
       \le C|k|^{-3/2}|\lambda|^{-3}M.
\end{equation}
For $|k|q^2\ge c_*$ and $|\lambda|\sim|k|$, the characteristic
of $D_\lambda$ with $v(\lambda_1)=\log q$ satisfies
\[
 |v'(\mu)|\le C_{c_*},
 \qquad
 \frac{d}{d\mu}g_{f,\pm,k}(e^{v(\mu)},\mu)
       =(\Lambda_{g,\pm}+R_{g,\pm})_k(e^{v(\mu)},\mu).
\]
On a logarithmic support, $V=0$ and its critical index $n$ satisfies
\[
 \left|\partial_\lambda\log
       \frac{q}{\sqrt{8(1-\lambda/n)}}\right|
 =\frac{4}{|n|\,8(1-\lambda/n)}\sim_{c_*}1,
 \qquad |n|\sim|k|.
\]
Using $\int_{-1}^1|\log|u||^2du<\infty$ for $\Lambda_g$ and
the trace inequality for $R_g$, it follows from
\eqref{r2:high-jet-coefficient} and \eqref{r2:high-regular-g} that
on each fixed-width spectral interval $I$,
\[
 \int_I|(\Lambda_{g,\pm}+R_{g,\pm})_k
                (e^{v(\mu)},\mu)|^2\,d\mu
       \le C_{c_*}k^{2N}M^2.
\]
It follows from Cauchy--Schwarz and
$|g_k(v_2,\lambda)-g_k(v_1,\lambda)|
\le |v_2-v_1|^{1/2}\|\partial_vg_k(\lambda)\|_2$
that
\[
 |g_{f,\pm,k}(q,\lambda_2)-g_{f,\pm,k}(q,\lambda_1)|
 \le C_{c_*}|k|^N M|\lambda_2-\lambda_1|^{1/2}.
\]
Choose a fixed $K$ sufficiently large that the preceding
estimates apply for $|k|>2K$. If
$0<\lambda_2-\lambda_1\le1$, divide
$[\lambda_1,\lambda_2]$ at the partition endpoints:
\[
 \lambda_1=\mu_0<\mu_1<\cdots<\mu_m=\lambda_2,
 \qquad m\le C.
\]
By Proposition~\ref{global:LAP}, the boundary values agree
at each endpoint. The preceding local estimate therefore implies
\[
\begin{aligned}
 |g_{f,\pm,k}(q,\lambda_2)-g_{f,\pm,k}(q,\lambda_1)|
 &\le C_{c_*}|k|^NM
       \sum_{j=1}^m(\mu_j-\mu_{j-1})^{1/2}\\
 &\le C_{c_*}|k|^NM|\lambda_2-\lambda_1|^{1/2}.
\end{aligned}
\]
For $|\lambda_2-\lambda_1|>1$, use
\eqref{cubic:tail-traces}.
For $1\le|k|\le2K$, the condition $|k|q^2\ge c_*$
implies $q^2\ge c_*/(2K)$; hence
\eqref{r2:away-holder} and its version for $\psi_{f,\pm}$,
together with $g_{f,\pm}=f+i\lambda\psi_{f,\pm}$,
apply on this fixed annulus.
Consequently, for some fixed $\gamma>0$ and $N$,
\[
 \|kb(q)\Gamma_{f,k}(q,\cdot)\|_
       {C^\gamma([k/2,2k])}
 \le C_{c_*}k^NM,\qquad k\ge1.
\]

We treat $k>0$. Choose $0\le\chi\le1$ in
$C_c^\infty((1/2,2))$, equal to one on a fixed neighborhood
of the range of $b$ in the center disk.
For
\[
 H(\lambda)=kb(q)\chi(\lambda/k)\Gamma_{f,k}(q,\lambda),
\]
the preceding bound and \eqref{cubic:tail-traces} imply
\[
 \|H\|_\infty\le CMk^{-7/2},
 \qquad [H]_{C^\gamma(\mathbb R)}
       \le C_{c_*}k^NM,
\]
after increasing $N$ if necessary.
Applying \eqref{r2:small-amplitude-PV} with $L=k$ yields
\[
 \sup_{t\ge0}
 \left|\operatorname{PV}\int_{\mathbb R}
       e^{-it\lambda}\frac{H(\lambda)}
                            {\lambda-kb(q)}\,d\lambda\right|
 \le C_{c_*}Mk^{-7/2}(1+\log(2+k)).
\]
The critical-value term in \eqref{sd:jump} is bounded by
$CMk^{-7/2}$, by \eqref{cubic:tail-traces}.

Fix $0<\lambda_0<\min(\lambda_*,1/4)$.
On the support of $1-\chi(\lambda/k)$ with
$|\lambda|\le2k$, we have $|\lambda-kb(q)|\ge ck$.
Thus \eqref{cubic:tail-traces} implies
\[
 \int_{\lambda_0\le|\lambda|\le2k}
 \bigl(1-\chi(\lambda/k)\bigr)
 \left|\frac{kb(q)\Gamma_{f,k}(q,\lambda)}
              {\lambda-kb(q)}\right|\,d\lambda
 \le CMk^{-7/2}\log(2+k).
\]
This includes the transition regions of $\chi$ and the
negative spectral interval. On $|\lambda|<\lambda_0$,
the bound
$|\Gamma_{f,k}|\le Ck^{-4}e^{-c/|\lambda|}M$,
which follows from \eqref{r2:tail-zero}, implies
\[
 \int_{|\lambda|<\lambda_0}
 \left|\frac{kb(q)\Gamma_{f,k}(q,\lambda)}
              {\lambda-kb(q)}\right|\,d\lambda
 \le CMk^{-4}.
\]
Finally, \eqref{cubic:tail-far} implies
\[
 \int_{|\lambda|>2k}
 \left|\frac{kb(q)\Gamma_{f,k}(q,\lambda)}
              {\lambda-kb(q)}\right|\,d\lambda
 \le CMk^{-1/2}\int_{2k}^\infty\lambda^{-4}\,d\lambda
 \le CMk^{-7/2}.
\]
The same argument with $(k,\lambda)$ replaced by
$(-k,-\lambda)$ treats negative indices. Summing these
estimates proves
\[
 \sum_{|k|q^2\ge c_*}|(\omega-p_\infty)_k(t,q)|
 \le C_{c_*}M
       \sum_{|k|\ge c_*/q^2}|k|^{-7/2}(1+\log(2+|k|))
 \le C_{c_*}Mq^5(1+|\log q|)
 \le C_{c_*}Mq^3.
\]

For the far spectral contribution, assume
$0<c_{\rm low}\le c_{\rm spec}/2$.
On the region
\[
 0<|k|q^2<c_{\rm low},
 \qquad |\lambda|q^2\ge c_{\rm spec},
\]
we have $|\lambda-kb(q)|\ge|\lambda|/2$.
Moreover, if the angular sum is nonempty, then
$q^2<c_{\rm low}$ and hence $|\lambda|>2$.
It follows from \eqref{cubic:tail-far} that
\[
\begin{aligned}
 &\sum_{0<|k|q^2<c_{\rm low}}
 \int_{|\lambda|q^2\ge c_{\rm spec}}
 \left|\frac{kb(q)\Gamma_{f,k}(q,\lambda)}
              {\lambda-kb(q)}\right|\,d\lambda\\
 &\qquad\le CM
       \sum_{1\le|k|<c_{\rm low}/q^2}|k|^{-1/2}
       \int_{c_{\rm spec}/q^2}^{\infty}\lambda^{-4}\,d\lambda\\
 &\qquad\le CMc_{\rm spec}^{-3}c_{\rm low}^{1/2}q^5
 \le C_{c_{\rm spec}}Mq^3.
\end{aligned}
\]
\end{proof}

\subsection{The noncritical spectral contribution}

After removing a neighborhood of $\lambda=k$, we estimate one
spectral derivative of the vorticity integrand in $L^1_\lambda$.
The factor $q$ in this estimate follows from the radial Green
formula for the nonzero angular modes. One integration by parts
then proves the required $q/(1+t)$ bound.

\begin{lemma}\label{cubic:noncritical}
Let $f\in H^5\cap Y$ satisfy the moment conditions of
Theorem~\ref{th1:Int}. Fix $0<\kappa_1<\kappa/8$, where
$\kappa$ is the parameter in \eqref{good-der}, and choose an even
$\chi\in C_c^\infty((-\kappa_1,\kappa_1))$ such that
$0\le\chi\le1$ and $\chi=1$ on $[-\kappa_1/2,\kappa_1/2]$.
Set $\rho_k(\lambda)=1-\chi(\lambda-k)$. Choose
\[
 0<c_{\rm spec}<\kappa_1/4096,\qquad
 0<c_{\rm low}<c_{\rm spec}/2,
\]
and $\eta\in C_c^\infty((-2c_{\rm spec},2c_{\rm spec}))$ with
$\eta=1$ on $[-c_{\rm spec},c_{\rm spec}]$.
In either center chart, for $0<q<q_*/4$,
\begin{equation}\label{cubic:noncritical-W11}
 \sum_{0<|k|q^2<c_{\rm low}}\sum_{j=0}^1
 \left\|\partial_\lambda^j\left[
   \eta(q^2\lambda)\rho_k(\lambda)
   \frac{kb(q)\Gamma_{f,k}(q,\lambda)}{\lambda-kb(q)}
                   \right]\right\|_{L^1(\mathbb R_\lambda)}
 \le Cq\|f\|_{H^5}.
\end{equation}
The derivatives are distributional. The expression in brackets
is defined to be zero near $\lambda=kb(q)$, where $\rho_k=0$.
Consequently, for $t\ge0$,
\begin{equation}\label{r2:nonlocal-integral}
 \left|\sum_{0<|k|q^2<c_{\rm low}}e^{ik\theta}
 \int_{\mathbb R}e^{-it\lambda}\eta(q^2\lambda)\rho_k(\lambda)
       \frac{kb(q)\Gamma_{f,k}(q,\lambda)}{\lambda-kb(q)}\,d\lambda
 \right|
 \le \frac{Cq}{1+t}\|f\|_{H^5}.
\end{equation}
The constant depends on the fixed cutoffs, but not on $q,t,f$.
\end{lemma}

\begin{proof}
Set $M=\|f\|_{H^5}$ and fix $0<\alpha<1$.
All cylinder gradients below are in $(v,\theta)$, where
$v=\log q$, and $[h]_k=(2\pi)^{-1}\int_{\mathbb T}h e^{-ik\theta}\,d\theta$.

\smallskip\noindent
\emph{Step 1: Spectral intervals and denominator bounds.}
First take $\lambda\ge1-\kappa_1/2$. On each of the intervals
\[
 |\lambda-n|<\kappa_1/2,\qquad
 n+\kappa_1/2\le\lambda\le n+1-\kappa_1/2,
 \qquad n\ge1,
\]
choose, respectively,
\[
 r_n^2=\frac{32\kappa_1}{n},\qquad
 r_n^2=\frac{\kappa_1}{16(n+1)}.
\]
These disks lie inside the spatial cutoffs in \eqref{good-der};
thus $D_\lambda=\partial_\lambda$ on $q<r_n$.
On the first interval, $\rho_n=0$, and the only possible
logarithmic correction is the one with index $n$. Its support
satisfies
\[
 q^2<e\delta^2\le4e\kappa_1/n<r_n^2,
 \qquad \delta^2=8(1-\lambda/n).
\]
On the second interval, every possible logarithmic support
satisfies
\[
 q^2>e^{-1}\delta^2\ge\frac{4\kappa_1}{e(n+1)}>r_n^2.
\]
There is therefore no logarithmic correction or critical circle
in the second disk.
Since $b(q)=1-q^2/8$, for $q<r_n$ and $k\ne0$ with
$\rho_k(\lambda)\ne0$ we have
\begin{equation}\label{cubic:noncritical-denominators}
 \begin{aligned}
 |kb(q)-\lambda|&\ge c(1+|k-n|),\\
 |T_k(q,\lambda)|&\le\frac{C|k|}{1+|k-n|},\qquad
 |\partial_\lambda T_k(q,\lambda)|
       \le\frac{C|k|}{(1+|k-n|)^2},\qquad
 T_k=\frac{kb(q)}{kb(q)-\lambda}.
 \end{aligned}
\end{equation}
Indeed, on the first interval $k\ne n$; on the second,
$|k-\lambda|\ge\kappa_1/2$ and
$|k|q^2/8\le\kappa_1/64$ when $0<k\le2(n+1)$.
For $k>2(n+1)$ or $k<0$, the denominator is bounded below
by a fixed multiple of $|k|+n$.
Our choice of $c_{\rm spec}$ ensures $q<r_n/4$ whenever
$\eta(q^2\lambda)\ne0$.

\smallskip\noindent
\emph{Step 2: From cylinder energy to pointwise bounds.}
For $v_0=\log r_n$, the Dirichlet Green kernel of
$-\partial_v^2+k^2$ on $(-\infty,v_0)$ is
\[
 \mathcal K_k(v,u)
 =\frac{e^{-|k||v-u|}-e^{-|k|(2v_0-v-u)}}{2|k|}.
\]
It satisfies $|\partial_u\mathcal K_k|+|k||\mathcal K_k|
\le Ce^{-|k||v-u|}$ and, for $q=e^v<r_n$,
\[
 q^{-2}\int_{-\infty}^{v_0}e^{4u-2|k||v-u|}\,du
 =\frac{q^2}{2|k|+4}
   +q^{2|k|-2}\int_v^{v_0}e^{(4-2|k|)u}\,du
 \le \frac{Cr_n^2}{|k|}.
\]
For $|k|=2$, use $q^2\log(r_n/q)\le Cr_n^2$.
Thus an energy solution of
$(-\partial_v^2+k^2)u_k
=\partial_v(q^2\mathcal F)+ikq^2\mathcal G+q^2\mathcal H$
satisfies, for $q<r_n/4$,
\begin{equation}\label{cubic:radial-green}
 \frac{|u_k(q)|}{q}
 \le Cr_n^{-1}|k|^{-1/2}
        \|u_k\|_{H^1_k(v_0-1,v_0)}
 +Cr_n|k|^{-1/2}
        (\|\mathcal F\|_2+\|\mathcal G\|_2
                       +|k|^{-1}\|\mathcal H\|_2).
\end{equation}
Here $\|u\|_{H^1_k(J)}^2=\int_J(|u'|^2+k^2|u|^2)\,dv$,
and the source norms are on $(-\infty,v_0)$ with measure $dv$.
The boundary contribution is $u_k(v_0)(q/r_n)^{|k|}$;
the bound follows from
$|u_k(v_0)|\le C|k|^{-1/2}\|u_k\|_{H^1_k(v_0-1,v_0)}$
and Cauchy--Schwarz for the Green integral.

Recall explicitly the energy norm used in \eqref{r2:angular-g}:
\[
 \|h\|_{\mathscr E(q<r_n)}^2
 =\int_{-\infty}^{\log r_n}\int_{\mathbb T}
       (|\partial_vh|^2+|\partial_\theta h|^2)\,d\theta\,dv.
\]
It follows from \eqref{r2:tangential-high},
$\Gamma_f=(\psi_{f,+}-\psi_{f,-})/(2\pi)$,
\eqref{r2:high-regular-density}, and
\eqref{r2:compact-regular-density} that
\begin{equation}\label{cubic:noncritical-energy}
 \sum_{j=0}^2\|\partial_\theta^j\Gamma_f\|_{\mathscr E(q<r_n)}
       \le Cn^{-1}M,\qquad
 \sum_{j=0}^1\|\partial_\theta^jR_f\|_{\mathscr E(q<r_n)}
       \le C_\alpha n^{-1-\alpha/2}M.
\end{equation}

For the indices in \eqref{cubic:noncritical-denominators},
the boundary difference equation is
\begin{equation}\label{cubic:noncritical-density-equation}
 (-\partial_v^2+k^2)\Gamma_{f,k}
 =[\operatorname{div}((\mathsf G-I)\nabla\Gamma_f)]_k
       +q^2\xi_0T_k\Gamma_{f,k}.
\end{equation}
The point source in the index $n$ is excluded by $\rho_n=0$
on the first interval and lies outside the disk on the second.
By Lemma~\ref{geom:metric} and angular integration by parts,
\begin{equation}\label{cubic:noncritical-flux-bounds}
 \begin{aligned}
 \|[q^{-2}(\mathsf G-I)\nabla\Gamma_f]_k\|_2
       &\le C|k|^{-2}n^{-1}M,&
 \|\Gamma_{f,k}\|_2&\le C|k|^{-3}n^{-1}M,\\
 \|[q^{-2}(\mathsf G-I)\nabla R_f]_k\|_2
       &\le C_\alpha|k|^{-1}n^{-1-\alpha/2}M,&
 \|R_{f,k}\|_2&\le C_\alpha|k|^{-2}n^{-1-\alpha/2}M.
 \end{aligned}
\end{equation}
The zero angular coefficients enter only through their radial
derivatives, which are included in \eqref{cubic:noncritical-energy}.
Applying \eqref{cubic:radial-green} to
\eqref{cubic:noncritical-density-equation}, and using
$r_n\sim n^{-1/2}$, we obtain
\[
 |\Gamma_{f,k}|
 \le CqM\left(|k|^{-5/2}n^{-1/2}
       +|k|^{-5/2}n^{-3/2}
       +\frac{|k|^{-7/2}n^{-3/2}}{1+|k-n|}\right).
\]
In particular,
\begin{equation}\label{cubic:noncritical-pointwise}
 |\Gamma_{f,k}(q,\lambda)|
       \le Cq|k|^{-5/2}n^{-1/2}M.
\end{equation}

\smallskip\noindent
\emph{Step 3: The logarithmic source and the regular derivative.}
By \eqref{spectral:decomposition},
$\partial_\lambda\Gamma_f=R_f+\Lambda_f$ in the chosen disk.
Differentiating \eqref{cubic:noncritical-density-equation}
distributionally, we obtain
\begin{equation}\label{cubic:noncritical-remainder-equation}
 \begin{aligned}
 (-\partial_v^2+k^2)R_{f,k}
 ={}&[\operatorname{div}((\mathsf G-I)\nabla R_f)]_k
          +q^2\xi_0T_kR_{f,k}\\
 &+q^2\xi_0(\partial_\lambda T_k)\Gamma_{f,k}
  +[\operatorname{div}(\mathsf G\nabla\Lambda_f)]_k
                          +q^2\xi_0T_k(\Lambda_f)_k.
 \end{aligned}
\end{equation}
The two terms containing $\Lambda_f$ vanish on the intervening
intervals. We estimate their sum on a resonant threshold interval
by explicitly cancelling the radial pole. Consider
one sign in the correction defined before
Proposition~\ref{prop:spectral-bounds}. Write
\[
 \begin{gathered}
 \Lambda_{f,n}^{\pm}=a_{n,\pm}A(\theta)\ell(v),\qquad
 A(\theta)=\frac{e^{in\theta}}{\mathsf G_{vv}(\delta,\theta)},
 \qquad \delta^2=8(1-\lambda/n),\\
 a_{n,\pm}=\frac{\pm8\xi_0(\delta)}{\pi i n^2\delta^2}
                     g_{f,\pm,n}(\delta,\lambda),\qquad
 \ell(v)=\chi_{\log}(v-\log\delta)
                   \log(v-\log\delta\mp i0).
 \end{gathered}
\]
Here $\chi_{\log}$ denotes the fixed cutoff in that definition,
supported in $(-1/2,1/2)$; it is distinct from the spectral
cutoff $\chi$ in this lemma. It follows from
\eqref{spectral:critical-amplitude} that
$|a_{n,\pm}|\le Cn^{-9/2}M$.

For either sign, suppress the sign on $a_n$ and set
\begin{equation}\label{cubic:noncritical-log-fluxes}
 \begin{aligned}
 q^2\mathcal F_{n,k}
 &=a_n\left([\mathsf G_{vv}A]_k\ell'
       +([\mathsf G_{v\theta}\partial_\theta A]_k
                     +ik[\mathsf G_{v\theta}A]_k)\ell\right),\\
 q^2\mathcal G_{n,k}
 &=a_n\left([\mathsf G_{\theta\theta}\partial_\theta A]_k
                  -[(\partial_v\mathsf G_{v\theta})A]_k\right)\ell,\\
 \mathcal H_{n,k}&=a_n\xi_0T_k[A]_k\ell.
 \end{aligned}
\end{equation}
Symmetry of $\mathsf G$ and the product rule imply
\begin{equation}\label{cubic:noncritical-log-source}
 [\operatorname{div}(\mathsf G\nabla\Lambda_{f,n}^{\pm})]_k
                  +q^2\xi_0T_k(\Lambda_{f,n}^{\pm})_k
 =\partial_v(q^2\mathcal F_{n,k})
                 +ikq^2\mathcal G_{n,k}+q^2\mathcal H_{n,k}.
\end{equation}
For $k\ne n$, the radial pole cancels since
$[\mathsf G_{vv}(\delta,\theta)A]_k=0$.
More explicitly, with $u=v-\log\delta$,
\[
 [\mathsf G_{vv}(q,\theta)A]_k\ell'
 =\left(\int_0^1
     [(\partial_v\mathsf G_{vv})(\delta e^{tu},\theta)A]_k\,dt\right)
   \bigl(\chi_{\log}(u)+u\chi_{\log}'(u)\log(u\mp i0)\bigr).
\]
This uses the distributional identity $u/(u\mp i0)=1$.
The mixed pole is included in the radial flux through
$ik[\mathsf G_{v\theta}A]_k\ell'
=\partial_v(ik[\mathsf G_{v\theta}A]_k\ell)
-ik[(\partial_v\mathsf G_{v\theta})A]_k\ell$.
On this support $q\sim\delta$. Lemma~\ref{geom:metric}
and angular integration by parts bound the coefficients
without $\partial_\theta A$ by
$C_N\delta^2(1+|k-n|)^{-N}$, and those with
$\partial_\theta A$ by $C_Nn\delta^2(1+|k-n|)^{-N}$.
Since the displayed logarithmic factors belong to $L^2(du)$,
\begin{equation}\label{cubic:noncritical-log-flux-bound}
 \|\mathcal F_{n,k}\|_2+\|\mathcal G_{n,k}\|_2
             +|k|^{-1}\|\mathcal H_{n,k}\|_2
 \le C_Nn^2(1+|k-n|)^{-N}|a_n|.
\end{equation}

Apply \eqref{cubic:radial-green} to
\eqref{cubic:noncritical-remainder-equation}, using
\eqref{cubic:noncritical-flux-bounds} and
\eqref{cubic:noncritical-log-flux-bound}.
The boundary and metric terms are bounded by
$C_\alpha q|k|^{-3/2}n^{-1/2-\alpha/2}M$.
The two potential terms and the logarithmic source are bounded by
\[
 C_\alpha qM\left(
 \frac{|k|^{-5/2}n^{-3/2-\alpha/2}}{1+|k-n|}
 +\frac{|k|^{-7/2}n^{-3/2}}{(1+|k-n|)^2}
 +|k|^{-1/2}n^{-3}(1+|k-n|)^{-N}\right).
\]
Since $|k|\le n+|k-n|$, choosing $N\ge2$ we get that
\begin{equation}\label{cubic:noncritical-R}
 |R_{f,k}(q,\lambda)|
       \le C_\alpha q|k|^{-3/2}n^{-1/2-\alpha/2}M.
\end{equation}

\smallskip\noindent
\emph{Step 4: Integration of the logarithms and summation.}
For $k\ne n$, $\mathsf G_{vv}(\delta,\theta)^{-1}=1+O(\delta^2)$
implies
\begin{equation}\label{r2:offdiagonal-jet}
 |(\Lambda_{f,n})_k(q,\lambda)|
 \le C_N\delta^2n^{-9/2}(1+|k-n|)^{-N}
  \mathbf1_{\{|\log(q/\delta)|<1/2\}}
       (1+|\log|\log(q/\delta)||)M.
\end{equation}
At fixed $q$, set $u=\log(q/\delta)$. Then
$\lambda=n(1-q^2e^{-2u}/8)$ and
$d\lambda=(nq^2/4)e^{-2u}\,du$. Hence
\begin{equation}\label{cubic:offdiagonal-log-integral}
 \int|(\Lambda_{f,n})_k(q,\lambda)|\,d\lambda
 \le C_Nq^4n^{-7/2}(1+|k-n|)^{-N}M.
\end{equation}
The integral is over the threshold interval with index $n$.

Splitting the sums into $0<|k|<n/2$, $n/2\le|k|\le2n$,
and $|k|>2n$, we obtain
\[
 \sum_{k\ne0}\frac{|k|^{-3/2}}{1+|k-n|}\le Cn^{-1},\qquad
 \sum_{k\ne0}\frac{|k|^{-1/2}}{1+|k-n|}
       \le Cn^{-1/2}\log(2+n).
\]
Thus \eqref{cubic:noncritical-pointwise} and
\eqref{cubic:noncritical-R} imply, on either interval,
\[
 \begin{aligned}
 \sum_{\substack{k\ne0\\\rho_k(\lambda)\ne0}}
       \frac{|k||\Gamma_{f,k}|}{1+|k-n|}
       &\le Cq n^{-3/2}M,\\
 \sum_{\substack{k\ne0\\\rho_k(\lambda)\ne0}}
       \frac{|k||R_{f,k}|}{1+|k-n|}
       &\le C_\alpha q n^{-1-\alpha/2}\log(2+n)M.
 \end{aligned}
\]
Both right-hand sides are summable in $n$. For the logarithms,
\eqref{cubic:offdiagonal-log-integral} implies, when $N>2$,
\[
 \sum_{n\ge1}\sum_{k\ne0,n}\frac{|k|}{1+|k-n|}
       \int|(\Lambda_{f,n})_k|\,d\lambda
 \le C_Nq^4M\sum_{n\ge1}n^{-7/2}
       \sum_{k\ne0}\frac{|k|}{(1+|k-n|)^{N+1}}
 \le Cq^4M.
\]
Indeed, the inner sum is at most $Cn$.
These estimates bound the integrand in
\eqref{cubic:noncritical-W11} and its derivative: differentiating
the denominator replaces $(1+|k-n|)^{-1}$ by
$(1+|k-n|)^{-2}$; the cutoff derivatives contribute bounded
factors $\rho_k'$ and $q^2\eta'(q^2\lambda)$.

\smallskip\noindent
\emph{Step 5: Bounded spectral parameters and integration by parts.}
Choose $0<\lambda_0<\kappa/2$.
On $q<\sqrt{\kappa_1}/4$, $T$ and $\partial_\lambda T$ are bounded
on each physical Sobolev space for $|\lambda|\le1-\kappa_1/2$.
Using the equations
\[
 (\Delta+T)\Gamma_f=0,\qquad
 (\Delta+T)\partial_\lambda\Gamma_f
              =-(\partial_\lambda T)\Gamma_f,
\]
interior elliptic estimates, \eqref{tailD:F-W14}, and
\eqref{tailD:zero-derivative} imply, on $q<\sqrt{\kappa_1}/16$,
\[
 \|\Gamma_f(\lambda)\|_{C^4}
 +\|\partial_\lambda\Gamma_f(\lambda)\|_{C^4}
       \le Ce^{-c/|\lambda|}M,\qquad 0<|\lambda|<\lambda_0.
\]
For $\lambda_0\le|\lambda|\le1-\kappa_1/2$, the same norms
are bounded by $CM$: use Proposition~\ref{global:LAP},
\eqref{tailD:zero-derivative} below
$\kappa$, and \eqref{r2:compact-regular-density} above
$\kappa$, where $D_\lambda=\partial_\lambda$ and $\Lambda_f=0$
on the larger disk. The choice $c_{\rm low}<\kappa_1/10000$
ensures that every nonempty angular sum in the lemma is
evaluated in $q<\sqrt{\kappa_1}/16$.
By the coordinate expansion in Lemma~\ref{geom:metric},
$\|\partial_\theta^3h(q,\cdot)\|_\infty\le Cq\|h\|_{C^4}$.
Consequently,
\[
 \sum_{k\ne0}|k|
       (|\Gamma_{f,k}|+|\partial_\lambda\Gamma_{f,k}|)
 \le Cq\bigl(\|\Gamma_f\|_{C^4}
                  +\|\partial_\lambda\Gamma_f\|_{C^4}\bigr)
             \sum_{k\ne0}|k|^{-2}.
\]
This proves the remaining $L^1_\lambda$ bounds. The bound \eqref{r2:nonlocal-integral} follows from integration by parts.
\end{proof}

\subsection{The pointwise estimate away from the elliptic points}

For a spectral cutoff $\chi$, write
\[
 \omega_\chi(t)=\Delta\int_{\mathbb R}e^{-it\lambda}
                      \chi(\lambda)\Gamma_f(\lambda)\,d\lambda.
\]
On the fixed exterior $K$ of \eqref{r2:away-W1q}, it follows
from \eqref{sd:jump} that, for $k\ne0$,
\begin{equation}\label{sd:localized-vorticity}
\begin{aligned}
 (\omega_\chi)_k(t,r)
 &=\frac12\chi(kb(r))e^{-itkb(r)}
   \bigl[g_{f,+,k}(r,kb(r))
        +g_{f,-,k}(r,kb(r))\bigr]\\
 &\quad+\operatorname{PV}\int_{\mathbb R}e^{-it\lambda}
   \frac{kb(r)\chi(\lambda)\Gamma_{f,k}(r,\lambda)}
        {\lambda-kb(r)}\,d\lambda.
\end{aligned}
\end{equation}

\emph{Compact intervals away from zero.}
Fix $\chi\in C_c^\infty((c,C))$. For the indices
$c/2\le |kb(r)|\le2C+1$, the identity $b(r)=2\pi/S(r)$ implies
$|k|^{-1}\le4\pi/(cS(r))$.
Apply Lemma~\ref{r2:orbit-coefficients}\textup{(i)} to
$g_{f,\pm}=f+i\lambda\psi_{f,\pm}$, to $\Gamma_f$, and to
its spectral differences. It follows from \eqref{r2:away-holder}
that
\[
 |\chi(kb(r))g_{f,\pm,k}(r,kb(r))|
 +\bigl\|kb(r)\chi(\cdot)\Gamma_{f,k}(r,\cdot)
                                      \bigr\|_{C^\gamma(\mathbb R)}
 \le \frac{C_\chi M}{S(r)}.
\]
Here the amplitude is extended by zero outside $(c,C)$.
Applying \eqref{r2:small-amplitude-PV} with
$A=B=C_\chi M/S(r)$ and fixed $L$ bounds each principal-value
integral by $C_\chi M/S(r)$, uniformly in $t$. Since
\[
 \#\{k:c/2\le |kb(r)|\le2C+1\}
 \le \frac{2C+1}{\pi}S(r),
\]
the sum of the absolute values of these integrals and the
transported terms is at most $C_\chi M$, uniformly in $t$
and in the regular streamlines contained in $K$.
For $|kb|<c/2$ use $|(L_0\Gamma_f)_k|\le CM/S$; for
$|kb|>2C+1$ use $|\Gamma_{f,k}|\le CM S^{-1}|kb|^{-1-\beta}$.
Indeed, it follows from the spacing $2\pi/S$ that
\[
 S^{-1}\#\{k:0<|kb|<c/2\}\le C_c,\qquad
 S^{-1}\sum_{|kb|>2C+1}|kb|^{-1-\beta}\le C_{C,\beta}.
\]
Negative spectral intervals follow by reflection. Hence
\begin{equation}\label{r2:compact-vorticity-Linfty}
 \sup_t\|\omega_\chi(t)\|_{L^\infty(K)}\le C_\chi M,
       \qquad \chi\in C_c^\infty(\mathbb R\setminus\{0\}).
\end{equation}

\emph{Large spectral parameters.}
Choose an even $\chi_{\rm low}\in C_c^\infty(\mathbb R)$ equal
to one on $[-4\Lambda_0,4\Lambda_0]$, and set
$\chi_{\rm high}=1-\chi_{\rm low}$.

We first estimate the first term on the right-hand side of
\eqref{sd:localized-vorticity}, with $\chi=\chi_{\rm high}$. Its coefficients are
the boundary values $g_{f,\pm,k}(r,kb(r))$, multiplied by
$\chi_{\rm high}(kb(r))e^{-itkb(r)}$. The resolvent identity
implies
\begin{equation}\label{hp:amplitude-expansion}
 f+i\lambda\psi_{f,\pm}
 =f+\Delta^{-1}f+(i\lambda)^{-1}\Delta^{-1}Lf+r_\pm,
 \qquad r_\pm=(i\lambda)^{-1}\psi_{L^2f,\pm}.
\end{equation}
For $0<\beta<1/8$, interpolation in \eqref{hp:base-H1}, applied
with datum $L^2f$, implies
\[
 \|r_\pm\|_{C^\beta(K)}
 \le C\|r_\pm\|_{H^1(K)}^{3/4}
       \|r_\pm\|_{W^{1,4}(K)}^{1/4}
 \le C|\lambda|^{-7/4}M.
\]
Since $b=2\pi/S$, Lemma~\ref{r2:orbit-coefficients}\textup{(i)}
therefore implies
\[
 \sum_{k\ne0}
 \left|\chi_{\rm high}(kb(r))(r_\pm)_k(r,kb(r))\right|
 \le \frac{CM}{S(r)}
       \sum_{|kb(r)|\ge4\Lambda_0}|kb(r)|^{-7/4}
 \le CM.
\]
The remaining contribution from the first term in
\eqref{hp:amplitude-expansion} are
\[
 e^{-tL_0}\left\{
 \chi_{\rm high}(-iL_0)(f+\Delta^{-1}f)
 +\left[\frac{\chi_{\rm high}(\lambda)}{i\lambda}\right]
       (-iL_0)\Delta^{-1}Lf\right\}.
\]
Their $L^\infty$ norm is bounded by $CM$: the transport group
is an $L^\infty$ isometry,
$\chi_{\rm high}(-iL_0)=I-\chi_{\rm low}(-iL_0)$, and
\[
 \widehat{\chi_{\rm low}}\in L^1,\qquad
 \frac{\chi_{\rm high}(\lambda)}{\lambda}\in H^1(\mathbb R)
 \quad\Longrightarrow\quad
 \mathcal F\!\left(\frac{\chi_{\rm high}(\lambda)}{\lambda}\right)
 \in L^1.
\]

For the principal-value integral, decompose the positive spectral
region using smooth cutoffs $\chi_\Lambda$, supported in
$(\Lambda-1,\Lambda+1)$ and uniformly bounded in $C^1$, such that
\[
 \sum_{\substack{\Lambda\in\mathbb N\\\Lambda\ge2\Lambda_0}}
       \chi_\Lambda(\lambda)=\chi_{\rm high}(\lambda),
 \qquad \lambda>0.
\]
Fix one such $\Lambda$. For
$\Lambda/2\le|kb(r)|\le2\Lambda$, it follows from
Lemma~\ref{r2:orbit-coefficients}\textup{(i)},
\eqref{hp:G-Holder}, and \eqref{hp:spectral-holder} that
\[
 \begin{aligned}
 \|kb(r)\chi_\Lambda\Gamma_{f,k}(r,\cdot)\|_{L^\infty_\lambda}
 &\le \frac{CM}{S(r)}\Lambda^{-3/2},\\
 \|kb(r)\chi_\Lambda\Gamma_{f,k}(r,\cdot)\|_{C^\gamma_\lambda}
 &\le \frac{CM}{S(r)}\Lambda.
 \end{aligned}
\]
For $|kb(r)-\Lambda|\le2$, apply
\eqref{r2:small-amplitude-PV} with these upper bounds for
$A,B$ and fixed $L$. Each integral is bounded by
$CM S(r)^{-1}\Lambda^{-3/2}\log(2+\Lambda)$ and there are at most $CS(r)$ such indices.

For $\Lambda/2\le|kb(r)|\le2\Lambda$ and
$|kb(r)-\Lambda|>2$, the denominator is separated from zero.
Absolute integration and the spacing $2\pi/S(r)$ imply
\[
 \frac{CM\Lambda^{-3/2}}{S(r)}
 \sum_{\substack{\Lambda/2\le|kb(r)|\le2\Lambda\\
                 |kb(r)-\Lambda|>2}}
       \frac1{|kb(r)-\Lambda|}
 \le CM\Lambda^{-3/2}\log(2+\Lambda).
\]

It remains to consider $0<|kb(r)|<\Lambda/2$ and
$|kb(r)|>2\Lambda$. By
Lemma~\ref{r2:orbit-coefficients}\textup{(iii)},
$L_0(i\lambda\Gamma_f)$ vanishes at the hyperbolic fixed points.
Applying part~\textup{(ii)} of that lemma and
\eqref{hp:G-Holder}, we obtain, for $|\lambda-\Lambda|<1$,
\[
 |(i\lambda L_0\Gamma_f)_k(r,\lambda)|
 \le \frac{CM}{S(r)}\Lambda^{-3/4}
             \min\{1,|kb(r)|^{-\beta}\}.
\]
Since
$|kb(r)\Gamma_{f,k}|=|\lambda|^{-1}
|(i\lambda L_0\Gamma_f)_k|$, the sum of the absolute integrals
over these indices is bounded by
\[
 \begin{aligned}
 CM\Lambda^{-3/4}\left\{
 \Lambda^{-2}\int_0^{\Lambda/2}\min(1,u^{-\beta})\,du
 +\Lambda^{-1}\int_{2\Lambda}^\infty
                  \frac{u^{-\beta}}{\Lambda+u}\,du\right\}
 \le C_\beta M\Lambda^{-7/4-\beta}.
 \end{aligned}
\]
The negative spectral region satisfies the same estimates by
complex conjugation, applied to the datum $\overline f$.
Finally,
\[
 \sum_{\Lambda\ge2\Lambda_0}
 \left[\Lambda^{-3/2}\log(2+\Lambda)
             +\Lambda^{-7/4-\beta}\right]<\infty.
\]
Combining this sum with the transported-term estimates proves
\begin{equation}\label{hp:conclusion}
 \sup_t\|\omega_{\chi_{\rm high}}(t)\|_{L^\infty(K)}\le CM.
\end{equation}
\emph{A neighborhood of zero.}
Set $f_*=f-p_\infty$ and choose
$\chi_0\in C_c^\infty((-\lambda_*/4,\lambda_*/4))$ equal to
one near zero, where $\lambda_*$ is as in
Lemma~\ref{pointwise:elliptic-inputs}.
Since
\[
 g_{f,\pm}=f_*+i\lambda\psi_{f_*,\pm},
\]
the contribution of $f_*$ to the first term in
\eqref{sd:localized-vorticity} is
$\chi_0(-iL_0)e^{-tL_0}f_*$. The transport isometry and
$\widehat{\chi_0}\in L^1$ imply
\begin{equation}\label{zp:free-bound}
 \sup_t\|\chi_0(-iL_0)e^{-tL_0}f_*\|_\infty
 \le C\|\widehat{\chi_0}\|_1\|f_*\|_\infty\le CM.
\end{equation}
We estimate the remaining contribution
\[
 \omega_{\chi_0,\mathrm{reg}}
 =\omega_{\chi_0}-\chi_0(-iL_0)e^{-tL_0}f_*.
\]

Decompose the positive spectral interval using smooth cutoffs
$\chi_\delta$, where $\delta=2^{-j}\lambda_*/4$, $j\ge0$, such that
\[
 \operatorname{supp}\chi_\delta\subset(\delta/2,2\delta),
 \qquad
 \|\chi_\delta(\delta\,\cdot)\|_{C^1}\le C,
 \qquad
 \sum_\delta\chi_\delta(\lambda)=\chi_0(\lambda)
 \quad(\lambda>0).
\]
Use the reflected cutoffs on the negative interval.
For $\delta/2\le|kb(r)|\le2\delta$, it follows from
\eqref{zp:individual-holder} and
Lemma~\ref{r2:orbit-coefficients}\textup{(i)} that
\[
 |kb(r)\psi_{f_*,\pm,k}(r,kb(r))|
 \le C\delta\bigl(S(r)^{-1}+|k|^{-1}\bigr)M
 \le \frac{CM}{S(r)}.
\]
Since $b(r)=2\pi/S(r)$, there are at most $C\delta S(r)$
such indices. Thus the correction to the explicit sum in
\eqref{sd:localized-vorticity} is bounded by $C\delta M$
at this scale.

For the principal-value integral with cutoff $\chi_\delta$,
first consider $\delta/4\le|kb(r)|\le4\delta$.
Lemma~\ref{r2:orbit-coefficients}\textup{(i)} and
\eqref{zp:scaled-holder} imply
\[
 \bigl\|kb(r)\chi_\delta(\delta\,\cdot)
          \Gamma_{f,k}(r,\delta\,\cdot)\bigr\|_{C^\gamma(\mathbb R)}
 \le \frac{Ce^{-c/\delta}M}{S(r)}.
\]
After $\lambda=\delta u$, the integral becomes
\[
 \operatorname{PV}\int_{\mathbb R}e^{-it\delta u}
 \frac{kb(r)\chi_\delta(\delta u)\Gamma_{f,k}(r,\delta u)}
      {u-kb(r)/\delta}\,du.
\]
It follows from \eqref{r2:small-amplitude-PV} that its absolute
value is at most $Ce^{-c/\delta}M/S(r)$.
Summation over the at most $C\delta S(r)$ indices therefore
produces the bound $C\delta e^{-c/\delta}M$.

For the remaining indices, parts~\textup{(ii)--(iii)} of
Lemma~\ref{r2:orbit-coefficients} and
\eqref{zp:scaled-holder} imply
\[
 |kb(r)\Gamma_{f,k}(r,\lambda)|
 =|(L_0\Gamma_f)_k(r,\lambda)|
 \le \frac{Ce^{-c/\delta}M}{S(r)}
        \min\{1,|kb(r)|^{-\beta}\},
 \qquad \delta/2<\lambda<2\delta.
\]
Here part~\textup{(iii)} ensures that $L_0\Gamma_f$ vanishes
at the hyperbolic fixed points, as required in part~\textup{(ii)}.
If $0<|kb(r)|<\delta/4$, then
$|\lambda-kb(r)|\ge\delta/4$ on $\operatorname{supp}\chi_\delta$.
Consequently,
\[
 \sum_{0<|kb(r)|<\delta/4}
 \int_{\mathbb R}
 \left|\frac{kb(r)\chi_\delta(\lambda)\Gamma_{f,k}(r,\lambda)}
              {\lambda-kb(r)}\right|\,d\lambda
 \le C\delta e^{-c/\delta}M.
\]
If $|kb(r)|>4\delta$, then
$|\lambda-kb(r)|\ge |kb(r)|/2$. Absolute integration and the
spacing $2\pi/S(r)$ bound the sum over these indices by
\[
 CM\delta e^{-c/\delta}
 \left(\int_{4\delta}^1\frac{du}{u}
       +\int_1^\infty u^{-1-\beta}\,du\right)
 \le C_\beta M\delta e^{-c/\delta}(1+|\log\delta|).
\]
The negative spectral interval satisfies the same estimates
by complex conjugation, applied to $\overline f$.

The angular sums converge absolutely on each regular streamline.
Their uniform bounds identify the resulting functions with the
$L^2$ spectral representation by dominated convergence.
The bounds $C\delta M$ and
$C_\beta\delta e^{-c/\delta}(1+|\log\delta|)M$ are summable
over the dyadic scales, so
\begin{equation}\label{zp:zero-vorticity-Linfty}
 \sup_t\|\omega_{\chi_0,\mathrm{reg}}(t)\|_{L^\infty(K)}
 \le CM.
\end{equation}
Combining \eqref{zp:free-bound},
\eqref{zp:zero-vorticity-Linfty},
\eqref{r2:compact-vorticity-Linfty}, and \eqref{hp:conclusion},
we conclude that
\[
 \sup_t\|\omega(t)-p_\infty\|_{L^\infty(K)}\le CM.
\]
\subsection{Proof of Theorem~\ref{th1:Int}}
\begin{proof}
It follows from Propositions~\ref{zero:constraints} and
\ref{qual:theorem} that $p_\infty\in E$ exists. By
Theorem~\ref{r2:global-damping},
$\psi_c=\Delta^{-1}(\omega-p_\infty)$ satisfies
\begin{equation}\label{synthesis:full-damping}
 \|\psi_c(t)\|_2\le C(1+t)^{-2}M,\qquad
 \|\psi_c(t)\|_{H^1}\le C(1+t)^{-1}M.
\end{equation}
Since $Q$ commutes with $\Delta^{-1}$ and annihilates $p_\infty$,
it follows from \eqref{synthesis:full-damping} that
\eqref{Intm7} and \eqref{Intquadratic} hold.

Since $\psi_{p_\infty}(z)=-p_\infty/z$, subtracting $p_\infty$
does not change $g_f$ or $\Gamma_f$. In a center disk,
it follows from \eqref{sd:jump} that, for $k\ne0$,
\begin{equation}\label{cubic:vorticity-formula}
 \begin{split}
 (\omega-p_\infty)_k(t,q)
 &=\frac12e^{-itkb(q)}
      [g_{f,+,k}+g_{f,-,k}](q,kb(q))\\
 &\quad+\operatorname{PV}\int e^{-it\lambda}
       \frac{kb(q)\Gamma_{f,k}(q,\lambda)}{\lambda-kb(q)}\,d\lambda.
 \end{split}
\end{equation}
The zeroth coefficient vanishes because $\omega-p_\infty\in Y$.
With $c_0,c_1,q_0$ as in Lemmas~\ref{cubic:profiles} and
\ref{cubic:first-mode}, choose the threshold cutoff $\chi$
in Lemma~\ref{cubic:noncritical} to be even, with
\[
 \operatorname{supp}\chi\subset(-\kappa,\kappa),\qquad
 0<\kappa<\tfrac12,\qquad
 \sqrt{8\kappa}<\min(c_1/2,q_0/8).
\]
Next fix $\eta$ as in that lemma and choose $c_{\rm spec}>0$
such that $\eta(s)=1$ for $|s|\le c_{\rm spec}$. Choose
\[
 0<c_{\rm low}<\min\{c_0^2/4,c_1^2,q_0^2/64,c_{\rm spec}/2\}.
\]
It follows from
\eqref{cubic:first-trace} and \eqref{cubic:higher-trace} that
\begin{equation}\label{synthesis:depleted-profile}
 \sum_{0<|k|q^2<c_{\rm low}}
       |g_{f,\pm,k}(q,kb(q))|
 \le CMq^3\left(1+\sum_{|k|\ge2}|k|^{-2}\right)
 \le CMq^3.
\end{equation}
It follows from Proposition~\ref{r2:high-output-tail} that
the modes $|k|q^2\ge c_{\rm low}$ and the spectral contribution
containing $1-\eta(q^2\lambda)$ are bounded by $CMq^3$.

It remains to estimate the part of the integral in
\eqref{cubic:vorticity-formula} near $\lambda=k$.
We have $\eta(q^2\lambda)=1$ whenever
$0<|k|q^2<c_{\rm low}$ and $\chi(\lambda-k)\ne0$.
For $0<|k|q^2<c_{\rm low}$ and $t\ge0$, we prove
\begin{equation}\label{synthesis:threshold-integral}
 \left|\operatorname{PV}\int_{\mathbb R}
       e^{-it\lambda}\chi(\lambda-k)
       \frac{kb(q)\Gamma_{f,k}(q,\lambda)}
            {\lambda-kb(q)}\,d\lambda\right|
 \le CM\left(|k|^{-3/2}q^3
               +\frac{|k|^{-2}q}{1+t}\right).
\end{equation}

We first take $k\ge1$ and $\lambda=k(1-\delta^2/8)$.
The cutoffs defining $G_k$ in
Lemmas~\ref{cubic:profiles} and \ref{cubic:first-profiles}
equal one on the region under consideration. Consequently,
\begin{equation}\label{synthesis:threshold-numerator}
 \chi(\lambda-k)kb(q)\Gamma_{f,k}(q,\lambda)
 =\frac{b(q)\chi(-k\delta^2/8)}
        {2\pi i(1-\delta^2/8)}G_k(q,\delta).
\end{equation}
The factor multiplying $G_k$ has bounded
$(\delta\partial_\delta)^j$ derivatives for $j\le2$:
indeed, $\delta\partial_\delta\chi(-k\delta^2/8)
=-(k\delta^2/4)\chi'(-k\delta^2/8)$ and
$k\delta^2\le8\kappa$ on its support.

Choose a fixed smooth function $\varphi$, equal to one on
$[0,1/3]$ and zero on $[1/2,\infty)$, and use
\[
 \varphi(\delta/q)+[1-\varphi(\delta/q)-\varphi(q/\delta)]
                      +\varphi(q/\delta)=1.
\]

\smallskip\noindent
\emph{The region $\delta<q/2$.}
Since $|d\lambda|/|\lambda-kb(q)|
=2\delta\,d\delta/|q^2-\delta^2|$, it follows from
\eqref{cubic:first-profile-outer},
\eqref{cubic:profile-outer}, and
\eqref{synthesis:threshold-numerator} that this contribution
is bounded by
\[
 C\int_0^{q/2}
       \frac{\delta|G_k(q,\delta)|}{q^2}\,d\delta
 \le CMk^{-2}q^3,\qquad k\ge1.
\]

\smallskip\noindent
\emph{The region $q/3<\delta<3q$.}
Set
\[
 v=\log(q/\delta),\qquad
 \xi=\frac{\lambda-kb(q)}{kq^2}
     =\frac{1-e^{-2v}}8.
\]
Recall that
$G_k(q,qe^{-v})=q^3e^{-3v}U_k(\log q-v,v)$.
The product rule, the one-dimensional Sobolev embedding, and
integration in $w$ imply, by \eqref{cubic:profile-core} and
\eqref{cubic:first-profile-core},
\[
 \begin{aligned}
 &\left\|\chi(-kq^2e^{-2v}/8)U_k(\log q-v,v)
                   \right\|_{C^{1/2}(-\log3,\log3)}\\
 &\qquad\le C\sup_w\sum_{j=0}^1
    \left\|\partial_w^j\bigl[
      \chi(-ke^{2w}/8)U_k(w,\cdot)\bigr]\right\|_{H^1_v(-2,2)}
 \le CMk^{-3/2}.
 \end{aligned}
\]
The cutoff product is extended by zero in $w$ as above.
Since $v\mapsto\xi$ is a smooth diffeomorphism on
$|v|<\log3$, it follows from
\eqref{synthesis:threshold-numerator} that
\begin{equation}\label{synthesis:threshold-holder}
 \left\|[1-\varphi(e^{-v})-\varphi(e^v)]
       \chi(\lambda-k)kb(q)\Gamma_{f,k}(q,\lambda)
                         \right\|_{C^{1/2}(\mathbb R_\xi)}
 \le CMk^{-3/2}q^3.
\end{equation}
In this norm $\lambda=kb(q)+kq^2\xi$,
$v=-\frac12\log(1-8\xi)$, and the expression is extended
by zero outside $|v|<\log3$.
Its support lies in a fixed interval. Applying
\eqref{r2:small-amplitude-PV} with
$\gamma=1/2$, $L=1$, and $A=B=CMk^{-3/2}q^3$
bounds the middle contribution by $CMk^{-3/2}q^3$.
The principal value is preserved by the affine change
from $\lambda$ to $\xi$.

\smallskip\noindent
\emph{The region $\delta>2q$.}
Choose $\rho\in C_c^\infty((1/2,2))$ such that
$\sum_{m\in\mathbb Z}\rho(\delta/(2^mq))=1$, and set $d=2^mq$.
Only $q<d<2\sqrt{8\kappa/k}$ contributes. In the measure
$d\delta/\delta$, the rational factor is
$2\delta^2/(q^2-\delta^2)$, whose first two logarithmic
derivatives are bounded on $\delta>2q$. For $t>0$, use
\[
 e^{-itk(1-\delta^2/8)}
 =\frac{4}{itk\delta^2}\,
       \delta\partial_\delta e^{-itk(1-\delta^2/8)}.
\]
Two integrations by parts for $k=1$, and one for $k\ge2$,
together with \eqref{cubic:first-profile-inner},
\eqref{cubic:profile-inner}, and the absolute estimates imply
\begin{equation}\label{synthesis:threshold-dyadic}
 \begin{aligned}
 &\left|\int_{k-\kappa}^{k}e^{-it\lambda}
       \rho(\delta/d)\varphi(q/\delta)\chi(\lambda-k)
       \frac{kb(q)\Gamma_{f,k}(q,\lambda)}
            {\lambda-kb(q)}\,d\lambda\right|\\
 &\qquad\le CM
 \begin{cases}
 qd^2\min\{1,(td^2)^{-2}\},&k=1,\\
 k^{-1}q^{5/4}d^{7/4}\min\{1,(tkd^2)^{-1}\},&k\ge2.
 \end{cases}
 \end{aligned}
\end{equation}
The dyadic cutoff and $\chi$ eliminate the endpoint terms.
For $t\ge1$,
\[
 \begin{aligned}
 \sum_{d>q}d^2\min\{1,(td^2)^{-2}\}
   &\le Ct^{-1},\\
 \sum_{d>q}d^{7/4}\min\{1,(tkd^2)^{-1}\}
   &\le (tk)^{-1}\sum_{d>q}d^{-1/4}
    \le C(tk)^{-1}q^{-1/4}.
 \end{aligned}
\]
The first inequality follows by splitting at $d=t^{-1/2}$.
For $0\le t\le1$, summing the absolute estimates instead
implies the bound $CMqk^{-2}$, since
$\delta^2<8\kappa/k$ and this region is empty unless
$q<\sqrt{8\kappa/k}$.
Thus the contribution is at most $CMqk^{-2}/(1+t)$.

For $\lambda=k(1+\delta^2/8)$, replace $q^2-\delta^2$ by
$q^2+\delta^2$. On $q/3<\delta<3q$,
$|G_k|\le CMk^{-3/2}q^3$ and
$|d\lambda|/|\lambda-kb(q)|\le C\,d\delta/\delta$,
so absolute integration suffices; the other two estimates
are unchanged. Negative modes follow by applying the estimates
to $\overline f$ and using
$\Gamma_{f,-k}(q,-\lambda)
=\overline{\Gamma_{\overline f,k}(q,\lambda)}$ for $k>0$.
This proves \eqref{synthesis:threshold-integral}.
Summing it in $k$ and using
\eqref{r2:nonlocal-integral},
\eqref{synthesis:depleted-profile}, and
Proposition~\ref{r2:high-output-tail}, we obtain
\begin{equation}\label{synthesis:center-depletion}
 |\omega(t,x)-p_\infty(x)|
 \le CM\left(q^3+\frac{q}{1+t}\right).
\end{equation}

Recall that
$E=\operatorname{span}\{\sin x,\sin y,\cos x-\cos y\}$.
Each basis function vanishes at both elliptic points, so
\eqref{synthesis:full-damping} implies
\[
 |P\psi_c(t,x)|
 \le Cd(x)\|\psi_c(t)\|_2
 \le Cd(x)(1+t)^{-2}M.
\]
Since $\omega-p_\infty=\Delta\psi_c$, $P\Delta=-P$, and
$Qp_\infty=0$, we have
\[
 Q\omega=(\omega-p_\infty)-P(\omega-p_\infty)
        =(\omega-p_\infty)+P\psi_c.
\]
It follows from \eqref{synthesis:center-depletion} and
$q\sim d(x)$ that, in the center disks,
\[
 |Q\omega(t,x)|
 \le CM\left(q^3+\frac{q}{1+t}
                    +\frac{d(x)}{(1+t)^2}\right)
 \le CM\left(d(x)^3+\frac{d(x)}{1+t}\right).
\]

On $K$, it follows from
\eqref{r2:compact-vorticity-Linfty}, \eqref{hp:conclusion},
\eqref{zp:zero-vorticity-Linfty}, and \eqref{zp:free-bound}
that
\[
 \|Q\omega(t)\|_{L^\infty(K)}
 \le \|\omega(t)-p_\infty\|_{L^\infty(K)}
       +C\|\psi_c(t)\|_2
 \le CM.
\]
Since $d(x)\ge c_K>0$ on $K$, this also implies
\eqref{Intm9} there, after enlarging the constant.
For each fixed $t$, spatial continuity extends the estimate
to the elliptic points and the separatrix.
\end{proof}

\section*{Declaration on the use of AI}

The project began in early 2025 and is a continuation of the second author's research program on understanding the precise dynamics of smooth solutions to the two-dimensional incompressible Euler and slightly viscous Navier--Stokes equations. Many of the ideas in the paper, such as the identification of the main equation near the elliptic centers, the use of integration by parts to treat multiplier singularities and the method of adapted spectral derivatives, can be traced to earlier papers of the second author and collaborators. The main new difficulty for the dipole case is the coupling of the angular modes and the presence of the separatrix.

The authors introduced the streamline coordinates, analyzed the behavior of the multiplier singularity near the separatrix and the summability of the angular modes. After substantial work, we obtained sufficient bounds to show that the separatrix and angular coupling would not be an obstacle to either sharp inviscid damping or vorticity depletion estimates. We drafted a manuscript with preliminary results, including the expected main theorem, many detailed estimates on the streamline coordinates and the lower-order operators. However, it was clear to us that to obtain the sharp estimates in the paper, a significant amount of additional technical work was still needed. This was expected based on our experience working with shear flows and vortices, where the control of higher-order spectral derivatives, although not difficult, nonetheless requires heavy computations and the treatment of many strongly singular terms. Moreover, the dipole case is significantly more complicated on the technical level.

After GPT-6~Astra was released, we uploaded our manuscript to it and provided additional suggestions on how to handle spectral singularities. After several sessions, GPT-6~Astra produced a paper completing all the required detailed estimates in early October 2026. Under our detailed instructions and after several iterations we obtained a version where the proof was along the lines that we had in mind. We then carefully went through and revised the entire paper, with the assistance of GPT-6~Astra. Most of the arguments that GPT-6~Astra produced are expected, with the notable exception of the bounds on the vector field $V$ in Lemma~\ref{geom:common-field} and \eqref{quad:saddle-coefficients}, which played an important role when commuting the adapted spectral derivative with $\Delta$. The use of the explicit homogeneous solution in the proof of Lemma~\ref{cubic:first-mode} is also worth mentioning. Without the special solution, standard perturbative arguments would result in a logarithmic factor in the sharp cubic depletion rates for the vorticity field near the elliptic centers. However, this was already discovered in~\cite{Bed2} and also used in~\cite{JiaVortex} in the case of strictly monotonic vortices, and is therefore not completely new.

The proof of Theorem~\ref{Intm6} is a different story. In our preliminary manuscript, the absence of nonzero (embedded) eigenvalues was stated as an assumption. We carried out extensive numerical calculations both on the spectrum directly and on the long-time dynamics of the solution to the linearized flow. These numerical results provided convincing evidence, but not a mathematical proof, for the absence of embedded eigenvalues. It came as a surprise to us (at the time) that after several prompts, GPT-6~Astra was able to produce a complete proof of their absence. The AI-generated paper and a Lean verification can be found at~\cite{NoEmbeddedEigenvalues}.

We take full responsibility for the results and proofs in the paper.

\printbibliography

@article{Bassom,
  author       = {Bassom, Andrew P. and Gilbert, Andrew D.},
  title        = {The spiral wind-up of vorticity in an inviscid planar vortex},
  journal      = {Journal of Fluid Mechanics},
  shortjournal = {J. Fluid Mech.},
  volume       = {371},
  year         = {1998},
  pages        = {109--140},
  mrnumber     = {1650153},
  doi          = {10.1017/S0022112098001955},
}

@article{Bed2,
  author       = {Bedrossian, Jacob and Coti Zelati, Michele and Vicol, Vlad},
  title        = {Vortex axisymmetrization, inviscid damping, and vorticity depletion in the linearized 2{D} {E}uler equations},
  journal      = {Annals of PDE},
  shortjournal = {Ann. PDE},
  volume       = {5},
  year         = {2019},
  number       = {1},
  pages        = {Paper No. 4, 192 pp.},
  mrnumber     = {3987441},
  doi          = {10.1007/s40818-019-0061-8},
}

@article{BM,
  author       = {Bedrossian, Jacob and Masmoudi, Nader},
  title        = {Inviscid damping and the asymptotic stability of planar shear flows in the 2{D} {E}uler equations},
  journal      = {Publications Math\'{e}matiques de l'IH\'{E}S},
  shortjournal = {Publ. Math. Inst. Hautes \'{E}tudes Sci.},
  volume       = {122},
  year         = {2015},
  pages        = {195--300},
  mrnumber     = {3415068},
  doi          = {10.1007/s10240-015-0070-4},
}

@article{BeekieChenJia2026,
  author       = {Beekie, Rajendra and Chen, Shan and Jia, Hao},
  title        = {Uniform vorticity depletion and inviscid damping for periodic shear flows in the high {R}eynolds number regime},
  journal      = {Archive for Rational Mechanics and Analysis},
  shortjournal = {Arch. Ration. Mech. Anal.},
  volume       = {250},
  year         = {2026},
  number       = {1},
  pages        = {Paper No. 7, 90 pp.},
  mrnumber     = {5006587},
  doi          = {10.1007/s00205-025-02162-4},
}

@misc{BCJ26,
  author        = {Beekie, Rajendra and Chen, Shan and Jia, Hao},
  title         = {Optimal linear inviscid damping and vorticity depletion for non-monotonic shear flows},
  year          = {2026},
  eprint        = {2610.08912},
  archiveprefix = {arXiv},
  primaryclass  = {math.AP},
}

@article{BG26,
  author       = {Bian, Dongfen and Grenier, Emmanuel},
  title        = {Asymptotic behavior of solutions of the linearized {E}uler equations near a shear layer},
  journal      = {Peking Mathematical Journal},
  shortjournal = {Peking Math. J.},
  year         = {2026},
  note         = {Online first},
  doi          = {10.1007/s42543-026-00124-7},
}

@article{Bouchet,
  author       = {Bouchet, Freddy and Morita, Hidetoshi},
  title        = {Large time behavior and asymptotic stability of the 2{D} {E}uler and linearized {E}uler equations},
  journal      = {Physica D: Nonlinear Phenomena},
  shortjournal = {Phys. D},
  volume       = {239},
  year         = {2010},
  number       = {12},
  pages        = {948--966},
  mrnumber     = {2639613},
  doi          = {10.1016/j.physd.2010.01.020},
}

@misc{CJWZ,
  author        = {Chen, Qi and Jia, Hao and Wei, Dongyi and Zhang, Zhifei},
  title         = {Asymptotic stability of the {K}olmogorov flow at high {R}eynolds numbers},
  year          = {2025},
  eprint        = {2510.13181},
  archiveprefix = {arXiv},
  primaryclass  = {math.AP},
}

@article{CWZ2,
  author       = {Chen, Qi and Wei, Dongyi and Zhang, Zhifei},
  title        = {Linear inviscid damping and enhanced dissipation for monotone shear flows},
  journal      = {Communications in Mathematical Physics},
  shortjournal = {Comm. Math. Phys.},
  volume       = {400},
  year         = {2023},
  number       = {1},
  pages        = {215--276},
  mrnumber     = {4581474},
  doi          = {10.1007/s00220-022-04597-2},
}

@article{CKRZ,
  author       = {Constantin, P. and Kiselev, A. and Ryzhik, L. and Zlato\v{s}, A.},
  title        = {Diffusion and mixing in fluid flow},
  journal      = {Annals of Mathematics},
  shortjournal = {Ann. of Math. (2)},
  volume       = {168},
  year         = {2008},
  number       = {2},
  pages        = {643--674},
  mrnumber     = {2434887},
  doi          = {10.4007/annals.2008.168.643},
}

@software{NoEmbeddedEigenvalues,
  author       = {Dahne, Joel and Jia, Hao},
  title        = {{Joel-Dahne/EulerDipoleEigenvalues}},
  year         = {2026},
  howpublished = {GitHub},
  url          = {https://github.com/Joel-Dahne/EulerDipoleEigenvalues},
}

@article{Elgindi,
  author       = {Elgindi, Tarek M.},
  title        = {Remark on the stability of energy maximizers for the 2d {E}uler equation on {$\mathbb T^2$}},
  journal      = {Communications on Pure and Applied Analysis},
  shortjournal = {Commun. Pure Appl. Anal.},
  volume       = {23},
  year         = {2024},
  number       = {10},
  pages        = {1562--1568},
  mrnumber     = {4799455},
  doi          = {10.3934/cpaa.2024019},
}

@article{Iyer,
  author       = {Ionescu, Alexandru D. and Iyer, Sameer and Jia, Hao},
  title        = {On the stability of shear flows in bounded channels, {II}: non-monotonic shear flows},
  journal      = {Vietnam Journal of Mathematics},
  shortjournal = {Vietnam J. Math.},
  volume       = {52},
  year         = {2024},
  number       = {4},
  pages        = {851--882},
  mrnumber     = {4772726},
  doi          = {10.1007/s10013-023-00661-z},
}

@article{JiaVortex,
  author       = {Ionescu, Alexandru D. and Jia, Hao},
  title        = {Linear vortex symmetrization: the spectral density function},
  journal      = {Archive for Rational Mechanics and Analysis},
  shortjournal = {Arch. Ration. Mech. Anal.},
  volume       = {246},
  year         = {2022},
  number       = {1},
  pages        = {61--137},
  mrnumber     = {4487510},
  doi          = {10.1007/s00205-022-01815-y},
}

@article{IJacta,
  author       = {Ionescu, Alexandru D. and Jia, Hao},
  title        = {Non-linear inviscid damping near monotonic shear flows},
  journal      = {Acta Mathematica},
  shortjournal = {Acta Math.},
  volume       = {230},
  year         = {2023},
  number       = {2},
  pages        = {321--399},
  mrnumber     = {4628607},
  doi          = {10.4310/acta.2023.v230.n2.a2},
}

@article{JiaSIMA,
  author       = {Jia, Hao},
  title        = {Linear inviscid damping near monotone shear flows},
  journal      = {SIAM Journal on Mathematical Analysis},
  shortjournal = {SIAM J. Math. Anal.},
  volume       = {52},
  year         = {2020},
  number       = {1},
  pages        = {623--652},
  mrnumber     = {4062802},
  doi          = {10.1137/19M1273232},
}

@article{JiaG,
  author       = {Jia, Hao},
  title        = {Linear inviscid damping in {G}evrey spaces},
  journal      = {Archive for Rational Mechanics and Analysis},
  shortjournal = {Arch. Ration. Mech. Anal.},
  volume       = {235},
  year         = {2020},
  number       = {2},
  pages        = {1327--1355},
  mrnumber     = {4064200},
  doi          = {10.1007/s00205-019-01445-x},
}

@article{JiaUM,
  author       = {Jia, Hao},
  title        = {Uniform linear inviscid damping and enhanced dissipation near monotonic shear flows in high {R}eynolds number regime ({I}): {T}he whole space case},
  journal      = {Journal of Mathematical Fluid Mechanics},
  shortjournal = {J. Math. Fluid Mech.},
  volume       = {25},
  year         = {2023},
  number       = {3},
  pages        = {Paper No. 42, 38 pp.},
  mrnumber     = {4589745},
  doi          = {10.1007/s00021-023-00794-8},
}

@article{Zhiwu3,
  author       = {Lin, Zhiwu and Xu, Ming},
  title        = {Metastability of {K}olmogorov flows and inviscid damping of shear flows},
  journal      = {Archive for Rational Mechanics and Analysis},
  shortjournal = {Arch. Ration. Mech. Anal.},
  volume       = {231},
  year         = {2019},
  number       = {3},
  pages        = {1811--1852},
  mrnumber     = {3902476},
  doi          = {10.1007/s00205-018-1311-8},
}

@article{NaZh,
  author       = {Masmoudi, Nader and Zhao, Weiren},
  title        = {Nonlinear inviscid damping for a class of monotone shear flows in a finite channel},
  journal      = {Annals of Mathematics},
  shortjournal = {Ann. of Math. (2)},
  volume       = {199},
  year         = {2024},
  number       = {3},
  pages        = {1093--1175},
  mrnumber     = {4740211},
  doi          = {10.4007/annals.2024.199.3.3},
}

@article{Orr,
  author       = {Orr, W. M.},
  title        = {The stability or instability of the steady motions of a perfect liquid and of a viscous liquid. {P}art {I}: {A} perfect liquid},
  journal      = {Proceedings of the Royal Irish Academy. Section A: Mathematical and Physical Sciences},
  shortjournal = {Proc. Roy. Irish Acad. Sect. A},
  volume       = {27},
  year         = {1907},
  pages        = {9--68},
}

@article{dongyi,
  author       = {Wei, Dongyi and Zhang, Zhifei and Zhao, Weiren},
  title        = {Linear inviscid damping for a class of monotone shear flow in {S}obolev spaces},
  journal      = {Communications on Pure and Applied Mathematics},
  shortjournal = {Comm. Pure Appl. Math.},
  volume       = {71},
  year         = {2018},
  number       = {4},
  pages        = {617--687},
  mrnumber     = {3772399},
  doi          = {10.1002/cpa.21672},
}

@article{Dongyi2,
  author       = {Wei, Dongyi and Zhang, Zhifei and Zhao, Weiren},
  title        = {Linear inviscid damping and vorticity depletion for shear flows},
  journal      = {Annals of PDE},
  shortjournal = {Ann. PDE},
  volume       = {5},
  year         = {2019},
  number       = {1},
  pages        = {Paper No. 3, 101 pp.},
  mrnumber     = {3919496},
  doi          = {10.1007/s40818-019-0060-9},
}

@article{Dongyi3,
  author       = {Wei, Dongyi and Zhang, Zhifei and Zhao, Weiren},
  title        = {Linear inviscid damping and enhanced dissipation for the {K}olmogorov flow},
  journal      = {Advances in Mathematics},
  shortjournal = {Adv. Math.},
  volume       = {362},
  year         = {2020},
  pages        = {Paper No. 106963, 103 pp.},
  mrnumber     = {4050586},
  doi          = {10.1016/j.aim.2019.106963},
}

@article{Yudovich1,
  author       = {Yudovich, V. I.},
  title        = {Non-stationary flow of an ideal incompressible liquid},
  journal      = {USSR Computational Mathematics and Mathematical Physics},
  shortjournal = {USSR Comput. Math. Math. Phys.},
  volume       = {3},
  year         = {1963},
  number       = {6},
  pages        = {1407--1456},
  note         = {Translated from \v{Z}. Vy\v{c}isl. Mat. i Mat. Fiz. \textbf{3} (1963), no.~6, 1032--1066},
  mrnumber     = {158189},
  doi          = {10.1016/0041-5553(63)90247-7},
}

@article{Zillinger2,
  author       = {Zillinger, Christian},
  title        = {Linear inviscid damping for monotone shear flows in a finite periodic channel, boundary effects, blow-up and critical {S}obolev regularity},
  journal      = {Archive for Rational Mechanics and Analysis},
  shortjournal = {Arch. Ration. Mech. Anal.},
  volume       = {221},
  year         = {2016},
  number       = {3},
  pages        = {1449--1509},
  mrnumber     = {3509006},
  doi          = {10.1007/s00205-016-0991-1},
}

@article{Zillinger1,
  author       = {Zillinger, Christian},
  title        = {Linear inviscid damping for monotone shear flows},
  journal      = {Transactions of the American Mathematical Society},
  shortjournal = {Trans. Amer. Math. Soc.},
  volume       = {369},
  year         = {2017},
  number       = {12},
  pages        = {8799--8855},
  mrnumber     = {3710645},
  doi          = {10.1090/tran/6942},
}

\end{document}